\documentclass[11pt,letterpaper]{article}

\usepackage[T1]{fontenc}
\usepackage{float,graphicx,comment,color,lmodern,cite}
\usepackage[margin=1.1in]{geometry}
\usepackage{microtype}  
\usepackage{amsmath,amssymb,amsthm,mathtools}
\usepackage{microtype,enumitem,booktabs,array,longtable,needspace}
\usepackage[hidelinks]{hyperref}
\hypersetup{pdftitle={Maximal Independent Sets in Set Systems},pdfauthor={}}
\setlist[enumerate]{label=\textup{(\arabic*)},ref=(\arabic*),leftmargin=*,itemsep=.3em}
\setlist[itemize]{leftmargin=*,itemsep=.25em}
\numberwithin{equation}{section}
\newtheorem{theorem}{Theorem}[section]
\newtheorem{lemma}[theorem]{Lemma}
\newtheorem{proposition}[theorem]{Proposition}
\newtheorem{corollary}[theorem]{Corollary}
\newtheorem{conjecture}[theorem]{Conjecture}
\newtheorem*{conjecture*}{Conjecture}
\theoremstyle{definition}
\newtheorem{definition}{Definition}[section]
\theoremstyle{remark}
\newtheorem{remark}[theorem]{Remark}

\newcommand{\HH}{\mathcal H}

\newcommand{\PG}{\operatorname{PG}}

\newcommand{\vb}[1]{\boldsymbol{#1}}
\newcommand{\cS}{\mathcal S}
\newcommand{\cL}{\mathcal L}
\newcommand{\Zq}{\mathbb Z_q}
\newcommand{\F}{\mathbb F}

\author{}
\date{}

\begin{document}

\title{Maximal Independent Sets \\ {\normalsize in} \\ Set Systems}

\author{Jeroen Schillewaert\thanks{Department of Mathematics, University of Auckland, New Zealand.
       Supported by an FRDF grant of the University of Auckland. E-mail: j.schillewaert@auckland.ac.nz.} \and Jacques Verstra\"ete\thanks{Department of Mathematics, University of California at San Diego, CA, United States. Research supported by NSF Award DMS-2347832. E-mail: jverstraete@math.ucsd.edu.} }

\maketitle

\begin{abstract}
For a set system $\mathcal{S}$ of subsets of a set $X$, an \emph{independent set} is a subset $I$ of $X$
containing no members of $\mathcal{S}$. The set system $\mathcal{S}$ is an \emph{antichain} if no element of $\mathcal{S}$ is contained in any other.
In this paper, we prove a general theorem on the existence of small maximal independent
sets in set systems based on natural local density conditions. We suppose that
$\mathcal S$ is an antichain and that, for some sequence $d_2,d_3,\dots,d_r$, every $x \in X$ and $2 \leq i \leq r$ satisfy
\[
 N_i(x) = \{S\setminus\{x\}\colon S \in \mathcal S,\ x \in S,\ |S| = i\},
 \qquad |N_i(x)| = d_i,
\]
and let
$D(\vb{d}) = \max\{d_i^{1/(i - 1)} \colon 2 \leq i \leq r\}$.
We say $\mathcal{S}$ is \emph{$\delta$-sparse} if, for $2 \leq i \leq r$,
every set $A \subseteq X$ with $2 \leq |A| < i$ is contained in at most
$D(\vb{d})^{i - |A| - \delta}$ blocks of size $i$, and every pair of distinct
elements $x,y \in X$ satisfies
$|N_i(x)\cap N_i(y)| \leq D(\vb{d})^{i - 1 - \delta}$.
When $D(\vb{d}) > 1$, let $\rho(\vb{d})$ denote the unique positive
real number $\kappa$ satisfying
\begin{equation*}
\sum_{i = 2}^r d_i \kappa^{i - 1} = \log D(\vb{d}).
\end{equation*}
We prove that for $\delta > 0$, as $D(\vb{d}) \to \infty$, there is a maximal
independent set $I$ in $\mathcal{S}$ satisfying
\[ |I| = (1 + o(1)) \rho(\vb{d})|X|.\]
This improves on the known results for the size of maximum independent sets in uniform hypergraphs~\cite{AKPSS,BB,Dhawan,Dhawan2,DLMVY, DMV,RS,VWilson}, and the proof is accompanied by a polynomial-time algorithm for finding maximal independent sets of the prescribed size. In particular, in a linear $(r + 1)$-uniform hypergraph with $n$ vertices and maximum degree $\triangle \geq 2$, an independent set of size at least $(1 - o(1))n((\log \triangle)/r\triangle)^{1/r}$ is found. The theorem has many applications, for example when $r \geq 2$, we obtain a complete arc in $\PG(r,q)$ of size roughly
\[\bigl((r + 1)(r - 1)(r - 1)!\,q\log q\bigr)^{\frac{1}{r}},\]
improving the known results on Segre's complete arc problem to within a factor roughly $(\log q)^{1/r}$ of the counting lower bound. 
We also show that, for a wide class of multivariable polynomials
$f \in \mathbb Z[x_1,x_2,\ldots,x_r]$, there exists a maximal $f$-free set
$I \subseteq \mathbb F_q$ of size roughly $(q\log q)^{1/(r - 1)}$ and no such set of size smaller than roughly $q^{1/(r - 1)}$. 
\end{abstract}

\section{Introduction}

The problem of finding small maximal independent objects arises in a wide
range of combinatorial settings. The general question is to find a
collection which avoids a prescribed family of forbidden configurations, but
which cannot be enlarged while preserving this property. Examples include
maximal independent sets in graphs and hypergraphs, partial ovoids and
spreads in finite geometry, and sets avoiding polynomial configurations. Large independent sets have been studied extensively, but the problem of
finding small maximal independent sets is considerably more subtle. The
additional requirement of maximality imposes a global condition, and many
known results rely on arguments that exploit special features of the
particular problem under consideration. In this paper, we prove a general
theorem which gives bounds for the size of small maximal independent sets in families of sets under natural local degree conditions. We first introduce the general set-system framework in which our main
result is stated.

\textbf{Set systems.} We start with a set $X$ whose elements shall be referred to as \emph{atoms}, and a family $\mathcal{S}$ of
subsets of $X$ referred to as \emph{blocks}. Here $\mathcal{S}$ is referred to generically as a \emph{set system}, and we denote
by $X(\mathcal{S}) = X$ the set of atoms of $\mathcal{S}$.
If every atom in $X$ is contained in exactly $d$ blocks of $\mathcal{S}$ for some positive integer $d$, and every block in $\mathcal{S}$ has size $r$ and $|X| = n$, we accordingly refer to $\mathcal{S}$ as an \emph{$(n,d,r)$-system}. For $A \subseteq X$, let $d_r(A)$ denote the number of blocks in $\mathcal{S}$ containing $A$, and
for $x,y \in X$, let $d(x,y)$ denote the number of pairs $\{S,S'\} \subseteq \mathcal{S}$ of distinct blocks
such that $x \in S$ and $y \in S'$ and $S \setminus \{x\} = S' \setminus \{y\}$. A \emph{maximal independent set} in a set system $\mathcal{S}$ is a set $I$ of atoms containing no block but such that the addition of any atom to $I$ results in a set containing some block of $\mathcal{S}$. In this paper, we give an estimate for the smallest size of a maximal independent set when $\mathcal{S}$ is an $(n,d,r)$-system satisfying certain sparseness conditions:

\begin{definition}\label{def:sparse}
Let $\delta > 0$, $n \geq r \geq 2$ and $d \geq 1$. An $(n,d,r)$-system $\mathcal{S}$ is \emph{$\delta$-sparse} if
\begin{enumerate}[topsep=0pt,leftmargin=0.4in,noitemsep]
\item for any set $A \subseteq X$ with $2 \le |A| < r$, $d_r(A) \leq d^{(r - |A| - \delta)/(r - 1)}$ and
\item for any pair of distinct atoms $x,y \in X$, $d(x,y) \leq d^{1 - \delta/(r - 1)}$.
\end{enumerate}
\end{definition}

Ajtai, Koml\'{o}s, Pintz, Spencer and Szemer\'{e}di~\cite{AKPSS}
proved that so-called \emph{uncrowded} $(n,d,r + 1)$-systems have
independent sets of size at least
$Kn\left((\log d)/d\right)^{1/r}$ for some
$K = K(r) > 0$, but their proof implies these are not
necessarily maximal independent sets. A short proof of their theorem
was found by Verstra\"ete and Wilson~\cite{VWilson}. Bennett and
Bohman~\cite{BB} studied a randomized greedy algorithm in dense
$\delta$-sparse $(n,d,r + 1)$-systems --- where
$d = \Omega(n^{\epsilon})$ for some $\epsilon > 0$ --- and showed that
the maximal independent set produced has size at least of order
$n\left((\log d)/d\right)^{1/r}$. They leave the open problem of
proving that their maximal independent set produced has size of order
at most $n\left((\log d)/d\right)^{1/r}$. Our main theorem shows that
$\delta$-sparse set systems indeed have small maximal independent sets:

\begin{theorem}\label{main}
Let $\delta>0$, and let $r\ge1$ be an integer. There exist constants
$K = K(\delta,r) \geq 1$ and $\epsilon = \epsilon(\delta,r) > 0$ such that,
for every integer $d > 1$, every $\delta$-sparse $(n,d,r + 1)$-system
$\mathcal{S}$ contains a maximal independent set $I$ satisfying
\[ \Bigl||I| - n \cdot \Bigl(\frac{\log d}{rd}\Bigr)^{\frac{1}{r}}\Bigr| \leq \frac{K n}{(\log d)^{\epsilon}} \cdot \Bigl(\frac{\log d}{rd}\Bigr)^{\frac{1}{r}}.\]
\end{theorem}

It is straightforward to see that both conditions in the definition of $\delta$-sparseness are necessary for the conclusion of Theorem~\ref{main}.
A short calculation shows the lower bound in Theorem~\ref{main} does not hold for complete set systems (the maximum independent sets are all small) and the upper bound in Theorem~\ref{main}
does not hold for complete $r$-partite set systems (the maximal independent sets are all very large). For both of these systems, one verifies that condition (1) of $\delta$-sparseness
does not hold. Condition (2) is also necessary, as all maximal independent sets in appropriate blowups of Steiner $r$-tuple systems -- see Section \ref{sec:blowups} -- are much larger than the bound in Theorem~\ref{main}. We also remark that for ``almost all'' $\delta$-sparse $(n,d,r + 1)$-systems, no asymptotically smaller maximal independent sets can be found than those provided by Theorem \ref{main}, as shown by the configuration model of random regular hypergraphs -- see Section \ref{sec:rrh}.

{\bf Maximum independent sets.}
Theorem \ref{main} gives the current best lower bound for the size of a maximum independent set
in linear hypergraphs with maximum degree $\triangle$. Let us say that a set system $\mathcal{S}$ with blocks of size $(r + 1)$ is $(\delta,\triangle)$-sparse if (1) and (2) in Definition 1.1 hold with $d$ replaced by $\triangle$ and, in addition, $\max\{d_{r+1}(x) : x \in X\} \leq \triangle$. 
Then the following holds:

\begin{theorem}\label{thm:maxdegree}
Let $r\geq 1$ be an integer, let $\delta > 0$ be fixed and let $\triangle \geq 2$, and let $\mathcal{S}$ be a $(\delta,\triangle)$-sparse set system with $n$ atoms and blocks of size $r + 1$. Then $\mathcal{S}$ contains an independent set of size at least 
\[ (1 - o(1))n \cdot \Bigl(\frac{\log \triangle}{r\triangle}\Bigr)^{\frac{1}{r}} \quad \mathrm{ as } \quad \Delta\to\infty.\]
\end{theorem}
 
 We prove this theorem via a reduction to Theorem \ref{main} in Section \ref{sec:maxdegree}. 
Theorem \ref{thm:maxdegree} gives the current best lower bound for the size of a maximum independent set in $(\delta,\triangle)$-sparse hypergraphs -- see~\cite{AKPSS,Dhawan,Dhawan2,NV,DMV,VWilson}, together with a deterministic polynomial-time algorithm for finding such an independent set -- see Section \ref{sec:alg}. It is plausible that for a fixed constant $c > 1$, no independent set of size $cn((\log \triangle)/r\triangle)^{1/r}$ can be found efficiently, even in the context of random hypergraphs with average degree $\triangle$ -- see Dhawan et al.~\cite{Dhawan2} and the references therein. This particular problem for graphs -- efficiently finding for a fixed $c > 1$ an independent set of size $cn(\log d)/d$ in the random graph $G_{n,d/n}$ -- is a well-studied problem in combinatorics and theoretical computer science, and a major question of Karp~\cite{Karp1976}.

\textbf{Non-uniform set systems.} Theorem~\ref{main} is a special case of a substantially more general theorem on set systems with mixed block sizes. To describe the general theorem,
we require some notation and definitions. For $r \geq 2$,
a set system $\mathcal{S}$ of distinct subsets of a set $X$ of atoms has \emph{upper rank} at most $r$ if $\max\{|S| \colon S \in \mathcal{S}\} \leq r$.
We call $\mathcal{S}$ an \emph{antichain} if no block in $\mathcal{S}$ is a subset of any other block in $\mathcal{S}$.
Let $\vb{d} = (d_2,d_3,\dots,d_r)$
be a sequence of non-negative integers, which we call a {\em degree vector}. An $(n,\vb{d})$\emph{-antichain} is an antichain $\mathcal{S}$ of upper rank at most $r$ on a set $X$ of $n$ atoms
such that for every atom $x \in X$ and $2 \leq i \leq r$, the number $d_i(x)$ of blocks in $\mathcal{S}$ of size $i$
containing $x$ is exactly $d_i$. We sometimes refer to $d_i(x)$ as the {\em degree} of $x$. In particular, an $(n,d,r)$-system is an $(n,\vb{d})$-antichain with $\vb{d} = (0,0,\dots,0,d_r)$ where $d_r = d$.
For a set $A \subseteq X$, let
\begin{equation*}
N_i(A) = \{S \setminus A \colon |S| = i, A \subseteq S \in \mathcal{S}\} \qquad\text{and}\qquad d_i(A) = |N_i(A)|.
\end{equation*}
If $A = \{x\}$, we write $N_i(A) = N_i(x)$ and $d_i(A) = d_i(x)$. For distinct atoms $x,y \in X$, let
\begin{equation*}
N_i(x,y) = N_i(x) \cap N_i(y) \qquad\text{and}\qquad d_i(x,y) = |N_i(x,y)|.
\end{equation*}
For a sequence $\vb{d} = (d_2,d_3,\dots,d_r)$ of non-negative integers, it is convenient to define
\begin{equation*}
D(\vb{d}) = \max\{d_i^{\frac{1}{i - 1}} \colon 2 \leq i \leq r\}.
\end{equation*}
When $D(\vb{d}) > 1$, let $\rho(\vb{d})$ denote the unique positive
real number $\kappa$ satisfying
\begin{equation*}
\sum_{i = 2}^r d_i \kappa^{i - 1} = \log D(\vb{d}).
\end{equation*}
Existence and uniqueness follow because the sum increases strictly from zero to infinity for
$\kappa \ge 0$. We now come to the definition of
$\delta$-sparse set systems for systems with mixed block sizes.

\begin{definition}
Let $\delta > 0$, $n \geq r \geq 2$, and let $\vb{d} = (d_2,d_3,\dots,d_r)$
be a sequence of non-negative integers with $D(\vb{d}) > 1$.
An $(n,\vb{d})$-antichain is \emph{$\delta$-sparse} if for $2 \leq i \leq r$,
\begin{enumerate}[topsep=0pt,leftmargin=0.4in,noitemsep]
\item for any set $A \subseteq X$ with $2 \le |A| < i$, $d_i(A) \leq D(\vb{d})^{i - |A| - \delta}$ and
\item for any pair of distinct atoms $x,y \in X$, $d_i(x,y) \leq D(\vb{d})^{i - 1 - \delta}$.
\end{enumerate}
\end{definition}

A \emph{maximal independent set} in an $(n,\vb{d})$-antichain $\mathcal{S}$ is a set $I \subseteq X$ containing no blocks from $\mathcal{S}$
such that for any $x \in X \setminus I$, $I \cup \{x\}$ contains at least one block from $\mathcal{S}$. We now state the main theorem of the paper:

\begin{theorem}\label{thm:main}
Let $\delta>0$, and let $r\ge2$ be an integer. There exist constants
$K = K(\delta,r) \geq 1$ and $\epsilon = \epsilon(\delta,r) > 0$ such that,
for every sequence $\vb{d} = (d_2,d_3,\dots,d_r)$ of non-negative integers
with $D(\vb{d}) > 1$, every $\delta$-sparse $(n,\vb{d})$-antichain has a
maximal independent set $I$ satisfying
\begin{equation*}
\bigl||I| - \rho(\vb{d})n\bigr| \leq K(\log D(\vb{d}))^{-\epsilon} \cdot \rho(\vb{d})n.
\end{equation*}
\end{theorem}

It is not hard to see that Theorem \ref{main} follows from Theorem \ref{thm:main} by setting parameters appropriately. We do not make a concerted attempt at optimizing the constants $K$ and $\epsilon$ in the proof of Theorem \ref{thm:main}. The antichain condition in this theorem is necessary, without which the conclusion of the theorem may not hold -- see Section \ref{subsec:antichain} for the details.
We note that if $D(\vb{d}) \geq e$ and $\max_{2 \leq i \leq r}d_i \geq 3$, then
\begin{equation*}
\rho(\vb{d}) \leq \min\Bigl\{\Bigl(\frac{\log d_{i + 1}}{i d_{i + 1}}\Bigr)^{\frac{1}{i}}
\colon 1 \leq i < r,\ d_{i + 1} \geq 3\Bigr\}.
\end{equation*}
This comparison, together with the bounded-degree argument in
Subsection~\ref{subsec:process-performance}, gives the following
corollary for every degree vector with $\max_i d_i \geq 3$:
\begin{corollary}\label{cor:main}
For any $\delta > 0$ and integer $r \geq 2$, there exist constants
$K = K(\delta,r) \geq 1$ and $\epsilon = \epsilon(\delta,r) > 0$ such that,
for every sequence $\vb{d} = (d_2,d_3,\dots,d_r)$ of non-negative integers
with $\max_{2 \leq i \leq r}d_i \geq 3$, every $\delta$-sparse
$(n,\vb{d})$-antichain has a maximal independent set $I$ satisfying
\begin{equation*}
|I| \leq (1 + K(\log D(\vb{d}))^{-\epsilon}) \cdot n \cdot \min_{1 \leq i < r} \Bigl\{\Bigl(\frac{\log d_{i + 1}}{id_{i + 1}}\Bigr)^{\frac{1}{i}} \colon d_{i + 1} \geq 3\Bigr\}.
\end{equation*}
\end{corollary}
This corollary is more explicit than Theorem~\ref{thm:main}, and therefore more convenient to apply. Without the antichain condition, the conclusions of these results do not hold -- we give a construction in Section \ref{subsec:antichain}. 

{\bf Deterministic Polynomial-Time Algorithm.} Theorems \ref{main} and \ref{thm:main} and Corollary~\ref{cor:main} have accompanying polynomial-time deterministic algorithms for producing a maximal independent set of the given size. In the following theorem, we assume a set system $\mathcal{S}$ is given by its $X(\mathcal{S}) \times \mathcal{S}$ incidence matrix, whose entry $(x,S)$ is 1 if $x \in S$ and zero otherwise, and polynomial time refers to a polynomial function of the dimensions of the incidence matrix.

\begin{theorem}
Let $\delta>0$, and let $r\ge2$ be an integer. There exist constants
$C = C(\delta,r)$ and $\epsilon_0 = \epsilon_0(\delta,r) > 0$ such that for every sequence $\vb{d} = (d_2,d_3,\dots,d_r)$ of non-negative integers
with $D(\vb{d}) > 1$ and every $\delta$-sparse $(n,\vb{d})$-antichain, there is a polynomial time algorithm  which produces a
maximal independent set $I$ satisfying
\begin{equation*}
\bigl||I| - \rho(\vb{d})n\bigr| \leq C(\log D(\vb{d}))^{-\epsilon_0} \cdot \rho(\vb{d})n.
\end{equation*} 
\end{theorem}

The algorithm is given in Section \ref{sec:alg}, and also gives an independent set of size at least $(1 - o(1))n((\log \triangle)/r\triangle)^{1/r}$ in any $n$-vertex $(r + 1)$-uniform $(\delta,\triangle)$-sparse hypergraph -- see Section \ref{sec:maxdegree}.

Next, we give some select applications of Theorem \ref{thm:main}. We anticipate that 
Theorem \ref{thm:main} has many further applications other than those which we select to present here.

\textbf{Zero sets of polynomials.} Let $\mathbb F_q$ denote the
finite field of order $q$, and let
$\mathbb F_q[x_1,x_2,\dots,x_r]$ and $\mathbb Z[x_1,x_2,\dots,x_r]$ denote the rings of polynomials in
variables $x_1,x_2,\dots,x_r$ over the finite field $\mathbb F_q$ and integers $\mathbb Z$, respectively. The \emph{total
degree} of a polynomial $f$ in these rings, denoted $\deg(f)$, is the largest value of $d_1 + d_2 + \dots + d_r$ over all
monomials $x_1^{d_1}x_2^{d_2} \dots x_r^{d_r}$ with a nonzero
coefficient in the expansion of $f$ as a finite sum of monomials. 
A set $I \subseteq \mathbb F_q$ is \emph{$f$-free} if
$f(x_1,\ldots,x_r) \neq 0$ in $\mathbb F_q$ for every 
$(x_1,\ldots,x_r) \in I^r$, and \emph{maximal $f$-free} if it is
$f$-free while $I \cup \{x\}$ is not $f$-free for all
$x \in \mathbb F_q \setminus I$. 
If the diagonal polynomial $f(x,x,\ldots,x)$ is nonzero over $\mathbb F_q$
and $d=\deg(f)$, then, for $q>d$, every maximal $f$-free set has size
at least $c_{r,d}q^{1/(r-1)}$, where
$c_{r,d}=[1+d(2^r-1)]^{-1/(r-1)}$.
Indeed, if $I$ is maximal and $s=|I|$, then $s\ge1$ and
\[
 q\le s+d\bigl((s+1)^r-s^r\bigr)
   \le[1+d(2^r-1)]s^{r-1}.
\]
For each element outside $I$, a forbidden tuple specifies its nonempty
set of positions and the entries from $I$ in the other positions. The
resulting polynomial in that element is nonzero, since otherwise its
value at any element of $I$ would contradict freeness, and has at most
$d$ roots. This proves the displayed inequality.
The order $\sqrt q$ is attained for $f(x,y,z)=x+y-z$ over prime fields
and for $f(x,y,z)=xy-z$ over finite fields, by
Haviv and Levy~\cite[Theorem 1.5]{HavivLevy}. For the second assertion,
the cyclic group $\mathbb F_q^\times$ has order $q-1$~\cite{LidlNiederreiter}, and an exponent
isomorphism transfers a maximal sum-free set to a maximal product-free
set; adjoining $0$ is forbidden by $0\cdot0=0$.
These respectively correspond to maximal sum-free and product-free
subsets of $\mathbb F_q$.
In this section, we give sufficient conditions on 
polynomials $f \in \mathbb Z[x_1,\ldots,x_r]$ for which this lower bound is almost tight. A polynomial $f \in \mathbb Z[x_1,x_2,\dots,x_r]$ is \emph{symmetric}
if, for every permutation $\sigma$ of $\{1,2,\dots,r\}$,
\[ f(x_{\sigma(1)},\ldots,x_{\sigma(r)}) = f(x_1,\ldots,x_r). 
\]
We impose the following 
hypotheses:

\begin{definition}
For a fixed prime power $q$, a polynomial $f \in \mathbb F_q[x_1,x_2,\dots,x_r]$ is \emph{eligible} if:
\begin{enumerate}[topsep=0pt,leftmargin=0.4in,noitemsep]
\item  $f$ is symmetric in its $r$ variables and is homogeneous and
\item $f$ is absolutely irreducible, that is, irreducible over the algebraic closure
       $\overline{\mathbb F}_q$ of $\mathbb F_q$.
\end{enumerate}
\end{definition}
An integer polynomial is eligible over $\mathbb F_q$ if its reduction
$f_q$ has the same total degree and is eligible in this sense. When
working over $\mathbb F_q$, we write $f$ for this reduction; evaluations
and polynomial identities are then taken in $\mathbb F_q$.
It is convenient for a polynomial $f \in \mathbb Z[x_1,x_2,\dots,x_r]$ to let $\gamma_q(f)$ denote the smallest size of a maximal $f$-free subset of $\mathbb F_q$. The following theorem
is an application of Theorem~\ref{thm:main} (or Corollary~\ref{cor:main}), together with some mild algebra to check the required conditions of that theorem:

\begin{theorem}\label{thm:poly-maximal-sets}
Let $r\ge3$ and $d\ge1$ be integers. Then there exist 
$\gamma_1 = \gamma_1(r,d) \ge 4$ and
$\gamma_2 = \gamma_2(r,d) > 0$ such that for every prime power
$q \ge \gamma_1$ and every polynomial
$f \in \mathbb Z[x_1,x_2,\dots,x_r]$ of degree $d$ that is eligible over $\mathbb F_q$,
\begin{equation}\label{eq:poly-two-sided}
\gamma_q(f) \le 
4d\bigl((r - 2)(r - 2)!\bigr)^{\frac{1}{r - 1}}
(q\log q)^{\frac{1}{r - 1}}.
\end{equation}
If $f(1,1,\dots,1) \ne 0$ in $\mathbb F_q$, then
$\gamma_q(f) \ge \gamma_2q^{1/(r - 1)}$.
\end{theorem}

The lower bound in Theorem~\ref{thm:poly-maximal-sets} is a simple
counting argument. For homogeneous
but asymmetric polynomials, the conclusion
(\ref{eq:poly-two-sided}) may also not hold: for example, if
$f(x,y,z) = x$, then $\gamma_q(f) = q - 1$. Condition (2) of eligibility
is necessary, as the example $f(x,y,z) = xyz$ over $\mathbb F_q$
shows, where $\gamma_q(f) = q - 1$.
For an example which is symmetric and homogeneous and irreducible over $\mathbb F_q$, but not absolutely irreducible, let $q$ be odd with
$q \equiv 2 \pmod{3}$, and let $\kappa_0 \in \mathbb F_q$ be a nonsquare. We set
\[
 f(x,y,z) = (x + y + z)^4 - \kappa_0(xy + xz + yz)^2.
\]
If $\theta^2 = \kappa_0$ in $\mathbb F_{q^2}$, then
$f = Q_+Q_-$, where
$Q_\pm = (x + y + z)^2\pm\theta(xy + xz + yz)$. The matrices of $Q_\pm$ have
nonzero eigenvalues $3\pm\theta$ and $\mp\theta/2$ with multiplicity two,
so both quadratics are absolutely irreducible. They are interchanged by
$\operatorname{Gal}(\mathbb F_{q^2}/\mathbb F_q)$, and hence $f$ is
irreducible over $\mathbb F_q$. If $f(x,y,z) = 0$ with $x,y,z \in \mathbb F_q$,
the nonsquareness of $\kappa_0$ forces $x + y + z = xy + xz + yz = 0$. If $x \ne 0$, then
$t = y/x$ satisfies $t^2 + t + 1 = 0$, contrary to $3\nmid q - 1$; if $x = 0$, then
$y = z = 0$. Thus the origin is the only $\mathbb F_q$-zero and
$\gamma_q(f) = q - 1$.

Asymmetric homogeneous examples are also interesting, and it is possible to show that 
Theorem \ref{thm:main} applies effectively to such cases. We consider the special case of {\em Pythagorean triple free} sets, namely, 
$f$-free subsets of $\mathbb F_q$ where $f(x,y,z) = x^2 + y^2 - z^2$.  In Section \ref{sec:Pythag}, we prove the following theorem for Pythagorean triples:

\begin{theorem}\label{thm:Pythag}
Let $f(x,y,z) = x^2 + y^2 - z^2$. Then for odd values of $q \neq 3^{2m}$, $\gamma_q(f) \leq (3 + o(1))\sqrt{q\log q}$ as $q \rightarrow \infty$, 
whereas if $q = 3^{2m}$ then $\gamma_q(f) \leq 2(\sqrt{q} - 1)$ and if $q$ is even then $\gamma_q(f) \leq 3\sqrt{q/2}$. In all cases, $\gamma_q(f) \geq (1 - o(1))\sqrt{q/3}$ as $q \rightarrow \infty$. 
\end{theorem}

If $q$ is even, then Pythagorean triple free sets are simply sum-free sets. More generally, the problem of $f$-free sets is well studied when $f$ is linear~\cite{GiudiciHart,HavivLevy,Plagne2002,Ruzsa,
RedmanRoseWalker,HancockTreglown1,HancockTreglown2}, and has interest in the more general setting of abelian groups. We discuss this in Section \ref{sec:Pythag}. One of the main open questions is as follows:

\begin{conjecture}
For each $r \geq 3$, there exists a symmetric homogeneous polynomial $f \in \mathbb Z[x_1,x_2,\dots,x_r]$ that is absolutely irreducible over $\mathbb Q$, satisfies $f(1,\ldots,1)\ne0$ over $\mathbb Z$, and has
$\{\gamma_q(f)q^{-1/(r - 1)}\}$ unbounded over primes $q$. 
\end{conjecture}

By Theorem \ref{thm:poly-maximal-sets}, $\gamma_q(f)q^{-1/(r - 1)}$ cannot grow faster than roughly $(\log q)^{1/(r - 1)}$. Perhaps a good candidate is $f(x_1,x_2,\dots,x_r) = x_1^2 + x_2^2 + \dots + x_r^2$; we do not know the order of magnitude of $\gamma_q(f)$ for primes $q$ even when $r = 3$.
Another very interesting question is to determine, for each fixed $r\ge3$, whether there exists a polynomial $f \in \mathbb Z[x_1,x_2,\dots,x_r]$ with $f(x,x,\ldots,x)\not\equiv0$ such that every $f$-free subset of $\mathbb F_q$ has size at most $q^{1/(r - 1) + o(1)}$ as $q \rightarrow \infty$ through primes for which the diagonal polynomial remains nonzero.

\textbf{Geometric Applications.} We discuss applications in finite
geometry. A \emph{complete arc} in the $r$-dimensional projective space
$\PG(r,q)$ over $\mathbb F_q$ is a maximal set of points every subset
of at most $r + 1$ of which is projectively independent. Segre~\cite{S1}
raised the fundamental problems of studying the minimum and maximum sizes
of complete arcs, as well as their structure. More generally, for
$2 \le k \le r$, a \emph{$k$-complete arc} is a maximal set of points
every subset of at most $k + 1$ of which is projectively independent.
For a set with at least $k + 1$ points, this independence condition is
equivalent to having no $k + 1$ points in a $(k - 1)$-dimensional
subspace. In the case $k = r = 2$, Kim and Vu~\cite{KV} proved that
there exist complete arcs in $\PG(2,q)$ of size
$O(\sqrt{q}\log^{300}q)$. This comes within polylogarithmic factors of
the counting lower bound observed by Lunelli and Sce~\cite{LS}.
Bartoli, Faina, Marcugini and Pambianco~\cite{BMP} subsequently extended this
result to $2$-complete arcs in higher-dimensional projective spaces
$\PG(r,q)$.
The independent sets in the set system whose atoms are the points of
$\PG(r,q)$ and whose blocks are the inclusion-minimal dependent sets of
at most $k + 1$ points are precisely the sets satisfying this
independence condition. Each such block $S$ spans an
$(|S| - 2)$-dimensional subspace. Applying Corollary~\ref{cor:main} gives the following result:

\begin{theorem}\label{thm:k-complete}
Let $r\ge k\ge2$ be integers, and let $q$ be a prime power. Then there is a $k$-complete arc $I$ of points of $\PG(r,q)$ such that
\begin{equation}\label{eq:pg-circuit-size}
 |I| = (1 + o_{r,k}(1))
 ((r + 1)(k - 1)(k - 1)!\,
 q^{r - k + 1}\log q)^{\frac{1}{k}}.
\end{equation}
\end{theorem}

Fixing $r$ and $k$, the term $o_{r,k}(1)$ denotes a function of $q$ which decays to zero as $q \rightarrow \infty$.
For $r = k = 2$, this sharply improves the aforementioned results, giving $|I| = (1 + o(1))\sqrt{3q\log q}$.
An elementary lower bound from counting shows that no $k$-complete arc can have size less than
\begin{equation*}
(1 - o_{r,k}(1))(k! \cdot q^{r - k + 1})^{\frac{1}{k}},
\end{equation*}
so the theorem is tight up to a factor of order $(\log q)^{1/k}$. We may also apply Theorem~\ref{thm:main} to the problem of small maximal
subsets of $\PG(2,q)$ without $k+1$ collinear points. The parameter calculation in Section~\ref{sec:pg-circuits} gives the following:

\begin{theorem}\label{thm:plane-collinear}
Let $k\ge2$ be an integer, and let $q$ be a prime power. Then there is a maximal set $I$ of points of $\PG(2,q)$ with no $k + 1$ points in $I$ collinear such that
\begin{equation}\label{eq:plane-collinear-size}
 |I| = (1 + o_{k}(1))
 ((k + 1)(k - 1)!\,
 q^{k - 1}\log q)^{\frac{1}{k}}.
\end{equation}
\end{theorem}

The problem of $k$-complete arcs may be studied in other classical geometries such as \emph{generalized quadrangles} (Van Maldeghem~\cite{vanMaldeghem}) and \emph{polar spaces}; the extensive literature on this problem is discussed by De Beule, Klein, Metsch and Storme~\cite{deBKMS}, Thas~\cite{Thas}, and the references therein. In this context, $1$-complete arcs are called \emph{maximal partial ovoids}. Ovoids have a long history going back to the seminal work of Tits~\cite{Tits1962}, and their construction and classification have been studied extensively. From Theorem~\ref{main}, we obtain substantially stronger bounds in almost all instances, although sharper bounds and exact values are known for several individual classical geometries~\cite{BKMS,CDWS,Metsch,MS,AEL,GacSzo,GSVM}. By duality, bounds on maximal partial ovoids yield bounds on maximal partial line spreads in the dual geometries; for example, the elliptic quadric $Q^-(5,q)$ is dual to the Hermitian generalized quadrangle $H(3,q^2)$. In the following theorem, for $q,t \geq 2$, let $GQ(q,t)$ denote any generalized quadrangle where the lines have size $q + 1$ and each point lies in $t + 1$ lines. The standard point counts, residue counts, and generalized-polygon parameters  may be found in~\cite{HT,Thas,vanMaldeghem,deBKMS}. The following theorem gives sharply stronger estimates on the sizes of maximal partial ovoids in a number of classical geometries (the list here is by no means exhaustive). In all cases here, the appropriate set system $\mathcal{S}$ is the graph 
whose vertex set is the set of points of the geometry, and whose edges are pairs of collinear points. We record a value of $\delta$ for which $\mathcal{S}$ is $\delta$-sparse in each case, after which it is a matter of substituting the appropriate parameters in Theorem \ref{main}.

\begin{theorem}\label{thm:classical-geometries}
Let $n\ge2$ be an integer. For the classical polar spaces in the table,
let $q$ be a prime power; for $GQ(q,t)$, let $q,t\ge2$ be integers.
For each geometry $\mathcal Q$ in the table, let $\gamma(\mathcal Q)$
denote the minimum size of a maximal partial ovoid in $\mathcal Q$.
Then, as $q\to\infty$, the following hold:
\begingroup
\renewcommand{\arraystretch}{1.2}
\begin{center}
\begin{tabular}{l|c|c}
\hline
{\rm Geometry} $\mathcal{Q}$ & $\delta$ & {\rm Asymptotic upper bound on} $\gamma(\mathcal{Q})$\\
\hline
$Q(2n,q),q$ odd & $\frac{1}{4n-4}$ & $(2n - 2)q\log q$ \\
\hline
$Q^+(2n + 1,q)$ & $\frac{1}{4n-2}$ & $(2n - 1)q\log q$ \\
\hline
$Q^-(2n + 1,q)$ & $\frac{1}{4n-2}$ & $(2n - 1)q\log q$ \\
\hline
$H(2n,q^2)$ & $\frac{1}{4n - 3}$ & $(4n - 3)q^2\log q$ \\
\hline
$H(2n - 1,q^2)$ & $\frac{1}{4n - 3}$ & $(4n - 5)q^2\log q$ \\
\hline
$GQ(q,t)$ & $\frac{1}{4}$ & $3q\log q$.\\
\hline
\end{tabular}
\end{center}
\endgroup
\end{theorem}

For the polar space $Q(2n,q)$ of rank $n$ with $n \geq 3$ the lower bound is generally $q + 1$ and the preceding best upper bound is asymptotic to $q^n$ as $q \rightarrow \infty$. The only lower bound on maximal partial ovoids known for polar spaces $Q^+(2n + 1,q)$ and $Q^-(2n + 1,q)$ of odd projective dimension with $n \geq 2$ is $2q + 1$,
whereas the upper bounds are $q^n + 1$ and $2 + ((q^n - 1)(q^2 + q + 2))/(2(q + 1))$ -- as discussed by De Beule, Klein, Metsch and Storme~\cite{deBKMS}. Schillewaert and Verstra\"ete~\cite{SchillewaertV} give a short probabilistic
proof of the upper bound $\gamma(\mathcal{Q}) = O(q(\log q)^5)$ when $\mathcal{Q}$ is the elliptic quadric $Q^-(5,q)$, and Theorem~\ref{thm:classical-geometries}
improves this to  $(3 + o(1))q\log q$. Ebert and Hirschfeld~\cite{EH} proved that any maximal partial ovoid in $Q^-(5,q)$ has at least $2q + 1$ points for $q \geq 2$ and at least $2q + 2$ points when $q \geq 5$. The following question was posed in~\cite[Problem 2]{SchillewaertV}; we state it as a conjecture:

\begin{conjecture}
Every maximal partial ovoid in the elliptic quadric $Q^-(5,q)$ has size at least $q\omega(q)$ for some function $\omega(q)\to\infty$ as $q\to\infty$.  
\end{conjecture}

More generally, we would like to know for which of the geometries in the table the logarithmic terms in the
upper bound on $\gamma(\mathcal{Q})$ are necessary.

\textbf{Independent sets in Steiner systems.} A \emph{Steiner
$r$-tuple system} is a set system of $r$-element subsets of some finite
point set $X$ such that every pair of elements of $X$ is contained in
exactly one of the $r$-element sets. For a Steiner $(r + 1)$-tuple system with
$n$ points, it is necessary that $r\mid n - 1$ and
$r(r + 1)\mid n(n - 1)$. For $r \geq 2$, an $n$-point
Steiner $(r + 1)$-tuple system is regular of degree $d = (n - 1)/r$ and
is $1$-sparse, so Theorem~\ref{main} gives the following consequence.

\begin{theorem}\label{thm:steiner-maximal}
As $n \rightarrow \infty$ through orders for which an $n$-point Steiner triple system exists, every such system contains a maximal independent set of
size $(1 + o(1))\sqrt{n\log n}$.
\end{theorem}

This theorem follows directly from Theorem~\ref{main}, and the same substitution applies to Steiner $(r+1)$-tuple systems. Section~\ref{sec:Steiner} constructs Steiner triple systems showing that its bound on the minimum maximal size is asymptotically sharp. An important open question, raised in~\cite{EustisV,KMV}, is the following:

\begin{conjecture}
Every $n$-vertex Steiner triple system contains an independent set of size $(1 - o(1))\sqrt{3n\log n}$ as $n \rightarrow \infty$.
\end{conjecture}

This would be asymptotically best possible by the construction in~\cite{EustisV}. The conjecture could be extended to Steiner $(r + 1)$-tuple systems, where the plausible size of a maximum independent set is $(1 - o(1))n^{1 - 1/r}((r + 1)\log n)^{1/r}$. Bennett and Frieze~\cite[Theorem 1]{BennettFrieze} show that in an $n$-vertex random $d$-regular $(r+1)$-uniform hypergraph, the size of a largest independent set is $(1+o_d(1))n((r+1)(\log d)/rd)^{1/r}$ with high probability as $n\to\infty$, with $r$ fixed and then $d\to\infty$.

\medskip

\textbf{Organization.} Section~\ref{sec:process} describes the idealized
and actual processes and states the quantitative conclusion.
Section~\ref{sec:tools} collects the probabilistic tools.
The proof is organized around the two random choices in each round:
Section~\ref{sec:general-focus} establishes the focus estimates and
collects the required focus inequalities, while Section~\ref{sec:construction-analysis}
analyzes the construction on a fixed good focus and proves that a good
construction round exists. Section~\ref{sec:completing-process} proves
termination and the required bounds on the final maximal independent set.
In Section \ref{sec:maxdegree}, we give the reduction of Theorem~\ref{thm:maxdegree} to Theorem~\ref{main}.
Sections~\ref{sec:polynomial-zero-sets} --~\ref{sec:pg-circuits} 
give the applications to complete arcs and polynomial zero sets, as well as Theorem \ref{thm:Pythag}.  In Section \ref{sec:blowups}, we show that  condition (1) alone of $\delta$-sparseness is not sufficient to guarantee the conclusion of Theorem \ref{main}. In Section \ref{subsec:antichain} we construct a $\delta$-sparse set system that is not an antichain for which the conclusion of Theorem \ref{thm:main} does not hold. In the remainder of  Section~\ref{sec:Steiner}, we construct $n$-point Steiner triple systems whose smallest maximal independent sets have size $(1 + o(1))\sqrt{n\log n}$, showing the asymptotic sharpness of Theorem~\ref{thm:steiner-maximal}. In Section \ref{sec:rrh}, we show that Theorem \ref{main} is asymptotically tight by appealing to the configuration model of random regular hypergraphs.  In Section \ref{sec:alg}, we give a deterministic algorithm which runs in polynomial time to give maximal independent sets as guaranteed by Theorem \ref{thm:main}.

\section{The random process}\label{sec:process}

In this section, we give the random process for constructing the maximal
independent set required for Theorem~\ref{thm:main}.
First, we give an idealized process whose parameters will accurately
describe the tracking of the actual process, which is described in
Sections~2.3 and~2.4.

\Needspace{16\baselineskip}

\subsection{Definition of an idealized process}\label{subsec:idealized-process}

Let $T \geq 1$ and $q_t \in [0,1]$ for $0 \leq t < T$. Suppose
at a stage $t$ with $0 \leq t < T$, we have an antichain $\mathcal{S}^t$ with atom set $X^t$.
For a set $Y \subseteq X^t$, we define
\begin{equation*}
\partial Y = \{x \in X^t \colon \exists S \subseteq Y, S \cup \{x\} \in \mathcal{S}^t\}.\end{equation*}
The following is the process we follow for $0 \leq t < T$: having defined $\mathcal{S}^t$ and $X^t$ for some $t < T$,

\begin{enumerate}[leftmargin=*,itemsep=.35em]
\item Let the elements of $Y^t \subseteq X^t$ be chosen independently with probability $q_t$.
\item We define the independent set
\[
I^t = Y^t\setminus\bigcup\{S \in \mathcal S^t\colon S \subseteq Y^t\}.
\]
\item Let $\mathcal S^{t + 1}$ denote the set system whose atom set is 
\[ X^{t + 1} = X^t\setminus(I^t\cup\partial I^t)\]
and whose block set is the distinct inclusion-minimal members of
\[
\{S\setminus I^t : S \in \mathcal S^t, S \setminus I^t\subseteq X^{t + 1}\}.
\]
\end{enumerate}

Once the process hits time $t = T$, we gather the independent set
$J^T = I^0 \cup I^1 \cup \dots \cup I^{T - 1}$ so far, and add to it
a maximal independent set $J'$ in $X^T$ to obtain a maximal independent
set $I = J^T \cup J'$ in the entire set system. We have deliberately
not specified $T$ and $q_t$ in this process. The quantities $T$ and $q_t$ have to be
chosen very carefully as a function of the numbers of blocks of
different sizes containing each atom in $\mathcal{S}^t$, so that the
process arrives at the independent set required by
Theorem~\ref{thm:main} with positive probability. In the next subsection,
we describe these choices, and indicate for the reader why the process
should terminate with a maximal independent set of the required size.

\subsection{Parameters of the idealized process}

A set system $\mathcal{S}$ is \emph{locally sparse} if for $2 \leq k \leq 5$ and any set of atoms $\{x_0,x_1,\dots,x_{k - 1}\}$,
there are no distinct blocks $S_0,S_1,\dots,S_{k - 1}$ of $\mathcal{S}$ such that $\{x_i,x_{i + 1}\} \subseteq S_i$ with
subscripts modulo $k$. Then $d_i(A) \leq 1$ for any $A \subseteq X$ of size two or more, and $d_i(x,y) \leq 1$ for any distinct atoms $x,y \in X$.
In particular, $\mathcal{S}$ is $\delta$-sparse in a very strong sense: if the upper rank of $\mathcal{S}$ is $r$, then one may take $\delta = 1/(r - 1)$. We describe one step of the process executed on such an $(n,\vb{d})$-antichain $\mathcal{S}$ with atom set $X$. It is convenient to let
\begin{equation*}
\triangle = \max\{d_i \colon 2 \leq i \leq r\} \qquad\text{and}\qquad \Sigma(\vb{d}) = \sum_{i = 2}^r d_i^{\frac{1}{i - 1}}.
\end{equation*}
Throughout this subsection, $r$ is fixed, $D(\vb{d}) \to \infty$, and
$0 < \eta \leq 1/4$ is chosen so that $\eta \to 0$ in this limit.
The notation $\lambda \sim \mu$ means that $\lambda/\mu \to 1$.
We distinguish the expected values in one random round from the
deterministic reference values used in the multiple-round calculation.

\textbf{One step of the process.} We let $q = \eta/(\Sigma(\vb{d}))$ and apply one step of the process to track expected values. We temporarily suppress the parameter $t$,
since we are only discussing one step, and let $\mathcal{S}'$ be the set system after one step, with atom set $X'$. We write $d'_i(A)$ and $d'_i(x,y)$
for the quantities in $\mathcal{S}'$ corresponding to $d_i(A)$ and $d_i(x,y)$ in $\mathcal{S}$.
The expected size of the set $Y$ in step (1) is exactly $qn$.
Conditional on $x \in Y$, the sets $S\setminus\{x\}$, over blocks $S$
containing $x$, are disjoint. Thus, defining $\sigma$ by the product below,
\begin{equation*}
\mathbb P(x \in I) = q \cdot \prod_{i = 2}^r (1 - q^{i - 1})^{d_i} = q\sigma.
\end{equation*}
Therefore $\mathbb E(|I|) = q\sigma n$ is also exact.
An atom that is not contained in $I \cup \partial I$ is said to \emph{survive}. These are precisely the atoms of $X'$.
Acceptance events at different atoms are dependent, so the following
product gives an asymptotic comparison and defines the reference value $\omega$:
\begin{equation*}
\omega := (1 - q\sigma) \cdot \prod_{i = 2}^r
 (1 - (q\sigma)^{i - 1})^{d_i},
\qquad
\mathbb P(x \in X') \sim \omega.
\end{equation*}
Consequently $\mathbb E(|X'|) \sim \omega n$.
For $2 \leq i \leq r$, we consider the degree of $x$ conditional on
$x \in X'$. For $S\in\mathcal S$ of size $j$ containing $x$, the set
$S\setminus I$ has size $i$ and lies in $X'$ when $i-1$ atoms of
$S\setminus\{x\}$ survive and its remaining $j-i$ atoms lie in $I$. The absence of the configurations forbidden by local sparseness gives the
one-step comparison
\begin{equation*}
 \bar d_i' := \sum_{j = i}^r
 \binom{j - 1}{i - 1}
 (q\sigma)^{j - i}
 \omega^{i - 1} d_j,
 \qquad
 \mathbb E(d_i'(x)\mid x \in X') \sim \bar d_i'.
\end{equation*}
Here the comparison applies when $\bar d_i' > 0$; if $\bar d_i' = 0$, the
conditional degree is identically zero.
These calculations suggest an $(\omega n,\bar{\vb{d}}')$-antichain as a
reference for $\mathcal{S}'$, where $\bar{\vb{d}}' = (\bar d_2',\dots,\bar d_r')$.
They do not assert that the random degrees are exactly regular.

\textbf{Multiple steps in the process.} We set $n_0 = n$ and $d_{i,0} = d_i$.
We now define deterministic reference values $n_t$ and
$\vb{d}_t = (d_{2,t},\dots,d_{r,t})$ by iterating the preceding formulas.
At each round, we put $\Sigma_t = \Sigma(\vb{d}_t)$ and $q_t = \eta/\Sigma_t$, and define
\begin{equation*}
\sigma_t := \prod_{i = 2}^r (1 - q_t^{i - 1})^{d_{i,t}} \qquad\text{and}\qquad \omega_t := (1 - q_t\sigma_t) \cdot \prod_{i = 2}^r (1 - (q_t \sigma_t)^{i - 1})^{d_{i,t}}.
\end{equation*}
We set $n_{t + 1} = \omega_t n_t$, and for $2 \leq i \leq r$ we define
\begin{equation}\label{eq:degreechange}
d_{i,t + 1} := \omega_t^{i - 1} \sum_{j = i}^r
 \binom{j - 1}{i - 1}
 (q_t\sigma_t)^{j - i} d_{j,t}.
\end{equation}
This recurrence is exact by definition of the reference degrees.
Its interpretation as an expected-degree update is asymptotic.
Let $\pi_0 = 1$ and, for $t \geq 1$, let
$\pi_t = \prod_{s = 0}^{t - 1}\omega_s$, so $n_t = \pi_t n$ exactly.
Let $T$ be the first value of $t$ such that
$\Sigma_t \leq \sqrt{\log D(\vb{d})}$. The set
\begin{equation*}
J^T = I^0 \cup I^1 \cup \dots \cup I^{T - 1}
\end{equation*}
is independent in $\mathcal{S}^0$. We choose a maximal independent set
$J'$ in the residual system $\mathcal{S}^T$ on $X^T$; then
$I = J^T \cup J'$ is maximal in $\mathcal{S}^0$.
If the random populations and degrees follow the reference values
uniformly up to $T$, the expected-size calculation is
\begin{equation*}
\frac{\mathbb E(|J^T|)}{n} = \frac{1}{n}\sum_{t = 0}^{T - 1}\mathbb E(|I^t|) \sim \sum_{t = 0}^{T - 1}\pi_tq_t\sigma_t = \eta\sum_{t = 0}^{T - 1}\frac{\pi_t\sigma_t}{\Sigma_t}.
\end{equation*}
The later sections justify the required comparisons for the actual
process. Here we calculate the size predicted by the reference sequence.

Let $\iota_0 = 0$ and $\iota_{t + 1} = \iota_t + \pi_tq_t\sigma_t$. We define
\begin{equation*}
f_t(\kappa) = \sum_{i = 2}^r d_{i,t} \kappa^{i - 1} \qquad\text{and}\qquad g(\kappa) = \sum_{i = 2}^r \Bigl(\sum_{j = i}^r \binom{j - 1}{i - 1} d_{j,0} \kappa^{j - i}\Bigr)^{\frac{1}{i - 1}}.
\end{equation*}
By definition, $f_0(\rho(\vb{d})) = \log D(\vb{d})$.
Equation~(\ref{eq:degreechange}) and the binomial theorem give the exact identities
\begin{align*}
 f_{t+1}(\kappa)
 &=\sum_{i=2}^r\sum_{j=i}^r\binom{j-1}{i-1}
       (\omega_t\kappa)^{i-1}(q_t\sigma_t)^{j-i}d_{j,t}\\
 &=\sum_{j=2}^r d_{j,t}\sum_{i=2}^j\binom{j-1}{i-1}
       (\omega_t\kappa)^{i-1}(q_t\sigma_t)^{j-i}\\
 &=f_t(q_t\sigma_t+\omega_t\kappa)-f_t(q_t\sigma_t).
\end{align*}
Induction, using $\pi_{t + 1} = \omega_t\pi_t$ and the definition of
$\iota_{t + 1}$, now gives
\begin{equation}\label{eq:fidentity}
f_t(\kappa) = \sum_{i = 2}^r d_{i,t} \kappa^{i - 1} = f_0(\iota_t + \pi_t \kappa) - f_0(\iota_t).
\end{equation}
For $2 \leq i \leq r$, comparison of coefficients gives
\begin{equation*}
d_{i,t} = \pi_t^{i - 1} \sum_{j = i}^r \binom{j - 1}{i - 1} \iota_t^{j - i} d_{j,0}.
\end{equation*}
Taking $(i - 1)$st roots and summing yields
\begin{equation*}
\Sigma_t = \sum_{i = 2}^r d_{i,t}^{\frac{1}{i - 1}} = \pi_t \sum_{i = 2}^r \Bigl(\sum_{j = i}^r \binom{j - 1}{i - 1} \iota_t^{j - i} d_{j,0}\Bigr)^{\frac{1}{i - 1}} = \pi_t g(\iota_t).
\end{equation*}
Consequently
\begin{equation*}
\iota_{t + 1} - \iota_t = \pi_t q_t \sigma_t = \frac{\eta \pi_t \sigma_t}{\Sigma_t} = \frac{\eta \sigma_t}{g(\iota_t)}.
\end{equation*}
Setting $\kappa = q_t\sigma_t$ in (\ref{eq:fidentity}) gives
\begin{equation*}
\sum_{i = 2}^r d_{i,t} (q_t \sigma_t)^{i - 1} = f_0(\iota_t + \pi_t q_t \sigma_t) - f_0(\iota_t) = f_0(\iota_{t + 1}) - f_0(\iota_t).
\end{equation*}
All the identities in this paragraph are exact for the reference sequence.

To pass from the product defining $\pi_t$ to an exponential, we need
an additive estimate for its logarithm. For $s < T$ and $D(\vb{d})$
sufficiently large, expansion of the product defining $\omega_s$ gives
\begin{equation*}
0 \leq - \log\omega_s - \sum_{i = 2}^r d_{i,s}(q_s\sigma_s)^{i - 1} \leq \frac{3\eta}{\Sigma_s}.
\end{equation*}
The same scale calculation as in the proof of
Lemma~\ref{lem:general-scale}, followed by geometric summation,
gives $\sum_{s < T}\eta/\Sigma_s \to 0$.
Thus, uniformly for $0 \leq t \leq T$,
\begin{equation*}
\log\pi_t + f_0(\iota_t) = \sum_{s = 0}^{t - 1}\left(\log\omega_s+
 \sum_{i = 2}^r d_{i,s}(q_s\sigma_s)^{i - 1}\right)\longrightarrow0.
\end{equation*}
This stronger additive statement, rather than the ratio of logarithms,
justifies exponentiation. We obtain
\begin{equation*}
\pi_t \sim \exp(-f_0(\iota_t))
\qquad\text{and}\qquad
\Sigma_t \sim \exp(-f_0(\iota_t))g(\iota_t).
\end{equation*}
Since $\Sigma_{t + 1} \geq \omega_t\Sigma_t \geq e^{-2\eta}\Sigma_t$,
the definition of $T$ gives
\begin{equation*}
e^{-2\eta}\sqrt{\log D(\vb{d})} \leq \Sigma_T \leq \sqrt{\log D(\vb{d})}.
\end{equation*}
As $\eta \to 0$, it follows that
$\Sigma_T \sim \sqrt{\log D(\vb{d})}$.
For $\kappa \geq 0$ with $f_0(\kappa) \geq 1$, we also have
\begin{equation*}
D(\vb{d}) \leq g(\kappa) \leq (r - 1)2^r(1 + f_0(\kappa))^{r - 2}D(\vb{d}).
\end{equation*}
Combining these bounds with
$f_0(\iota_T) \sim \log g(\iota_T) - \log\Sigma_T$ gives
\begin{equation*}
f_0(\iota_T) \sim \log D(\vb{d}) = f_0(\rho(\vb{d})).
\end{equation*}
The nonnegative coefficients and fixed degree of $f_0$ then imply
$\iota_T \sim \rho(\vb{d})$.
The predicted accumulated size is therefore
\begin{equation*}
\frac{\mathbb E(|J^T|)}{n} \sim \sum_{t = 0}^{T - 1}\pi_tq_t\sigma_t = \sum_{t = 0}^{T - 1}(\iota_{t + 1} - \iota_t) = \iota_T \sim \rho(\vb{d}).
\end{equation*}

It remains to compare the completion with this size.
The exact inequality is $\mathbb E(|J'|) \leq \mathbb E(|X^T|)$,
and the reference prediction is $\mathbb E(|X^T|) \sim n\pi_T$.
The $i = 2$ term of $g(\kappa)$ is $f_0'(\kappa)$; the remaining terms have
sum negligible relative to $f_0'(\kappa)$ as $f_0(\kappa) \to \infty$.
Hence $g(\iota_T) \sim f_0'(\iota_T)$.
Also $\iota_T \sim \rho(\vb{d})$ and the nonnegative coefficients of
$f_0'$ give $f_0'(\iota_T) \sim f_0'(\rho(\vb{d}))$. Thus
\begin{equation*} \pi_T = \frac{\Sigma_T}{g(\iota_T)} \sim \frac{\sqrt{\log D(\vb{d})}}{f_0'(\iota_T)} \sim \frac{\sqrt{\log D(\vb{d})}}{f_0'(\rho(\vb{d}))}.
\end{equation*}
By the definition of $\rho(\vb{d})$,
\begin{equation*}
\rho(\vb{d})f_0'(\rho(\vb{d})) \geq f_0(\rho(\vb{d})) = \log D(\vb{d}).
\end{equation*}
Consequently the expected-size comparison gives
\begin{equation*}
\frac{\mathbb E(|X^T|)}{\rho(\vb{d})n} \sim \frac{\sqrt{\log D(\vb{d})}}
 {\rho(\vb{d})f_0'(\rho(\vb{d}))} \leq \frac{1}{\sqrt{\log D(\vb{d})}}\longrightarrow0.
\end{equation*}
Since $|J^T| \leq |I| \leq |J^T| + |J'|$, the reference calculations above therefore predict
$\mathbb E(|I|) \sim \rho(\vb{d})n$. The subsequent analysis establishes the required asymptotic size of
the output for the actual process.

\subsection{The actual process}

The idealized calculations do not control random deviations: one
round need not produce an exactly regular set system.
Sections~\ref{sec:general-focus} and~\ref{sec:construction-analysis}
bound the probabilities of local failures, and
Lemma~\ref{lem:conditional-lll} makes the required inequalities hold
simultaneously with positive probability. Under $\delta$-sparseness,
blocks through nearby atoms can overlap, and survival probabilities
can differ between atoms. We must therefore control both random
deviations and accumulated errors in the expected values.

To address these two issues, each iteration begins with a
\emph{focusing round}, followed by a \emph{construction round}.
At the start of a round, let $\mathsf{W} \subseteq X$ consist of the atoms
still present that were made ineligible for sampling in earlier rounds.
Initially $\mathsf{W} = \varnothing$. In the focusing round, we sample from
$X\setminus \mathsf{W}$ with probability $q = \Sigma(\vb{d})^{-1 + \delta_2}$ to
obtain a set $Z$, called a \emph{focus}, where $\delta_2 > 0$ is very
small relative to $\delta$. We choose a focus with the local
properties required below and freeze it. In the construction round,
each atom of $Z$ is sampled with probability
$p = \eta\Sigma(\vb{d})^{-\delta_2}$ to obtain $Y$.
The \emph{accepted set} is $I = Y\setminus\partial Y$.
The next full atom set is $X' = X\setminus(Y\cup\partial I)$;
the unused sampled atoms are kept aside for final completion.

The focus controls the effect of local intersections. To handle
the different survival probabilities at different atoms, we also
perform an \emph{equalization}, following Alon, Kim and
Spencer~\cite[Section~2]{AlonKimSpencer}. We write $s_x$ for the
probability that $x \in X'$ and $\omega = \min_{x \in X}s_x$.
For each $x \in X\setminus \mathsf{W}$, an independent retention choice
with probability $\omega/s_x$ makes its probability of surviving
and remaining eligible exactly $\omega$. A surviving atom that
fails this choice stays in the residual system but is added to
the ineligible set. Thus $\mathsf{W}'$ consists of the surviving atoms of
$\mathsf{W}$ together with these new exclusions, and sampling in the next
round is from $X'\setminus \mathsf{W}'$. Previously excluded atoms never
become eligible again. They leave $\mathsf{W}$ only when they leave the
full atom set.

\Needspace{18\baselineskip}
\subsection{Definition of the actual process}\label{subsec:actual-process}

We now give the counterpart of Subsection~\ref{subsec:idealized-process}.
We fix an integer $r \ge 2$, $c \ge 1$, and $0 < \delta_0 \le 1/2$.
Let $\mathcal S^0$ be the given system on $X^0$, with $n = |X^0|$
and initial degree vector $\vb{d} = (d_2,\ldots,d_r)$.
Assume $D(\vb{d}) \ge e$.
Initially, we set
\[
 J^0 = \mathsf{W}^0 = \varnothing,\qquad \pi_0 = 1,\qquad \iota_0 = 0,
\]
and, for $2 \le i \le r$, set $d_{i,0} = d_i$.
We write $\vb{d}_t = (d_{2,t},\ldots,d_{r,t})$, and put
\[
 D = D(\vb{d}_0),\qquad
 \triangle = \max_{2 \le i \le r}d_{i,0}.
\]
The initial $D$ stays fixed throughout the run.

We fix the error exponents and sampling rate by
\[
 \delta_1=\frac{\delta_0}{16r},\qquad
 \delta_2=\frac{\delta_1}{64r},\qquad
 \delta_3=\frac{\delta_2}{64},\qquad
 \eta=(\log D)^{-\delta_3/8}.
\]
For $2\le i\le r$ and $3\le i\le r$, respectively, let
\[
 \delta_4(i)=(i-1)\delta_1,\qquad \delta_5(i)=(r+i-2)\delta_1.
\]
The final relative error is
\[
 \epsilon=\frac{\delta_3}{32}=\frac{\delta_0}{2^{21}r^2},\qquad
 \varepsilon(\triangle)=(\log\triangle)^{-\epsilon}.
\]

\textbf{State and one-round notation.}
At the beginning of round $t$, the state is
$(J^t,X^t,\mathsf{W}^t,\mathcal S^t)$, with nonnegative reference vector
$\vb{d}_t$. The sets $J^t,X^t$ are disjoint: $J^t$ contains the
points already accepted, while $X^t$ contains all points still
present in the residual antichain $\mathcal S^t$.
The set $\mathsf{W}^t \subseteq X^t$ consists of points made permanently
ineligible in earlier rounds. Sampling uses $X^t\setminus \mathsf{W}^t$.
We put
\[
 \Sigma_t = \Sigma(\vb{d}_t).
\]
The focus and construction probabilities are
$q_t = \Sigma_t^{-1 + \delta_2}$ and
$p_t = \eta\Sigma_t^{-\delta_2}$.
The reference degrees are distinct from the actual degrees in
$\mathcal S^t$. The maximum scale in Section~1 gives
\[
 D(\vb{d}_t) \le \Sigma_t \le (r - 1)D(\vb{d}_t);
\]
in particular, $D \le \Sigma_0 \le (r - 1)D$.

In one-round formulas we suppress $t$; thus $X=X^t$, $\mathsf{W}=\mathsf{W}^t$,
$\mathcal S=\mathcal S^t$, $\Sigma=\Sigma_t$, $q=q_t$ and $p=p_t$.

We stop before round $t$ if $\Sigma_t \le \sqrt{\log D}$; thus
\[
 T = \min\{t \in \mathbb Z_{ \ge 0}\colon \Sigma_t \le \sqrt{\log D}\},\qquad
 \min\varnothing = \infty.
\]
The analysis proves that $T$ is finite. For $0 \le t < T$
with $X^t \ne \varnothing$, we perform the following steps.
The proposed next state is adopted only in step~\ref{step:process-good}.

\begin{enumerate}
\item \textbf{The focus.}
Let $\mathcal G_{{\rm foc},t}\subseteq 2^{X^t\setminus \mathsf{W}^t}$ be a
prescribed nonempty family. We form the \emph{focus} $Z^t$ by including
each atom independently with probability $q_t$, and repeat until
$Z^t\in\mathcal G_{{\rm foc},t}$. Such a focus is \emph{good}.
We fix $Z=Z^t$ for the construction round.

\item \textbf{Sampling, acceptance, and full survival.}
\label{step:process-sample}
We take independent Bernoulli variables $(\chi_x)_{x \in Z}$ of
parameter $p = p_t = \eta\Sigma_t^{-\delta_2}$, and set $\chi_x = 0$ outside $Z$.
The \emph{sample} is
\[
 Y^t = Y = \{x \in Z\colon \chi_x = 1\}.
\]
For $B \subseteq X$, we set
\[
 \partial B = \{v \in X\colon \exists S \in \mathcal S,\ v \in S,\
 S\setminus\{v\} \subseteq B\}.
\]
We set
\[
 I^t = I = Y\setminus\partial Y,\qquad
 J^{t + 1} = J^t\cup I^t,\qquad
 X^{t + 1} = X' = X\setminus(Y\cup\partial I).
\]
The set $I$ is the \emph{accepted set}. This is a single deletion
step, and $I$ is independent: every atom of a block contained in
$Y$ lies in $\partial Y$. For $x \in X$, we define the acceptance
and survival indicators by
\[
 J_x = \vb1_{\{x \in I\}},\qquad
 \psi_x = \vb1_{\{x \notin Y\cup\partial I\}}.
\]
An atom \emph{survives} if it lies in
$X' = \{x \in X\colon \psi_x = 1\}$.
Every point of the sample leaves the residual system, including
the points of $Y^t\setminus I^t$. We keep those unused sampled points
from an adopted round for final completion. Points outside $Z^t$
also leave the system if they belong to $\partial I^t$.

\item \textbf{Eligibility for the next round.}
\label{step:process-equalization}
With the current state and $Z^t$ fixed, we use the independent
sampling law in step~\ref{step:process-sample} to define, for $x\in X^t$,
\begin{equation}\label{eq:process-survival-probabilities}
 s_{x,t}=s_x=\mathbb P(\psi_x=1),\qquad
 \omega_t=\omega=\min_{x\in X^t}s_{x,t}.
\end{equation}
These values are fixed before the sample is realized, and the
possibility $Y^t=\varnothing$ gives $s_{x,t}>0$.
For $x\in X^t\setminus\mathsf{W}^t$, let
$\mathsf{w}_{x,t}=\mathsf{w}_x$ be Bernoulli random variables,
mutually independent and independent of $(\chi_y:y\in Z^t)$, with
\begin{equation}\label{eq:process-w-law}
 \mathbb P(\mathsf{w}_{x,t}=1)=1-\frac{\omega_t}{s_{x,t}}.
\end{equation}
For $x\in\mathsf{W}^t$, we set $\mathsf{w}_{x,t}=\mathsf{w}_x=1$.
The value $\mathsf{w}_{x,t}=1$ records the decision to exclude $x$
from future sampling. If $x$ survives round $t$, it belongs to $\mathsf{W}^{t+1}$
exactly when $\mathsf{w}_{x,t}=1$. Thus
\begin{equation}\label{eq:process-eligible-atoms}
 \begin{aligned}
 \mathsf{W}^{t+1}
 &=\mathsf{W}'=\{x\in X^{t+1}:\mathsf{w}_{x,t}=1\},\\
 X^{t+1}\setminus\mathsf{W}^{t+1}
 &=\{x\in X^t:\psi_x(1-\mathsf{w}_x)=1\}.
 \end{aligned}
\end{equation}
For $x \in X^t\setminus \mathsf{W}^t$, these independent choices give
\[
 \mathbb P(x \in X^{t + 1}\setminus \mathsf{W}^{t + 1}) = s_{x,t}\frac{\omega_t}{s_{x,t}} = \omega_t.
\]
The set $\mathsf{W}^{t + 1}$ carries forward all earlier exclusions still
present and adds the new exclusions from this round. In particular,
\[
 \mathsf{W}^{t + 1}\cap \mathsf{W}^t = \mathsf{W}^t\cap X^{t + 1},\qquad X^{t + 1}\setminus \mathsf{W}^{t + 1} \subseteq X^t\setminus \mathsf{W}^t.
\]
A previously ineligible point never becomes eligible again.
The set $\mathsf{W}^{t + 1}$ need not contain $\mathsf{W}^t$, because atoms of
$\mathsf{W}^t$ may not survive. Retention affects eligibility only: the points of
$\mathsf{W}^{t + 1}$ remain in $X^{t + 1}$ with their incident residual blocks.
They may leave the system in a later round or be considered during final completion.

\item \textbf{The residual system.}\label{step:process-residual}
Let $\mathcal S^{t + 1}$ consist of the distinct inclusion-minimal
members of
\[
 \{S\setminus I^t\colon S \in \mathcal S^t,\
 S\setminus I^t \subseteq X^{t + 1},\ |S\setminus I^t| \ge 2\}.
\]
The residual system has atom set $X^{t + 1}$; its blocks may
contain points of $\mathsf{W}^{t + 1}$. The process is illustrated in Figure~\ref{fig:theprocess}.

\begin{figure}[H]
\centering
\includegraphics[width=6.4in]{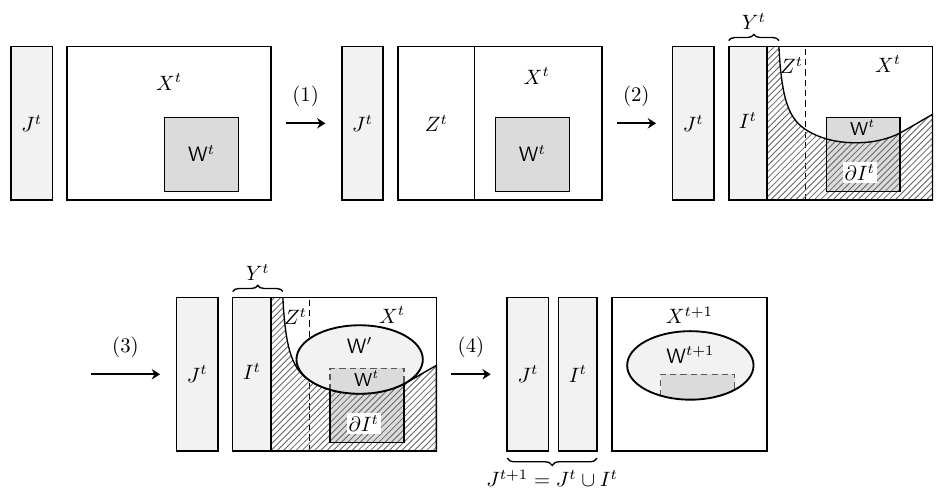}
\caption{The actual process.}
\label{fig:theprocess}
\end{figure}

\item \textbf{A good construction round.}\label{step:process-good}
Let $\mathcal G_{{\rm con},t}(Z)\subseteq
\{0,1\}^{Z}\times\{0,1\}^{X\setminus \mathsf{W}}$ be a prescribed nonempty
set of assignments to $(\chi_x:x\in Z)$ and $(\mathsf{w}_x:x\in X\setminus \mathsf{W})$,
each having positive probability under their independent law.
A round is \emph{good} when its assignment belongs to this set.
We condition the sampling and retention choices jointly on this event.
Equivalently, we repeat steps~\ref{step:process-sample}--\ref{step:process-residual},
with the same state and focus, until this event occurs, and then
adopt the proposed next state.
Unsuccessful trials leave the current state, including $\mathsf{W}^t$, fixed;
only their proposed next state is discarded. The values
$s_{x,t}$ and $\omega_t$ stay fixed throughout these trials.
The exact eligibility marginal in step~\ref{step:process-equalization}
refers to the law before this conditioning.

\item \textbf{The reference parameters.}\label{step:process-reference}
We put
\[
 \sigma_t = \prod_{i = 2}^r
 (1 - p_t^{i - 1})^{q_t^{i - 1}d_{i,t}},\qquad
 \vartheta_t = q_t p_t\sigma_t = \frac{\eta\sigma_t}{\Sigma_t}.
\]
For $2 \le i \le r$, we use the reference recurrence
\[
 d_{i,t + 1} = \omega_t^{i - 1}
 \sum_{j = i}^r\binom{j - 1}{i - 1}\vartheta_t^{j - i}d_{j,t}.
\]
We also set
\[
 \pi_{t + 1} = \omega_t\pi_t,\qquad
 \iota_{t + 1} = \iota_t + \pi_t\vartheta_t.
\]
These are reference quantities; the actual degrees are those of
$\mathcal S^{t + 1}$.

\item \textbf{The next round.}
We increase $t$ by one, check the stopping condition, and choose a
fresh focus from $X^t\setminus \mathsf{W}^t$. There is one successful
construction round for each focus. The preceding focus has no
further role: its unselected points can be sampled in later rounds
if they survive and remain eligible.
\end{enumerate}

\subsection{Process terminology}\label{subsec:process-terminology}

{\bf Roots, links and degrees.}
Let $2\le i\le r$, and let $\mathcal S$ and $\mathcal S'$ be the
set systems before and after a round, with atom sets $X$ and $X'$,
respectively. We write $\mathcal S_i=\{S\in\mathcal S:|S|=i\}$.
A \emph{root} is a specified atom $x\in X$.
An \emph{$i$-link} at $x$ is a set $S\setminus\{x\}$ with
$x\in S\in\mathcal S_i$. For distinct roots $x,y$, a
\emph{common $i$-link} is an $(i-1)$-set $F$ such that
$F\cup\{x\},F\cup\{y\}\in\mathcal S_i$. Their families are
\[
 N_i(x)=\{S\setminus\{x\}:x\in S\in\mathcal S_i\},\qquad
 N_i(x,y)=N_i(x)\cap N_i(y).
\]
The \emph{$i$-degree} of $x$ is $d_i(x)=|N_i(x)|$.
For $A\subseteq X$, $d_i(A)$ counts the blocks of size $i$ containing
$A$, and $d_i(x,y)$ counts the common $i$-links of $x$ and $y$.

{\bf Acceptance, survival, and retention.}
An atom $x$ is \emph{accepted} if $x\in I=Y\setminus\partial Y$,
and it \emph{survives} if $x\in X'=X\setminus(Y\cup\partial I)$.
A \emph{surviving root} is a root in $X'$.
For $x\in X'$, \emph{retention} means that $x$ remains eligible
for sampling: $x\notin\mathsf W'$, equivalently $\mathsf{w}_x=0$.

{\bf Contractions and repetitions.}
A block $S\in\mathcal S$ \emph{contracts to} $B\in\mathcal S'$
when $B=S\setminus I$; $S$ is then a \emph{parent} of $B$.
The phrase \emph{created by contraction} additionally requires
$S\cap I\ne\varnothing$.
A \emph{repetition} occurs when distinct blocks $S,P\in\mathcal S$
give the same set $S\setminus I=P\setminus I$.
Step~\ref{step:process-residual} keeps each distinct remaining set
once and discards any such set that properly contains another
remaining set.

{\bf Old and new blocks.}
An \emph{old block} is a member of $\mathcal S$; a \emph{new block}
is a member of $\mathcal S'\setminus\mathcal S$.

\subsection{Parameters and conditioning}
\label{subsec:process-focus-parameters}

For $2\le i\le r$, $A,U\subseteq X$ and distinct $x,y\in X$, we define
the restricted families and their cardinalities by
\begin{equation}\label{eq:restricted-neighborhoods}
 \begin{aligned}
 N_i(A,U)&=\{S\in N_i(A):S\subseteq U\},\\
 d_i(A,U)&=|N_i(A,U)|.
 \end{aligned}
\end{equation}
\begin{equation}\label{eq:restricted-common-links}
 \begin{aligned}
 N_i(x,y,U)&=\{S\in N_i(x,y):S\subseteq U\},\\
 d_i(x,y,U)&=|N_i(x,y,U)|.
 \end{aligned}
\end{equation}
For $i=2$, singleton completions are identified with their atoms.

Primes denote the proposed next state. The current eligible-degree
and full-degree errors are
\[
 E_i=\max_{v\in X}|d_i(v,X\setminus \mathsf{W})-d_i|,\qquad
 \widehat E_i=\max_{v\in X}(d_i(v,X)-d_i)_+.
\]
We use the same definitions for $E_i',\widehat E_i'$ in the proposed
state, relative to the fixed-focus reference values
$\bar d_i'=d_{i,t+1}$ of step~\ref{step:process-reference}, and put
$\Sigma'=\Sigma(\bar{\vb d}')$. Empty maxima are zero. We also put
\[
 \mathsf E_{\rm acc}=\sum_{i=2}^r
       \eta^{i-1}\frac{E_i}{\Sigma^{i-1}},\qquad
 \mathcal B_k=\sum_{j=k+1}^r\binom{j-1}{k-1}
       \vartheta^{j-k}(d_j+E_j),
\]
where $\vartheta=qp\sigma$ is as in step~\ref{step:process-reference}.
The scale $\Sigma$ and the probabilities $q,p,s_x,\omega$ always
have the meanings assigned in the seven-step process.

All constants are fixed for $r,\delta_0,c$ before a run.

If $X^t=\varnothing$ before the stopping time, we keep all subsequent
point sets empty except for the fixed accepted set $J^t$.
We continue the reference updates in step~\ref{step:process-reference}, taking
\[
 \omega_t=\exp\left(-\sum_{i=2}^r d_{i,t}\vartheta_t^{i-1}\right).
\]
After $T$, we freeze the state and reference parameters and take
later samples and increments of $\iota_t$ to be zero.

\textbf{Selection of a run.}
We put $\varepsilon=\varepsilon(\triangle)$.
We repeat entire runs independently until
\[
 |J^T|\ge(1-\varepsilon)n\rho(\vb d_0),\qquad
 |X^T|+\sum_{t=0}^{T-1}|Y^t|
 \le(1+\varepsilon)n\rho(\vb d_0).
\]
We use the first such run for the following completion.

\textbf{Maximal completion.}
Starting with $J^T$, we consider the points of
\[
 X^T\ \cup\ \bigcup_{t=0}^{T-1}(Y^t\setminus I^t)
\]
in any fixed order, adding a point whenever independence in the
original system $\mathcal S^0$ is preserved. Every other unused
point $x$ satisfies $S\setminus\{x\}\subseteq J^T$ for some
$S\in\mathcal S^0$ containing $x$, so the resulting set is maximal and
\[
 |J^T|\le|I|\le|X^T|+\sum_{t=0}^{T-1}|Y^t|.
\]
The points of $\mathsf{W}^T\subseteq X^T$ are included in this completion.

\Needspace{10\baselineskip}
\subsection{Performance of the process}\label{subsec:process-performance}

\Needspace{8\baselineskip}
The next theorem states the performance of the process in
Subsection~\ref{subsec:actual-process} and is the large-degree
quantitative form of Theorem~\ref{thm:main}.

\nopagebreak[4]
\begin{theorem}\label{prop:quantitative-main}
Let $r \ge 2$ be an integer, and let $c \ge 1$ and $0 < \delta_0 \le 1/2$.
There is a threshold $C = C(r,\delta_0,c) \ge e$ with the following property.
Let $\mathcal S$ be an $(n,\vb{d})$-antichain on $X$, and let
$D = D(\vb{d})$ and $\triangle = \max_{2 \le i \le r}d_i$.
Suppose $D \ge C$ and, for $B \subseteq X$ and $2 \le |B| < i \le r$,
\begin{equation}\label{eq:initialP}
 d_i(B) \le c D^{i - |B| - \delta_0},
\end{equation}
and that, for distinct $x,y \in X$ and $2 \le i \le r$,
\begin{equation}\label{eq:initialCommon}
 d_i(x,y) \le c D^{i - 1 - \delta_0}.
\end{equation}
There exist conditioning events for which the process in
Subsection~\ref{subsec:actual-process}, including the selection of a run
and maximal completion, terminates almost surely and produces a maximal
independent set $I$ such that
\begin{equation*}
 \bigl||I| - n\rho(\vb{d})\bigr|
 \le \varepsilon(\triangle)n\rho(\vb{d}).
\end{equation*}
\end{theorem}

The parameter $c \ge 1$ multiplies both initial sparsity bounds
\eqref{eq:initialP}--\eqref{eq:initialCommon} and is independent
of $n$ and the degree vector. Allowing this multiplier lets the
polynomial applications use the same estimates; it affects the
threshold $C(r,\delta_0,c)$.
For the theorem in Section~1, we take $c = 1$ and
$\delta_0 = \min\{\delta,1/2\}$.
The proof is completed in Subsection~\ref{subsec:quantitative-assembly}.

\begin{proof}[Deduction of Theorem~\ref{thm:main}]
We fix $\delta > 0$ and $r \geq 2$. We put $s = r - 1$,
$\delta_0 = \min\{\delta,1/2\}$, and take $c = 1$ in
Theorem~\ref{prop:quantitative-main}. We define
\[
 \epsilon = \frac{\delta_0}{2^{21}r^2},\qquad
 D_0 = C(r,\delta_0,1),\qquad
 K = \frac{s^2D_0(\log D_0)^\epsilon}{\log 2}.
\]
These constants depend only on $\delta,r$, and $K \geq 1$ since
$D_0 \geq e$. We write $D = D(\vb{d})$ and
$\triangle = \max_{2 \leq i \leq r}d_i$.
If $D \geq D_0$, Theorem~\ref{prop:quantitative-main} supplies the required maximal
independent set, because $\log\triangle \geq \log D$ and hence
\[
 \varepsilon(\triangle) = (\log\triangle)^{-\epsilon} \leq (\log D)^{-\epsilon} \leq K(\log D)^{-\epsilon}.
\]

Suppose instead that $1 < D < D_0$. Some integer $d_i$ is at least two,
so $\log D \geq (\log 2)/s$. We put
\[
 \kappa = \frac{1}{D}\min\left\{1,\frac{\log D}{s}\right\}.
\]
Since $d_i \leq D^{i - 1}$ and $0 < D\kappa \leq 1$, we have
\[
 \sum_{i = 2}^r d_i \kappa^{i - 1} \leq \sum_{j = 1}^s(D\kappa)^j \leq sD\kappa \leq \log D.
\]
Since $\sum_{i=2}^r d_i\rho^{i-1}=\log D$, we obtain
\[
 \rho(\vb{d}) \geq \kappa \geq \frac{\log 2}{s^2D_0}.
\]
We choose any maximal independent set $I$. Since
$0 < \rho(\vb{d}) < 1$ and $0 \leq |I| \leq n$,
\[
 \bigl||I| - \rho(\vb{d})n\bigr| \leq n \leq K(\log D)^{-\epsilon}\rho(\vb{d})n.
\]
This proves the theorem for every $D > 1$.
\end{proof}

\begin{proof}[Deduction of Corollary~\ref{cor:main}]
We use $s,\epsilon,D_0,K$ from the proof of Theorem~\ref{thm:main}, and put
\[
 \iota = \min\left\{\left(\frac{\log d_{k + 1}}{k d_{k + 1}}\right)^{\frac{1}{k}}
 \colon 1 \leq k \leq s,\ d_{k + 1} \geq 3\right\}.
\]
We first verify $\rho(\vb{d}) \leq \iota$ when $D \geq e$.
Let $f(\kappa) = \sum_{i = 2}^r d_i \kappa^{i - 1}$, and let $\ell$ satisfy
$d_{\ell + 1} = D^\ell$. For any $k$ with $d_{k + 1} \geq 3$, we put
$\vartheta = d_{k + 1}^{1/k}$ and $\kappa = ((\log \vartheta)^{1/k})/\vartheta$.
If $\vartheta = D$, then $f(\kappa) \geq \log D$.
If $\vartheta \leq e$, then $3 \leq d_{k + 1} \leq e^k$; the function
$(\log \lambda)/\lambda$ decreases for $\lambda \geq 3$, so $\kappa \geq e^{-1}$ and
\[
 f(\kappa) \geq (D\kappa)^\ell \geq \left(\frac{D}{e}\right)^\ell \geq \frac{D}{e} \geq \log D.
\]
If $e < \vartheta < D$, then $k \ne \ell$ and
\[
 f(\kappa) \geq \log \vartheta + (D\kappa)^\ell \geq \log \vartheta + \frac{D}{\vartheta} \geq \log D.
\]
In all cases $\rho(\vb{d}) \leq \kappa$, proving the comparison.
For $D \geq D_0$, Theorem~\ref{thm:main} now gives the corollary.

For $1 < D < D_0$ and each $k$ occurring in the minimum defining $\iota$,
we have $d_{k + 1}^{1/k} < D_0$, $\log d_{k + 1} > 1$, and $k \leq s$.
Consequently,
\[
 \iota \geq \frac{1}{sD_0} \geq \frac{\log 2}{s^2D_0}.
\]
Any maximal independent set therefore satisfies
\[
 |I| \leq n \leq K(\log D)^{-\epsilon}\iota n \leq \bigl(1 + K(\log D)^{-\epsilon}\bigr)\iota n.
\]
This proves the corollary throughout its stated range, with the
same constants as Theorem~\ref{thm:main}.
\end{proof}

We use $b_1,\ldots,b_6$ for the constants in Section~4 and
$c_{1},c_{2},c_{3}$ for those in
Sections~\ref{sec:construction-analysis} and~\ref{sec:completing-process}.
In the deduction of Theorem~\ref{thm:main} from Theorem~\ref{prop:quantitative-main}, we take $c = 1$, so $K$ depends only on $\delta$ and $r$. An explicit choice of $K$ that works is
$\log\log\log\log K = \left((4r)/(\min\{1/2,\delta\})\right)^{16}$.

\section{Probabilistic tools}\label{sec:tools}

The following proposition is given by Alon and Spencer~\cite[Appendix~A]{AS}.
\begin{proposition}\label{prop:bernstein}
Let $\iota,\tau>0$, $h\ge2$ and $a\ge0$.
Let $I$ be finite and let $0\le\beta_x\le\iota$ for every $x\in I$.
Let $(U_x:x\in I)$ be independent Bernoulli random variables.
Let $\Phi=\sum_{x\in I}\beta_xU_x$ and $\mu=\mathbb E(\Phi)$. Then
\[
 \mathbb P(|\Phi-\mu|\ge\tau)
 \le2\exp\left(-\frac{\tau^2}{2\iota(\mu+\tau)}\right)
\]
and
\[
 \mathbb P\!\left(|\Phi-\mu|>
 4(a+1)\bigl(\sqrt{\iota\mu\log h}+\iota\log h\bigr)\right)
 \le h^{-a}.
\]
\end{proposition}

\medskip
\begin{definition}\label{def:coordinate-supports}
Let $I$ and $X$ be finite sets, and let $(Y_i:i\in I)$ be independent
Bernoulli random variables. For each $x\in X$, let
$T_x\subseteq I$ and let $f_x:\{0,1\}^{T_x}\to\{0,1\}$ be a function, and let
\[
 U_x=f_x((Y_i)_{i\in T_x}).
\]
The set $T_x$ is a \emph{coordinate support} for $U_x$.
For $k\ge1$, the family $(U_x:x\in X)$ with these sets has
\emph{read bound} $k$ if, for every $i\in I$,
\[
 |\{x\in X:i\in T_x\}|\le k.
\]
\end{definition}

\emph{Remark.} The read bound limits how many $U_x$ can use the same input
$Y_i$. Proposition~\ref{prop:readk} uses it to bound deviations of
$\sum_{x\in X}U_x$, even when the $U_x$ are dependent.

The following is a consequence of Gavinsky, Lovett, Saks and
Srinivasan~\cite[Theorem~1.1]{GLSS}.
\begin{proposition}\label{prop:readk}
Let $k\ge1$ and $\tau>0$. Let $(U_x:x\in X)$ be a family of
indicators with read bound $k$ as in Definition~\ref{def:coordinate-supports},
and let $\mu=\mathbb E(\sum_{x\in X}U_x)$. Then
\[
 \mathbb P\!\left(\left|\sum_{x\in X}U_x-\mu\right|\ge\tau\right)
 \le2\exp\left(-\frac{\tau^2}{2k(\mu+\tau)}\right).
\]
\end{proposition}

\subsection{Polynomial moments and concentration}

\begin{definition}\label{def:polynomial}
Let $(U_x:x\in I)$ be independent Bernoulli random variables, where
$I$ is finite, and let $s\ge0$ be an integer. For each $F\subseteq I$,
let $z_F\ge0$. For $T\subseteq I$, we define
\[
 \Phi=\sum_{F\subseteq I:\,|F|\le s}
 z_F\prod_{x\in F}U_x,\qquad \partial_T\Phi=\sum_{T\subseteq F\subseteq I:\,|F|\le s}
 z_F\prod_{x\in F\setminus T}U_x.
\]
Let $E_0=\mathbb E(\Phi)$, and for $j\ge1$ let
$E_j=\max_{T\subseteq I:\,|T|=j}\mathbb E(\partial_T\Phi)$,
where an empty maximum is zero.
\end{definition}

\emph{Remark.} The variable $\partial_T\Phi$ sums the terms containing every
index in $T$, with those Bernoulli factors removed. The quantities $E_j$
control the moment and deviation bounds in
Propositions~\ref{prop:positive-moments} and~\ref{prop:polynomial}.

\begin{proposition}\label{prop:positive-moments}
Let $s\ge0$ and $m\ge1$ be integers, and let $a,\iota\ge0$ and $h\ge2$.
Let $\Phi$ and $E_j$ be as in Definition~\ref{def:polynomial}. Then
\[
 \mathbb E(\Phi^m)\le
 \left(\sum_{j=0}^s\binom{sm}{j}E_j\right)^m.
\]
If $s\ge1$, then
\[
 \mathbb P\!\left(\Phi>
 e(s+1)[2s(a+4)]^s\sum_{j=0}^s(\log h)^jE_j\right)
 \le h^{-a}.
\]
If $s\ge1$ and $E_j\le\iota$ for $0\le j\le s$, then
\[
 \mathbb P\!\left(\Phi>
 e(s+1)[2s(a+4)]^s\iota(\log h)^s\right)
 \le h^{-a}.
\]
\end{proposition}
\begin{proof}
We expand $\Phi^m$ as a sum over ordered tuples $(F_1,\ldots,F_m)$,
where $F_\ell\subseteq I$ and $|F_\ell|\le s$ for $1\le\ell\le m$.
We fix $F_1,\ldots,F_{m-1}$ and put $R=F_1\cup\cdots\cup F_{m-1}$,
so $|R|\le s(m-1)$. By independence, summing over the last set
$F=F_m$ multiplies the expected contribution of the first $m-1$ factors by
\[
 \sum_{\substack{F\subseteq I\\|F|\le s}}
 z_F\prod_{x\in F\setminus R}\mathbb P(U_x=1)
 \le\sum_{\substack{T\subseteq R\\|T|\le s}}
 \mathbb E(\partial_T\Phi)
 \le\sum_{j=0}^s\binom{sm}{j}E_j.
\]
The first inequality groups the sets $F$ according to $T=F\cap R$.
Iterating this bound proves the moment estimate.
For $s\ge1$, we put $m=\lceil(a+1)\log h\rceil$. Markov's inequality gives
\[
 \mathbb P\!\left(\Phi>e(\mathbb E(\Phi^m))^{1/m}\right)
 \le e^{-m}\le h^{-a}.
\]
Since $m\le(a+3)\log h$, $\binom{sm}{j}\le(sm)^j$, and
$\sum_{j=0}^s(sm)^j\le(s+1)(sm)^s$, the moment estimate gives
both probability bounds with coefficient $e(s+1)[2s(a+4)]^s$.
The case $s=0$ is deterministic.
\end{proof}

\begin{corollary}\label{cor:polynomial-envelope}
Let $s\ge1$ be an integer, $a,\iota\ge0$ and $h\ge2$.
Let $\Phi$ and $E_j$ be as in Definition~\ref{def:polynomial},
and let $\mu=\mathbb E(\Phi)$.
\begin{enumerate}
\item If $E_j\le\iota$ for $1\le j\le s$, then
\[
 \mathbb P\!\left(\Phi>
 e(s+1)[2s(a+4)]^s\bigl(\mu+\iota(\log h)^s\bigr)\right)
 \le h^{-a}.
\]
\item Let $\alpha>0$. If $h\ge\max\{e,[4s(a+3)/(e\alpha)]^{2/\alpha}\}$
and $E_j\le\iota h^{(s-j)\alpha}$ for $0\le j\le s$, then
\[
 \mathbb P\!\left(\Phi>2e\iota h^{s\alpha}\right)\le h^{-a}.
\]
\item Let $\alpha>0$ and let $m\ge1$ be an integer.
Suppose $\Phi$ is homogeneous of degree $s$, its coefficients
$z_F$ satisfy $z_F\le\iota$ for every $F\subseteq I$ with $|F|=s$,
and $E_j\le\iota h^{-\alpha}$ for $0\le j<s$. Then
\[
 \mathbb P\bigl(\Phi>\iota(sm)^s\bigr)\le h^{-\alpha m}.
\]
\end{enumerate}
\end{corollary}
\begin{proof}
We use Proposition~\ref{prop:positive-moments} with
$m=\lceil(a+1)\log h\rceil$. For (1), we separate $E_0$ in the moment
bound and use $(\log h)^j\le2^s(\log h)^s$ for $1\le j\le s$.
For (2), $\log h\le2h^{\alpha/2}/(e\alpha)$ and the stated lower bound
on $h$ give $2s(a+3)\log h\le h^\alpha$, so $smh^{-\alpha}\le1/2$.
Thus $(\mathbb E(\Phi^m))^{1/m}$ is at most
$\iota h^{s\alpha }\sum_{j=0}^s(sm h^{-\alpha })^j\le2\iota h^{s\alpha }$.
The inequality $\mathbb P(\Phi>e(\mathbb E(\Phi^m))^{1/m})\le e^{-m}$ completes (1) and (2).

For (3), the case $\iota=0$ is deterministic. Dividing $\Phi$ and
all its coefficients $z_F$ by $\iota$, we may assume $\iota=1$.
For the integer $m$ in (3), we define
\[
 W_m=
 \sum_{\substack{F_1,\ldots,F_m\in\binom Is\\
 F_j\nsubseteq\bigcup_{\ell=1}^{j-1}F_\ell\ (1\le j\le m)}}
 \prod_{j=1}^m\left(z_{F_j}\prod_{x\in F_j}U_x\right).
\]
We put $n=s(m-1)$, $N=\sum_{j=0}^{s-1}\binom nj$ and $C=\binom ns$.
For a fixed initial tuple $F_1,\ldots,F_{j-1}$, we put
$R=F_1\cup\cdots\cup F_{j-1}$. Summing over the next set $F=F_j$
multiplies the expected contribution of the preceding factors by
\[
 \sum_{\substack{F\in\binom Is\\F\nsubseteq R}}
 z_F\prod_{x\in F\setminus R}\mathbb P(U_x=1)
 \le\sum_{\substack{T\subseteq R\\|T|<s}}
 \mathbb E(\partial_T\Phi)\le Nh^{-\alpha}.
\]
Here the first inequality groups the sets $F$ by $T=F\cap R$.
Iterating gives $\mathbb E(W_m)\le(Nh^{-\alpha})^m$.
On $\Phi>C+N$, the terms indexed by $F\subseteq R$ have total weight
at most $C$, since $z_F\le1$ and $|R|\le n$.
The remaining terms therefore have total weight greater than $N$
at every choice, so $W_m>N^m$. Hence
\[
 \mathbb P(\Phi>C+N)
 \le\mathbb P(W_m>N^m)
 \le\frac{\mathbb E(W_m)}{N^m}\le h^{-\alpha m}.
\]
Since $C+N\le(n+1)^s\le(sm)^s$, this proves (3).
\end{proof}

The next proposition follows from Propositions~\ref{prop:positive-moments} and~\ref{prop:bernstein}.
\begin{proposition}\label{prop:polynomial}
Let $s\ge1$ be an integer, $a,K,\iota\ge0$ and $h\ge2$.
Let
\[
 \begin{aligned}
 \lambda_{1}&=a+K+4s+8,\\
 \lambda_{2}&=1+(s+1)(K+1)s^s/s!,\\
 \lambda_{3}&=1+e(s+1)[2s(\lambda_1+4)]^s.
 \end{aligned}
\]
Suppose $h\ge(1+s\lambda_2)(64\lambda_1\lambda_3)^s$.
Let $\Phi$ and $E_j$ be as in Definition~\ref{def:polynomial},
with $|I|\le h^K$ and $E_j\le\iota$ for $1\le j\le s$.
Let $\mu=\mathbb E(\Phi)$. Then
\[
 \mathbb P\Bigl(|\Phi-\mu|>(1+s\lambda_2)(64\lambda_1\lambda_3)^s \cdot\bigl(\sqrt{\iota\mu}(\log h)^{(s+1)/2} +\iota(\log h)^{s+1}\bigr)\Bigr)\le h^{-a}.
\]
\end{proposition}
\begin{proof}
The case $\iota=0$ is deterministic. We remove unused indices from $I$,
so $|I|\le h^K$. Since
$h\ge(1+s\lambda_2)(64\lambda_1\lambda_3)^s
\ge\max\{e,2s\lambda_2\}$, we put
$L=\log h$ and $z=64\lambda_{1}\lambda_{3}$.
First suppose $I$ is partitioned into $j$ classes $I_1,\ldots,I_j$,
where $1\le j\le s$, and each monomial uses one coordinate from each class.
Induction on $j$ gives
\begin{equation}\label{eq:partite-polynomial}
 \mathbb P\!\left(|\Phi-\mu|>
 z^j\bigl(\sqrt{\iota\mu}L^{j/2}+\iota L^j\bigr)\right)
 \le h^{-a-4-4(s-j)}.
\end{equation}
For $j=1$, we use Proposition~\ref{prop:bernstein}. For $j>1$, we write
$\Phi=\sum_{x\in I_j}\beta_xU_x$. Each $\beta_x$ depends only on the first $j-1$
classes and has all derivative expectations at most $\iota$.
Proposition~\ref{prop:positive-moments}, with exponent
$a+K+8+4(s-j)\le\lambda_1$, and a union bound give
\[
 \mathbb P\!\left(\max_x \beta_x>\lambda_3\iota L^{j-1}\right)
 \le h^{-a-8-4(s-j)}.
\]
The conditional mean $m=\sum_{x\in I_j}\mathbb P(U_x=1)\beta_x$ is a polynomial
in the first $j-1$ classes. Applying \eqref{eq:partite-polynomial} to it gives
\[
 \mathbb P\!\left(|m-\mu|>
 z^{j-1}\bigl(\sqrt{\iota\mu}L^{(j-1)/2}+\iota L^{j-1}\bigr)\right)
 \le h^{-a-8-4(s-j)}.
\]
The complementary event implies
$m\le2\mu+(1+z^{j-1})^2\iota L^{j-1}$.
We fix the first $j-1$ classes with $\max_x\beta_x\le\lambda_3\iota L^{j-1}$.
The last-class variables remain independent, and Proposition~\ref{prop:bernstein},
with exponent $a+8+4(s-j)$, gives
\[
 \mathbb P\!\left(
 |\Phi-m|>4\lambda_1
 \bigl(\sqrt{\lambda_3\iota m}L^{j/2}+\lambda_3\iota L^j\bigr)
 \,\middle|\, (U_x)_{x\in I\setminus I_j}\right)
 \le h^{-a-8-4(s-j)}.
\]
When none of these three tail events occurs, adding the bounds for
$|\Phi-m|$ and $|m-\mu|$ gives the required deviation bound, since the
coefficients are bounded by $64\lambda_1\lambda_3z^{j-1}=z^j$.
We integrate the conditional probability bound over choices with
$\max_x\beta_x\le\lambda_3\iota L^{j-1}$ and add the preceding two tail
probabilities. The resulting bound
$3h^{-a-8-4(s-j)}\le h^{-a-4-4(s-j)}$ proves \eqref{eq:partite-polynomial}.

For a homogeneous part of degree $j$, at most
\[
 \left\lceil\frac{(j+1)j^j\log(|I|+1)}{j!}\right\rceil
 \le\lambda_2L
\]
colourings of $I$ with $j$ colours suffice so that every $j$-element
subset receives distinct colours in at least one colouring.
Indeed, in a uniformly random colouring a fixed $j$-set receives
distinct colours with probability $j!/j^j$. For independent
colourings, the probability that some $j$-set never receives
distinct colours is at most
$|I|^j\exp[-(j+1)\log(|I|+1)]<1$.
We assign each monomial to the first colouring in which its variables
receive distinct colours. Each resulting polynomial has one variable
from each colour class in every monomial, so \eqref{eq:partite-polynomial}
applies. There are at most $s\lambda_2L\le h^2$ such polynomials,
with derivative means still at most $\iota$. The sum of the square roots of their means is at most
$\sqrt{s\lambda_{2}L\mu}$. Summing the deviation bounds in
\eqref{eq:partite-polynomial} over the pieces gives the coefficient
$(1+s\lambda_2)z^s$ in the claimed bound; a union bound gives the stated
probability. The constant part has no fluctuation.
\end{proof}

\subsection{Correlation inequalities and conditional probabilities}

The following is Harris' Inequality, as stated by Harris~\cite{Harris1960}
and Alon and Spencer~\cite[Chapter~6]{AS}.
\begin{proposition}\label{prop:harris}
Let $I$ be finite and let $Y=(Y_i:i\in I)$ be a vector of
independent Bernoulli random variables. Let
$f,g:\{0,1\}^I\to\mathbb R$.
When $f$ and $g$ are both increasing or both decreasing in the
coordinatewise order,
\[
 \mathbb E(f(Y)g(Y))\ge\mathbb E(f(Y))\mathbb E(g(Y)).
\]
The inequality is reversed when one of $f,g$ is increasing and the
other is decreasing.
\end{proposition}

\begin{definition}\label{def:coordinate-dependency}
Let $I,J$ be finite sets and let $(Y_i:i\in I)$ be independent
Bernoulli random variables on a probability space
$(\Omega,\mathcal F,\mathbb P)$. For each $j\in J$, let
$T_j\subseteq I$ and $\mathcal A_j\subseteq\{0,1\}^{T_j}$. We define
\[
 \mathcal B_j=\{(Y_i)_{i\in T_j}\in\mathcal A_j\}.
\]
The \emph{coordinate dependency graph} of these events and sets $T_j$
is the simple graph with vertex set $J$ and edge set
\[
 \bigl\{\{j,\ell\}\in\tbinom J2:T_j\cap T_\ell\ne\varnothing\bigr\}.
\]
\end{definition}

\emph{Remark.} An edge records a possible shared Bernoulli variable.
Lemma~\ref{lem:conditional-lll} gives conditions under which conditioning
on $\bigcap_{j\in J}\overline{\mathcal B_j}$ is possible and gives upper bounds
for conditional probabilities and expectations.

The following lemma follows from Alon and Spencer~\cite[Lemma~5.1.1 and its proof]{AS}.
\Needspace{12\baselineskip}
\begin{lemma}\label{lem:conditional-lll}
Let $\iota,K\ge0$, $a\ge K+\iota+3$ and $h\ge2$.
Let $(\mathcal B_j:j\in J)$ and $(T_j:j\in J)$ be as in
Definition~\ref{def:coordinate-dependency}. Suppose
$\mathbb P(\mathcal B_j)\le h^{-a}$ for every $j\in J$ and the
coordinate dependency graph has maximum degree at most $h^K$.
Let $\mathcal G=\bigcap_{j\in J}\overline{\mathcal B_j}$. Then
\begin{enumerate}
\item $\mathbb P(\mathcal G)>0$.
\item Let $S\subseteq I$, and let $F\subseteq\Omega$ be an event
whose indicator is a function of $(Y_i)_{i\in S}$. If
$|\{j\in J:T_j\cap S\ne\varnothing\}|\le h^K$, then
$\mathbb P(F\mid\mathcal G)\le(1+h^{-\iota})\mathbb P(F)$.
\item Let $S\subseteq I$, and let $\Phi\ge0$ be a random variable
that is a function of $(Y_i)_{i\in S}$. If
$|\{j\in J:T_j\cap S\ne\varnothing\}|\le h^K$, then
$\mathbb E(\Phi\mid\mathcal G)\le(1+h^{-\iota})\mathbb E(\Phi)$.
\end{enumerate}
\end{lemma}
\begin{proof}
We put $\kappa=2h^{-a}<1$. Since $\kappa h^K\le1/4$,
\[
 (1-\kappa)^{h^K}\ge1-\kappa h^K\ge\frac34,
 \qquad h^{-a}\le\kappa(1-\kappa)^{h^K}.
\]
For $T\subseteq J$, we put
$A_T=\bigcap_{\ell\in T}\overline{\mathcal B_\ell}$, with
$A_\varnothing=\Omega$. We prove by induction on $|T|$ that
$\mathbb P(A_T)>0$ and, for every $j\in J\setminus T$,
\[
 \mathbb P(\mathcal B_j\mid A_T)\le\kappa.
\]
Both assertions hold when $T=\varnothing$. For nonempty $T$, we choose
$\ell\in T$. The induction hypothesis gives
\[
 \mathbb P(A_T)
 =\mathbb P(A_{T\setminus\{\ell\}})
 \bigl(1-\mathbb P(\mathcal B_\ell\mid A_{T\setminus\{\ell\}})\bigr)
 \ge(1-\kappa)\mathbb P(A_{T\setminus\{\ell\}})>0.
\]
We fix $j\notin T$, and partition $T$ into
$J_1=\{\ell\in T:T_\ell\cap T_j\ne\varnothing\}$ and
$J_0=T\setminus J_1$. The event $\mathcal B_j$ is independent of
$A_{J_0}$, because their indicators are functions of disjoint sets of
independent variables. We write $J_1=\{\ell_1,\ldots,\ell_d\}$.
The chain rule and the induction hypothesis give
\[
 \mathbb P(A_{J_1}\mid A_{J_0}) =\prod_{q=1}^d \bigl(1-\mathbb P(\mathcal B_{\ell_q}\mid A_{J_0\cup\{\ell_1,\ldots,\ell_{q-1}\}})\bigr) \ge(1-\kappa)^d,
\]
since each conditioning set in this product has fewer than $|T|$
elements. The empty product is one. As $d\le h^K$,
\[
 \mathbb P(\mathcal B_j\mid A_T)
 \le\frac{\mathbb P(\mathcal B_j\mid A_{J_0})}
 {\mathbb P(A_{J_1}\mid A_{J_0})}
 \le\frac{h^{-a}}{(1-\kappa)^d}\le\kappa.
\]
This completes the induction and proves (1).

For (2), we partition $J$ into the at most $h^K$ indices with
$T_j\cap S\ne\varnothing$ and the remaining indices. Applying the
same chain-rule calculation and independence gives
\[
 \mathbb P(F\mid\mathcal G) \le\mathbb P(F)(1-\kappa)^{-h^K} \le\mathbb P(F)\exp(4h^{K-a}) \le(1+h^{-\iota})\mathbb P(F).
\]
Here $-\log(1-\kappa)\le2\kappa$ and
$\exp(4h^{K-a})\le1+8h^{K-a}\le1+h^{-\iota}$, since
$4h^{K-a}\le1/2$.

For (3), the indicator of $\{\Phi>\lambda\}$, where $\lambda\ge0$, is a
function of the variables indexed by $S$, so (2) applies. Integrating
the resulting inequality over $\lambda\in[0,\infty)$ gives
\[
 \mathbb E(\Phi\mid\mathcal G) =\int_0^\infty\mathbb P(\Phi>\lambda\mid\mathcal G)\,d\lambda \le(1+h^{-\iota})\int_0^\infty\mathbb P(\Phi>\lambda)\,d\lambda =(1+h^{-\iota})\mathbb E(\Phi).
\]
This completes the proof of Lemma~\ref{lem:conditional-lll}.
\end{proof}

\Needspace{8\baselineskip}
\section{Focus analysis}\label{sec:general-focus}

In this section, we use the parameters of Subsection~\ref{subsec:actual-process}
and the one-round notation of Subsection~\ref{subsec:process-focus-parameters}.
Let $t\ge0$ be the integer round under consideration. We write
$\Sigma=\Sigma(\vb d_t)$, $D=D(\vb d_0)$,
$U=X\setminus\mathsf W$, $q=\Sigma^{-1+\delta_2}$ and
$p=\eta\Sigma^{-\delta_2}$.
Assume $X\ne\varnothing$, $D,\Sigma\ge2$ and $0<\eta\le1/4$.
Here $\mathcal S\subseteq\bigcup_{i=2}^r\binom Xi$ is the current
set system. Empty maxima are zero.

For the focus $Z\subseteq U$ and $x\in X$, we write
$\xi_x=\vb1_{\{x\in Z\}}$. Before conditioning
$Z$ to be good, $(\xi_x:x\in U)$ are independent Bernoulli$(q)$
variables, and $\xi_x=0$ on $\mathsf W$. With $Z$ fixed, the variables
$(\chi_x:x\in Z)$ and $(\mathsf{w}_x:x\in U)$ are mutually independent
Bernoulli variables, with respective parameters $p$ and $1-\omega/s_x$
from \eqref{eq:process-survival-probabilities}--\eqref{eq:process-w-law}.
Here $\chi_x=0$ outside $Z$ and $\mathsf{w}_x=1$ on $\mathsf W$.
The variables $\chi_x$ and $J_x$ record membership in $Y$ and $I$,
respectively, while $\mathsf{w}_x=1$ puts a surviving atom in $\mathsf W'$.
By
\eqref{eq:process-eligible-atoms}, for $x\in X$,
\[
 \psi_x=\vb1_{\{x\in X'\}},\qquad
 \psi_x(1-\mathsf{w}_x)=\vb1_{\{x\in X'\setminus\mathsf W'\}}.
\]
The laws of $\chi_x$ and $\mathsf{w}_x$ here precede conditioning on
a good construction round.

The initial phase consists of the rounds $0\le t<T_{\mathrm{init}}$,
where the deterministic endpoint is
\begin{equation}\label{eq:general-initial-period-length}
 T_{\mathrm{init}}=\left\lfloor\frac{2^{8r^2}}{\eta}\right\rfloor.
\end{equation}
This phase is needed because $d_{2,t}$ may initially be too small
for deletion through $2$-blocks to control $d_i'(A)$ and $d_i'(x,y)$.
Lemma~\ref{lem:general-scale} ensures that $d_{2,t}\ge\Sigma_t/2$
and $\Sigma_t$ decreases geometrically by time $T_{\mathrm{init}}$.
Before then, the stronger initial bounds in terms of $D$ bound
$|N_i'(A)\setminus N_i(A)|$ and $|N_i'(x,y)\setminus N_i(x,y)|$
without needing a lower bound on the number of members of $N_i(A)$
or $N_i(x,y)$ removed in a round, as proved in
Lemma~\ref{lem:initial-focus-estimates} and
Subsection~\ref{subsec:general-initial-period}.

\begin{definition}[Good focus]\label{def:good-focus}
Let $a>0$, $b_{1}\ge1$ and $k_t>0$. We put
\[
 b_{2}=\left[\frac{4r(a+2^{10}r^3+1)}{\delta_2}\right]^r,
 \qquad b_{3}=(2^{16}r(b_{2}+1))^{2^{10}r^4},
 \qquad b_{4}=e^{-2^{8r^2+1}}.
\]
A set $Z\subseteq U$ is a \emph{good focus} in round $t$ if the
following conditions hold.
\begin{enumerate}[label=\textup{(GF\arabic*)},ref=GF\arabic*]
\item\label{gf:small-counts} For $3\le i\le r$ and $A\subseteq X$ with $2\le|A|<i$,
and for $2\le i\le r$ and distinct $x,y\in X$, respectively,
\begin{equation}\label{eq:general-constant-focused-completions}
 d_i(A,Z)\le b_{2},\qquad d_i(x,y,Z)\le b_{2}.
\end{equation}
\item\label{gf:degrees} For $v\in X$ and $2\le i\le r$,
\begin{equation}\label{eq:general-focused-degrees}
 |d_i(v,Z)-q^{i-1}d_i(v,U)|
 \le\Sigma^{(i-1)\delta_2-4\delta_3}.
\end{equation}
For $v\in X$ and integers $2\le k<j\le r$,
\begin{align}
 \left|
 \sum_{\substack{v\in S\in\mathcal S_j\\S\setminus\{v\}\subseteq U}}
 \binom{|(S\setminus\{v\})\cap Z|}{j-k}
 -\binom{j-1}{k-1}q^{j-k}d_j(v,U)
 \right|
 &\le\Sigma^{k-1+(j-k)\delta_2-4\delta_3},
 \label{eq:general-focused-contractors}\\
 \sum_{v\in S\in\mathcal S_j}
 \binom{|(S\setminus\{v\})\cap Z|}{j-k}
 -\binom{j-1}{k-1}q^{j-k}d_j(v,X)
 &\le\Sigma^{k-1+(j-k)\delta_2-4\delta_3}.
 \label{eq:general-full-contractors}
\end{align}
\item\label{gf:degree-dependence} For $v\in X$ and an integer $2\le k\le r$, let
\begin{equation}\label{eq:degree-candidate-sets}
 \mathcal S_k(v)=\left\{B\in\binom Xk:v\in B,\
 B\subseteq P,\ P\setminus B\subseteq Z
 \text{ for some }P\in\mathcal S\right\}.
\end{equation}
For $B\in\mathcal S_k(v)$, let
\begin{equation}\label{eq:degree-atoms}
 \mathcal A_B=\bigcup_{\substack{P\in\mathcal S,\ P\setminus B\subseteq Z\\
                  P\cap B\notin\{\varnothing,\{v\}\}}}(P\setminus B),
\end{equation}
\begin{equation}\label{eq:degree-support}
 T_B={}\{1\}\times\left[ ((B\setminus\{v\})\cap Z)\cup\mathcal A_B \cup\bigcup_{\substack{P\in\mathcal S,\ P\subseteq Z\setminus\{v\}\\ P\cap\mathcal A_B\ne\varnothing}}P\right] {}\cup\bigl(\{2\}\times(B\setminus(\mathsf{W}\cup\{v\}))\bigr).
\end{equation}
These sets are required to satisfy
\begin{equation}\label{eq:general-degree-read}
 \max\left\{1,\max_{\substack{j\in\{1,2\}\\x\in X}}
 |\{B\in\mathcal S_k(v):(j,x)\in T_B\}|\right\}
 \le b_{3}\Sigma^{k-1-\delta_1+3r\delta_2}.
\end{equation}
\item\label{gf:coverage} For integers $3\le i\le r$ and $A\subseteq X$ with
$2\le|A|<i$, or integers $2\le i\le r$ and distinct $x,y\in X$,
respectively, let
\[
 (\mathcal B,F_0)=(N_i(A),A)\quad\text{or}\quad
 (\mathcal B,F_0)=(N_i(x,y),\{x,y\}).
\]
For each choice of $i,A$ or $i,x,y$, the quantity $\ell$, the sets $C_u$
for $u\in\bigcup\mathcal B$, and $C_B,\mathcal T_B$ for
$B\in\mathcal B$ are defined by
\begingroup
\allowdisplaybreaks[4]
\begin{align}
 \ell&=\frac{2b_{3}(b_{1}+r)}{\delta_3}
 \left(\max_{v\in X}|\{B\in\mathcal B:v\in B\}|
                 +|\mathcal B|\Sigma^{-\delta_4(2)}\right),\nonumber\\
 C_u&=\left\{w\in N_2(u,U)\setminus F_0:
           \sum_{B\in\mathcal B}|B\cap N_2(w)|\le\ell\right\},\nonumber\\
 C_B&=\left(\bigcup_{u\in B}C_u\right)\setminus B,\nonumber\\
 \mathcal T_B&=(C_B\cap Z)\ \cup\!
 \bigcup_{\substack{w\in S\in\mathcal S,\ w\in C_B\cap Z\\
             S\subseteq Z,\ S\cap(F_0\cup B)=\varnothing}}S.
 \label{eq:coverage-candidate-definitions}
\end{align}
\endgroup
If $t\ge T_{\mathrm{init}}$ and
\begin{equation}\label{eq:coverage-family-activation}
 d_2\ge\frac12\Sigma,\qquad
 |\mathcal B|\ge\frac12
 \begin{cases}
  \Sigma^{i-|A|-\delta_5(i)},&\text{for the choice }i,A,\\
  \Sigma^{i-1-\delta_4(i)},&\text{for the choice }i,x,y,
 \end{cases}
\end{equation}
then, for every $u\in\bigcup\mathcal B$,
\begin{equation}\label{eq:coverage-focus-candidates}
 \bigl||C_u\cap Z|-q|C_u|\bigr|\le \frac{\delta_3}{b_{3}}q|C_u|,
\end{equation}
and
\begin{equation}\label{eq:coverage-focus-read}
 \max\left\{1,\max_{u\in Z}|\{B\in\mathcal B:u\in \mathcal T_B\}|\right\}
 \le r b_{1}b_{3}^2|\mathcal B|\Sigma^{2\delta_2-\delta_1}.
\end{equation}
\item\label{gf:new-counts} Let $(\xi_x:x\in U)$ be variables, with $\xi_x=0$ for
$x\in\mathsf W$. For $F\subseteq X$, we write
$\xi_F=\prod_{x\in F}\xi_x$, where $\xi_\varnothing=1$.
For a polynomial $\Phi$ in these variables, $\Phi(Z)$ denotes its
value at $\xi_x=\vb1_{\{x\in Z\}}$. For $T\subseteq U$,
$\partial_T\Phi$ is the coefficient derivative of
Definition~\ref{def:polynomial}.
Let $3\le i<k\le r$ be integers, and let $A\subseteq X$ satisfy
$2\le|A|<i$. We define
\begin{subequations}\label{eq:general-birth-witness-polynomials}
\begin{equation}\label{eq:proper-birth-witness-polynomial}
 \mathsf Y_{i,A;k}(\xi)
 =\sum_{A\subseteq P\in\mathcal S_k}
   \sum_{\substack{F\subseteq P\setminus A\\|F|=k-i}}\xi_F.
\end{equation}
Let $2\le i<r$ and $i\le k,l\le r$ be integers with $(k,l)\ne(i,i)$,
let $x,y\in X$ be distinct, and let $h$ be an integer with
$0\le h\le\min\{k-i,l-i\}$. Let $(S,R,P,Q)$ range over ordered
quadruples of pairwise disjoint subsets of $X\setminus\{x,y\}$ of
sizes $i-1,h,k-i-h,l-i-h$, respectively. We define
\begin{equation}\label{eq:common-birth-witness-polynomial}
 \mathsf Y_{i,x,y;k,l,h}(\xi)
 =\sum_{\substack{\{x\}\cup S\cup R\cup P\in\mathcal S_k\\
                  \{y\}\cup S\cup R\cup Q\in\mathcal S_l}}
       \xi_{R\cup P\cup Q}.
\end{equation}
\end{subequations}
For either polynomial $\mathsf Y$, let
\begin{equation}\label{eq:birth-polynomial-parameters}
 (s,m,\iota)=
 \begin{cases}
 (k-i,\ i-|A|,\ \Sigma^{i-|A|-\delta_5(k)}),
                 &\text{for }\mathsf Y_{i,A;k},\\
 (k+l-2i-h,\ i-1,\ \Sigma^{i-1-\delta_4(i+1)}),
                 &\text{for }\mathsf Y_{i,x,y;k,l,h}.
 \end{cases}
\end{equation}
For $T\subseteq U$ with $|T|\le s$, the following bounds hold
when $t<T_{\mathrm{init}}$ and when $t\ge T_{\mathrm{init}}$, respectively:
\begin{align}
 (\partial_T\mathsf Y)(Z)
 &\le b_{1}^2b_{3}^2b_{4}^{-(s-|T|)}k_t
 D^{m-\delta_0+(s-|T|)\delta_2}(\log D)^{s-|T|},
 \label{eq:general-initial-focus-derivatives}\\
 (\partial_T\mathsf Y)(Z)
 &\le b_{1}^2b_{3}^2\iota\Sigma^{(s-|T|)\delta_2}
                         (\log\Sigma)^{s-|T|}.
 \label{eq:general-birth-focus-derivatives}
\end{align}
\end{enumerate}
Let $\mathcal G_{{\rm foc},t}$ be the family of good focus sets in round $t$.
\end{definition}

\emph{Remark.} In \eqref{eq:degree-support}, $(1,x)$ represents $\chi_x$ and
$(2,x)$ represents $\mathsf{w}_x$. After fixing $\chi_v=0$, the
variables indexed by $T_B$ suffice to determine the contribution of
$B$ to the degree sums in Lemma~\ref{lem:degree-indicator-identities}.
When $v$ survives, those sums give $d_k'(v,X')$ and
$d_k'(v,X'\setminus\mathsf W')$. The bound
\eqref{eq:general-degree-read} limits how many contributions use any
one Bernoulli variable, which gives the concentration of these degrees
in Lemma~\ref{lem:general-degree-round}.

In \eqref{eq:coverage-candidate-definitions}, $C_u$ consists of atoms
in $U\setminus F_0$ lying in a $2$-block with $u$, after excluding those
used too often among the choices for members of $\mathcal B$;
$C_B$ combines these choices for $u\in B$. The variables
$(\chi_x:x\in\mathcal T_B)$ determine whether $C_B\cap Z$
contains an accepted atom after the variables at $F_0\cup B$ are
set to zero. Lemma~\ref{lem:general-coverage} uses
\eqref{eq:coverage-focus-candidates} and \eqref{eq:coverage-focus-read}
to bound the probability that too few members $B\in\mathcal B$ have
this property. If the atoms of $F_0$ survive, each such $B$ loses an atom.
Lemma~\ref{lem:focus-coverage} establishes the required bounds on
$|C_u\cap Z|$ and on how often one $\chi_x$ is used.

Evaluated at $Z$, \eqref{eq:proper-birth-witness-polynomial} counts
ways to obtain an $i$-set containing $A$ by deleting $k-i$ focus
atoms from a $k$-block. Equation~\eqref{eq:common-birth-witness-polynomial}
counts ordered choices of two blocks and deleted sets leaving the same
$(i-1)$-set beside $x$ and $y$, respectively. Distinct choices are
counted separately even when they leave the same sets. The derivatives count
these choices with the deleted atoms in $T$ prescribed and all other
deleted atoms required to lie in $Z$.
The bounds \eqref{eq:general-initial-focus-derivatives} and
\eqref{eq:general-birth-focus-derivatives} control terms from blocks
$S\in\mathcal S$ with $|S|>i$ when bounding
$d_i'(A)$ and $d_i'(x,y)$, as used in
Lemma~\ref{lem:general-births} and
Subsection~\ref{subsec:general-initial-period}.

For $\ell\in\mathbb Z_{\ge0}$, we define
\begin{equation}\label{eq:general-initial-envelope}
 k_\ell=c(1+b_{1}^2b_{3}^4b_{4}^{-4r})^\ell(\log D)^{4r\ell}
 .
\end{equation}
We also use $\mathcal G_{{\rm foc},t}$ for the event that $Z$ is good.

\subsection{\texorpdfstring{Concentration of $d_i(A,Z)$ and $d_i(x,y,Z)$}{\detokenize{Concentration of d_i(A,Z) and d_i(x,y,Z)}}}
\label{subsec:focus-basic-counts}

Assume $X\ne\varnothing$, $\Sigma\ge2$, $0<\eta\le1/4$, and
for $2\le i\le r$, $E_i,\widehat E_i\le\Sigma^{i-1}$. Thus
$d_i(v,X)\le2\Sigma^{i-1}$. For $A\subseteq X$ with
$2\le|A|<i\le r$, and for $2\le i\le r$ and distinct $x,y\in X$,
respectively, assume
\begin{align}
 d_i(A)&\le\Sigma^{i-|A|-\delta_5(i)},
 \label{eq:local-absolute-codegrees}\\
 d_i(x,y)&\le\Sigma^{i-1-\delta_4(i)}.
 \label{eq:local-absolute-common}
\end{align}

The exponents fixed in Subsection~\ref{subsec:actual-process} satisfy
$0<\delta_0\le1/2$ and
\[
 \delta_1=\frac{\delta_0}{16r},\qquad
 \delta_2=\frac{\delta_1}{64r},\qquad
 \delta_3=\frac{\delta_2}{64}.
\]
For $2\le i\le r$ and $3\le i\le r$, respectively,
\[
 \delta_4(i)=(i-1)\delta_1,\qquad
 \delta_5(i)=(r+i-2)\delta_1.
\]
In particular, $\delta_2<\delta_1/r$ and
$\delta_1-r\delta_2=63\delta_1/64>0$.

\begin{lemma}\label{lem:general-constant-focus}
Let $a>0$, and let
\[
 b_{2}=\left[\frac{4r(a+2^{10}r^3+1)}{\delta_2}\right]^r.
\]
For $3\le i\le r$ and $A\subseteq X$ with $2\le|A|<i$,
\begin{equation*}
 \mathbb P\bigl(d_i(A,Z)>b_{2}\bigr)\le\Sigma^{-a-2^{10}r^3}.
\end{equation*}
For $2\le i\le r$ and distinct $x,y\in X$,
\begin{equation*}
 \mathbb P\bigl(d_i(x,y,Z)>b_{2}\bigr)\le\Sigma^{-a-2^{10}r^3}.
\end{equation*}
\end{lemma}
\begin{proof}
We put
\[
 m=1+\left\lceil\frac{a+2^{10}r^3}{\delta_1-r\delta_2}\right\rceil.
\]
Since $\delta_1-r\delta_2=63r\delta_2$, we have
$rm\le2r+(a+2^{10}r^3)/(63\delta_2)$, and hence
$b_{2}\ge\max\{(rm)^{r-1},4(a+2^{10}r^3+1)/\delta_2\}$.
By \eqref{eq:restricted-neighborhoods} and
\eqref{eq:restricted-common-links},
\begin{equation}\label{eq:focus-completion-polynomials}
 d_i(A,Z)=\sum_{F\in N_i(A,U)}\xi_F,\qquad d_i(x,y,Z)=\sum_{F\in N_i(x,y,U)}\xi_F.
\end{equation}
Let $\Phi$ be either polynomial in
\eqref{eq:focus-completion-polynomials}, with
$\mathfrak d=i-|A|$ or $\mathfrak d=i-1$, respectively.
Its coefficients are zero or one. For $T\subseteq U$ with
$j=|T|<\mathfrak d$, the derivative of the first polynomial has at
most $d_i(A\cup T)$ terms when $T\cap A=\varnothing$, and is zero otherwise.
For the second polynomial, we use $d_i(x,y)$ when $T=\varnothing$
and $d_i(\{x\}\cup T)$ when $T\ne\varnothing$ and $x\notin T$;
the derivative is zero if $x\in T$.
Thus \eqref{eq:local-absolute-codegrees} and
\eqref{eq:local-absolute-common} give
\begin{equation}\label{eq:focus-completion-derivative-mean}
 \mathbb E(\partial_T\Phi)
 \le\Sigma^{(\mathfrak d-j)\delta_2-\delta_1}
 \le\Sigma^{-(\delta_1-r\delta_2)}.
\end{equation}
Since $1\le\mathfrak d\le r-1$, Corollary~\ref{cor:polynomial-envelope}~(3)
implies
\[
 \mathbb P(\Phi>b_{2})
 \le\mathbb P\bigl(\Phi>(\mathfrak d m)^{\mathfrak d}\bigr)
 \le\Sigma^{-a-2^{10}r^3}.
\]
This completes the proof of Lemma~\ref{lem:general-constant-focus}.
\end{proof}

\subsection{Set systems of indicators}
\label{subsec:set-system-indicators}

\begin{definition}\label{def:shadow-balls}
For a set system $\mathcal S$ on $X$, the \emph{shadow graph} is the
simple undirected graph on $X$ in which distinct atoms $u,v$ are
adjacent if and only if $\{u,v\}\subseteq S$ for some $S\in\mathcal S$.
For $v\in X$ and an integer $j\ge0$, the \emph{closed ball}
$\mathbb B_j(v,\mathcal S)$ consists of all atoms reachable from $v$ by a path
of at most $j$ edges in this graph; in particular,
$\mathbb B_0(v,\mathcal S)=\{v\}$.
\end{definition}

\emph{Remark.} The radius counts transitions between atoms lying in a common
block. Lemma~\ref{lem:conditioning-locality} uses these balls to bound
the number of random variables determining an event and the number of
events that can share any of those variables.

\begin{lemma}\label{lem:conditioning-locality}
Let $r\ge2$ and $\ell,s,\mathfrak d\ge0$ be integers, let
$\lambda,m\ge1$, and let
\[
  h\ge\max\left\{
 4[\lambda(r-1)^2]^{2\ell},\sqrt{(s+1)m},\mathfrak d+1\right\}.
\]
Let $\mathcal S$ be a set system on a finite set $X$, of upper rank
at most $r$, such that $d_i(v,X)\le\lambda h^{i-1}$ for every
$v\in X$ and $2\le i\le r$. For each $x\in X$, let $U_{x,1}$
and $U_{x,2}$ be Bernoulli random variables, so that
$\{U_{x,j}:x\in X,\ j\in\{1,2\}\}$ is a family of independent
random variables. Let $\mathcal E$ be a finite family of events in
the joint probability space of the variables $U_{x,j}$.
For each $E\in\mathcal E$, let
$v_E\in X$, let $S_E\subseteq\mathbb B_\ell(v_E,\mathcal S)$
with $|S_E|\le s$, and let
$T_E\subseteq\mathbb B_\ell(v_E,\mathcal S)\times\{1,2\}$
be such that the indicator of $E$ is a function of
$(U_{x,j})_{(x,j)\in T_E}$. Suppose that, for every $v\in X$ and
$S\subseteq\mathbb B_\ell(v,\mathcal S)$ with $|S|\le s$,
\[
 |\{E\in\mathcal E:v_E=v,\ S_E=S\}|\le m.
\]
Then, for every $v\in X$ and every
$T\subseteq\mathbb B_\ell(v,\mathcal S)\times\{1,2\}$,
\begin{enumerate}
\item $|T|\le h^{\ell(r-1)+1}$.
\item $|\{E\in\mathcal E:T_E\cap T\ne\varnothing\}|
      \le h^{\ell(s+2)(r-1)+s+3}$.
\item Every multilinear polynomial of degree at most $\mathfrak d$
in $(U_{x,j})_{(x,j)\in T}$ has at most
$ h^{\mathfrak d[\ell(r-1)+1]+1}$ nonzero coefficient derivatives.
\end{enumerate}
\end{lemma}
\begin{proof}[Proof of \textup{(1)}]
For every $v\in X$, the number of atoms $w\ne v$ lying in a block
with $v$ is at most
\[
 \sum_{i=2}^r(i-1)d_i(v,X)
 \le\lambda(r-1)^2 h^{r-1}.
\]
This upper bound is at least two. For $0\le j\le2\ell$, counting
paths from $v$ gives
\[
 |\mathbb B_j(v,\mathcal S)| \le\sum_{i=0}^j[\lambda(r-1)^2 h^{r-1}]^i \le2[\lambda(r-1)^2]^j h^{j(r-1)} \le\frac12 h^{j(r-1)+1}.
\]
The factor two bounds the geometric sum after division by its
last term; the final inequality uses
$ h\ge4[\lambda(r-1)^2]^{2\ell}$. There are two variables at
each atom, so $|T|\le2|\mathbb B_\ell(v,\mathcal S)|$ proves (1).

\emph{Proof of \textup{(2)}.} For each fixed $w\in X$, there are at most
$(s+1)|\mathbb B_\ell(w,\mathcal S)|^s$ possible sets $S_E$.
Consequently,
\[
 |\{E\in\mathcal E:v_E=w\}|
 \le m(s+1) h^{s[\ell(r-1)+1]}
 \le h^{2+s[\ell(r-1)+1]}.
\]
If $T_E\cap T\ne\varnothing$, then
$v_E\in\mathbb B_{2\ell}(v,\mathcal S)$. There are at most
$ h^{2\ell(r-1)+1}$ such atoms, giving (2).

\emph{Proof of \textup{(3)}.} A nonzero coefficient derivative is specified
by at most $\mathfrak d$ variables. By (1) and
$ h\ge\mathfrak d+1$, their number is at most
\[
 \sum_{j=0}^{\mathfrak d}\binom{|T|}{j}
 \le(\mathfrak d+1) h^{\mathfrak d[\ell(r-1)+1]}
 \le h^{\mathfrak d[\ell(r-1)+1]+1}.
\]
This completes the proof of Lemma~\ref{lem:conditioning-locality}.
\end{proof}

\Needspace{8\baselineskip}
\subsection{\texorpdfstring{Concentration of $d_i(v,Z)$}{\detokenize{Concentration of d_i(v,Z)}}}
\label{subsec:focus-degree-concentrations}

We retain the assumptions of Subsection~\ref{subsec:focus-basic-counts}.

Let $2\le k<j\le r$ be integers and $v\in X$.
We denote the sums in \eqref{eq:general-focused-contractors} and
\eqref{eq:general-full-contractors} by $M_{j,k}(v)$ and
$\widehat M_{j,k}(v)$, respectively. For $2\le k\le r$, we put
$M_{k,k}(v)=d_k(v,U)$ and $\widehat M_{k,k}(v)=d_k(v,X)$.
The hat removes the restriction $S\setminus\{v\}\subseteq U$;
the $j-k$ atoms selected in $Z$ still lie outside $\mathsf W$.
It is straightforward to check that
\begin{equation}\label{eq:focused-contraction-expectations}
 \mathbb E(M_{j,k}(v))
 =\binom{j-1}{k-1}q^{j-k}d_j(v,U),\qquad \mathbb E(\widehat M_{j,k}(v))
 \le\binom{j-1}{k-1}q^{j-k}d_j(v,X).
\end{equation}

We put
\[
 b_{5}=(b_{1}b_{3})^{16/\delta_3}.
\]

\begin{lemma}\label{lem:general-focus-estimates}
Let $a>0$ and $b_{1}\ge1$, and assume $\Sigma\ge b_{5}$.
For $v\in X$ and $2\le i\le r$,
\begin{equation}\label{eq:focus-degree-tail}
 \mathbb P\left(
 |d_i(v,Z)-q^{i-1}d_i(v,U)|
 >\Sigma^{(i-1)\delta_2-4\delta_3}\right)
 \le\frac{\Sigma^{-a}}{(r-1)^2}.
\end{equation}
For $v\in X$ and $2\le k<j\le r$,
\begin{align}
 \mathbb P\left(
 \left|M_{j,k}(v)-\binom{j-1}{k-1}q^{j-k}d_j(v,U)\right|
 >\Sigma^{k-1+(j-k)\delta_2-4\delta_3}\right)
 &\le\frac{\Sigma^{-a}}{(r-1)^2},
 \\
 \mathbb P\left(
 \widehat M_{j,k}(v)-\binom{j-1}{k-1}q^{j-k}d_j(v,X)
 >\Sigma^{k-1+(j-k)\delta_2-4\delta_3}\right)
 &\le\frac{\Sigma^{-a}}{(r-1)^2}.
 \label{eq:focus-full-contraction-tail}
\end{align}
\end{lemma}
\begin{proof}
We set $z=a+2^{10}r^3$. The choices of $b_3,b_5$ give
\[
 \left(\frac{32r}{7e\delta_2}\right)^r\le b_3,
 \qquad b_5\ge b_3^{64/(7\delta_2)}.
\]
Thus $\Sigma\ge b_5$ implies
\[
 \Sigma\ge(8r)^{2^{10}r^2},\qquad
 b_3(\log\Sigma)^r\le\Sigma^{\delta_2/2-4\delta_3}.
\]
We apply Lemma~\ref{lem:conditioning-locality}~(1) with
$h=\Sigma$, $\lambda=2$, $\ell=1$, $s=\mathfrak d=0$, $m=1$,
$U_{x,1}=\xi_x$, $U_{x,2}=0$ and $\mathcal E=\varnothing$.
The variables in $d_i(v,Z)$, $M_{j,k}(v)$ and $\widehat M_{j,k}(v)$
are indexed by atoms in $\mathbb B_1(v,\mathcal S)$, so their number
is at most $\Sigma^r\le\Sigma^{2^9r^3}$.
Since $b_2\ge4(z+1)/\delta_2$, the numerical factors supplied by
Proposition~\ref{prop:polynomial} in the following applications are at most
\[
 [2^{16}r(b_2+1)]^{8r^2}\le\frac{b_3}{2}.
\]
For $d_i(v,Z)$, a nonzero derivative on
$0<|T|=\ell_0<i-1$ coordinates has expectation at most
\[
 q^{i-1-\ell_0}d_i(\{v\}\cup T)
 \le\Sigma^{(i-1-\ell_0)\delta_2-\delta_5(i)}\le1.
\]
The full derivative is at most one, and
$\mathbb E(d_i(v,Z))\le2\Sigma^{(i-1)\delta_2}$. Proposition~\ref{prop:polynomial}
(Proposition~\ref{prop:bernstein} for $i=2$) gives deviation at most
$b_3\Sigma^{(i-1)\delta_2/2}(\log\Sigma)^r$.
The deviation is therefore at most
$\Sigma^{(i-1)\delta_2-4\delta_3}$, proving
\eqref{eq:focus-degree-tail}.

For $\ell=j-k\ge1$, a nonzero derivative of $M_{j,k}(v)$ on
$1\le|T|=\ell_0\le\ell$ coordinates has expectation at most
\[
 \binom{j-1-\ell_0}{\ell-\ell_0}q^{\ell-\ell_0}d_j(\{v\}\cup T,X\setminus \mathsf{W})
 \le2^r\Sigma^{k-1-\delta_5(j)+(\ell-\ell_0)\delta_2}.
\]
Here $|\{v\}\cup T|\le j-k+1\le j-1$, so
\eqref{eq:local-absolute-codegrees} applies also when $\ell_0=\ell$. The mean is at most
$2^{r+1}\Sigma^{k-1+\ell\delta_2}$. Proposition~\ref{prop:polynomial}
bounds the deviation by
\[
 \frac{b_3}{2}\Sigma^{k-1+\ell\delta_2} \left(\Sigma^{-\delta_5(j)/2-\delta_2/2} +\Sigma^{-\delta_5(j)-\delta_2}\right)(\log\Sigma)^r \le\Sigma^{k-1+\ell\delta_2-4\delta_3}.
\]
For $\widehat M_{j,k}(v)$, we use the same derivative bounds and its
mean bound in \eqref{eq:focused-contraction-expectations}.
Each probability in
\eqref{eq:focus-degree-tail}--\eqref{eq:focus-full-contraction-tail}
is at most $\Sigma^{-z}\le\Sigma^{-a}/(r-1)^2$.
\end{proof}

\subsection{\texorpdfstring{The degrees $d_k'(v,X')$ and $d_k'(v,X'\setminus \mathsf{W}')$}{The residual degrees}}
\label{subsec:general-degree-indicators}
\label{subsec:general-read-counts}

Assume $X\ne\varnothing$, $\Sigma\ge2$, $0<\eta\le1/4$,
\eqref{eq:local-absolute-codegrees}--\eqref{eq:local-absolute-common},
and $E_i,\widehat E_i\le\Sigma^{i-1}$ for $2\le i\le r$.

Let $2\le k\le r$ be an integer, $Z\subseteq U$ and $v\in X$.
We use $\mathcal S_k(v)$ and $T_B$ from
\eqref{eq:degree-candidate-sets} and \eqref{eq:degree-support}
for these $Z,v,k$. We define $\mathcal C_Z(A)$ for $A\subseteq X$
with $2\le|A|\le r$, $\widehat T_B$ for $B\in\mathcal S_k(v)$,
and $\mu_k(v)$ by
\begin{align}
 \mathcal C_Z(A)&=\{C\subseteq Z\setminus A:A\cup C\in\mathcal S\},
 \label{eq:degree-completions}\\
 \widehat T_B&=T_B\cap(\{1\}\times X),
 \label{eq:degree-hat-support}\\
 \mu_k(v)&=\max\left\{1,\max_{\substack{j\in\{1,2\}\\x\in X}}
 |\{B\in\mathcal S_k(v):(j,x)\in T_B\}|\right\}.
 \label{eq:degree-multiplicity}
\end{align}

\begin{definition}\label{def:degree-original-variables}
Let $Z\subseteq X\setminus \mathsf{W}$.
Let $(\chi_x:x\in Z)$ be independent Bernoulli$(p)$ random variables,
and let $\chi_x=0$ for $x\in X\setminus Z$. For $x\in X$, we recall
\begin{equation}\label{eq:degree-original-variables}
 J_x=\chi_x\prod_{x\in S\in\mathcal S}
                \left(1-\prod_{y\in S\setminus\{x\}}\chi_y\right),\qquad \psi_x=(1-\chi_x)\prod_{x\in S\in\mathcal S}
                \left(1-\prod_{y\in S\setminus\{x\}}J_y\right).
\end{equation}
With expectations under this law, for $x\in X$ we define
\begin{equation}\label{eq:degree-survival-probabilities}
 s_x=\mathbb E(\psi_x),\qquad
 \omega=\min_{x\in X}s_x.
\end{equation}
Independently of $(\chi_x:x\in Z)$, let
$(\mathsf{w}_x:x\in X\setminus \mathsf{W})$ be independent Bernoulli variables.
Let $\mathbb P_Z$ denote the joint law. For $x\in Z$ and
$y\in X\setminus\mathsf W$, its probabilities are
\begin{equation*}
 \mathbb P_Z(\chi_x=1)=p,\qquad
 \mathbb P_Z(\mathsf w_y=1)=1-\frac{\omega}{s_y}.
\end{equation*}
Let $\mathsf{w}_x=1$ for $x\in \mathsf{W}$.
\end{definition}

\emph{Remark.} Here $J_x=1$ records acceptance of $x$, and $\psi_x=1$ records
its survival in $X'$. The product $\psi_x(1-\mathsf{w}_x)$ records membership
in $X'\setminus\mathsf W'$, with probability $\omega$ under $\mathbb P_Z$ for every
$x\in X\setminus\mathsf W$. These variables are used in
Lemmas~\ref{lem:general-accepted-sets} and~\ref{lem:general-degree-round}.

\begin{definition}\label{def:degree-modified-variables}
Let $k\in\{2,\ldots,r\}$, $Z\subseteq X\setminus \mathsf{W}$ and $v\in X$.
We use the variables and law in Definition~\ref{def:degree-original-variables}
for this $Z$, and $\mathcal S_k(v),\mathcal C_Z(A)$ from
\eqref{eq:degree-candidate-sets}, \eqref{eq:degree-completions}.
Let
\begin{equation}\label{eq:degree-conditioned-law}
 \mathbb P_{Z,v}^0=\mathbb P_Z(\,\cdot\mid\chi_v=0).
\end{equation}
We write $\mathbb E_{Z,v}^0$ for expectation under $\mathbb P_{Z,v}^0$.
We define $\chi_x^{v},J_x^{v}$ for $x\in X$, $J_C^{v}$ for $C\subseteq X$,
$R_A^{v}$ for $A\subseteq X$ with $2\le|A|\le r$,
$\psi_u^{B,v}$ for $B\in\mathcal S_k(v)$ and $u\in B\setminus\{v\}$,
and $R_F^{B,v}$ for $B\in\mathcal S_k(v)$ and
$F\subsetneq B$ with $|F|\ge2$ by
\begingroup
\allowdisplaybreaks[4]
\begin{align}
 \chi_x^{v}&=\begin{cases}0,&x=v,\\\chi_x,&x\ne v,\end{cases}
 \notag\\
 J_x^{v}&=\chi_x^{v}
  \prod_{x\in S\in\mathcal S}
       \left(1-\prod_{y\in S\setminus\{x\}}\chi_y^{v}\right),
 \label{eq:degree-modified-J}\\
 J_C^{v}&=\prod_{x\in C}J_x^{v},
 \notag\\
 R_A^{v}&=1-\prod_{C\in\mathcal C_Z(A)}(1-J_C^{v}),
 \label{eq:degree-modified-R}\\
 \psi_u^{B,v}
 &=(1-\chi_u^{v})
 \prod_{\substack{u\in S\in\mathcal S,\ S\cap B=\{u\}\\
                  S\setminus\{u\}\subseteq Z}}
 (1-J_{S\setminus\{u\}}^{v}),
 \label{eq:degree-restricted-psi}\\
 R_F^{B,v}
 &=1-\prod_{\substack{C\in\mathcal C_Z(F)\\C\cap B=\varnothing}}
 (1-J_C^{v}).
 \label{eq:degree-restricted-R}
\end{align}
\endgroup
\end{definition}

\emph{Remark.} The superscript $v$ means that $\chi_v$ is set to zero.
Then $J_C^v=1$ means that every atom of $C$ is accepted, and $R_A^v=1$
means that some block containing $A$ has all its atoms outside $A$
accepted. The products defining $\psi_u^{B,v}$ and $R_F^{B,v}$ use
only blocks meeting $B$ in $\{u\}$ and $F$, respectively.
Lemma~\ref{lem:degree-indicator-identities} uses these variables to
identify the residual degrees.

\begin{definition}\label{def:degree-sums}
Let $k\in\{2,\ldots,r\}$, $Z\subseteq X\setminus \mathsf{W}$ and $v\in X$.
We use the variables in Definition~\ref{def:degree-modified-variables}
for these $k,Z,v$. For $B\in\mathcal S_k(v)$, we define
$\widehat{\mathsf Z}^{\deg}_B,\mathsf Z^{\deg}_B$, and then
$F_k(v),\widehat F_k(v)$, by
\begingroup
\allowdisplaybreaks[4]
\begin{align}
 \widehat{\mathsf Z}^{\deg}_B
 &=\left(\prod_{u\in B\setminus\{v\}}\psi_u^{B,v}\right)
 R_B^{v}\prod_{\substack{F\subsetneq B\\|F|\ge2}}
                         (1-R_F^{B,v}),
 \label{eq:degree-hat-indicator}\\
 \mathsf Z^{\deg}_B
 &=\widehat{\mathsf Z}^{\deg}_B\prod_{u\in B\setminus\{v\}}(1-\mathsf{w}_u),
 \label{eq:degree-indicator}\\
 F_k(v)&=\sum_{B\in\mathcal S_k(v)}\mathsf Z_B^{\deg},\qquad
 \widehat F_k(v)=\sum_{B\in\mathcal S_k(v)}\widehat{\mathsf Z}_B^{\deg}.
 \notag
\end{align}
\endgroup
\end{definition}

\emph{Remark.} When $v$ survives, $\widehat{\mathsf Z}^{\deg}_B=1$ records
that $B$ is a residual block; $\mathsf Z^{\deg}_B=1$ also requires
every atom of $B\setminus\{v\}$ to lie outside $\mathsf W'$.
Lemma~\ref{lem:degree-indicator-identities} identifies their sums with
$d_k'(v,X')$ and $d_k'(v,X'\setminus\mathsf W')$.
Lemma~\ref{lem:general-degree-round} bounds deviations of these sums.

\begin{lemma}\label{lem:degree-indicator-identities}
Let $2\le k\le r$, $v\in X$ and $Z\subseteq X\setminus \mathsf{W}$.
\begin{enumerate}
\item Under $\mathbb P_{Z,v}^0$, for every $B\in\mathcal S_k(v)$,
      $\mathsf Z_B^{\deg}$ is determined by
      $(\chi_x:(1,x)\in T_B)$ and $(\mathsf{w}_x:(2,x)\in T_B)$, and
      $\widehat{\mathsf Z}_B^{\deg}$ is determined by
      $(\chi_x:(1,x)\in\widehat T_B)$.
      For every $j\in\{1,2\}$ and $x\in X$,
\begin{equation}\label{eq:degree-variable-use}
 |\{B\in\mathcal S_k(v):(j,x)\in\widehat T_B\}|
 \le|\{B\in\mathcal S_k(v):(j,x)\in T_B\}|\le\mu_k(v).
\end{equation}
\item Conditional on $\psi_v=1$, the laws $\mathbb P_Z$ and
      $\mathbb P_{Z,v}^0$ agree, and
\begin{equation}\label{eq:degree-residual-identities}
 F_k(v)=d_k'(v,X'\setminus \mathsf{W}'),\qquad \widehat F_k(v)=d_k'(v,X').
\end{equation}
\end{enumerate}
\end{lemma}
\begin{proof}[Proof of \textup{(1)}]
Equation~\eqref{eq:degree-atoms} gives $\mathcal A_B\subseteq Z\setminus B$.
For $x\in Z\setminus\{v\}$, let
\begin{equation}\label{eq:degree-variable-neighborhood}
 U_{Z,v}(x)=\{x\}\cup
 \bigcup_{\substack{x\in P\in\mathcal S\\P\subseteq Z\setminus\{v\}}}P.
\end{equation}
For $x\in\mathcal A_B$, $J_x^{v}$ is determined by
$(\chi_y:y\in U_{Z,v}(x))$ by \eqref{eq:degree-modified-J}.
Equations~\eqref{eq:degree-support}, \eqref{eq:degree-hat-support} and
\eqref{eq:degree-restricted-psi}--\eqref{eq:degree-indicator}
give the two asserted sets of variables.
Since $\widehat T_B\subseteq T_B$, \eqref{eq:degree-multiplicity}
proves \eqref{eq:degree-variable-use}.

\emph{Proof of \textup{(2)}.}
The event $\psi_v=1$ implies $\chi_v=0$, so
\eqref{eq:degree-conditioned-law} gives the equality of conditional laws.
On $\chi_v=0$, the expression
\begin{equation}\label{eq:degree-unrestricted-product}
 \left(\prod_{u\in B\setminus\{v\}}\psi_u(1-\mathsf{w}_u)\right)
 R_B^{v}\prod_{\substack{F\subsetneq B\\|F|\ge2}}(1-R_F^{v})
\end{equation}
equals $\mathsf Z_B^{\deg}$. Indeed, if $\chi_u=1$ for some
$u\in B\setminus\{v\}$, both expressions are zero. Otherwise
$J_x^{v}=0$ for every $x\in B$, so every factor omitted in
\eqref{eq:degree-restricted-psi} and \eqref{eq:degree-restricted-R}
is one. Replacing $\psi_u(1-\mathsf{w}_u)$ by $\psi_u$ in
\eqref{eq:degree-unrestricted-product} gives
$\widehat{\mathsf Z}_B^{\deg}$.
For such a $B$, $R_B^{v}=1$ precisely when $B=P\setminus I$
for some $P\in\mathcal S$. The factors $1-R_F^{v}$ exclude
smaller such sets $F\subsetneq B$ with $|F|\ge2$, and the factors
$\psi_u(1-\mathsf{w}_u)$ or $\psi_u$ require every $u\in B\setminus\{v\}$ to lie
in $X'\setminus \mathsf{W}'$ or $X'$, respectively.
Together with $\psi_v=1$, these are exactly the conditions for
$B$ to contribute to the corresponding degree in
\eqref{eq:degree-residual-identities}. Summing over $B$ proves
\eqref{eq:degree-residual-identities}.
\end{proof}

\begin{lemma}\label{lem:two-parent-count}
Let $2\le i,j\le r$, and let $A,T\subseteq X$ be nonempty sets
containing distinct atoms $v\in A$, $x\in T$. Let $\ell$ be an
integer with $|A\cap T|<\ell\le\min\{i,j\}$, and let
$f=i+j-\ell-|A\cup T|$. There are at most
$2(r!)^2\Sigma^{f-\delta_1}$ ordered pairs
$P\in\mathcal S_i$, $Q\in\mathcal S_j$ satisfying
$|P\cap Q|=\ell$, $P\cap(A\cup T)=A$ and $Q\cap(A\cup T)=T$.
\end{lemma}
\begin{proof}
A common atom lies outside $A\cup T$, so $|A|<i$ and $|T|<j$.
There are at most $2\Sigma^{i-|A|}$ choices for $P$;
when $|A|\ge2$, there are at most
$\Sigma^{i-|A|-\delta_5(i)}$. There are at most $2^r$ choices for
$P\cap Q$. The block $Q$ then has
$\iota=|T|+\ell-|A\cap T|\ge2$ prescribed atoms.
If $\iota<j$, \eqref{eq:local-absolute-codegrees} for $Q$ reduces
the exponent $(i-|A|)+(j-\iota)=f$ by $\delta_5(j)$.
If $\iota=j$, then $Q$ is determined; we use the bound for $P$
when $|A|\ge2$, or reverse the order when $|T|\ge2$.
For $A=\{v\}$ and $T=\{x\}$, either $\ell+1<j$ or
$\ell+1<i$ allows one of these two orders, or
$i=j=\ell+1$. In the latter case,
\eqref{eq:local-absolute-common} bounds the number of pairs by
$\Sigma^{f-\delta_4(i)}$.
We now use $2^r\le(r!)^2$ and $\delta_4(i),\delta_5(i)\ge\delta_1$.
\end{proof}

\begin{lemma}\label{lem:degree-indicator-read}
Let $2\le k\le r$ and $v\in X$, and suppose
\begin{equation}\label{eq:degree-read-scale}
 \Sigma\ge\max\left\{(8r)^{2^{10}r^2},
 \left(\frac{8r(2^{12}r^3+3)}{e\delta_2}\right)^{2/\delta_2}\right\}.
\end{equation}
Let $Z=\{x\in X:\xi_x=1\}$.
\begin{enumerate}
\item $\mathbb P\!\left(
 \max_{x\in Z\setminus\{v\}}
 |\{B\in\mathcal S_k(v):x\in\mathcal A_B\}|
 >e2^{2^{13}r^2+2}(r!)^2\Sigma^{k-1-\delta_1+2r\delta_2}
 \right)\le\Sigma^{-2^{10}r^3}$.
\item Let $\mathcal D_{v,k}$ be the event that
      $\mu_k(v)>b_{3}\Sigma^{k-1-\delta_1+3r\delta_2}$ and
      \eqref{eq:general-focused-degrees} holds for every
      $u\in\mathbb B_3(v,\mathcal S)$ and $2\le i\le r$,
      with $u$ in place of $v$. Then
\begin{equation*}
 \mathbb P(\mathcal D_{v,k})\le\Sigma^{-2^{10}r^3},
\end{equation*}
      and $\vb1_{\mathcal D_{v,k}}$ is a function of
      $(\xi_x:x\in\mathbb B_4(v,\mathcal S))$.
\end{enumerate}
\end{lemma}
\begin{proof}[Proof of \textup{(1)}]
We fix $x\in X\setminus\{v\}$. In each tuple below,
$B\subseteq P\in\mathcal S$, $v\in B$ and $|B|=k$.
In cases \textup{(2)} and \textup{(3)}, we also require
$Q\in\mathcal S$ and $x\in Q\setminus B$.
\begin{enumerate}
\item $\mathbf Q=(P,B)$, where $x\in P\setminus B$.
\item $\mathbf Q=(P,Q,B,u)$, where $u\in B\setminus\{v\}$ and $Q\cap B=\{u\}$.
\item $\mathbf Q=(P,Q,B,F)$, where $F=Q\cap B$ and $2\le|F|<k$.
\end{enumerate}
For every such tuple, we define
\[
 R(\mathbf Q)=
 \begin{cases}
 P\setminus(B\cup\{x\}),&\text{in case \textup{(1)}},\\
 (P\cup Q)\setminus(B\cup\{x\}),&\text{in cases \textup{(2)} and \textup{(3)}}.
 \end{cases}
\]
We write $R=R(\mathbf Q)$ when the tuple is fixed.
By \eqref{eq:degree-atoms}, each membership $x\in\mathcal A_B$
is represented by one of these tuples with $R\subseteq Z$.
We partition the tuples by their case, their block sizes, and
$\mathfrak d=|R(\mathbf Q)|\le2r$. For each resulting family
$\mathcal Q$, we define
\begin{equation}\label{eq:degree-occurrence-polynomial}
 \Phi=\sum_{\mathbf Q\in\mathcal Q}
             \prod_{y\in R(\mathbf Q)}\xi_y.
\end{equation}
For $x\in Z\setminus\{v\}$, the sum of these polynomials bounds
$|\{B\in\mathcal S_k(v):x\in\mathcal A_B\}|$.

For $0\le s\le\mathfrak d$, we fix $s$ atoms of $R$.
In case \textup{(1)}, with $|P|=i$ and $\mathfrak d=i-k-1$,
\eqref{eq:local-absolute-codegrees} bounds the number of tuples by
$2^{2^{12}r^2+1}(r!)^2\Sigma^{k-1+\mathfrak d-s-\delta_5(i)}$.
In cases \textup{(2)} and \textup{(3)}, we assign each fixed atom to
$P$, $Q$ or both, and apply Lemma~\ref{lem:two-parent-count}.
The shared atom $u$ in case \textup{(2)}, or an atom of
$F\setminus\{v\}$ in case \textup{(3)}, lies outside
$R\cup\{v,x\}$ and is not fixed by the derivative.
Since $|P\cup Q|=k+\mathfrak d+1$, the exponent in
that lemma is $k-1+\mathfrak d-s-\delta_1$.
There are at most $2^{2^{12}r^2}$ choices for these assignments
and for the subsets of the blocks. Consequently, for every
$T\subseteq X$ with $|T|=s$,
\begin{equation}\label{eq:degree-occurrence-derivatives}
 \mathbb E(\partial_T\Phi)
 \le2^{2^{12}r^2+1}(r!)^2
       \Sigma^{k-1-\delta_1+(\mathfrak d-s)\delta_2}.
\end{equation}
Each atom in $R(\mathbf Q)$ occurs once in the product
\eqref{eq:degree-occurrence-polynomial}.
For $\mathfrak d=0$, \eqref{eq:degree-occurrence-derivatives}
is a deterministic bound. For $\mathfrak d>0$,
Corollary~\ref{cor:polynomial-envelope}~(2), with $a=2^{12}r^3$, gives
\begin{equation}\label{eq:degree-occurrence-tail}
 \mathbb P\!\left(\Phi>
 e2^{2^{12}r^2+2}(r!)^2
       \Sigma^{k-1-\delta_1+\mathfrak d\delta_2}\right)
 \le\Sigma^{-2^{12}r^3}.
\end{equation}
Condition~\eqref{eq:degree-read-scale} supplies the hypothesis of
the corollary for $\mathfrak d\le2r$.
There are at most $2^{2^{12}r^2}$ families $\mathcal Q$ for each
$x$, and every relevant $x$ lies in $\mathbb B_2(v,\mathcal S)$.
The inequalities $d_i(u,X)\le2\Sigma^{i-1}$ give
$|\mathbb B_2(v,\mathcal S)|\le\Sigma^{2r-1}$ by counting paths
of length at most two. Since $2^{2^{12}r^2}\le\Sigma^2$, a union
bound over at most $\Sigma^{2r+3}$ polynomials in
\eqref{eq:degree-occurrence-tail} gives
\begin{equation}\label{eq:degree-occurrence-union-tail}
 \mathbb P\!\left(
 \max_{x\in Z\setminus\{v\}}
 |\{B\in\mathcal S_k(v):x\in\mathcal A_B\}|
 >e2^{2^{13}r^2+2}(r!)^2\Sigma^{k-1-\delta_1+2r\delta_2}
 \right)\le\Sigma^{2r+3-2^{12}r^3}\le\Sigma^{-2^{10}r^3}.
\end{equation}
This proves part~\textup{(1)}.

\emph{Proof of \textup{(2)}.}
Assume \eqref{eq:general-focused-degrees} holds for every
$u\in\mathbb B_3(v,\mathcal S)$ and $2\le i\le r$.
For each such $u$, \eqref{eq:degree-variable-neighborhood} gives
\begin{equation}\label{eq:degree-J-variable-count}
 |\{x\in Z\setminus\{v\}:u\in U_{Z,v}(x)\}|
 \le1+\sum_{i=2}^r(i-1)d_i(u,Z)
 \le4r^2\Sigma^{(r-1)\delta_2}.
\end{equation}
For $k=2$ and $u\ne v$, there is at most one
$B\in\mathcal S_k(v)$ containing $u$.
For $k\ge3$, fixed $k\le i\le r$ and $u\ne v$, the polynomial
\begin{equation}\label{eq:degree-direct-occurrence-polynomial}
 \sum_{\substack{\{v,u\}\subseteq B\subseteq P\in\mathcal S_i\\|B|=k}}
       \prod_{y\in P\setminus B}\xi_y
\end{equation}
counts the pairs $B\subseteq P\in\mathcal S_i$ with
$\{v,u\}\subseteq B$, $|B|=k$ and $P\setminus B\subseteq Z$.
We apply \eqref{eq:local-absolute-codegrees} and, when $i>k$,
Corollary~\ref{cor:polynomial-envelope}~(2) to
\eqref{eq:degree-direct-occurrence-polynomial}.
A union bound over $i$ and $u$ gives
\begin{equation}\label{eq:degree-direct-occurrence-tail}
 \mathbb P\!\left(
 \max_{u\in\mathbb B_1(v,\mathcal S)\setminus\{v\}}
 |\{B\in\mathcal S_k(v):u\in B\}|>
 e2^{2^{13}r^2+2}(r!)^2\Sigma^{k-2-\delta_1+r\delta_2}
 \right)\le\Sigma^{2r+3-2^{12}r^3}.
\end{equation}
Outside the events in \eqref{eq:degree-occurrence-union-tail} and
\eqref{eq:degree-direct-occurrence-tail},
\eqref{eq:degree-J-variable-count}, \eqref{eq:degree-support} and
\eqref{eq:degree-multiplicity} give \eqref{eq:general-degree-read}: we multiply the bound on
$|\{B:x\in\mathcal A_B\}|$ by $4r^2\Sigma^{(r-1)\delta_2}$
and add the bound on $|\{B:u\in B\}|$.
For $k=2$, we use $1$ for this last bound. Hence
\begin{equation}\label{eq:degree-read-union-tail}
 \mathbb P(\mathcal D_{v,k})
 \le2\Sigma^{2r+3-2^{12}r^3}\le\Sigma^{-2^{10}r^3}.
\end{equation}
Every $B\in\mathcal S_k(v)$ lies in $\mathbb B_1(v,\mathcal S)$,
and \eqref{eq:degree-atoms} gives
$\mathcal A_B\subseteq\mathbb B_2(v,\mathcal S)$.
Equations~\eqref{eq:degree-candidate-sets}, \eqref{eq:degree-support}
and \eqref{eq:degree-multiplicity} therefore determine $\mu_k(v)$ from
$(\xi_x:x\in\mathbb B_3(v,\mathcal S))$.
For $u\in\mathbb B_3(v,\mathcal S)$, the variables in $d_i(u,Z)$
are indexed by $\mathbb B_4(v,\mathcal S)$.
This proves the claimed dependence of $\vb1_{\mathcal D_{v,k}}$.
\end{proof}

\subsection{\texorpdfstring{Concentration of $|C_x\cap Z|$}{\detokenize{Concentration of the sets C_x intersect Z}}}
\label{subsec:focus-coverage}

For each choice of $i,A$ or $i,x,y$ in
Definition~\ref{def:good-focus}~\eqref{gf:coverage} with $\mathcal B\ne\varnothing$,
and each $Z\subseteq U$, we use the associated
$\mathcal B,F_0,\ell,C_u,C_B,\mathcal T_B$. We put
\begin{equation*}
 k^{\rm cov}=\max\left\{1,\max_{u\in Z}
          |\{B\in\mathcal B:u\in\mathcal T_B\}|\right\}.
\end{equation*}
This bounds the number of sets $\mathcal T_B$ containing any one atom;
the letter $k$ matches the parameter in Proposition~\ref{prop:readk}.
We define
$h,m,\zeta,m_0$ and $r(w)$ for $w\in U$ by
\begin{equation}\label{eq:coverage-proof-parameters}
 \begin{aligned}
 h&=|\mathcal B|,\qquad
 m_0=\max_{v\in X}|\{B\in\mathcal B:v\in B\}|,\\
 (m,\zeta)&=\begin{cases}
 (i-|A|,\delta_5(i)),&\text{for the choice }i,A,\\
 (i-1,\delta_4(i)),&\text{for the choice }i,x,y,
 \end{cases}\\
 r(w)&=\sum_{B\in\mathcal B}|B\cap N_2(w)|.
 \end{aligned}
\end{equation}

For $2\le i\le r$ and $v\in X$, assume
$d_i(v,X)\le b_{1}\Sigma^{i-1}$ and $E_i\le\Sigma^{i-1-\delta_3}$.
Also assume $\frac12\le\sigma\le1$ and
$\Sigma^{-\delta_3/2}\le\eta\le1/4$.
For $A\subseteq X$ with $2\le|A|<i\le r$, and for
$2\le i\le r$ and distinct $x,y\in X$, respectively, assume
\begin{align}
 d_i(A)&\le\Sigma^{i-|A|-\delta_5(i)},
 \label{eq:general-codegree-barriers}\\
 d_i(x,y)&\le\Sigma^{i-1-\delta_4(i)}.
 \label{eq:general-common-barriers}
\end{align}
\begin{definition}\label{def:coverage-random-variables}
Let $i,A$ or $i,x,y$ be as in Definition~\ref{def:good-focus}~\eqref{gf:coverage},
with $\mathcal B\ne\varnothing$, and let $Z\subseteq U$.
We use the associated $\mathcal B,F_0$ and the sets $C_B$ in
\eqref{eq:coverage-candidate-definitions} for this choice.
Let $(\chi_x:x\in Z)$ be independent Bernoulli$(p)$ random variables.
For $B\in\mathcal B$, we define
\begin{equation}\label{eq:coverage-variable-formula}
 \mathsf Z_B^{\rm cov}
 =1-\prod_{w\in C_B\cap Z}
 \left(1-\chi_w
 \prod_{\substack{w\in S\in\mathcal S,\ S\subseteq Z\\
                         S\cap(F_0\cup B)=\varnothing}}
 \left(1-\prod_{u\in S\setminus\{w\}}\chi_u\right)\right).
\end{equation}
\end{definition}

\emph{Remark.} The equality $\mathsf Z_B^{\rm cov}=1$ means that an atom of $C_B\cap Z$
is accepted after setting $\chi_x=0$ for $x\in(F_0\cup B)\cap Z$.
Formula~\eqref{eq:coverage-variable-formula} uses only the variables
$(\chi_x:x\in \mathcal T_B)$. Lemma~\ref{lem:general-coverage}
bounds how many $B\in\mathcal B$ have this property. If the atoms
of $F_0$ survive, each such $B$ loses an atom.

Assume \eqref{eq:coverage-family-activation} and $\Sigma\ge b_{5}$
for the choice in Definition~\ref{def:coverage-random-variables},
with the notation in \eqref{eq:coverage-proof-parameters}.
For $x\in\bigcup\mathcal B$, double counting gives the following bound,
with the summand for
$y=x$ interpreted as $d_2(x,X\setminus \mathsf{W})$:
\begin{equation}\label{eq:general-trimming}
 \sum_{w\in N_2(x,X\setminus \mathsf{W})}r(w) =\sum_{B\in\mathcal B}\sum_{y\in B}d_2(x,y,X\setminus \mathsf{W}) \le b_{1}\Sigma m_0+mh\Sigma^{1-\delta_1}.
\end{equation}
By \eqref{eq:general-trimming}, discarding neighbors with
$r(w)>\ell$ removes at most $\frac{\delta_3}{b_{3}}d_2$ choices.
Since $|F_0|\le r$, $E_2\le\Sigma^{1-\delta_3}$ and
$\Sigma\ge b_{5}$, deleting $F_0$ and allowing for the degree error
reduces the lower bound by at most another $\frac{\delta_3}{b_{3}}d_2$. Thus
\begin{equation*}
 (1-\frac{2\delta_3}{b_{3}})d_2\le|C_x|\le(1+\frac{\delta_3}{2b_{3}})d_2.
\end{equation*}

\begin{lemma}\label{lem:focus-coverage}
Let $i,A$ or $i,x,y$ be as in Definition~\ref{def:good-focus}~\eqref{gf:coverage}
and satisfy \eqref{eq:coverage-family-activation}.
We use the associated $\mathcal B,F_0$ and the notation in
\eqref{eq:coverage-proof-parameters}.
Assume \eqref{eq:general-codegree-barriers}--\eqref{eq:general-common-barriers},
$\Sigma\ge b_{5}$, $\frac12\le\sigma\le1$ and
$\Sigma^{-\delta_3/2}\le\eta\le \frac{\delta_3}{b_{3}}$.
For $2\le i\le r$ and $v\in X$, assume
$d_i(v,X)\le b_{1}\Sigma^{i-1}$ and $E_i\le\Sigma^{i-1-\delta_3}$.
Then the following hold.
\begin{enumerate}
\item For every $x\in\bigcup\mathcal B$,
\begin{align}
 \mathbb P\!\left(\bigl||C_x\cap Z|-q|C_x|\bigr|>\frac{\delta_3}{b_{3}}q|C_x|\right)
 &\le\Sigma^{-a},\notag\\
 \mathbb P\!\left(k^{\rm cov}>rb_{1}b_{3}^2h\Sigma^{2\delta_2-\delta_1}\right)
 &\le\Sigma^{-a}.\label{eq:coverage-read-probability}
\end{align}
\item If \eqref{eq:general-constant-focused-completions} and
\eqref{eq:coverage-focus-candidates} hold, then for every $B\in\mathcal B$,
\begin{equation}\label{eq:coverage-focused-candidate-size}
 (m-\frac{3r\delta_3}{b_{3}}-2r^2(b_{2}+1)\Sigma^{-\delta_2})qd_2
 \le |C_B\cap Z|\le2rqd_2.
\end{equation}
\end{enumerate}
\end{lemma}
\begin{proof}
We put $M=a+2^{10}r^3$. Since $\Sigma\ge b_{5}$,
$\frac{\delta_3^2}{8b_{3}^2}\Sigma^{\delta_2}\ge(M+1)\log\Sigma$.
Proposition~\ref{prop:bernstein} gives
\begin{equation}\label{eq:coverage-candidate-bernstein}
 \begin{aligned}
 \mathbb P\!\left(\bigl||C_x\cap Z|-q|C_x|\bigr|>\frac{\delta_3}{b_{3}}q|C_x|\right)
 &\le2\exp\left(-\frac{\delta_3^2q|C_x|}{2b_{3}(b_{3}+\delta_3)}\right)\\
 &\le2\exp(-\frac{\delta_3^2\Sigma^{\delta_2}}{8b_{3}^2})
 \le2\Sigma^{-M-1}\le\Sigma^{-M}.
 \end{aligned}
\end{equation}
By \eqref{eq:coverage-focus-candidates}, $|C_x\cap Z|$ lies between
$(1-\frac{3\delta_3}{b_{3}})qd_2$ and $2qd_2$.
Equation~\eqref{eq:general-constant-focused-completions} with $i=2$
bounds pairwise intersections by $b_{2}$; deleting $B$ removes at most $r$ atoms. This proves \eqref{eq:coverage-focused-candidate-size}.

\emph{Proof of \textup{\eqref{eq:coverage-read-probability}}.}
Under \eqref{eq:coverage-family-activation}, for
$u\in X\setminus(\mathsf{W}\cup F_0)$ and $w\in X\setminus \mathsf{W}$, we put
\[
 m_u(w)=|\{B\in\mathcal B:u\notin B,\ w\in C_B\}|.
\]
For $2\le s\le r$, we put
\[
 \begin{aligned}
 \mathcal M_{u,s}&=
 \left\{f\in N_s(u)\cap\binom{X\setminus \mathsf{W}}{s-1}:
                   \sum_{w\in f}m_u(w)>0\right\},\\
 {\mathsf Q}_{u,s}(\xi)&=
 \sum_{f\in\mathcal M_{u,s}}
       \left(\sum_{w\in f}m_u(w)\right)\prod_{x\in f}\xi_x.
 \end{aligned}
\]
If $\mathcal M_{u,s}=\varnothing$, then $\mathsf Q_{u,s}=0$.

The restriction $r(w)\le\ell$ gives $m_u(w)\le\ell$. Since $u\notin B$ in this weight,
there is no diagonal term, and
\[
 \mathbb E(\mathsf Q_{u,2}) \le q\sum_{\substack{B\in\mathcal B\\u\notin B}} \sum_{x\in B}d_2(x,u) \le rh\Sigma^{\delta_2-\delta_1}.
\]
For $\lambda=(\delta_2/2)\log\Sigma$, independence and
$0\le m_u(w)\le\ell$ give
\[
 \mathbb E(\exp(\lambda\mathsf Q_{u,2}/\ell))
 \le\exp\bigl((e^\lambda-1)\mathbb E(\mathsf Q_{u,2})/\ell\bigr).
\]
Since $b_{2}\lambda\ge2(M+1)\log\Sigma\ge1$, this gives
\[
 \mathbb P\bigl(\mathsf Q_{u,2}>
 b_{2}(\ell+rh\Sigma^{2\delta_2-\delta_1})\bigr)
 \le e^{-b_{2}\lambda}\le\Sigma^{-M}.
\]

For $s\ge3$, we set $\Phi=\mathsf Q_{u,s}/((s-1)\ell)$.
Its coefficients are at most one. For a nonzero derivative with
$1\le|T|=\ell_0<s-1$,
\[
 \mathbb E(\partial_T\Phi) \le q^{s-1-\ell_0}d_s(\{u\}\cup T) \le\Sigma^{(s-1-\ell_0)\delta_2-\delta_5(s)}.
\]
Also
\[
 \sum_f\sum_{w\in f}m_u(w) \le\Sigma^{s-2-\delta_5(s)} \sum_{\substack{B\in\mathcal B\\u\notin B}}|C_B| \le2rh\Sigma^{s-1-\delta_5(s)}.
\]
Since $\ell\ge8rh\Sigma^{-\delta_1}$, the mean of $\Phi$ is at
most $\Sigma^{(s-1)\delta_2-\delta_5(s)+\delta_1}$.
Thus all proper expected derivatives of $\Phi$, including its mean,
are at most $\Sigma^{-(\delta_1-r\delta_2)}$.
Corollary~\ref{cor:polynomial-envelope}~(3), with $b_{2}$ from
Lemma~\ref{lem:general-constant-focus} and exponent $M$, gives
\begin{equation}\label{eq:coverage-weighted-polynomial-probability}
 \mathbb P(\Phi>b_{2})\le\Sigma^{-M}.
\end{equation}
For these events, we apply Lemma~\ref{lem:conditioning-locality}
with parameters
\begin{equation}\label{eq:indicator-locality-parameters}
 (h,\lambda,\ell,s,m,\mathfrak d)
 =\bigl(\Sigma,b_{1},16r,16r,2^{2^{12}r^2},2r\bigr).
\end{equation}
The bound $\Sigma\ge(4rb_{1})^{2^{10}r^2}$ implies
$4[b_{1}(r-1)^2]^{32r}\le\Sigma$,
$(16r+1)2^{2^{12}r^2}\le\Sigma^2$ and $2r+1\le\Sigma$.
The exponents in parts (2) and (3) of that lemma satisfy
\[
 (256r^2+32r)(r-1)+16r+3\le2^9r^3,\qquad 2r[16r(r-1)+1]+1\le2^9r^3.
\]
Thus both conclusions give the bound $\Sigma^{2^9r^3}$.

For $u\in Z\setminus F_0$, the number of sets $B$ with $u\in C_B$
is at most $\ell$. Each further occurrence of $u$ in $\mathcal T_B$
comes from a block $S\subseteq Z$ containing $u$ and an atom of $C_B$;
it contributes to $\mathsf Q_{u,|S|}$.
Thus
\[
 |\{B:u\in \mathcal T_B\}|
 \le\ell+\sum_{s=2}^r\mathsf Q_{u,s}.
\]
For $u\in F_0$, no set $\mathcal T_B$ contains $u$. Using $\ell\ge1$,
$\mathsf Q_{u,2}\le b_{2}(\ell+rh\Sigma^{2\delta_2-\delta_1})$
and, for $3\le s\le r$, $\mathsf Q_{u,s}\le(s-1)\ell b_{2}$ gives
\[
 k^{\rm cov}\le
 (1+b_{2}\tbinom r2)\ell
 +rb_{2}h\Sigma^{2\delta_2-\delta_1}.
\]
Dividing by $h$ and substituting $\ell$ from
\eqref{eq:coverage-candidate-definitions} gives
\[
 \frac{k^{\rm cov}}h
 \le b_{1}b_{3}^2\left(\frac{m_0}h+\Sigma^{2\delta_2-\delta_1}\right).
\]
By \eqref{eq:general-codegree-barriers} and
\eqref{eq:coverage-family-activation}, proper codegrees with
$m\ge2$ satisfy $m_0/h\le2\Sigma^{-1}$; for $m=1$, multiplicity
one gives $m_0/h\le2\Sigma^{-1+\delta_5(i)}$.
For common links with $i\ge3$,
\eqref{eq:general-codegree-barriers} and
\eqref{eq:coverage-family-activation}, applied to
$d_i(\{x,u\})$, give $m_0/h\le2\Sigma^{-1-\delta_5(i)+\delta_4(i)}$;
for $i=2$, multiplicity one gives $m_0/h\le2\Sigma^{-1+\delta_1}$.
Each bound is at most $2\Sigma^{-3/4}\le\Sigma^{-1/2}$,
proving \eqref{eq:coverage-focus-read}.
For each family, Lemma~\ref{lem:conditioning-locality}~(2), with
\eqref{eq:indicator-locality-parameters}, bounds the number of failure
events in the union over $u,s$ by
$4\Sigma^{2^9r^3}$. Since $M=a+2^{10}r^3$, their union has
probability at most
$4\Sigma^{2^9r^3-M}\le\Sigma^{-a}$.
This proves \eqref{eq:coverage-read-probability}.
\end{proof}

\subsection{\texorpdfstring{Polynomials for $d_i^{\prime}(A)$ and $d_i^{\prime}(x,y)$}{\detokenize{Polynomials for d_i'(A) and d_i'(x,y)}}}
\label{subsec:focus-births}

\begin{lemma}\label{lem:focus-birth-derivatives}
Let $b_{1}\ge1$ and $a>0$. Assume
\eqref{eq:local-absolute-codegrees}--\eqref{eq:local-absolute-common}
and $d_i(v,X)\le b_{1}\Sigma^{i-1}$ for $2\le i\le r$ and $v\in X$.
Let $\mathsf Y$ be either polynomial in
\eqref{eq:proper-birth-witness-polynomial} or
\eqref{eq:common-birth-witness-polynomial}, with its indices in the
stated ranges and parameters $s,m,\iota$ from
\eqref{eq:birth-polynomial-parameters}.
For every $T\subseteq U$ with $|T|\le s$,
\begin{equation}\label{eq:birth-derivative-mean}
 \mathbb E(\partial_T\mathsf Y)
 \le b_{1}^2b_{3}\iota\Sigma^{(s-|T|)\delta_2}.
\end{equation}
If $\Sigma\ge(4rb_{1})^{2^{10}r^2}$ and $0<\eta\le1/4$, then
\begin{equation}\label{eq:birth-derivative-tail}
 \mathbb P\left((\partial_T\mathsf Y)(Z)
 >b_{1}^2b_{3}^2\iota\Sigma^{(s-|T|)\delta_2}
             (\log\Sigma)^{s-|T|}\right)
 \le\Sigma^{-(a+2^{11}r^3)}.
\end{equation}
\end{lemma}
\begin{proof}
For $\mathsf Y_{i,A;k}$, we put $j=|A|$, $m=i-j$ and $s=k-i$.
Each summand chooses $A\subseteq P\in\mathcal S_k$ and
$F\subseteq P\setminus A$ with $|F|=s$. For a nonzero derivative with
$T\subseteq U$ and $\ell_0=|T|\le s$,
\eqref{eq:local-absolute-codegrees} bounds the coefficient sum by
\[
 \binom{k-j-\ell_0}{s-\ell_0}d_k(A\cup T)
 \le b_{1}^2b_{3}\Sigma^{k-j-\ell_0-\delta_5(k)}.
\]
The set $A\cup T$ is proper because $m\ge1$ atoms remain free.
Multiplying by $q^{s-\ell_0}$ gives the expected-derivative bound
with $\iota=\Sigma^{m-\delta_5(k)}$.

For $\mathsf Y_{i,x,y;k,l,h}$, let $(S,R,P_x,P_y)$ be an ordered
quadruple summed over in \eqref{eq:common-birth-witness-polynomial},
with $P=P_x$ and $Q=P_y$.
We put $h_0=|R|$, $h_1=|P_x|$, $h_2=|P_y|$ and
$s=h_0+h_1+h_2\ge1$. For $T\subseteq U$ with $\ell_0=|T|\le s$,
we partition $T$ among $R,P_x,P_y$ in at most $3^{\ell_0}$ ways, with
$h_f=|T\cap R|$, $\ell_1=|T\cap P_x|$ and
$\ell_2=|T\cap P_y|$.
There are at most $\Sigma^{l-1-h_f-\ell_2-\delta_5(l)}$ choices
for the block through $y$ if $h_f+\ell_2>0$, and at most
$b_1\Sigma^{l-1}$ otherwise. We choose $S,R$ inside it in
at most $3^r$ ways. If $h_1-\ell_1>0$, there are at most
$\Sigma^{h_1-\ell_1-\delta_5(k)}$ completions of the block through $x$;
otherwise that block is determined.
Apart from $-\delta_5(k),-\delta_5(l)$, the exponents sum to $i-1+s-\ell_0$.

Equation~\eqref{eq:local-absolute-codegrees} applies unless
$h_1=\ell_1$ and $h_f=\ell_2=0$. In the remaining cases:
\begin{enumerate}
\item If $h_1>0$, we reverse the order: the proper prescribed set
$\{x\}\cup(P_x\cap T)$ gives at most
$\Sigma^{i+h_0-1-\delta_5(k)}$ choices, and there are at most
$\Sigma^{h_2}$ completions of the block through $y$.
\item If $h_1=0<h_2$, then $\ell_0=0$. We choose the block through $x$ in at
most $b_{1}\Sigma^{i+h_0-1}$ ways; the block through $y$ then has a
proper prescribed set and at most $\Sigma^{h_2-\delta_5(l)}$ completions.
\item If $h_1=h_2=\ell_0=0$, then $h_0\ge1$ and $k=l=i+h_0$.
The blocks share their whole link; \eqref{eq:local-absolute-common}
gives at most
$\binom{k-1}{h_0}d_k(x,y)
\le b_{1}^2b_{3}\Sigma^{i-1+s-\delta_4(i+1)}$ choices.
\end{enumerate}
Since every $\delta_5(j)>\delta_4(r)$, all cases give
coefficient sum at most
$b_{1}^2b_{3}\Sigma^{i-1+s-\ell_0-\delta_4(i+1)}$.
Multiplying by $q^{s-\ell_0}$ gives the expected-derivative bound
with $\iota=\Sigma^{i-1-\delta_4(i+1)}$.

If $|T|=s$, the derivative is deterministic and
\eqref{eq:birth-derivative-mean} makes the event in
\eqref{eq:birth-derivative-tail} empty.
For $|T|<s$, we apply Proposition~\ref{prop:positive-moments} to
$\partial_T\mathsf Y$, of degree $s-|T|$, with target $a+2^{11}r^3$. For $T'\supseteq T$,
$\Sigma^{(s-|T'|)\delta_2}\le\Sigma^{(s-|T|)\delta_2}$,
which gives
\eqref{eq:birth-derivative-tail}.
For fixed $A$ or $x,y$, to obtain \eqref{eq:general-birth-focus-derivatives}
for all the corresponding polynomials simultaneously, there are at most $2^{2^{12}r^2}$ choices
of $k,l$ and $h=|R|$, and at most
$\Sigma^{2^9r^3}$ derivatives for each by
Lemma~\ref{lem:conditioning-locality}~(3), with
\eqref{eq:indicator-locality-parameters} and $\mathcal E=\varnothing$.
Their variables are indexed
by atoms at distance at most two from a prescribed atom in the shadow graph.
Since $\Sigma\ge(4rb_{1})^{2^{10}r^2}$, we have
$2^{2^{12}r^2}\le\Sigma^2$. Each event has probability at most
$\Sigma^{-(a+2^{11}r^3)}$, so their union has probability at most
\[
 2^{2^{12}r^2}\Sigma^{2^9r^3-a-2^{11}r^3}
 \le\Sigma^{-a}.
\]
This completes the proof of Lemma~\ref{lem:focus-birth-derivatives}.
\end{proof}

\subsection{\texorpdfstring{Bounds for $(\partial_T\mathsf Y)(Z)$ when $t<T_{\mathrm{init}}$}{\detokenize{Bounds for (partial_T Y)(Z) when t < T_init}}}
\label{subsec:focus-initial-period}

We set
\[
 b_{6}=\exp\!\left(
 \frac{2^{16r^2}(a+1)(1+\log b_{1})}{\delta_0}\right).
\]

\begin{lemma}\label{lem:initial-focus-estimates}
Let $T_0\ge0$ be an integer. Assume the initial bounds
\eqref{eq:initialP}--\eqref{eq:initialCommon}, $D\ge b_{6}$, and
\[
 2b_{4}^{-(r-1)}k_{T_0}\le D^{\delta_0/2}.
\]
Let $0\le t\le T_0$. Assume $b_{4}D\le\Sigma_t\le 2^{4r^2}D$,
and, for $2\le i\le r$ and $v\in X^t$,
$d_i(v,X^t)\le b_{1}D^{i-1}$.
For $3\le i\le r$ and $A\subseteq X^t$ with $2\le|A|<i$,
and for $2\le i\le r$ and distinct $x,y\in X^t$, respectively,
assume
\begin{equation}\label{eq:initial-current-codegree-bounds}
 d_i(A)\le k_tD^{i-|A|-\delta_0},\qquad d_i(x,y)\le k_tD^{i-1-\delta_0}.
\end{equation}
Then \eqref{eq:local-absolute-codegrees}--\eqref{eq:local-absolute-common}
hold. For $t<T_0$, let $\mathsf Y$ be either polynomial in
\eqref{eq:proper-birth-witness-polynomial} or
\eqref{eq:common-birth-witness-polynomial}, with its indices in the
stated ranges and parameters $s,m,\iota$ from
\eqref{eq:birth-polynomial-parameters}.
For every $T\subseteq U$ with $|T|\le s$,
\begin{equation}\label{eq:initial-derivative-tail}
 \mathbb P\left((\partial_T\mathsf Y)(Z)
 >b_{1}^2b_{3}^2b_{4}^{-(s-|T|)}k_t
 D^{m-\delta_0+(s-|T|)\delta_2}(\log D)^{s-|T|}\right)
 \le D^{-(a+2^{11}r^3)}.
\end{equation}
For $3\le i\le r$ and $A\subseteq X^t$ with $2\le|A|<i$,
\begin{equation}\label{eq:initial-proper-completion-tail}
 \mathbb P\bigl(d_i(A,Z)>b_{2}\bigr)\le\Sigma_t^{-a}.
\end{equation}
For $2\le i\le r$ and distinct $x,y\in X^t$,
\begin{equation}\label{eq:initial-common-completion-tail}
 \mathbb P\bigl(d_i(x,y,Z)>b_{2}\bigr)\le\Sigma_t^{-a}.
\end{equation}
\end{lemma}
\begin{proof}
For $t<T_0$, \eqref{eq:initial-current-codegree-bounds} and
$k_t\le D^{\delta_0/2}$ give, in the stated ranges,
\[
 d_i(A)\le D^{i-|A|},\qquad d_i(x,y)\le D^{i-1}.
\]
In the proof of Lemma~\ref{lem:focus-birth-derivatives}, we use
$k_tD^{\iota-\delta_0}$ for one proper codegree or common link,
where $\iota=i-|A|$ or $i-1$, respectively, and $D^\iota$ for
the others, including cases~(1)--(3).
We use $b_{1}D^\iota$ for degrees and one for a determined block.
Consequently, prescribing $\ell_0$ atoms of $F$ or $R\cup P\cup Q$,
of total size $s\le2r$, leaves total coefficient at most
\[
 b_{1}^2b_{3}k_tD^{m+s-\ell_0-\delta_0}.
\]
There is only one factor $k_t$. Since $0<b_{4}\le1$ and
$0<\delta_2<1$, we have $q_t\le b_{4}^{-1}D^{-1+\delta_2}$.
The expected derivative
of the corresponding polynomial $\mathsf Y$ is at most
\[
 b_{1}^2b_{3}b_{4}^{-(s-\ell_0)}k_t
 D^{m-\delta_0+(s-\ell_0)\delta_2}.
\]
If $\ell_0=s$, the derivative is deterministic and the event in
\eqref{eq:initial-derivative-tail} is empty.
For $\ell_0<s$, Proposition~\ref{prop:positive-moments}, with target
$a+2^{11}r^3$, gives \eqref{eq:initial-derivative-tail}.
Both $b_{4}^{-(s-\ell_0)}$ and $D^{(s-\ell_0)\delta_2}$
decrease as $\ell_0$ increases.
The local union of derivative failures has probability at most
$D^{-a-3}\le\Sigma_t^{-a}$, since
$D\ge 2^{2ar^2}$. Lemma~\ref{lem:conditioning-locality}~(3),
with \eqref{eq:indicator-locality-parameters} and
$\mathcal E=\varnothing$, applies with $h=D$, $\lambda=b_{1}$,
and with $h=\Sigma_t$, $\lambda=b_{1}b_{4}^{-(r-1)}$:
$D\ge(4rb_{1})^{2^{10}r^2}$ and
$\Sigma_t\ge b_{4}D\ge(4rb_{1}b_{4}^{-(r-1)})^{2^{10}r^2}$.

Since $D\ge b_{6}$,
$b_{4}^{-(r-1)}\le D^{\delta_0/4-(r-1)\delta_2}$.
For either polynomial in \eqref{eq:focus-completion-polynomials},
let $\iota=i-|A|$ or $\iota=i-1$, respectively. Its mean and the
mean of every proper coefficient derivative are at most
\[
 b_{4}^{-\iota} k_tD^{-\delta_0+\iota\delta_2}
 \le b_{4}^{-\iota}D^{-\delta_0/2+\iota\delta_2}
 \le D^{-\delta_0/4}.
\]
Its coefficients are at most one. Since
$\delta_1-r\delta_2<\delta_0/4$,
Corollary~\ref{cor:polynomial-envelope}~(3), with $b_{2}$ from
Lemma~\ref{lem:general-constant-focus} and target $a+2^{10}r^3$,
gives \eqref{eq:initial-proper-completion-tail} and
\eqref{eq:initial-common-completion-tail}, since
$D^{-a-3}\le\Sigma_t^{-a}$.
Finally, for $t\le T_0$, we put $(m,\zeta)=(i-|A|,\delta_5(i))$
for $d_i(A)$ and $(m,\zeta)=(i-1,\delta_4(i))$ for $d_i(x,y)$.
In both cases $0<\zeta<\delta_0/2$ and $0<m-\zeta\le r-1$, so
\[
 \frac{k_tD^{m-\delta_0}}{\Sigma_t^{m-\zeta}}
 \le b_{4}^{-(r-1)}k_{T_0}D^{-\delta_0/2}\le\frac12.
\]
This proves \eqref{eq:local-absolute-codegrees}--\eqref{eq:local-absolute-common}.
\end{proof}

\subsection{Existence of a good focus}
\label{subsec:choosing-good-focus}

For the process, we set
\[
 b_{1}=2^{4r^3},\qquad a=2^{10}r^3.
\]

\begin{lemma}\label{lem:good-focus-existence}
Let $t\in\mathbb Z_{\ge0}$. Assume the hypotheses of
Lemmas~\ref{lem:general-constant-focus},
\ref{lem:general-focus-estimates} and \ref{lem:degree-indicator-read},
and let $\Sigma\ge(4rb_{1})^{2^{10}r^2}$.
If $t<T_{\mathrm{init}}$, also assume the hypotheses of
Lemma~\ref{lem:initial-focus-estimates} with $T_0=T_{\mathrm{init}}$.
If $t\ge T_{\mathrm{init}}$, also assume the hypotheses of
Lemmas~\ref{lem:focus-coverage} and \ref{lem:focus-birth-derivatives}.
Then
\begin{equation}\label{eq:good-focus-positive-probability}
 \mathbb P(\mathcal G_{{\rm foc},t})>0.
\end{equation}
For $v\in X$, let $F$ be an event and let $\Phi\ge0$ be a random
variable, each determined by $(\xi_x:x\in\mathbb B_{16r}(v,\mathcal S))$.
Then
\begin{align}
 \mathbb P(F\mid\mathcal G_{{\rm foc},t})
 &\le(1+\Sigma^{-128})\mathbb P(F),
 \\
 \mathbb E(\Phi\mid\mathcal G_{{\rm foc},t})
 &\le(1+\Sigma^{-128})\mathbb E(\Phi).
 \label{eq:good-focus-local-expectation}
\end{align}
\end{lemma}
\begin{proof}
Let $\mathcal E$ consist of the events that the inequalities in
Definition~\ref{def:good-focus} fail, except that for
\eqref{eq:general-degree-read} we use $\mathcal D_{v,k}$ from
Lemma~\ref{lem:degree-indicator-read}~(2). The degree inequalities
\eqref{eq:general-focused-degrees} hold everywhere when their failure
events are excluded, so
\begin{equation*}
 \mathcal G_{{\rm foc},t}=\bigcap_{E\in\mathcal E}\overline E.
\end{equation*}
The lemmas in the hypotheses give, for every $E\in\mathcal E$,
\begin{equation*}
 \mathbb P(E)\le\Sigma^{-2^{10}r^3}.
\end{equation*}

We discard events that are empty, and fix an ordering of $X$.
For an event at $v$, let $v_E=v$. For an event prescribed by a set
$A$ or a pair $x,y$, let $v_E$ be the first atom of that set or pair.
Let $S_E$ consist of the other prescribed atoms, including the atom
whose set $C_x$ is counted and, for a derivative event, the atoms
of its differentiation set. Thus
$|S_E|\le16r$. For each $x\in X$, we take $U_{x,1}=\xi_x$ and
$U_{x,2}=0$, and put
\[
 T_E=\mathbb B_{16r}(v_E,\mathcal S)\times\{1,2\}.
\]
The indicator of $E$ is a function of $(U_{x,j})_{(x,j)\in T_E}$:
\begin{enumerate}
\item The formulas for $d_i(A,Z)$, $d_i(x,y,Z)$, $d_i(v,Z)$,
$M_{j,k}(v)$ and $\widehat M_{j,k}(v)$ use one block, or two blocks
with a common set of atoms. Every atom in these blocks is at
distance at most two from a prescribed atom in the shadow graph.
\item By \eqref{eq:degree-candidate-sets}--\eqref{eq:degree-support}
and \eqref{eq:degree-multiplicity}, $\mu_k(v)$ is determined by $(\xi_x:x\in\mathbb B_4(v,\mathcal S))$.
The same is true of $\mathcal D_{v,k}$, since the degrees in its
definition are evaluated at atoms in $\mathbb B_3(v,\mathcal S)$.
\item Each choice of $i,A$ or $i,x,y$ in
Definition~\ref{def:good-focus}~\eqref{gf:coverage} specifies $\mathcal B,F_0$.
The atoms in $\bigcup\mathcal B$ are at distance
at most one from every atom of $F_0$,
the atoms of $C_x$ are at distance at most two, and those of
$\mathcal T_B$ are at distance at most three.
\item Each polynomial in \eqref{eq:proper-birth-witness-polynomial}
or \eqref{eq:common-birth-witness-polynomial} uses one or two blocks. The atoms of a nonzero derivative belong
to those blocks and are at distance at most two from a prescribed
atom.
\end{enumerate}
For fixed $v\in X$ and
$S\subseteq\mathbb B_{16r}(v,\mathcal S)$ with $|S|\le16r$,
the type of inequality, its ranks, and the roles of the atoms
of $S\cup\{v\}$ can be recorded with fewer than $2^{12}r^2$
binary choices. Consequently,
\[
 |\{E\in\mathcal E:v_E=v,\ S_E=S\}|\le2^{2^{12}r^2}.
\]
Lemma~\ref{lem:conditioning-locality}~(2), with its parameter
$h=\Sigma$ and \eqref{eq:indicator-locality-parameters}, bounds by
$\Sigma^{2^9r^3}$ both the degree of the coordinate dependency graph
and the number of events using variables at atoms of any fixed
ball $\mathbb B_{16r}(v,\mathcal S)$.
We apply Lemma~\ref{lem:conditional-lll} with
$h=\Sigma$, $K=2^9r^3$, $\iota=128$ and $a=2^{10}r^3$.
Its three conclusions give
\eqref{eq:good-focus-positive-probability}--\eqref{eq:good-focus-local-expectation}.
\end{proof}

\section{Construction analysis}\label{sec:construction-analysis}

In this section, we use the one-round notation of
Subsections~\ref{subsec:actual-process} and~\ref{subsec:process-focus-parameters}.
Primes denote the proposed next state.

\begin{definition}[Good construction round]
\label{def:good-construction}\label{def:process-conditioning-events}
Let $t\ge0$ be an integer with $X^t\ne\varnothing$, and let
$Z\subseteq X\setminus\mathsf W$ be a good focus in round $t$. We recall
\begin{equation}\label{eq:general-acceptance-reference}
\mathsf E_{\rm acc}=\sum_{i=2}^r
       \eta^{i-1}\frac{E_i}{\Sigma^{i-1}},\qquad \sigma=\prod_{i=2}^r(1-p^{i-1})^{q^{i-1}d_i},\qquad z=qp\sigma=\frac{\eta\sigma}{\Sigma}.
\end{equation}
An assignment to $(\chi_x:x\in Z)$ and
$(\mathsf{w}_x:x\in X\setminus\mathsf W)$ is a \emph{good construction}
if it has positive probability under the independent joint law
$\mathbb P_Z$ of Definition~\ref{def:degree-original-variables} and
satisfies the following conditions.
\begin{enumerate}[label=\textup{(GC\arabic*)},ref=GC\arabic*]
\item\label{gc:degrees} For every $2\le k\le r$,
\begin{align}
 E_k'&\le{}\omega^{k-1}\sum_{j=k}^r
 \binom{j-1}{k-1}z^{j-k}E_j
 +b_{3}^{r+4}\left(\mathsf E_{\rm acc}\mathcal B_k
 +\Sigma^{k-1-2\delta_3}\right),
 \label{eq:general-local-active-error-update}\\
 \widehat E_k'&\le{}
 \omega^{k-1}(1+b_{3}^{r+4}\mathsf E_{\rm acc})
 \sum_{j=k}^r\binom{j-1}{k-1}z^{j-k}\widehat E_j
 +b_{3}^{r+4}\Sigma^{k-1}
 (\mathsf E_{\rm acc}+\Sigma^{-2\delta_3}),
 \label{eq:general-local-full-error-update}
\end{align}
\item\label{gc:intersections} For $3\le i\le r$ and
$A\subseteq X'$ with $2\le|A|<i$, and for $2\le i\le r$ and
distinct $x,y\in X'$, respectively,
if $t<T_{\mathrm{init}}$, then
\begin{equation}\label{eq:process-initial-round-bounds}
 d_i'(A)\le k_{t+1}D^{i-|A|-\delta_0},\qquad d_i'(x,y)\le k_{t+1}D^{i-1-\delta_0}.
\end{equation}
If $t\ge T_{\mathrm{init}}$, then
\begin{equation}\label{eq:general-next-intersection-barriers}
 d_i'(A)\le(\Sigma')^{i-|A|-\delta_5(i)},\qquad d_i'(x,y)\le(\Sigma')^{i-1-\delta_4(i)}.
\end{equation}
\end{enumerate}
Let $\mathcal G_{{\rm con},t}(Z)$ be the set of good constructions.
\end{definition}

\emph{Remark.} In \eqref{gc:degrees}, $E_k'$ is the maximum, over
$v\in X'$, of the absolute difference between
$d_k'(v,X'\setminus\mathsf W')$ and its reference value.
The quantity $\widehat E_k'$ is the maximum positive excess of
$d_k'(v,X')$ over that value, also over all $v\in X'$.
Lemma~\ref{lem:general-degree-round} supplies the probability bounds
needed to impose these two inequalities.
Condition~\eqref{gc:intersections} controls blocks through a prescribed
set and pairs of blocks leaving the same set after deleting their
prescribed atoms. Before $T_{\mathrm{init}}$, it uses the initial
scale $D$; afterward, it uses the next-round scale $\Sigma'$.
Lemma~\ref{lem:general-overlap-round} and the initial-period argument
in Subsection~\ref{subsec:general-initial-period} establish these bounds.
Lemma~\ref{lem:general-error-closure} proves that the two conditions
can be imposed together with positive probability and remain valid
throughout the process.

\subsection{\texorpdfstring{Probabilities for $J_x$}{\detokenize{Probabilities for J_x}}}\label{sec:general-local-probabilities}

We fix a good focus as in Definition~\ref{def:good-focus}, and put
$h=\Sigma^{\delta_2}$.
For $2\le i\le r$ and $v\in X$, we have $d_i(v,Z)\le3h^{i-1}$;
proper codegrees and common-link counts in $Z$ are at most $b_{2}$.
We work before conditioning on $\mathcal G_{{\rm con},t}(Z)$, with
$0<p\le1/(4h)$. The probability law is $\mathbb P_Z$ unless a modified
experiment is specified.
We put $\tau_h=(\log h)^2/h$.
For $h\ge b_{3}^4$, we have
\begin{equation}\label{eq:local-probability-thresholds}
16r\log h\le h^{1/4},\qquad 2^{12}r^2b_{3}\tau_h\le1,\qquad h\ge b_{3}^2.
\end{equation}
Indeed, $\log h\le(8/e)h^{1/8}$,
$(\log h)^2\le16e^{-2}h^{1/2}$ and $b_{3}\ge2^{16}r^2$.

\Needspace{8\baselineskip}
\begin{lemma}\label{lem:bounded-incidence-avoidance}
Let $\iota\ge1$ and $0\le p\le1/2$. Let $\mathcal B$ be a finite
family of distinct nonempty sets of size at most $r$, with each atom
in at most $\iota$ members. Let $Y\subseteq\bigcup\mathcal B$ contain
each atom $x$ independently with probability $p_x\le p$, and let
$\mu=\sum_{B\in\mathcal B}\prod_{x\in B}p_x$.
If $\iota p\le(4r^2)^{-1}$, then
\[
 \left|\log\mathbb P(B\nsubseteq Y\text{ for every }B\in\mathcal B)+\mu\right|
 \le8r^2\iota p\mu.
\]
\end{lemma}
\begin{proof}
For $B\in\mathcal B$, let $\beta_B=\prod_{x\in B}p_x$.
If $\beta_B=0$, then $\mathbb P(B\subseteq Y)=0$, so this set
does not affect the probability of avoiding every member of $\mathcal B$.
Among the sets with positive weight, let $\mathcal D$ consist of the
inclusion-minimal members. Every positive-weight member of $\mathcal B$
contains a set in $\mathcal D$. Thus the probabilities of avoiding
every member of $\mathcal D$ and every member of $\mathcal B$ are equal.
We assign each omitted positive-weight set $B$ to one $C\in\mathcal D$
with $C\subsetneq B$. There are at most $\iota$ sets assigned to $C$,
since they all contain any fixed atom of $C$. Also
$\beta_B\le p\beta_C$, because $B\setminus C\ne\varnothing$.
Consequently,
\[
 \sum_{B\in\mathcal B\setminus\mathcal D}\beta_B
 \le\iota p\sum_{C\in\mathcal D}\beta_C\le\iota p\mu.
\]
We order the members of $\mathcal D$ arbitrarily.
For each $B\in\mathcal D$, let $\mathcal V$ be the event
that all earlier sets are avoided.
The union outside $B$ of earlier sets meeting $B$ has size at most
$r^2\iota$. Proposition~\ref{prop:harris} and independence imply
$\mathbb P(B\subseteq Y\mid \mathcal V)\ge
\beta_B(1-p)^{r^2\iota}$.
Proposition~\ref{prop:harris} also gives the lower bound in
\[
 \prod_{B\in\mathcal D}(1-\beta_B)
 \le\mathbb P(B\nsubseteq Y\text{ for every }B\in\mathcal B)
 \le\prod_{B\in\mathcal D}(1-(1-p)^{r^2\iota}\beta_B).
\]
The inequalities $1-(1-p)^{r^2\iota}\le r^2\iota p$,
$|\log(1-\beta)+\beta|\le2\beta^2$ and $\sum_B\beta_B^2\le p\mu$,
together with $\sum_{B\in\mathcal B\setminus\mathcal D}\beta_B\le\iota p\mu$,
give an error of at most $8r^2\iota p\mu$.
\end{proof}

\begin{lemma}\label{lem:general-accepted-sets}
Let $h\ge b_{3}^4$ and let $m$ be an integer with $0\le m\le h^{1/4}$.
Let $F\subseteq Z$ be an
independent set of size $m$. For every $x\in Z$,
\begin{equation}\label{eq:local-single-acceptance}
e^{-4}p\le \mathbb P(x\in I)\le p,\qquad \left|\log\frac{\mathbb P(x\in I)}
 {p\prod_{i=2}^r(1-p^{i-1})^{d_i(x,Z)}}\right|
 \le\frac{4r^3(b_{2}+1)}h.
\end{equation}
\begin{equation}\label{eq:local-joint-acceptance}
 \left|\log\frac{\mathbb P(F\subseteq I)}
                    {\prod_{x\in F}\mathbb P(x\in I)}\right|
 \le8r^3(b_{2}+1)\frac{m^2}{h}.
\end{equation}
For any $A\subseteq X$, the same estimates hold for $x\in Z\setminus A$
and independent $F\subseteq Z\setminus A$ in the system
$\{S\in\mathcal S:S\cap A=\varnothing\}$ on $X\setminus A$,
with $Z$ replaced by $Z\setminus A$ and $I,d_i$ and the probabilities
computed in that system.
\end{lemma}
\begin{proof}
The probability that every atom of $F$ is sampled is $p^m$. At most $rb_{2}m^2$ blocks meet $F$
twice; since $F$ is independent, each forbids a nonempty outside
set. Proposition~\ref{prop:harris} bounds their logarithmic contribution by $2rb_{2}m^2p$.
Blocks meeting $F$ once give links of incidence at most $rb_{2}m$
and total weight at most
$3m\sum_{j\ge1}(hp)^j\le m$.
Lemma~\ref{lem:bounded-incidence-avoidance} applies by
\eqref{eq:local-probability-thresholds}.
Passing from distinct links to separate products over $x\in F$
inserts at most $2r^2b_{2}m^2$ factors. These correspond to sets in
$N_i(x,Z)$ meeting $F\setminus\{x\}$, or belonging also to
$N_i(y,Z)$ for another $y\in F$. Each changes the logarithm by
at most $2p$. Including the product-to-exponential error gives
\[
 \left|\log\frac{\mathbb P(F\subseteq I)}
 {p^m\prod_{x\in F}\prod_{i=2}^r
                     (1-p^{i-1})^{d_i(x,Z)}}\right|
 \le4r^3(b_{2}+1)\frac{m^2}{h}.
\]
The case $m=1$ proves the logarithmic bound in
\eqref{eq:local-single-acceptance}; subtracting it for each $x\in F$
proves \eqref{eq:local-joint-acceptance}.
Proposition~\ref{prop:harris} gives $\mathbb P(x\in I)\ge p\exp(-2\sum_i d_i(x,Z)p^{i-1})
\ge e^{-2}p$, and $\mathbb P(x\in I)\le p$ follows from $I\subseteq Y$.
Finally, suppose that $\chi_x=0$ for every $x\in A$.
No block meeting $A$ can be contained in $Y$, so such a block cannot
cause an atom of $Z\setminus A$ to be excluded from $I$.
Acceptance on $Z\setminus A$ is therefore computed using precisely
the blocks in $\{S\in\mathcal S:S\cap A=\varnothing\}$.
The variables $(\chi_x:x\in Z\setminus A)$ remain independent
Bernoulli variables with parameter $p$. Deleting atoms and blocks
cannot increase any of the degree, codegree or common-link counts
bounded in this proof. The same argument therefore gives the claimed
estimates with $I$, $d_i$ and the probabilities computed in this
restricted system.
\end{proof}

\subsection{\texorpdfstring{Conditional probabilities for $I$}{Conditional probabilities for I}}\label{subsec:conditional-acceptance}

Let $A\subseteq X$ with $|A|\le r$, let $\chi_x=0$ for $x\in A$, and let
$\mathcal B=\bigcup_{v\in A}\bigcup_{i=2}^rN_i(v,Z\setminus A)$
be the family of distinct links at these atoms. For $1\le j\le r-1$,
at most $3rh^j$ members have size $j$, and each atom belongs to at most
$r^2b_{2}$ members. We put $\mathcal A=\{x:\{x\}\in\mathcal B\}$.
For every $u\in X\setminus A$,
\begin{equation}\label{eq:local-link-inputs}
\sum_{B\in\mathcal B}p^{|B|}\le r,\qquad |\mathcal A|\le3rh,\qquad |N_2(u,Z)\cap\mathcal A|\le rb_{2}.
\end{equation}
The last bound follows from \eqref{eq:general-constant-focused-completions}.

\begin{lemma}\label{lem:general-accepted-blocker-avoidance}
Let $h\ge b_{3}^4$, and let $Q\subseteq Z\setminus A$ be independent
with $|Q|\le r$. All probabilities in this lemma are in the experiment with
$\chi_x=0$ for $x\in A$. Let
$\mathcal B_Q=\{B\in\mathcal B:B\cap Q=\varnothing\}$ and
$\lambda_Q=\sum_{B\in\mathcal B_Q}\prod_{x\in B}\mathbb P(x\in I)$. Then
\begin{equation}\label{eq:local-conditional-avoidance}
 \left|\log\mathbb P(B\nsubseteq I\text{ for every }B\in\mathcal B_Q
              \mid Q\subseteq I)+\lambda_Q\right|
 \le\frac{b_{3}}{2^8r^2}\tau_h.
\end{equation}
\end{lemma}
\begin{proof}
Since $|Q|\le r\le h^{1/4}$, Lemma~\ref{lem:general-accepted-sets}
gives $\mathbb P(Q\subseteq I)>0$, so the conditional probability in
the statement is defined. We first record a bound for this conditional
law. Let $F\subseteq Z\setminus A$ be independent. Once every atom of
$F$ is sampled, all of them are accepted precisely when no block meeting
$F$ is entirely sampled. Since $F$ is independent, each such block has
an atom outside $F$. Changing an outside sampling variable from one to
zero cannot violate this requirement. Thus $F\subseteq I$, conditional
on $F\subseteq Y$, is a decreasing event in the remaining variables.
For $T\subseteq Z\setminus(A\cup F)$, the event $T\subseteq Y$ is
increasing, and Proposition~\ref{prop:harris} gives
\begin{equation}\label{eq:local-conditional-domination}
 \mathbb P(T\subseteq Y\mid F\subseteq I)\le p^{|T|}.
\end{equation}
The same bound holds when $T$ meets $A$, because then its left side is zero.

Let $\mathcal D$ consist of the inclusion-minimal members of
$\mathcal B_Q$, and put $\iota=r^2b_2$ and
$\mu=\sum_{B\in\mathcal D}p^{|B|}\le r$.
Every member of $\mathcal B_Q$ contains a member of $\mathcal D$.
Consequently a set contains no member of $\mathcal B_Q$ if and only if
it contains no member of $\mathcal D$.
We assign each $E\in\mathcal B_Q\setminus\mathcal D$ to one contained
$B\in\mathcal D$. There are at most $\iota$ sets assigned to any $B$,
because they all contain any fixed atom of $B$. Also $E\setminus B$ is
nonempty and $\mathbb P(x\in I)\le p$, so
\begin{equation}\label{eq:local-omitted-link-weights}
 0\le\lambda_Q-\sum_{B\in\mathcal D}\prod_{x\in B}\mathbb P(x\in I)
 \le\iota p\mu.
\end{equation}

For comparison, let $Y^\ast\subseteq\bigcup\mathcal D$ include its atoms
independently, with $\mathbb P(x\in Y^\ast)=\mathbb P(x\in I)$.
We take $Y^\ast$ independent of the original sampling variables.
Let $N$ count those $B\in\mathcal D$ for which $B\subseteq I$ and no
other $E\in\mathcal D$ meeting $B$ is contained in $Y$.
Let $N^\ast$ count those $B\in\mathcal D$ contained in $Y^\ast$ for which
no other member of $\mathcal D$ meeting $B$ is contained in $Y^\ast$.
These restrictions ensure that the links counted by either variable
are pairwise disjoint.
We also define the events
\[
 \mathcal U=\bigcap_{B\in\mathcal D}\{B\nsubseteq I\},\qquad
 \mathcal U^\ast=\bigcap_{B\in\mathcal D}\{B\nsubseteq Y^\ast\}.
\]

For distinct intersecting $B,E\in\mathcal D$, minimality gives
$E\setminus B\ne\varnothing$. For each $B$, at most $r\iota$ members
of $\mathcal D$ meet $B$. By \eqref{eq:local-conditional-domination}
with $F=Q$, the conditional probability that some such pair is sampled
is therefore at most
$\sum_{B\in\mathcal D}r\iota p^{|B|+1}=r\iota p\mu$.
The same estimate holds for $Y^\ast$, since its atom probabilities are
at most $p$. If no such pair is present, every member of $\mathcal D$ contained in
$I$ is counted by $N$; the corresponding statement holds for $Y^\ast$
and $N^\ast$. Hence
\begin{equation}\label{eq:local-counting-restriction}
\begin{aligned}
0&\le\mathbb P(N=0\mid Q\subseteq I)
 -\mathbb P(\mathcal U\mid Q\subseteq I)
 \le r\iota p\mu,\\
0&\le\mathbb P(N^\ast=0)
 -\mathbb P(\mathcal U^\ast)
 \le r\iota p\mu.
\end{aligned}
\end{equation}

We put $m=2\lceil2\log_2h\rceil$, so $m$ is even and
$4\log_2h\le m\le8\log h$.
For $\mathcal F\subseteq\mathcal D$, let $H_{\mathcal F}$ be the event
that every member of $\mathcal F$ is counted by $N$, and let
$H_{\mathcal F}^\ast$ be the corresponding event for $N^\ast$.
If two members of $\mathcal F$ intersect, both events are impossible.
Suppose therefore that its members are pairwise disjoint, and put
$U=\bigcup\mathcal F$. If some $E\in\mathcal D\setminus\mathcal F$
is contained in $U$, both events are again impossible: this $E$ is
present and meets a member of $\mathcal F$.

Suppose no such $E$ exists and $Q\cup U$ is independent, with
$|\mathcal F|\le m$. We have
$|Q\cup U|\le r(1+m)\le9r\log h\le16r\log h\le h^{1/4}$ by
\eqref{eq:local-probability-thresholds}. Applying
Lemma~\ref{lem:general-accepted-sets} to $Q\cup U$ and to $Q$ gives
\[
 \left|\log\frac{\mathbb P(U\subseteq I\mid Q\subseteq I)}
                   {\prod_{x\in U}\mathbb P(x\in I)}\right|
 \le\frac{8r^3(b_2+1)}h\bigl((|Q|+|U|)^2+|Q|^2\bigr)
 \le2^{10}r^5(b_2+1)\tau_h.
\]
At most $\iota|U|$ other members of $\mathcal D$ meet $U$.
Each has an atom outside $U$, and all are disjoint from $Q$.
Conditional on $Q\cup U\subseteq I$, the probability that any of these
members is sampled is at most $\iota|U|p$ by
\eqref{eq:local-conditional-domination}. Conditional on $U\subseteq Y^\ast$,
the same bound follows from independence of the remaining atoms.
Thus the additional restrictions defining $H_{\mathcal F}$ and
$H_{\mathcal F}^\ast$ each change the respective conditional probability
by a relative amount at most $\iota|U|p$.
Using $|e^\beta-1|\le2|\beta|$ for $|\beta|\le1/2$, with
\eqref{eq:local-probability-thresholds}, we obtain
\begin{equation}\label{eq:local-counted-joint-comparison}
 \left|\mathbb P(H_{\mathcal F}\mid Q\subseteq I)
             -\mathbb P(H_{\mathcal F}^\ast)\right|
 \le2^{14}r^5(b_2+1)\tau_h
           \prod_{x\in U}\mathbb P(x\in I).
\end{equation}

It remains to bound the contribution from pairwise disjoint families
$\mathcal F$ for which $Q\cup\bigcup\mathcal F$ contains a block.
Let $\mathfrak M_Q$ be the collection of subfamilies
$\mathcal F\subseteq\mathcal D$ whose members are pairwise disjoint
and for which $Q\cup\bigcup\mathcal F$ contains a block.
For every $B\in\mathcal D$, let $v_B\in A$ satisfy
$B\cup\{v_B\}\in\mathcal S$; such an atom exists by the definition
of $\mathcal B$.
For $\mathcal F\in\mathfrak M_Q$, let $S$ be the first block contained
in $Q\cup\bigcup\mathcal F$, in a fixed ordering of $\mathcal S$, and let
$\mathcal F_0=\{B\in\mathcal F:B\cap(S\setminus Q)\ne\varnothing\}$.
These links cover $S\setminus Q$, and each is needed because they are
pairwise disjoint. We write $\mathcal F_0=\{B_1,\ldots,B_k\}$ in a fixed
order. For $1\le i\le k$, let $j_i=|B_i|$, and put
$\ell=j_1+\cdots+j_k$.
Since $Q$ is independent, $1\le k\le r$ and $1\le\ell\le r(r-1)$.

We put $s=2r^2+3r+1$ and $R=2(s+1)(4+r^2b_2)$.
To count $S$ and $\mathcal F_0$, we list their atoms, those of $Q$ and
the roots $v_{B_i}$, and mark one atom of $S\setminus Q$ in each $B_i$.
At most $s$ positions are needed. There are at most
$2^{s^2}(s+1)^{5s}\le R^{s^2}$ ways to specify equal positions,
set membership, set sizes, the order of the links, their roots and
the marked positions. Once this information is fixed, for $1\le j\le r-1$, a link of size
$j$ has at most $3rh^j$ choices; a link containing a prescribed atom
has at most $\iota=r^2b_2$ choices; and a block containing two prescribed
atoms has at most $rb_2$ choices. We use these bounds as follows.
\begin{enumerate}[leftmargin=*,itemsep=.2em]
\item Suppose $S\cap Q=\varnothing$. We cannot have $k=1$, because
$S\subseteq B_1$ would give $S\subsetneq B_1\cup\{v_{B_1}\}$,
contrary to the antichain assumption.
If $k\ge3$, we first specify $B_1,B_2$, then $S$ through their marked
atoms, and finally every other link through its marked atom in $S$.
The only powers of $h$ arise from the first two links, and their
exponent is $j_1+j_2\le\ell-1$.
If $k=2$ and one link contains two atoms of $S$, we specify that link,
then $S$ through those two atoms, and then the other link. The exponent
is at most $\ell-1$ again.
Otherwise $S$ is a pair. If one link has size at least two, we specify
the other link, the second endpoint of $S$ in at most $3h$ ways,
and the larger link through that endpoint. If the first link has size
$j_1$, the exponent is $j_1+1\le\ell-1$.
If both links are singletons, their atoms lie in $\mathcal A$;
\eqref{eq:local-link-inputs} gives at most $3r^2b_2h$ choices for the pair.
\item Suppose $|S\cap Q|\ge2$. We specify $S$ through two of these fixed
atoms, then specify each link through its marked atom in $S$.
These choices contribute no power of $h$, and $\ell\ge1$.
\item Suppose $S\cap Q=\{u\}$. If $k\ge2$, we specify $B_1$, then $S$
through $u$ and the marked atom of $B_1$, and then the remaining links
through their marked atoms in $S$. The exponent is $j_1\le\ell-1$.
If $k=1$, we put $j_0=|S\setminus\{u\}|$. When $j_1>j_0$, we specify $S$
through $u$ in at most $3h^{j_0}$ ways and then $B_1$ through one of
its atoms in $S\setminus\{u\}$; here $j_0\le\ell-1$.
When $j_1=j_0$, we have $B_1=S\setminus\{u\}$, so $B_1$ is a common
link of $u$ and $v_{B_1}$. The two roots are distinct because
$u\notin A$, and the common-link bound gives at most $rb_2$ choices,
including the choice of $v_{B_1}$.
\end{enumerate}
Each case uses at most $s$ choices with constant factors at most $R$.
There are therefore at most $R^{s^2+s}h^{\ell-1}$ choices of
$S,\mathcal F_0$ for a fixed $\ell$.
Since $s\le4r^2$, the definition of $b_3$ gives
$R^{s^2+s}\le R^{2s^2}\le b_3^{1/8}$.

For fixed $\mathcal F_0$, its contribution to
$p^{|\bigcup\mathcal F|}$ is $p^\ell$.
Every additional link must be disjoint from $\bigcup\mathcal F_0$ and
from the other additional links. Dropping these restrictions and summing
over all subfamilies of $\mathcal D$ can only increase the sum.
Expanding the product shows that the total contribution of all such
extensions is at most
\[
 p^\ell\sum_{\mathcal H\subseteq\mathcal D}
              \prod_{B\in\mathcal H}p^{|B|}
 =p^\ell\prod_{B\in\mathcal D}(1+p^{|B|})
 \le p^\ell e^\mu\le p^\ell e^r.
\]
Summing over the possibilities for $S,\mathcal F_0$ therefore gives
\begin{equation}\label{eq:local-nonindependent-link-families}
 \sum_{\mathcal F\in\mathfrak M_Q}p^{|\bigcup\mathcal F|}
 \le e^r b_3^{1/8}\sum_{\ell=1}^{r(r-1)}h^{\ell-1}p^\ell
 \le\frac{e^r b_3^{1/8}}h,
\end{equation}
The final bound in \eqref{eq:local-nonindependent-link-families}
uses $hp\le1/4$.

For $j\ge0$, $\binom Nj$ counts the $j$-element subfamilies whose links
are all counted by $N$, and the corresponding statement holds for
$N^\ast$. Only pairwise disjoint subfamilies contribute.
By \eqref{eq:local-conditional-domination} with $F=Q$, each such
subfamily has conditional probability at most $p^{|\bigcup\mathcal F|}$;
the same upper bound holds in the independent model.
Expanding $\mu^m$ counts every $m$-element subfamily at least $m!$ times.
Thus
\[
\begin{aligned}
 \mathbb E\!\left(\binom Nm\,\middle|\,Q\subseteq I\right)
 &\le\frac{\mu^m}{m!}\le\frac{r^m}{m!}\le h^{-4},\\
 \mathbb E\!\left(\binom{N^\ast}m\right)
 &\le\frac{\mu^m}{m!}\le\frac{r^m}{m!}\le h^{-4}.
\end{aligned}
\]
Indeed, $h\ge b_3^4\ge e^{r^2}$ gives $m\ge2er$, and
$r^m/m!\le(er/m)^m\le2^{-m}\le h^{-4}$.
For every outcome, the difference between $\vb1_{\{N=0\}}$ and
$\sum_{j=0}^{m-1}(-1)^j\binom Nj$ has absolute value at most
$\binom Nm$. The same statement holds with $N^\ast$ in place of $N$.
Taking expectations gives, in both models,
\[
\begin{aligned}
 \left|\mathbb P(N=0\mid Q\subseteq I)
 -\sum_{j=0}^{m-1}(-1)^j
       \mathbb E\!\left(\binom Nj\,\middle|\,Q\subseteq I\right)\right|
 &\le\mathbb E\!\left(\binom Nm\,\middle|\,Q\subseteq I\right)\le h^{-4},\\
 \left|\mathbb P(N^\ast=0)
 -\sum_{j=0}^{m-1}(-1)^j\mathbb E\!\left(\binom{N^\ast}j\right)\right|
 &\le\mathbb E\!\left(\binom{N^\ast}m\right)\le h^{-4}.
\end{aligned}
\]

We now compare the two finite sums. Each expectation of $\binom Nj$
or $\binom{N^\ast}j$ is the sum of the probabilities that all links
in the corresponding $j$-element subfamily are counted.
For pairwise disjoint $\mathcal F$ with $Q\cup\bigcup\mathcal F$
independent, we use \eqref{eq:local-counted-joint-comparison}.
Because the links are disjoint, the product over their union is
$\prod_{B\in\mathcal F}\prod_{x\in B}\mathbb P(x\in I)$.
Summing the products in \eqref{eq:local-counted-joint-comparison}
over all sizes $j<m$ introduces no factor $m$:
each subfamily occurs once, and their total is at most
\[
 \prod_{B\in\mathcal D}
       \left(1+\prod_{x\in B}\mathbb P(x\in I)\right)
 \le\exp\!\left(\sum_{B\in\mathcal D}p^{|B|}\right)=e^\mu\le e^r.
\]
For $\mathcal F\in\mathfrak M_Q$, the actual joint acceptance probability
is zero, while the comparison probability is at most
$p^{|\bigcup\mathcal F|}$. We bound their total using
\eqref{eq:local-nonindependent-link-families}.
Together with the two remainders, this gives
\[
 \left|\mathbb P(N=0\mid Q\subseteq I)-\mathbb P(N^\ast=0)\right|
 \le2^{14}r^5(b_2+1)e^r\tau_h
       +\frac{e^r b_3^{1/8}}h+2h^{-4}.
\]
Applying \eqref{eq:local-counting-restriction} adds at most
$2r\iota p\mu$. Using $p\le1/(4h)$, $\mu\le r$ and $\log h\ge1$, we obtain
\begin{equation}\label{eq:local-avoidance-model-comparison}
\begin{aligned}
 \left|\mathbb P(\mathcal U\mid Q\subseteq I)-\mathbb P(\mathcal U^\ast)\right|
 &\le2r\iota p\mu+2^{14}r^5(b_2+1)e^r\tau_h
            +\frac{e^r b_3^{1/8}}h+2h^{-4}\\
 &\le e^{-r}b_3^{1/4}\tau_h.
\end{aligned}
\end{equation}
The last inequality in \eqref{eq:local-avoidance-model-comparison}
follows by substituting $\iota=r^2b_2$ and
$b_3=(2^{16}r(b_2+1))^{2^{10}r^4}$.

Since $\iota p\le r^2b_2/(4h)\le(4r^2)^{-1}$,
Lemma~\ref{lem:bounded-incidence-avoidance} applies to $Y^\ast$ and gives
\[
 \left|\log\mathbb P(\mathcal U^\ast)
       +\sum_{B\in\mathcal D}\prod_{x\in B}\mathbb P(x\in I)\right|
 \le8r^2\iota p\mu.
\]
Here $8r^2\iota p\mu\le2r^5b_2/h\le1/4$, so the comparison
avoidance probability is at least $e^{-r-1/4}>3e^{-r}/4$.
By \eqref{eq:local-probability-thresholds},
$b_3^{1/4}\tau_h\le1/4$. Thus the difference in
\eqref{eq:local-avoidance-model-comparison} is at most $e^{-r}/4$,
and the actual conditional avoidance probability is at least $e^{-r}/2$.
The logarithm has derivative at most $2e^r$ throughout the interval
between these two positive probabilities. Their logarithms therefore differ by at most
$2e^r$ times the bound in \eqref{eq:local-avoidance-model-comparison},
which is at most $2b_3^{1/4}\tau_h$.
Finally, \eqref{eq:local-omitted-link-weights} adds at most $\iota p\mu$
when the sum over $\mathcal D$ is replaced by $\lambda_Q$.
Together with the error from Lemma~\ref{lem:bounded-incidence-avoidance},
these remaining terms satisfy
\[
 (8r^2+1)\iota p\mu
 \le3r^5b_2\tau_h\le4b_3^{1/4}\tau_h.
\]
The total logarithmic error is consequently at most
\[
 6b_3^{1/4}\tau_h\le\frac{b_3}{2^8r^2}\tau_h.
\]
This proves \eqref{eq:local-conditional-avoidance}, since avoiding
$\mathcal D$ is equivalent to avoiding $\mathcal B_Q$.
This completes the proof of Lemma~\ref{lem:general-accepted-blocker-avoidance}.
\end{proof}

\subsection{\texorpdfstring{Estimates for $\mathbb P(x\in I)$ and $s_v$}{Estimates for acceptance and survival probabilities}}\label{subsec:survival-probabilities}

Probabilities and expectations in this subsection are under $\mathbb P_Z$,
unless a conditioning event is specified.

\begin{lemma}\label{lem:general-survival-estimates}
Let $Z$ be a good focus with $h=\Sigma^{\delta_2}\ge b_{3}^4$.
Suppose $0<\eta\le1/4$ and $d_i(v,Z)\le3h^{i-1}$ for
$v\in X$ and $2\le i\le r$. For each $v\in X$, let
\begin{equation}\label{eq:general-lambda}
 \lambda_v=\sum_{i=2}^r\sum_{B\in N_i(v,Z)}
                 \prod_{x\in B}\mathbb P(x\in I).
\end{equation}
\begin{enumerate}
\item For $v\in X$ and $A\subseteq X$ with $|A|\le r$,
\begin{align}
 \left|s_v-(1-p)^{\vb1_{\{v\in Z\}}}e^{-\lambda_v}\right|
 &\le b_{3}(1-p)^{\vb1_{\{v\in Z\}}}e^{-\lambda_v}\tau_h,
 \label{eq:general-survival-one}\\
 \left|\mathbb P(\psi_v=1\text{ for all }v\in A)
                    -\prod_{v\in A}s_v\right|
 &\le b_{3}\left(\prod_{v\in A}s_v\right)\tau_h.
 \label{eq:general-survival-joint}
\end{align}
\item Let $Q,A\subseteq X$ be disjoint, with $Q\cup A\in\mathcal S$,
$Q\subseteq Z$ and $|A|\ge2$. Then
\begin{equation}\label{eq:general-survival-mixed}
 \left|\mathbb P(Q\subseteq I,\ \psi_v=1\text{ for all }v\in A)
       -\prod_{x\in Q}\mathbb P(x\in I)\prod_{v\in A}s_v\right|
 \le b_{3}\prod_{x\in Q}\mathbb P(x\in I)\prod_{v\in A}s_v\,\tau_h.
\end{equation}
For $v\in A$ with $A\setminus\{v\}\subseteq X\setminus\mathsf W$,
\begin{equation}\label{eq:general-marked-mixed}
 \left|\mathbb E\left(
 \prod_{x\in Q}J_x\prod_{u\in A\setminus\{v\}}\psi_u(1-\mathsf{w}_u)
 \,\middle|\,\psi_v=1\right)
 -\omega^{|A|-1}\prod_{x\in Q}\mathbb P(x\in I)\right|
 \le b_{3}\omega^{|A|-1}\prod_{x\in Q}\mathbb P(x\in I)\,\tau_h.
\end{equation}
For every $v\in A$, the same bound holds with $\psi_u(1-\mathsf{w}_u)$ replaced by
$\psi_u$ and $\omega^{|A|-1}$ replaced by
$\prod_{u\in A\setminus\{v\}}s_u$.
\item Suppose additionally that $\Sigma\ge b_{5}$ and
$E_i\le\Sigma^{i-1}$ for every $2\le i\le r$.
For $x\in Z$ and $v\in X$,
\begin{equation}\label{eq:general-uniform-acceptance}
 |\mathbb P(x\in I)-p\sigma|
 \le b_{3}p\sigma(\mathsf E_{\rm acc}+\Sigma^{-4\delta_3}),
\end{equation}
\begin{equation}\label{eq:general-survival-log-reference}
 \left|\log s_v+\sum_{i=2}^r d_i z^{i-1}\right|
 \le b_{3}(\mathsf E_{\rm acc}+\Sigma^{-4\delta_3}),
\end{equation}
\begin{equation}\label{eq:general-survival-comparison}
s_v\ge\frac18,\qquad \omega\le s_v\le
 \omega\left(1+b_{3}
       (\mathsf E_{\rm acc}+\Sigma^{-4\delta_3})\right).
\end{equation}
Inequality~\eqref{eq:general-survival-log-reference} also holds with
$s_v$ replaced by $\omega$.
\end{enumerate}
\end{lemma}
\begin{proof}
For \textup{(1)}, setting $\chi_u=0$ for $u\in A$ deletes at most
$r^2b_{2}$ blocks through any $x\in Z\setminus A$.
Proposition~\ref{prop:harris} therefore gives
\begin{equation*}
 0\le\log\frac{\mathbb P(x\in I\mid\chi_u=0\text{ for all }u\in A)}
                    {\mathbb P(x\in I)}\le2r^2b_{2}p.
\end{equation*}
At most $2r^3b_{2}$ sets in
$\bigcup_{v\in A}\bigcup_iN_i(v,Z)$ meet $A$ or belong to the links
of two distinct atoms of $A$. Removing them and replacing each $\mathbb P(x\in I)$ by
$\mathbb P(x\in I\mid\chi_u=0\text{ for all }u\in A)$ alters their
total product weight by at most $8r^4b_{2}p$.
Exponentiating \eqref{eq:local-conditional-avoidance} and including
$(1-p)^{|A\cap Z|}$ and the weight error proves
\eqref{eq:general-survival-one}--\eqref{eq:general-survival-joint}.

For \textup{(2)}, no link at an atom of $A$ lies in $Q$:
adding its root would give a block properly contained in $Q\cup A$.
At most $r^3b_{2}$ root links meet $Q$; conditional on $Q\subseteq I$,
each is accepted with probability at most $p$ by
\eqref{eq:local-conditional-domination}. Deleting them changes
the probability and sum of link probabilities by at most $r^3b_{2}p$.
The remaining avoidance probability is at least $e^{-r}/2$,
so restoring these links changes the logarithm of the avoidance
probability by at most $4e^r r^3b_2p$.
We apply \eqref{eq:local-conditional-avoidance} and
\eqref{eq:local-joint-acceptance}.
Since $h\ge b_{3}^4$, the value of $b_{3}$ bounds each logarithmic error by
$b_{3}\tau_h/(32r)$; exponentiation proves
\eqref{eq:general-survival-mixed}.
The independence of the $\mathsf{w}_u$ from the $\chi_x$, followed by
division by $s_v$, gives \eqref{eq:general-marked-mixed}.
Omitting the factors $1-\mathsf{w}_u$ gives its stated version with $\psi_u$.

For \textup{(3)}, the additional hypothesis $E_i\le\Sigma^{i-1}$,
together with $\eta\le1/4$, gives
$\mathsf E_{\rm acc}\le1/3$ and
$\exp(-2\sum_{j\ge1}4^{-j})\le\sigma\le1$.
The condition $\Sigma\ge b_{5}$ gives
$16r\le\Sigma^{4\delta_3}$ and
$(4b_{3}+4)(\log\Sigma)^2\le\Sigma^{\delta_2-4\delta_3}$.
Putting $\varepsilon_0=\Sigma^{-4\delta_3}$ therefore gives
$b_{3}\tau_h\le \varepsilon_0/4$ and $\varepsilon_0\le1/(16r)$.
Equation~\eqref{eq:general-focused-degrees} and
Lemma~\ref{lem:general-accepted-sets} give, for every $x\in Z$,
\[
 \left|\log\frac{\mathbb P(x\in I)}{p\sigma}\right|
 \le2\mathsf E_{\rm acc}+2(r-1)\varepsilon_0+\frac{4r^3(b_{2}+1)}h
 \le2\mathsf E_{\rm acc}+2r\varepsilon_0<1.
\]
Exponentiating proves \eqref{eq:general-uniform-acceptance}.
To estimate $\lambda_v$, we put $L=2\mathsf E_{\rm acc}+2r\varepsilon_0<1$.
Since $\eta e^L\le3/4$, replacing the factors $\mathbb P(x\in I)$
by $p\sigma$ changes the sum in \eqref{eq:general-lambda} by at most
$3L\sum_{j\ge1}j(\eta e^L)^j\le36L$.
By \eqref{eq:general-focused-degrees} and \eqref{eq:general-acceptance-reference},
the degree errors add at most $\mathsf E_{\rm acc}+(r-1)\varepsilon_0$. Thus
$|\lambda_v-\sum_i d_i z^{i-1}|
\le128r(\mathsf E_{\rm acc}+\varepsilon_0)$.
In \eqref{eq:general-survival-one}, the factor for not sampling $v$
and the relative error contribute at most
$2p+2b_{3}\tau_h\le \varepsilon_0$. Hence
\[
 \left|\log s_v+\sum_{i=2}^r d_i z^{i-1}\right|
 \le256r(\mathsf E_{\rm acc}+\varepsilon_0),
\]
which proves \eqref{eq:general-survival-log-reference}; the same bound
holds for $\omega=\min_v s_v$.
Since $\lambda_v\le1$ and $b_{3}\tau_h\le1/2$,
$s_v\ge3/(8e)>1/8$. Comparing the logarithms, and using
$512r e^{512r}\le b_{3}$, proves
\eqref{eq:general-survival-comparison}.
Finally, $b_{3}\tau_h\le\Sigma^{-4\delta_3}/4$ bounds the relative
errors in \eqref{eq:general-survival-joint} and
\eqref{eq:general-marked-mixed} as required in later applications.
\end{proof}

\subsection{\texorpdfstring{Concentration of $d_k'(v,X')$ and $d_k'(v,X'\setminus\mathsf W')$}{Concentration of the residual degrees}}
\label{sec:general-degree-concentration}

We fix $Z$ and use $F_k(v),\widehat F_k(v)$ under $\mathbb P_{Z,v}^0$ as
defined in Subsection~\ref{subsec:general-degree-indicators}.
Throughout this subsection, we assume $E_i,\widehat E_i\le\Sigma^{i-1}$
for $2\le i\le r$, the intersection bounds
\eqref{eq:general-codegree-barriers}--\eqref{eq:general-common-barriers},
$\Sigma\ge b_{5}$ and $0<\eta\le1/4$.

We first estimate the conditional means. We assume that the focus satisfies
\eqref{eq:general-constant-focused-completions} and
\eqref{eq:general-focused-degrees}--\eqref{eq:general-full-contractors}.
For $2\le k\le r$ and $v\in X$, we count pairs $(S,C)$ with
$v\in S\in\mathcal S$, $C\subseteq(S\setminus\{v\})\cap I$,
$|S\setminus C|=k$ and
$S\setminus(C\cup\{v\})\subseteq X'\setminus\mathsf W'$
separately, even if distinct pairs leave the same set or one
remaining set contains another:
\[
 T_k(v)=\sum_{j=k}^r\sum_{v\in S\in\mathcal S_j}
 \sum_{\substack{C\subseteq(S\setminus\{v\})\cap Z\\|C|=j-k}}
 J_C\prod_{u\in S\setminus(C\cup\{v\})}\psi_u(1-\mathsf{w}_u).
\]
Replacing $\psi_u(1-\mathsf{w}_u)$ by $\psi_u$ defines $\widehat T_k(v)$.
Replacing each factor $\psi_u$ by $1$ and applying
\eqref{eq:general-full-contractors} gives
\begin{equation}\label{eq:general-oriented-mean-upper}
 \mathbb E_{Z,v}^0(\widehat T_k(v))
 \le\sum_{j=k}^rp^{j-k}\widehat M_{j,k}(v)
 \le 2^{r+2}\Sigma^{k-1}.
\end{equation}
Indeed, the sum is at most
$2\Sigma^{k-1}\sum_{j=k}^r\binom{j-1}{k-1}\eta^{j-k}
+r\Sigma^{k-1-4\delta_3}$ by
\eqref{eq:general-full-contractors}.

\begin{lemma}\label{lem:general-cleanup}
Let $2\le k\le r$ and $v\in X$. Under $\mathbb P_{Z,v}^0$, both
unconditionally and conditional on $\psi_v=1$,
\[
\begin{aligned}
 0&\le\mathbb E(T_k(v)-F_k(v))
 \le8r(b_{2}+1)(1+2^r)2^{r+2}\,p\Sigma^{k-1},\\
 0&\le\mathbb E(\widehat T_k(v)-\widehat F_k(v))
 \le8r(b_{2}+1)(1+2^r)2^{r+2}\,p\Sigma^{k-1}.
\end{aligned}
\]
\end{lemma}
\begin{proof}
For a fixed set $B$ of size at least two, at most $r b_{2}$ blocks
$S\in\mathcal S$ satisfy $S\setminus B\subseteq Z$ and $B\subseteq S$.
For distinct such blocks $B\cup C$ and $B\cup C'$, the antichain
property makes $C,C'$ incomparable, so $|C\cup C'|\ge|C|+1$.
Thus counting both when $C\cup C'\subseteq I$ adds at most
$r b_{2}p^{|C|+1}$ per choice of $S=B\cup C$.
If $B$ is deleted because $B'\subsetneq B$ also remains, its source
block has the form $B'\cup C'\in\mathcal S$, where
$C'\subseteq I$ and $C'\cap B=\varnothing$.
There are at most $2^r r b_{2}$ such source blocks, and the antichain
property gives $C'\setminus C\ne\varnothing$.
We sum these bounds using \eqref{eq:general-oriented-mean-upper}.
Finally, conditioning on survival multiplies a nonnegative
expectation by at most eight, by
\eqref{eq:general-survival-comparison}.
\end{proof}

In a summand of $T_k(v)$, write $S=C\cup B$, where $|B|=k$.
The estimate
\eqref{eq:general-marked-mixed} has relative error at most
$\Sigma^{-4\delta_3}$. For $C=\varnothing$, its reference value is $\omega^{k-1}$. For $C\ne\varnothing$, \eqref{eq:general-uniform-acceptance}
replaces $\prod_{x\in C}\mathbb P(x\in I)$ by $(p\sigma)^{|C|}$ with
relative error at most
\[
 [1+r(1+2b_{3})^r]
       (\mathsf E_{\rm acc}+\Sigma^{-4\delta_3}).
\]
This acceptance error occurs only in terms with $C\ne\varnothing$.
Their total reference value is bounded by
\begin{equation*}
 \mathcal B_k=
 \sum_{j=k+1}^r\binom{j-1}{k-1}z^{j-k}(d_j+E_j)
 \le2^{r+1}\eta\Sigma^{k-1}.
\end{equation*}
We sum the estimate over $S\in\mathcal S$, apply
\eqref{eq:general-focused-contractors}, and then
Lemma~\ref{lem:general-cleanup}. Their coefficients are bounded by
$b_{3}^{r+3}$, since
\[
 [1+r(1+2b_{3})^r]
 [2^{r+2}+8r(b_{2}+1)(1+2^r)2^{r+2}]
 \le b_{3}^{r+3}.
\]
Using $p\le\Sigma^{-4\delta_3}$ gives
\begin{equation}\label{eq:general-active-conditional-mean}
\left|\mathbb E(F_k(v)\mid\psi_v=1) -\omega^{k-1}\sum_{j=k}^r \binom{j-1}{k-1}z^{j-k}d_j(v,X\setminus \mathsf{W})\right| \le b_{3}^{r+3} (\mathsf E_{\rm acc}\mathcal B_k+\Sigma^{k-1-4\delta_3}).
\end{equation}
For $\widehat F_k(v)$, we replace $\psi_u(1-\mathsf{w}_u)$ by $\psi_u$, use
\eqref{eq:general-survival-comparison} and
\eqref{eq:general-full-contractors}, and omit deletions. This gives
\begin{equation}\label{eq:general-full-conditional-mean}
\mathbb E(\widehat F_k(v)\mid\psi_v=1) \le\omega^{k-1}(1+b_{3}^{r+3}\mathsf E_{\rm acc}) \sum_{j=k}^r\binom{j-1}{k-1}z^{j-k}d_j(v,X) {}+b_{3}^{r+3}\Sigma^{k-1-4\delta_3}.
\end{equation}

We next bound deviations from these conditional means.
We fix a good focus, $v\in X$, and $2\le k\le r$.
For $F\in\{F_k(v),\widehat F_k(v)\}$, we put
$\mu=\mathbb E_{Z,v}^0(F)$ and $\Delta_k=\Sigma^{k-1-2\delta_3}$.
Equation~\eqref{eq:general-oriented-mean-upper} gives
$\mu\le 2^{r+2}\Sigma^{k-1}$, and \eqref{eq:general-degree-read} bounds
$\mu_k(v)$. Since $\Sigma\ge b_{5}$,
$(2^{12}r^3+1)\log\Sigma\le
\Sigma^{\delta_1-4r\delta_2-4\delta_3}/[2b_{3}(2^{r+2}+1)]$.
Proposition~\ref{prop:readk} gives
\begin{equation}\label{eq:general-degree-unconditioned-tail}
 \mathbb P_{Z,v}^0(|F-\mu|>\Delta_k)
 \le2\exp\left(-\frac{\Sigma^{\delta_1-4r\delta_2-4\delta_3}}
 {2b_{3}(2^{r+2}+1)}\right)
 \le2\Sigma^{-2^{12}r^3-1}\le\Sigma^{-2^{12}r^3}.
\end{equation}
The bounds $F\le 2^{r+2}\Sigma^{r-1}$ and
$\mathbb P_{Z,v}^0(\psi_v=1)\ge1/8$ imply
\[
 |\mathbb E_{Z,v}^0(F\mid\psi_v=1)-\mu|
 \le \Delta_k+2^{r+5}\Sigma^{r-1-2^{12}r^3}\le2\Delta_k.
\]
Consequently
\begin{equation}\label{eq:general-degree-conditioned-tail}
 \mathbb P_{Z,v}^0\left(
 |F-\mathbb E_{Z,v}^0(F\mid\psi_v=1)|>4\Delta_k
 \,\middle|\,\psi_v=1\right)
 \le8\Sigma^{-2^{12}r^3}\le\Sigma^{-2^{10}r^3}.
\end{equation}

Step~\ref{step:process-reference} gives
\begin{equation}\label{eq:general-local-reference-update}
 \bar d_k'=\omega^{k-1}\sum_{j=k}^r
 \binom{j-1}{k-1}z^{j-k}d_j.
\end{equation}
\begin{lemma}\label{lem:general-degree-round}
Assume $\Sigma\ge b_{5}$ and let $Z$ be a good focus. For $2\le k\le r$, let
\begin{align*}
 R_k&=\omega^{k-1}\sum_{j=k}^r\binom{j-1}{k-1}z^{j-k}E_j
       +b_{3}^{r+4}\bigl(\mathsf E_{\rm acc}\mathcal B_k+\Sigma^{k-1-2\delta_3}\bigr),\\
 \widehat R_k&=\omega^{k-1}(1+b_{3}^{r+4}\mathsf E_{\rm acc})
       \sum_{j=k}^r\binom{j-1}{k-1}z^{j-k}\widehat E_j
       +b_{3}^{r+4}\Sigma^{k-1}\bigl(\mathsf E_{\rm acc}+\Sigma^{-2\delta_3}\bigr).
\end{align*}
\begin{enumerate}
\item For every $v\in X$ and $2\le k\le r$, the joint events below satisfy
\begin{align}
 \mathbb P_Z\bigl(\psi_v=1,\ |F_k(v)-\bar d_k'|>R_k\bigr)
 &\le\Sigma^{-2^{10}r^3},\label{eq:general-local-active-tail}\\
 \mathbb P_Z\bigl(\psi_v=1,\ \widehat F_k(v)-\bar d_k'>\widehat R_k\bigr)
 &\le\Sigma^{-2^{10}r^3}.\label{eq:general-local-full-tail}
\end{align}
\item Let $\mathcal E$ be a family of events containing both events in
\eqref{eq:general-local-active-tail} and \eqref{eq:general-local-full-tail}
for every $v,k$. Suppose $\mathcal E$ satisfies
Lemma~\ref{lem:conditioning-locality} with
\eqref{eq:indicator-locality-parameters} and $\lambda=2$, and that
$\mathbb P_Z(E)\le\Sigma^{-2^{10}r^3}$ for every $E\in\mathcal E$.
Let
\[
 \mathcal G=\bigcap_{E\in\mathcal E}\overline E.
\]
Then $\mathbb P_Z(\mathcal G)>0$. On $\mathcal G$, the degree bounds
\eqref{eq:general-local-active-error-update} and
\eqref{eq:general-local-full-error-update} hold. Moreover,
\begin{align}
 \mathbb E_Z(|X'|\mid\mathcal G)
 &\le\omega\bigl[1+b_{3}^{r+4}(\mathsf E_{\rm acc}+\Sigma^{-4\delta_3})\bigr]
             (1+\Sigma^{-128})|X|,\label{eq:general-one-round-population}\\
 \mathbb E_Z(|Y|\mid\mathcal G)&\le p(1+\Sigma^{-128})|Z|.
 \label{eq:general-one-round-raw}
\end{align}
\end{enumerate}
\end{lemma}
\begin{proof}[Proof of \textup{(1)}]
Let $v\in X$ and $2\le k\le r$. Subtracting
\eqref{eq:general-local-reference-update} from
\eqref{eq:general-active-conditional-mean} and using
$|d_j(v,X\setminus\mathsf W)-d_j|\le E_j$ bounds the difference
between the conditional mean of $F_k(v)$ and $\bar d_k'$.
Adding the deviation $4\Delta_k$, where
$\Delta_k=\Sigma^{k-1-2\delta_3}$, gives a bound at most $R_k$.
Thus \eqref{eq:general-degree-conditioned-tail} implies
\eqref{eq:general-local-active-tail}.

For the full degree, we substitute $d_j(v,X)\le d_j+\widehat E_j$
into \eqref{eq:general-full-conditional-mean}. The additional term
involving the reference degrees is at most
$2^r b_{3}^{r+3}\mathsf E_{\rm acc}\Sigma^{k-1}$.
The inequality $2^{r+4}(1+b_{3}^{r+3}+b_{3})\le b_{3}^{r+4}$
shows that this mean bound, together with $4\Delta_k$, is at most
$\widehat R_k$. Equation~\eqref{eq:general-degree-conditioned-tail}
therefore proves \eqref{eq:general-local-full-tail}.

\medskip
\noindent\emph{Proof of \textup{(2)}.}
On $\mathcal G$, every surviving atom $v\in X'$ satisfies
$|F_k(v)-\bar d_k'|\le R_k$ and
$\widehat F_k(v)-\bar d_k'\le\widehat R_k$.
By \eqref{eq:degree-residual-identities}, maximizing these bounds
over $v\in X'$ proves
\eqref{eq:general-local-active-error-update} and
\eqref{eq:general-local-full-error-update}.

With $Z$ fixed, let $U_{x,1}=\chi_x$ and $U_{x,2}=\mathsf w_x$;
these variables are independent. For either event in
\eqref{eq:general-local-active-tail} or
\eqref{eq:general-local-full-tail}, let $v_E=v$, $S_E=\varnothing$ and
\[
 T_E=\mathbb B_{16r}(v,\mathcal S)\times\{1,2\}.
\]
The formulas in Subsection~\ref{subsec:general-degree-indicators}
show that the event is determined by
$(U_{x,j}:(x,j)\in T_E)$.
Lemma~\ref{lem:conditioning-locality}~(2), with
\eqref{eq:indicator-locality-parameters} and $\lambda=2$, supplies
the dependence bounds for Lemma~\ref{lem:conditional-lll}.
We apply that lemma with $h=\Sigma$, $K=2^9r^3$, $\iota=128$
and $a=2^{10}r^3$. It gives $\mathbb P_Z(\mathcal G)>0$ and bounds
each conditional local expectation by $1+\Sigma^{-128}$ times
its unconditional value.
Applying part \textup{(3)} of the lemma to $\psi_x$ and $\chi_x$,
and summing over $x$, proves
\eqref{eq:general-one-round-population}--\eqref{eq:general-one-round-raw};
the bound for $\psi_x$ is \eqref{eq:general-survival-comparison}.

If every focus failure also has probability at most
$\Sigma^{-2^{10}r^3}$ and satisfies
Lemma~\ref{lem:conditioning-locality} with the same parameters,
Lemma~\ref{lem:conditional-lll}~(3), applied at the focus stage, gives
\[
 \mathbb E(|Y|)\le pq(1+\Sigma^{-128})^2|X\setminus\mathsf W|
                  \le2\eta|X|/\Sigma.
\]
\end{proof}

\subsection{\texorpdfstring{Estimates for $d_{i,t}$, $\Sigma_t$, $\pi_t$ and $\iota_t$}{\detokenize{Estimates for d_i,t, Sigma_t, pi_t and iota_t}}}\label{sec:general-iteration}

We use the reference variables of Subsection~\ref{subsec:actual-process},
with $\Sigma^{-\delta_3/2}\le\eta\le1/4$ at the stages considered.
We write $d_{i,0}=d_i$ and retain the initial value $D=D(\vb d)$
throughout. For $2\le i\le r$ and $t<T$, the reference update is
\begin{equation}\label{eq:general-reference-update}
d_{i,t+1}=\omega_t^{i-1}
 \sum_{j=i}^r\binom{j-1}{i-1}z_t^{j-i}d_{j,t},\qquad z_t=\frac{\eta\sigma_t}{\Sigma_t}.
\end{equation}
In particular, $\Sigma_{t+1}\ge\omega_t\Sigma_t$.
Whenever $E_{i,t}\le\Sigma_t^{i-1-\delta_3}$, we have
$\mathsf E_{\rm acc}\le2\eta\Sigma_t^{-\delta_3}$.
For $\Sigma_t\ge b_{5}$, the bound
$|\log\omega_t+\sum_i d_{i,t}z_t^{i-1}|
\le b_{3}(\mathsf E_{\rm acc}+\Sigma_t^{-4\delta_3})$ gives
\begin{equation}\label{eq:general-omega-estimate}
 \left|\log\omega_t+\sum_{i=2}^r d_{i,t}z_t^{i-1}\right|
 \le4b_{3}\eta\Sigma_t^{-\delta_3}.
\end{equation}
This also holds for the continuation when $X^t$ is empty.

To solve the reference recurrence, we put
\[
f_0(\kappa)=\sum_{i=2}^r d_i \kappa^{i-1},\qquad g_i(\kappa)=\frac{f_0^{(i-1)}(\kappa)}{(i-1)!},\qquad g(\kappa)=\sum_{i=2}^r g_i(\kappa)^{1/(i-1)}.
\]
We recall $\pi_0=1$, $\iota_0=0$, $\pi_{t+1}=\omega_t\pi_t$ and
$\iota_{t+1}=\iota_t+\pi_t z_t$ before stopping.

\begin{lemma}\label{lem:general-generating}
Let $0\le t\le T$, $2\le i\le r$ and $\kappa\in\mathbb R$.
For every reference trajectory,
\begin{align}
 \sum_{i=2}^r d_{i,t}\kappa^{i-1}
 &=f_0(\iota_t+\pi_t\kappa)-f_0(\iota_t),\label{eq:general-generating}\\
 d_{i,t}&=\pi_t^{i-1}g_i(\iota_t),\label{eq:general-degree-formula}\\
 \Sigma_t&=\pi_tg(\iota_t).\notag
\end{align}
For $0\le t<T$,
\begin{equation*}
 \Delta\iota_t=\frac{\eta\sigma_t}{g(\iota_t)}.
\end{equation*}
If \eqref{eq:general-omega-estimate} holds through time $t-1$, then
\begin{equation}\label{eq:general-survival-telescope}
 |\log\pi_t+f_0(\iota_t)|
 \le4b_{3}\sum_{s<t}\eta\Sigma_s^{-\delta_3}.
\end{equation}
\end{lemma}
\begin{proof}
For $0\le t\le T$, we put $f_t(\kappa)=\sum_{i=2}^rd_{i,t}\kappa^{i-1}$.
Expanding $(z_t+\omega_t\kappa)^{j-1}$ in
\eqref{eq:general-reference-update} gives
$f_{t+1}(\kappa)=f_t(z_t+\omega_t\kappa)-f_t(z_t)$.
Induction gives \eqref{eq:general-generating}; coefficient comparison
gives \eqref{eq:general-degree-formula}, hence $\Sigma_t=\pi_tg(\iota_t)$.
Finally $f_t(z_t)=f_0(\iota_{t+1})-f_0(\iota_t)$, so summing
\eqref{eq:general-omega-estimate} proves \eqref{eq:general-survival-telescope}.
\end{proof}

With $\beta_i=d_{i,t}^{1/(i-1)}/\Sigma_t$, we have $\sum_i \beta_i=1$.
Using $\sigma_t$ from \eqref{eq:general-acceptance-reference} and
$-2\kappa\le\log(1-\kappa)\le0$ for $0\le\kappa\le1/2$ gives
\begin{equation}\label{eq:general-sigma-bounds}
e^{-2\eta}\le\sigma_t\le1,\qquad 0\le f_0(\iota_{t+1})-f_0(\iota_t)
 =\sum_i(\eta\sigma_t \beta_i)^{i-1}\le\eta.
\end{equation}

The resulting identities also control the decrease of $\Sigma_t$.

\begin{lemma}\label{lem:general-scale}
Suppose \eqref{eq:general-omega-estimate} holds and
$\Sigma_t\ge b_{5}$ at the stages considered. Let
$t_0=\min\{t:f_0(\iota_t)\ge 2^{4r}\}$, with $\min\varnothing=\infty$.
For $0\le t\le\min\{t_0,T_{\mathrm{init}}\}$ at which these hypotheses hold,
$b_{4}D\le\Sigma_t\le 2^{4r^2}D$.
If they hold through the first $T_{\mathrm{init}}$ rounds, then
$t_0\le T_{\mathrm{init}}$. For $t\ge t_0$ at which they hold,
\begin{equation}\label{eq:general-geometric-D}
e^{-2\eta}\le\frac{\Sigma_{t+1}}{\Sigma_t}
 \le e^{-\eta/8},\qquad \omega_t\le e^{-\eta/8},
\end{equation}
\begin{equation}\label{eq:general-overlap-scale-margin}
\frac{d_{2,t}}{\Sigma_t}\ge\frac12,\qquad \frac{\Sigma_{t+1}}{\omega_t\Sigma_t}\le e^{\eta/2}.
\end{equation}
\end{lemma}
\begin{proof}
For $\iota>0$, $\kappa=f_0(\iota)\ge1$ and $2\le i\le r$,
coefficient comparison gives
\[
\kappa\le \iota f_0'(\iota)\le(r-1)\kappa,\qquad \iota^{i-1}g_i(\iota)\le2^{r-1}\kappa,
\]
and, for $h\ge0$,
\[
g_2(\iota)\le g(\iota)\le g_2(\iota)(1+2^{2r}\kappa^{-1/2}),\qquad \frac{g(\iota+h)}{g(\iota)}\le(1+h/\iota)^{r-2}.
\]
Here $\deg g_i\le r-i$; for $r=2$, $g=g_2$ is constant.
For $\kappa\ge 2^{4r}$ and $h=\Delta\iota_t$, we use
$hf_0'(\iota_t)=\eta\sigma_tf_0'(\iota_t)/g(\iota_t)$,
$g/g_2\le2$ and \eqref{eq:general-sigma-bounds} to obtain
\[
\eta/4\le hf_0'(\iota_t)\le f_0(\iota_t+h)-f_0(\iota_t)\le\eta,\qquad \log\frac{g(\iota_t+h)}{g(\iota_t)}\le\eta/16.
\]
The error in \eqref{eq:general-omega-estimate} is at most $\eta/16$.
Thus $-2\eta\le\log\omega_t\le-3\eta/16$, and
$\log(\Sigma_{t+1}/\Sigma_t)\le-\eta/8$.
Also $g_2/g\ge1/2$, proving the conclusions.

Let $k$ satisfy $d_k=D^{k-1}$. The level $f_0(\iota)=2^{4r}$ is
reached by $\iota=2^{4r}/D$. Until its first crossing, including the
overshoot, $\iota\le(2^{4r}+1)/D$. Coefficient comparison gives
\begin{align*}
 D&\le g(\iota)\le D\sum_{i=2}^r\left[
 \sum_{j=i}^r\binom{j-1}{i-1}(2^{4r}+1)^{j-i}
 \right]^{1/(i-1)}\le2^{4r^2}D,\\
 \frac{\eta e^{-1/2}}{2^{4r^2}D}&\le\Delta\iota_t\le\frac{\eta}{D}.
\end{align*}
Consequently, since $r\ge2$ and $\eta\le1/4$,
\[
 t_0\le\left\lceil\frac{e^{1/2}2^{4r^2+4r}}{\eta}\right\rceil+1
 \le T_{\mathrm{init}}.
\]
At every round satisfying the hypotheses,
\eqref{eq:general-omega-estimate} gives $\omega_t\ge e^{-2\eta}$;
therefore $\pi_t\ge e^{-2\eta t}\ge b_{4}$ for $t\le T_{\mathrm{init}}$.
Now $\Sigma_t=\pi_tg(\iota_t)$ gives the initial bounds, which exclude
a crossing below $b_{5}$ when $b_{4}D\ge b_{5}$.
\end{proof}

\subsection{\texorpdfstring{Bounds for $d_i^{\prime}(A)$ and $d_i^{\prime}(x,y)$}{\detokenize{Bounds for d_i'(A) and d_i'(x,y)}}}\label{sec:general-overlaps}\label{subsec:general-coverage}\label{subsec:general-initial-period}

For rounds $t\ge T_{\mathrm{init}}$, we use $b_{1}\ge1$, $a>0$ and the bounds
\eqref{eq:general-codegree-barriers}--\eqref{eq:general-common-barriers}
and, for $2\le i\le r$ and $v\in X$, assume
\[
d_i(v,X)\le b_{1}\Sigma^{i-1},\qquad E_i\le\Sigma^{i-1-\delta_3},\qquad \frac12\le\sigma\le1,\qquad \Sigma^{-\delta_3/2}\le\eta\le1/4.
\]
Probabilities for a fixed good focus are taken before conditioning
the round. For the estimates in rounds $t\ge T_{\mathrm{init}}$ below,
we also assume $d_2\ge\Sigma/2$ and $\eta\le \frac{\delta_3}{b_{3}}$.

We fix a choice of $i,A$ or $i,x,y$ from
Definition~\ref{def:good-focus}~\eqref{gf:coverage} with $\mathcal B\ne\varnothing$.
We use the associated $\mathcal B,F_0,C_B$, the values $h,m,\zeta$ in
\eqref{eq:coverage-proof-parameters}, and the variables
$\mathsf Z_B^{\rm cov}$ in Definition~\ref{def:coverage-random-variables}.
When $F_0\subseteq X'$, let $\mathcal B'=N_i'(A)$ or
$\mathcal B'=N_i'(x,y)$, respectively. We bound
$|\mathcal B'\cap\mathcal B|$ and $|\mathcal B'\setminus\mathcal B|$ separately.
For $\Sigma\ge b_{5}$ and $\eta\le \frac{\delta_3}{b_{3}}$,
\eqref{eq:general-uniform-acceptance} and \eqref{eq:general-omega-estimate} give
\[
\mathbb P(w\in I)\ge p\sigma(1-4b_{3}\Sigma^{-\delta_3}),\qquad \left|\log\omega+\sum_{j=2}^r z^{j-1}d_j\right|
 \le4b_{3}\eta\Sigma^{-\delta_3}.
\]
Indeed $\mathsf E_{\rm acc}\le2\eta\Sigma^{-\delta_3}$ and
$\Sigma^{-4\delta_3}\le\eta\Sigma^{-\delta_3}$.
The contribution from $j\ge3$ is at most $2\eta^2$, whereas
$z d_2\ge \eta/4$. With $q_0=1-\omega$, the
inequality $|1-e^{-\kappa}-\kappa|\le \kappa^2/2$ for $\kappa\ge0$ therefore gives
\begin{equation}\label{eq:general-graph-hazard}
|q_0-z d_2|\le 4(1+b_{3})\eta(\eta+\Sigma^{-\delta_3}),\qquad \eta/8\le q_0\le2\eta.
\end{equation}

\begin{lemma}\label{lem:general-coverage}
Assume the hypotheses of Lemma~\ref{lem:focus-coverage},
$\Sigma\ge b_{5}$ and $\eta\le \frac{\delta_3}{b_{3}}$.
Let $i,A$ or $i,x,y$ be as in Definition~\ref{def:good-focus}~\eqref{gf:coverage}
and satisfy \eqref{eq:coverage-family-activation}.
We use the associated $\mathcal B,F_0$ and the values $h,m,\zeta$ in
\eqref{eq:coverage-proof-parameters}, and let $Z$ be a focus satisfying
\eqref{eq:general-constant-focused-completions},
\eqref{eq:general-focused-degrees}, and
\eqref{eq:coverage-focus-candidates}--\eqref{eq:coverage-focus-read}.
For every $B\in\mathcal B$,
\begin{equation}\label{eq:general-coverage-mean}
 \mathbb E(\mathsf Z_B^{\rm cov})
 \ge\bigl(m-\frac{\delta_1}{16}+\frac{\delta_3}{b_{3}}-2^8r b_{3}(\eta+\Sigma^{-\delta_3})\bigr)q_0.
\end{equation}
\begin{equation}\label{eq:general-coverage-loss}
 \mathbb P_Z\!\left(\sum_{B\in\mathcal B}\mathsf Z_B^{\rm cov}
 <\bigl(m-\frac{\delta_1}{16}-2^8r b_{3}(\eta+\Sigma^{-\delta_3})\bigr)q_0h\right)
 \le\Sigma^{-a}.
\end{equation}
\end{lemma}
\begin{proof}
By \eqref{eq:coverage-focused-candidate-size}, $n_B=|C_B\cap Z|$ lies between
$(m-\frac{3r\delta_3}{b_{3}}-2r^2(b_{2}+1)\Sigma^{-\delta_2})qd_2$
and $2rqd_2$.
Setting $\chi_x=0$ for $x\in(F_0\cup B)\cap Z$ increases
$\mathbb P(w\in I)$ for $w\in C_B\cap Z$.
We sum these probabilities and subtract
$\mathbb P(w,w'\in I)\le p^2$ for distinct $w,w'\in C_B\cap Z$, to obtain
\[
 \mathbb E(\mathsf Z_B^{\rm cov})
 \ge m z d_2-\frac{3r\delta_3}{b_{3}} z d_2
 -(2r^2(b_{2}+1)+8r b_{3})\Sigma^{-\delta_3}z d_2-2r^2\eta^2.
\]
We use \eqref{eq:general-graph-hazard}, $z d_2\le\eta$ and
$q_0\ge \eta/8$. The coefficients of
$\frac{\delta_3}{b_{3}} q_0$, $\eta q_0$ and $\Sigma^{-\delta_3}q_0$ are at most
$24r$, $32r(1+b_{3})+16r^2$ and
$32r+96r b_{3}+16r^2(b_{2}+1)$, respectively.
These coefficients are bounded by $2^8r b_{3}-1$, $2^8r b_{3}$ and $2^8r b_{3}$,
respectively.

Since $\Sigma\ge b_{5}$,
\[
 \frac{\delta_3^2}{256r b_{1}b_{3}^4}
 \Sigma^{\delta_1-2\delta_2-\delta_3}
 \ge(a+3)\log\Sigma.
\]
For this fixed focus, let
\[
 V=\sum_{B\in\mathcal B}
       \bigl(\mathsf Z_B^{\rm cov}-\mathbb E(\mathsf Z_B^{\rm cov})\bigr),
 \qquad \Delta=\frac{\delta_3}{b_3}q_0h.
\]
Proposition~\ref{prop:readk}, with deviation $\Delta$, mean at most
$h$ and $\eta^2\ge\Sigma^{-\delta_3}$, gives
\begin{equation}\label{eq:general-coverage-deviation}
 \begin{aligned}
 \mathbb P_Z(V<-\Delta)
 &\le2\exp\left(-\frac{\delta_3^2\eta^2h}{256b_3^2k^{\rm cov}}\right)\\
 &\le2\exp\left(-\frac{\delta_3^2\Sigma^{\delta_1-2\delta_2-\delta_3}}
                         {256r b_1b_3^4}\right)\\
 &\le2\Sigma^{-a-3}\le\Sigma^{-a-2}.
 \end{aligned}
\end{equation}
This proves \eqref{eq:general-coverage-loss}.

Finally suppose $F_0\subseteq X'$ and $\mathsf Z_B^{\rm cov}=1$.
A sampled endpoint of $B$ is already outside $X'$. Otherwise
every coordinate forced to zero in its coverage experiment is
zero in the actual experiment: the roots survive and the endpoints
were not sampled. There is therefore an atom $w\in C_B\cap I$;
by the definition of $C_B$, some $u\in B$ satisfies
$\{u,w\}\in\mathcal S$. Hence $u\in\partial I$ and
$B\not\subseteq X'$. Whenever $F_0\subseteq X'$, at most
$h-\sum_B\mathsf Z_B^{\rm cov}$ members of $\mathcal B$ remain in $X'$.
\end{proof}

To bound $|N_i'(A)\setminus N_i(A)|$ and
$|N_i'(x,y)\setminus N_i(x,y)|$, we use the polynomials in
\eqref{eq:proper-birth-witness-polynomial} and
\eqref{eq:common-birth-witness-polynomial}.

\Needspace{10\baselineskip}
\begin{lemma}\label{lem:general-births}
Assume $d_i(v,X)\le b_{1}\Sigma^{i-1}$ for $2\le i\le r$, $v\in X$,
and \eqref{eq:general-codegree-barriers}--\eqref{eq:general-common-barriers}.
Suppose $\Sigma\ge (4rb_{1})^{2^{10}r^2}$ and $0<\eta\le1/4$.
We fix a focus satisfying \eqref{eq:general-birth-focus-derivatives}.
For $3\le i<r$ and $A\subseteq X$ with $2\le|A|<i$,
\begin{equation}\label{eq:general-new-codegree-tail}
 \mathbb P_Z\!\left(A\subseteq X',\
 |N_i'(A)\setminus N_i(A)|>
 b_{1}^2b_{3}^4\Sigma^{i-|A|-\delta_5(i+1)}(\log\Sigma)^{2r}\right)
 \le\Sigma^{-a}.
\end{equation}
For $2\le i<r$ and distinct $x,y\in X$,
\begin{equation}\label{eq:general-new-common-tail}
 \mathbb P_Z\!\left(x,y\in X',\
 |N_i'(x,y)\setminus N_i(x,y)|>
 b_{1}^2b_{3}^4\Sigma^{i-1-\delta_4(i+1)}(\log\Sigma)^{4r}\right)
 \le\Sigma^{-a}.
\end{equation}
\end{lemma}
\begin{proof}
Let $\mathsf Y$ be one of the polynomials in
\eqref{eq:general-birth-witness-polynomials}, of degree $s$ and with
the corresponding bound $\iota$. Let $\mathsf Y_Z$ be obtained by
setting $\xi_x=0$ for $x\notin Z$ and $\xi_x=\chi_x$ for $x\in Z$.
For $T\subseteq Z$ with $|T|\le s$,
\[
 \mathbb E_p(\partial_T\mathsf Y_Z)
 =p^{s-|T|}(\partial_T\mathsf Y)(Z)
 \le b_{1}^2b_{3}^2\iota(\log\Sigma)^s
\]
by \eqref{eq:general-birth-focus-derivatives}; the derivative is zero
if $T\not\subseteq Z$. Proposition~\ref{prop:positive-moments},
with $a+2^{10}r^3$ in place of $a$, gives
\begin{equation*}
 \mathbb P_Z\!\left(\mathsf Y_Z>
 b_{1}^2b_{3}^3\iota(\log\Sigma)^{2s}\right)
 \le\Sigma^{-(a+2^{10}r^3)}.
\end{equation*}
For the polynomials in \eqref{eq:general-birth-witness-polynomials}
with fixed $i,A$, or fixed $i,x,y$, there are at most $2^{2^{12}r^2}$
choices of $k,l$ and $h=|R|$.
Since $2^{2^{12}r^2}\le b_{3}$, the sum of these bounds has
coefficient at most $b_{1}^2b_{3}^4$. A union bound gives
\begin{equation}\label{eq:general-birth-joint-tail}
 \mathbb P_Z\!\left(\exists\mathsf Y:\ \mathsf Y_Z>
 b_{1}^2b_{3}^3\iota(\log\Sigma)^{2s}\right)
 \le2^{2^{12}r^2}\Sigma^{-(a+2^{10}r^3)}
 \le\Sigma^{-a-4}.
\end{equation}
Lemma~\ref{lem:focus-birth-derivatives} specifies the variables
$\chi_x$ occurring in these polynomials, so
Lemma~\ref{lem:conditioning-locality} applies with
\eqref{eq:indicator-locality-parameters}.

Each $B\in N_i'(A)\setminus N_i(A)$ comes from a block
$P\in\mathcal S_k$ with $k>i$ and $P\setminus I=A\cup B$.
For its polynomial, $\iota=\Sigma^{i-|A|-\delta_5(k)}$ and
$s=k-i\le r$, so \eqref{eq:general-birth-joint-tail} gives
\eqref{eq:general-new-codegree-tail}.
Each $B\in N_i'(x,y)\setminus N_i(x,y)$ comes from two blocks of
sizes $k,l$ with $(k,l)\ne(i,i)$; the corresponding bound is
$\iota=\Sigma^{i-1-\delta_4(i+1)}$ and $s\le2r$, giving
\eqref{eq:general-new-common-tail}.
The polynomials count every choice with $F\subseteq Y$ or
$R\cup P\cup Q\subseteq Y$, including choices for which these atoms
are not all in $I$ or the required remaining blocks are absent from
$\mathcal S'$. Distinct source choices are counted separately even
if they leave the same set, so the polynomial sums are upper bounds.
Finally, $\mathcal S_r'\subseteq\mathcal S_r$, and hence
$N_r'(A)\subseteq N_r(A)$ and $N_r'(x,y)\subseteq N_r(x,y)$
whenever the prescribed atoms lie in $X'$.
\end{proof}

Combining the bounds for $|\mathcal B'\cap\mathcal B|$ and
$|\mathcal B'\setminus\mathcal B|$ gives the following result.

\begin{lemma}\label{lem:general-overlap-round}
In a round $t\ge T_{\mathrm{init}}$, assume the hypotheses of
Lemma~\ref{lem:focus-coverage}, $\Sigma\ge b_{5}$,
$\Sigma^{-\delta_3/2}\le\eta\le \frac{\delta_3}{b_{3}}$ and $d_2\ge\Sigma/2$.
Assume $\Sigma'\ge\omega\Sigma$. For $3\le i\le r$ and
$A\subseteq X'$ with $2\le|A|<i$, and for $2\le i\le r$
and distinct $x,y\in X'$, respectively, the two conditioning steps
can impose \eqref{eq:general-next-intersection-barriers} simultaneously.
\end{lemma}
\begin{proof}
Let $(\mathcal B,\mathcal B')=(N_i(A),N_i'(A))$ and let
$(m,\zeta)=(i-|A|,\delta_5(i))$, or let
$(\mathcal B,\mathcal B')=(N_i(x,y),N_i'(x,y))$ and let
$(m,\zeta)=(i-1,\delta_4(i))$, respectively. We put $h=|\mathcal B|$.
If $h\ge\Sigma^{m-\zeta}/2$, then
\eqref{eq:coverage-family-activation} holds and
Lemma~\ref{lem:general-coverage} applies. Since
$m-\zeta\in[3/4,r]$ and $q_0\le1/2$,
\[
 (\Sigma')^{m-\zeta}
 \ge[1-(m-\zeta)q_0-2r^2q_0^2]\Sigma^{m-\zeta}.
\]
Since $\eta\le \frac{\delta_3}{b_{3}}$, $\Sigma^{-\delta_3}\le \frac{\delta_3}{4b_{3}}$ and
$2^{12}r\delta_3=\delta_1$, we have
$\delta_1/16+2^8r b_{3}(\eta+\Sigma^{-\delta_3})\le\delta_1/4$.
On the complement of the event in \eqref{eq:general-coverage-loss},
$|\mathcal B'\cap\mathcal B|\le[1-(m-\delta_1/4)q_0]h$.
The coefficient is nonnegative, so replacing $h$ by
$\Sigma^{m-\zeta}$ leaves a gap of at least
\[
 (\zeta-\delta_1/4-2r^2q_0)q_0\Sigma^{m-\zeta}
 \ge\frac12\zeta q_0\Sigma^{m-\zeta}.
\]
For $i<r$, Lemma~\ref{lem:general-births} gives
\[
 \frac{|\mathcal B'\setminus\mathcal B|}{\Sigma^{m-\zeta}}
 \le b_{1}^2b_{3}^4\Sigma^{-\delta_1}(\log\Sigma)^{4r}
 \le\frac{\delta_1}{32}\Sigma^{-\delta_3/2}
 \le\frac{\delta_1\eta}{32}.
\]
The second inequality uses $\Sigma\ge b_{5}$, which gives
$\frac{32b_{1}^2b_{3}^4}{\delta_1}(\log\Sigma)^{4r}
\le\Sigma^{\delta_1-\delta_3/2}$.
Since $q_0\ge\eta/8$ and $\zeta\ge\delta_1$, this leaves at least
half the gap. For $i=r$, $\mathcal B'\setminus\mathcal B=\varnothing$.

For smaller families, we put $z=2r+8r^2$.
Since $q_0\le2\eta$ and
$(1-q_0)^{m-\zeta}\ge1-rq_0-2r^2q_0^2$,
$(\Sigma')^{m-\zeta}\ge(1-z\eta)\Sigma^{m-\zeta}$.
If $h\le(1-2z\eta)\Sigma^{m-\zeta}$, then
$|\mathcal B'\cap\mathcal B|\le h$ and
$|\mathcal B'\setminus\mathcal B|\le z\eta\Sigma^{m-\zeta}$.
Thus $|\mathcal B'|\le(\Sigma')^{m-\zeta}$. Every other family has
$h\ge\Sigma^{m-\zeta}/2$, already treated, because
$2z\eta\le1/2$. This also covers empty families.
Keeping only the distinct inclusion-minimal sets $S\setminus I$
can only decrease these counts.
The events in \eqref{eq:general-coverage-deviation} and
\eqref{eq:general-birth-joint-tail} have probabilities at most
$\Sigma^{-a-2}$ and $\Sigma^{-a-4}$, respectively.
Their union has probability at most $\Sigma^{-a}$.
Lemma~\ref{lem:conditioning-locality}~(2), with
\eqref{eq:indicator-locality-parameters}, and
Lemma~\ref{lem:conditional-lll} combine them with the degree conditions.
The variables $\mathsf{w}_x$ do not affect $d_i'(A)$ or $d_i'(x,y)$.
\end{proof}

For the initial period, we use the following hypotheses in place of
those for rounds $t\ge T_{\mathrm{init}}$ above.
We fix an integer $T_0\ge0$ and a round $0\le t<T_0$ satisfying
the hypotheses of Lemma~\ref{lem:initial-focus-estimates}, with $k_t$ from
\eqref{eq:general-initial-envelope}.
We fix a focus satisfying \eqref{eq:general-initial-focus-derivatives}.
For $3\le i\le r$ and $A\subseteq X'$ with $2\le|A|<i$, and for
$2\le i\le r$ and distinct $x,y\in X'$, respectively, we establish
\eqref{eq:process-initial-round-bounds}.

For the fixed focus, the bounds
\eqref{eq:general-initial-focus-derivatives} control every derivative
of each polynomial in \eqref{eq:general-birth-witness-polynomials}.
Also $p_t\le b_{4}^{-\delta_2}D^{-\delta_2}
\le b_{4}^{-1}D^{-\delta_2}$.
The choices of $Z$ and $Y$ contribute a factor at most
$b_{4}^{-2(s-\ell_0)}\le b_{4}^{-4r}$.
We apply Proposition~\ref{prop:positive-moments} with $a+2^{10}r^3$ in place
of $a$, and sum over $k,l$ and $h=|R|$. For $3\le i\le r$ and $A\subseteq X$ with $2\le|A|<i$,
\begin{equation}\label{eq:general-initial-new-codegree-tail}
 \mathbb P_Z\!\left(A\subseteq X',\
 |N_i'(A)\setminus N_i(A)|>
 b_{1}^2b_{3}^4b_{4}^{-4r}k_tD^{i-|A|-\delta_0}(\log D)^{4r}\right)
 \le D^{-a-3}\le\Sigma_t^{-a}.
\end{equation}
For $2\le i\le r$ and distinct $x,y\in X$,
\begin{equation}\label{eq:general-initial-new-common-tail}
 \mathbb P_Z\!\left(x,y\in X',\
 |N_i'(x,y)\setminus N_i(x,y)|>
 b_{1}^2b_{3}^4b_{4}^{-4r}k_tD^{i-1-\delta_0}(\log D)^{4r}\right)
 \le D^{-a-3}\le\Sigma_t^{-a}.
\end{equation}
The last inequalities use $D\ge 2^{2ar^2}$.
The bounds $D\ge(4rb_{1})^{2^{10}r^2}$ and
$\Sigma_t\ge b_{4}D\ge(4rb_{1}b_{4}^{-(r-1)})^{2^{10}r^2}$
allow Lemma~\ref{lem:conditioning-locality} with
\eqref{eq:indicator-locality-parameters}, first with
$h=D$, $\lambda=b_{1}$, and then with
$h=\Sigma_t$, $\lambda=b_{1}b_{4}^{-(r-1)}$.
For $(\mathcal B,\mathcal B')=(N_i(A),N_i'(A))$ or
$(N_i(x,y),N_i'(x,y))$, let $m=i-|A|$ or $m=i-1$, respectively.
Using $|\mathcal B'\cap\mathcal B|\le|\mathcal B|$ and the bound on
$|\mathcal B'\setminus\mathcal B|$ gives
\[
 |\mathcal B'|\le k_tD^{m-\delta_0}
 [1+b_{1}^2b_{3}^4b_{4}^{-4r}(\log D)^{4r}]
 \le k_{t+1}D^{m-\delta_0}.
\]
If \eqref{eq:process-initial-round-bounds} fails at surviving atoms,
then the corresponding event in
\eqref{eq:general-initial-new-codegree-tail} or
\eqref{eq:general-initial-new-common-tail} occurs.
Starting with the initial bounds
\eqref{eq:initialP}--\eqref{eq:initialCommon} and $k_0=c$,
these one-step estimates propagate the stronger intersection bounds
whenever their preceding-round hypotheses hold. For the process, we take
$T_0=T_{\mathrm{init}}$: the bounds
\eqref{eq:process-initial-round-bounds} include the transition
$T_{\mathrm{init}}-1\to T_{\mathrm{init}}$.
When $r=2$, $\mathcal S'\subseteq\mathcal S$ and
$N_2'(x,y)\subseteq N_2(x,y)$, so the induction is deterministic.

\subsection{Existence of a good construction}\label{sec:general-error-closure}

We use $E_i,\widehat E_i,\mathsf E_{\rm acc}$ and $\mathcal B_i$
from Subsection~\ref{subsec:process-focus-parameters}. For $2\le i\le r$,
if $E_j\le\Sigma^{j-1}$ for every $2\le j\le r$, then
\begin{equation}\label{eq:general-born-mass}
 \mathcal B_i\le2^{r+1}\eta\Sigma^{i-1}.
\end{equation}
In the eligible-degree recurrence, acceptance errors multiply only
$\mathcal B_i$; in the full-degree recurrence, they multiply the whole degree.

We use $b_{1}=2^{4r^3}$ and $a=2^{10}r^3$ in both periods and set
\[
c_{1}=\frac{2^8b_{3}^{2r+10}}{\delta_3},\qquad c_{2}=2^{16}r\delta_3^{-1}\max\{c,b_{5},b_{6}\}.
\]
Here $b_{5}$ bounds $b_{3}^{2r+10}$, $c_{1}^{r-2}$ and $(2c_{1})^{r-1}$, while
$b_{6}\ge b_{1}^2b_{3}^4b_{4}^{-4r}$.
We assume $\log D\ge c_{2}^{64/\delta_3}$ and use
$\eta=(\log D)^{-\delta_3/8}$. Whenever $\Sigma\ge e^{-1/2}\sqrt{\log D}$,
$\Sigma\ge b_{5}$,
$\Sigma^{-\delta_3/2}\le\eta\le \frac{\delta_3}{b_{3}}$, and
\begin{equation}\label{eq:general-error-smallness}
b_{3}^{2r+10}\eta\le\frac{\delta_3c_{1}^{2-r}}{128},\qquad b_{3}^{2r+10}\Sigma^{-\delta_3}
 \le\frac{\delta_3\eta c_{1}^{2-r}}{128},\qquad (2c_{1})^{r-1}\Sigma^{-\delta_3}\le1.
\end{equation}
We use $\eta\le c_{2}^{-8}$, $\Sigma^{-\delta_3}/\eta\le2c_{2}^{-24}$
and $c_{2}\ge2^{16}r\delta_3^{-1}
\max\{b_{3}^{2r+10},c_{1}^{r-2},(2c_{1})^{r-1}\}$.

\begin{lemma}\label{lem:general-error-closure}
Suppose $\log D\ge c_{2}^{64/\delta_3}$. At every performed round $t$, and for every good
focus $Z$,
\[
\mathbb P(\mathcal G_{{\rm foc},t})>0,\qquad \mathbb P_Z(\mathcal G_{{\rm con},t}(Z))>0.
\]
Moreover, $t_0\le T_{\mathrm{init}}<T$.
For $0\le t\le T_{\mathrm{init}}$, intersection counts with $m$ free atoms
are at most $k_tD^{m-\delta_0}$, with $k_t$ from
\eqref{eq:general-initial-envelope}. Through the first stage with
$\Sigma_t\le\sqrt{\log D}$, the intersection bounds
\eqref{eq:general-codegree-barriers}--\eqref{eq:general-common-barriers}
hold. For $2\le i\le r$, the degree errors satisfy
\begin{equation}\label{eq:general-closed-errors}
E_{i,t}\le c_{1}^{2-i}\Sigma_t^{i-1-\delta_3},\qquad \widehat E_{i,t}\le (2c_{1})^{r-i+1}\Sigma_t^{i-1-\delta_3}.
\end{equation}
\end{lemma}
\begin{proof}
We work up to the first failure, so every local estimate has its hypotheses. The initial degree
errors vanish, and the initial intersections satisfy the stronger
bound with $k_0=c$ by \eqref{eq:initialP}--\eqref{eq:initialCommon}.

We first verify the error recurrences. After $t_0$, Lemma~\ref{lem:general-scale} gives the normalized diagonal
factor
\[
 \frac{\omega^{i-1}\Sigma^{i-1-\delta_3}}
 { (\Sigma')^{i-1-\delta_3}}
 \le\omega^{\delta_3}\le1-\delta_3\eta/16.
\]
The geometric sums satisfy
\begin{equation}\label{eq:general-full-weight-choice}
b_{3}^{2r+10}\sum_{j>i}c_{1}^{2-j}\le\frac{\delta_3c_{1}^{2-i}}{128},\qquad b_{3}^{2r+10}\left(1+\sum_{j>i}(2c_{1})^{r-j+1}\right)
 \le\frac{\delta_3(2c_{1})^{r-i+1}}{128}.
\end{equation}
For this calculation, let $\theta_i=(2c_1)^{r-i+1}$.
Equations~\eqref{eq:general-local-active-error-update}--\eqref{eq:general-local-full-error-update},
\eqref{eq:general-born-mass}, and
$\mathsf E_{\rm acc}\le2\eta\Sigma^{-\delta_3}$ yield
\begin{align}
 \frac{E_i'}{(\Sigma')^{i-1-\delta_3}}
 &\le(1-\delta_3\eta/16)c_1^{2-i}
 +b_3^{2r+10}\left(\eta\sum_{j>i}c_1^{2-j}+\eta^2+\Sigma^{-\delta_3}\right),
 \label{eq:general-active-error-update}\\
 \frac{\widehat E_i'}{(\Sigma')^{i-1-\delta_3}}
 &\le(1-\delta_3\eta/16)\theta_i
 +b_3^{2r+10}\left[\eta\left(1+\sum_{j>i}\theta_j+\Sigma^{-\delta_3}\theta_i\right)
                    +\Sigma^{-\delta_3}\right].
 \label{eq:general-full-error-update}
\end{align}
Equations~\eqref{eq:general-error-smallness} and
\eqref{eq:general-full-weight-choice} place every added term within
$\delta_3\eta/128$ times the corresponding coefficient in \eqref{eq:general-closed-errors}, so both
recurrences preserve \eqref{eq:general-closed-errors}.

During the initial period we use instead
$\mathsf E_t=\max_i(E_{i,t}+\widehat E_{i,t})/\Sigma_t^{i-1}$.
Equations~\eqref{eq:general-local-active-error-update}--\eqref{eq:general-local-full-error-update}
and $\Sigma_t\ge b_{4}D$ give
$\mathsf E_{t+1}\le c_{2}^4(\mathsf E_t+D^{-2\delta_3})$
as long as $\mathsf E_t\le2$. Hence, for $t\le T_{\mathrm{init}}$,
\[
 \mathsf E_t\le(c+t)c_{2}^{4(t+1)}D^{-2\delta_3}
 \le c_{1}^{2-r}2^{-4r^2\delta_3}D^{-\delta_3}.
\]
For the last inequality, the logarithm of the first coefficient
after multiplication by $2^{4r^2\delta_3}/c_{1}^{2-r}$ is at most
$c_{2}^4\sqrt{\log D}\le\delta_3\log D$: we use
$t\le 2^{8r^2}(\log D)^{\delta_3/8}$, $\log\log D\le\sqrt{\log D}$, and
$\log D\ge c_{2}^{64/\delta_3}$. Together with
$\Sigma_t\le 2^{4r^2}D$, this proves the error bounds throughout
the initial period as well.

We next check the initial intersection bounds and the hypotheses
needed to impose them. Lemma~\ref{lem:general-scale} gives $\Sigma_t\le 2^{4r^2}D$ before
$t_0$; \eqref{eq:general-geometric-D} preserves this upper bound afterward.
Also $\omega_t\ge e^{-2\eta}$ and $g(\iota_t)\ge D$, so for
$0\le t\le T_{\mathrm{init}}$,
\[
 b_{4}D\le\Sigma_t\le 2^{4r^2}D.
\]
Since $\log D\ge c_{2}^{64/\delta_3}$ and $c_{2}\ge b_{4}^{-1}$,
$b_{4}D>\sqrt{\log D}$, so stopping cannot intervene. The full-degree hypothesis of
Lemma~\ref{lem:initial-focus-estimates} follows from
\[
 d_i(v,X^t)\le2\Sigma_t^{i-1}
 \le2^{1+4r^2(i-1)}D^{i-1}
 \le b_{1}D^{i-1}.
\]
Moreover $D\ge b_{6}\ge (4rb_{1})^{2^{10}r^2}$ and
$\Sigma_t\ge b_{4}D\ge (4rb_{1}b_{4}^{-(r-1)})^{2^{10}r^2}$,
verifying Lemma~\ref{lem:conditioning-locality} with $h=D$ and $h=\Sigma_t$.
For $t\le T_{\mathrm{init}}$, \eqref{eq:general-initial-envelope} gives
\[
 \log(2b_{4}^{-(r-1)}k_t)
 \le c_{2}^4(\log D)^{3/4}\le\delta_0(\log D)/2.
\]
Here $\log\log D\le\sqrt{\log D}$ and $t\le 2^{8r^2}(\log D)^{\delta_3/8}$ suffice.
Thus Lemma~\ref{lem:initial-focus-estimates} applies. The bounds in
\eqref{eq:process-initial-round-bounds} propagate the stronger estimates
through $T_{\mathrm{init}}$, including
the last transition, and imply
$k_tD^{m-\delta_0}\le\Sigma_t^{m-\zeta}/2$.
Together with $\mathsf E_t\le c_{1}^{2-r}2^{-4r^2\delta_3}D^{-\delta_3}$,
this excludes initial failure and gives $t_0\le T_{\mathrm{init}}$.

For $t\ge T_{\mathrm{init}}$, \eqref{eq:general-overlap-scale-margin} permits
Lemma~\ref{lem:general-overlap-round}; it propagates the absolute
intersection bounds, while
\eqref{eq:general-active-error-update}--\eqref{eq:general-full-error-update}
preserve \eqref{eq:general-closed-errors}.
The bounds \eqref{eq:general-constant-focused-completions}--\eqref{eq:general-full-contractors}
follow from Lemmas~\ref{lem:general-constant-focus} and~\ref{lem:general-focus-estimates},
and Lemma~\ref{lem:degree-indicator-read}~(2) gives
\eqref{eq:general-degree-read}. For $t\ge T_{\mathrm{init}}$,
\eqref{eq:coverage-focus-candidates}--\eqref{eq:coverage-focus-read} and
\eqref{eq:general-birth-focus-derivatives} follow from
Lemmas~\ref{lem:focus-coverage} and~\ref{lem:focus-birth-derivatives};
before $T_{\mathrm{init}}$, we use Lemma~\ref{lem:initial-focus-estimates}.
These are the conditions of Definition~\ref{def:good-focus}.
For each fixed good focus, Lemma~\ref{lem:general-degree-round}
gives Definition~\ref{def:good-construction}~\eqref{gc:degrees}.
For \eqref{gc:intersections}, we use \eqref{eq:process-initial-round-bounds} before
$T_{\mathrm{init}}$ and \eqref{eq:general-next-intersection-barriers} afterward.
Lemma~\ref{lem:conditioning-locality}~(2), with
\eqref{eq:indicator-locality-parameters}, and
Lemma~\ref{lem:conditional-lll} give positive probability and multiply
local expectations by at most $1+\Sigma^{-128}$ at both stages.
Finally, the one-step overshoot is at least
$e^{-2\eta}\sqrt{\log D}\ge e^{-1/2}\sqrt{\log D}$, so the local hypotheses
persist through stopping. This completes the simultaneous induction.
\end{proof}

\section{Completing the process}\label{sec:completing-process}

Lemma~\ref{lem:general-error-closure} shows that, while
$\Sigma_t>\sqrt{\log D}$ and $X^t\ne\varnothing$, a good focus
and a good construction round exist, and the required degree and
intersection bounds remain valid. To complete the proof of
Theorem~\ref{prop:quantitative-main}, and hence
Theorem~\ref{thm:main}, we must show that
\[
 T=\min\{t\in\mathbb Z_{\ge0}:\Sigma_t\le\sqrt{\log D}\},
 \qquad \min\varnothing=\infty,
\]
is a finite stopping time. We then estimate $|J^T|$, $|X^T|$ and
$|X^T\setminus\mathsf W^T|$, with
$|\mathsf W^T|=|X^T|-|X^T\setminus\mathsf W^T|$.
Together with an upper bound for $\sum_{t<T}|Y^t|$, these estimates
allow us to select a run and complete $J^T$ to a maximal independent
set $I$ of the required size, as described in
Subsection~\ref{subsec:process-focus-parameters}.

\subsection{\texorpdfstring{Bounds for stopping time $T$}{Bounds for stopping time T}}\label{sec:general-completion}

We put $\varepsilon=(\log\triangle)^{-\epsilon}$, where $\epsilon=\delta_3/32$.
We put
\[
 c_{3}=c_{2}^{512/\delta_3}.
\]
We assume $\log D\ge c_{3}^{64/\epsilon}$. We set $C=e^{c_{3}^{64/\epsilon}}$.

\Needspace{8\baselineskip}
\begin{definition}\label{def:stopping}
Let $T$ be the first index with $\Sigma_T\le\sqrt{\log D}$, with
$\min\varnothing=\infty$. We use the continuation and freezing
rules of Subsection~\ref{subsec:process-focus-parameters}; the empty-state
continuation has zero survival-log error. Let $\mathcal F_t$
record the choices before round $t$, and let
$\mathcal G_t=\sigma(\mathcal F_t,Z^t)$. Thus
$\mathcal F_t\subseteq\mathcal G_t\subseteq\mathcal F_{t+1}$.
The reference vector and $\Delta \iota_t$ are $\mathcal F_t$-measurable;
$\omega_t$ and $1-\omega_t/s_{x,t}$ for $x\in X^t\setminus \mathsf{W}^t$
are $\mathcal G_t$-measurable.
After $T$, we take samples and $\Delta \iota_t$ to be zero.
\end{definition}

\emph{Remark.} Since $\Sigma_t$ is $\mathcal F_t$-measurable for every
integer $t\ge0$, we have
$\{T\le t\}=\bigcup_{s=0}^t\{\Sigma_s\le\sqrt{\log D}\}\in\mathcal F_t$;
hence $T$ is a stopping time for $(\mathcal F_t)_{t\ge0}$.
Conditioning on $\mathcal F_t$ fixes the earlier
rounds; conditioning on $\mathcal G_t$ also fixes the current focus.
Lemmas~\ref{lem:general-population-cost} and
\ref{lem:general-accepted-count-tracking} use these successive
conditionings to estimate sums over rounds.

Lemma~\ref{lem:general-scale} gives, for $t_0\le s<t\le T$,
\begin{equation}\label{eq:general-geometric-iteration}
 e^{-2\eta(t-s)}\Sigma_s\le\Sigma_t
 \le e^{-\eta(t-s)/8}\Sigma_s.
\end{equation}
Consequently $T\le n_0$, where the deterministic horizon is
\[
 n_0=T_{\mathrm{init}}+1+
 \max\left\{0,\left\lceil\frac8\eta
       \log\frac{2^{4r^2}D}{\sqrt{\log D}}\right\rceil\right\}
 \le\frac{c_{3}\log D}{\eta}.
\]
The lower bound in \eqref{eq:general-geometric-iteration} gives
\begin{equation}\label{eq:general-stopping-scale}
 e^{-2\eta}\sqrt{\log D}\le\Sigma_T\le\sqrt{\log D}.
\end{equation}
For every $\alpha>0$, reversing the upper ratio in
\eqref{eq:general-geometric-iteration} from $T-1$ gives
\begin{equation}\label{eq:general-geometric-sum}
 \sum_{t<T}\Sigma_t^{-\alpha}
 \le\frac{2^{8r^2}b_{4}^{-\alpha}}{\eta}D^{-\alpha}
 +\frac{(\log D)^{-\alpha/2}}{1-e^{-\eta\alpha/8}}.
\end{equation}
The first term covers the initial period, where
$\Sigma_t\ge b_{4}D$; the second is a geometric series.
Taking $\alpha=\delta_3,4\delta_3,128$, respectively, yields
\begin{equation}\label{eq:general-error-sum}
\begin{aligned}
 \sum_{t<T}\eta\Sigma_t^{-\delta_3}&\le c_{3}(\log D)^{-\delta_3/2},\\
 \sum_{t<T}\Sigma_t^{-4\delta_3}&\le c_{3}(\log D)^{-15\delta_3/8},\\
 \sum_{t<T}\Sigma_t^{-128}&\le c_{3}(\log D)^{\delta_3/8-64}.
\end{aligned}
\end{equation}
Here $1-e^{-\kappa}\ge \kappa/2$ for $0\le \kappa\le1$, and for the last
sum $\eta/(1-e^{-16\eta})\le2$ suffices.
The choice of $c_{3}$ also absorbs the coefficient in
\eqref{eq:general-survival-telescope}, giving, for every $t\le T$,
\begin{equation}\label{eq:general-w-exponential}
 |\log\pi_t+f_0(\iota_t)|\le c_{3}(\log D)^{-\delta_3/2}\le1.
\end{equation}

\subsection{\texorpdfstring{Expectations of $|X^T|$ and $\sum_{t<T}|Y^t|$}{Expectations of terminal population and total sample size}}\label{subsec:population-expectations}

\eqref{eq:general-one-round-population} and \eqref{eq:general-closed-errors} give
\[
 \mathbb E(|X^{t+1}|\mid\mathcal G_t)
 \le\omega_t(1+4b_{3}^{r+4}\eta\Sigma_t^{-\delta_3})
                  (1+\Sigma_t^{-128})|X^t|.
\]
We divide by $\omega_t$ before averaging over the focus, and set
\[
\Pi_t=\prod_{j<\min\{t,T\}}
       (1+4b_{3}^{r+4}\eta\Sigma_j^{-\delta_3})(1+\Sigma_j^{-128}),\qquad M_t=\frac{|X^t|}{\pi_t\Pi_t}.
\]
Since $\Pi_{t+1}$ is $\mathcal F_t$-measurable, $(M_t)$ is a
nonnegative supermartingale with $M_0=n$. We apply
\eqref{eq:general-geometric-sum} to
$\log\Pi_t\le4b_{3}^{r+4}\sum_{j<T}\eta\Sigma_j^{-\delta_3}
+\sum_{j<T}\Sigma_j^{-128}$ to obtain
\begin{equation}\label{eq:general-compensator-bounded}
 1\le\Pi_t\le\exp(c_{3}(\log D)^{-\delta_3/2}).
\end{equation}

\begin{lemma}\label{lem:general-population-cost}
If $\iota,\bar\pi$ are deterministic upper bounds for $\iota_T,\pi_T$ on every
history, then
\[
\mathbb E\left(\sum_{t<T}|Y^t|\right)\le e^{c_{3}\eta}n\iota,\qquad \mathbb E(|X^T|)\le e^{c_{3}(\log D)^{-\delta_3/2}}n\bar\pi.
\]
\end{lemma}
\begin{proof}
Applying Lemma~\ref{lem:conditional-lll}~(3) at both conditioning steps gives
\[
\mathbb E(|Y^t|\mid\mathcal F_t) \le\frac{\eta}{\Sigma_t}(1+\Sigma_t^{-128})^2|X^t| =\frac{\Pi_t}{\sigma_t}(1+\Sigma_t^{-128})^2M_t\Delta \iota_t \le e^{c_{3}\eta}M_t\Delta \iota_t.
\]
The last step uses \eqref{eq:general-sigma-bounds},
\eqref{eq:general-compensator-bounded} and $\log D\ge c_{3}^{64/\epsilon}$.
Because $\iota-\iota_{t+1}\ge0$ is $\mathcal F_t$-measurable,
\[
 \mathbb E\bigl((\iota-\iota_{t+1})M_{t+1}+M_t\Delta \iota_t
                  \mid\mathcal F_t\bigr)
 \le(\iota-\iota_t)M_t.
\]
Iterating to $n_0$ and discarding the nonnegative terminal term
proves the sampling bound. For the remaining population we use
$|X^T|\le\bar\pi e^{c_{3}(\log D)^{-\delta_3/2}}M_{n_0}$.
All expectations are finite: the horizon is bounded,
$|X^t|\le n$, and $\pi_t\ge e^{-2\eta n_0}>0$.
\end{proof}

\begin{lemma}\label{lem:general-stopping-cost}
Every history satisfies
\[
\left|\frac{\iota_T}{\rho(\vb d)}-1\right|
 \le c_{3}\frac{\log\log D}{\log D}\le\eta,\qquad \pi_T\le4(\log D)^{-1/2}\rho(\vb d),
\]
and
\begin{equation}\label{eq:general-total-cost}
 \mathbb E\left(\sum_{t<T}|Y^t|+|X^T|\right)
 \le(1+8c_{3}\eta)n\rho(\vb d).
\end{equation}
\end{lemma}
\begin{proof}
For $\kappa=f_0(\iota)\ge1$, a term $d_k=D^{k-1}$ gives $\iota D\le \kappa$.
The coefficient formula for $g_i$ therefore implies
$D\le g(\iota)\le(r-1)2^rD(1+\kappa)^{r-2}$.
We put $\kappa_T=f_0(\iota_T)$. Since $T>t_0$, $\kappa_T\ge1$, and
\eqref{eq:general-w-exponential} and \eqref{eq:general-stopping-scale} give
\[
 -1\le\kappa_T-\log D+\frac12\log\log D
 \le(r-2)\log(1+\kappa_T)+\log((r-1)2^r)+2.
\]
The function $\beta-(r-2)\log(1+\beta)$ is increasing beyond $2\log D$,
where its value exceeds $\log D-\tfrac12\log\log D+\log((r-1)2^r)+2$. Thus
$(\log D)/2\le \kappa_T\le2\log D$, and substitution gives
\begin{equation}\label{eq:general-terminal-hazard}
 |\kappa_T-\log D|\le c_{3}\log\log D.
\end{equation}
For $\beta\ge0$, $f_0(\lambda \beta)\ge\lambda f_0(\beta)$ when $\lambda\ge1$,
with the reverse inequality for $0\le\lambda\le1$.
As $f_0(\rho)=\log D$, this proves the bound for $\iota_T$.
Also $g(\iota_T)\ge f_0'(\iota_T)\ge \kappa_T/\iota_T$, so
$\pi_T=\Sigma_T/g(\iota_T)\le4(\log D)^{-1/2}\rho$.
We apply Lemma~\ref{lem:general-population-cost} with these bounds.
The inequalities $c_{3}\eta\le1$, $e^\beta\le1+2\beta$ for $0\le \beta\le1$,
and $(\log D)^{-1/2}\le\eta$ prove \eqref{eq:general-total-cost}.
\end{proof}

\subsection{\texorpdfstring{Concentration of $|X^t\setminus \mathsf{W}^t|$ and $|J^T|$}{\detokenize{Concentration of |X^t minus W^t| and |J^T|}}}

We put
\begin{align*}
 \varepsilon_1&=c_{3}(\log D)^{-15\delta_3/8},\\
 \varepsilon_2&=c_{3}(\log D)^{\delta_3/8-64},\\
 \Gamma_{\rm pop}&=c_{3}^2\left((\log D)^{-15\delta_3/8}
                         +(\log D)^{c_{3}+2}D^{-\delta_0}\right).
\end{align*}
For $v,w\in X$, let $F$ be an event determined by the variables
indexed by
$T_F=(\mathbb B_{16r}(v,\mathcal S)\cup\mathbb B_{16r}(w,\mathcal S))\times\{1,2\}$.
Lemma~\ref{lem:conditioning-locality}~(2), with
\eqref{eq:indicator-locality-parameters}, gives at most
$2\Sigma^{2^9r^3}\le\Sigma^{2^9r^3+1}$ events $E$ with
$T_E\cap T_F\ne\varnothing$.
We use $h=\Sigma$, $K=2^9r^3+1$, $\iota=128$ and $a=2^{10}r^3$
in Lemma~\ref{lem:conditional-lll}~(2), since
$2^{10}r^3\ge2^9r^3+132$. Applying it to $F$ and $\overline{F}$
at either conditioning gives
$\mathbb P(F)-\Sigma^{-128}\le\mathbb P(F\mid\mathcal G)
\le(1+\Sigma^{-128})\mathbb P(F)$.

\begin{lemma}\label{lem:general-active-population-tracking}
The conditioned process satisfies
\[
 \mathbb P\left(
 \max_{t\le T}\left|\frac{|X^t\setminus \mathsf{W}^t|}{n\pi_t}-1\right|
 >\varepsilon/16\right)\le\varepsilon/16.
\]
\end{lemma}
\begin{proof}
We choose $k$ with $d_k=D^{k-1}$. For $k\ge3$, summing proper
pair-codegrees gives
$(k-1)D^{k-1}\le(n-1)cD^{k-2-\delta_0}$.
For $k=2$, summing common-neighbour counts gives
$D(D-1)\le(n-1)cD^{1-\delta_0}$.
Consequently
\begin{equation}\label{eq:general-population-room}
 n/D\ge(2c)^{-1}D^{\delta_0}.
\end{equation}
Equations~\eqref{eq:general-terminal-hazard} and
\eqref{eq:general-w-exponential}, together with $\iota_t\le \iota_T$, give
$\pi_t\ge\pi_{\min}:=e^{-1}D^{-1}(\log D)^{-c_{3}}$.

We put $U_t=|X^t\setminus \mathsf{W}^t|/\pi_t$, frozen after $T$.
For $t<T$, before conditioning on $\mathcal G_{{\rm con},t}(Z^t)$,
$\mathbb P_{Z^t}(\psi_x(1-\mathsf{w}_x)=1)=\omega_t\ge e^{-1/2}$ for
$x\in X^t\setminus \mathsf{W}^t$. Dividing by
$\pi_t\omega_t$ before averaging over the focus gives
\begin{equation}\label{eq:general-active-first-moment}
 (1-2\Sigma_t^{-128})U_t\le\mathbb E(U_{t+1}\mid\mathcal F_t)
 \le(1+\Sigma_t^{-128})U_t.
\end{equation}
Equation~\eqref{eq:general-survival-joint}, the independence of the
$\mathsf{w}_{x,t}$ from the $\chi_{x,t}$, and
Lemma~\ref{lem:conditional-lll}~(2), applied to
$\mathbb B_{16r}(v,\mathcal S)\cup\mathbb B_{16r}(w,\mathcal S)$, give
\begin{equation}\label{eq:general-active-second-moment}
 \mathbb E(U_{t+1}^2\mid\mathcal F_t)
 \le(1+\Sigma_t^{-4\delta_3}+2\Sigma_t^{-128})U_t^2+2U_t/\pi_t.
\end{equation}
Since $\Sigma_t\ge b_{5}$, the relative error in
\eqref{eq:general-survival-joint} is at most $\Sigma_t^{-4\delta_3}$.

We set $\Lambda_0=1$ and
$\Lambda_{t+1}=\Lambda_t\mathbb E(U_{t+1}\mid\mathcal F_t)/U_t$ when
$U_t>0$ and $t<T$, keeping $\Lambda$ unchanged otherwise.
Then ${\mathsf M}_t=U_t/\Lambda_t$ is a martingale with ${\mathsf M}_0=n$.
Since $(\log D)^\epsilon\ge c_{3}^{64}$, we have
$\varepsilon_2=c_{3}(\log D)^{\delta_3/8-64}\le1/8$.
Equations~\eqref{eq:general-error-sum} and
\eqref{eq:general-active-first-moment} therefore give
\[
e^{-4\varepsilon_2}\le\Lambda_t\le e^{\varepsilon_2},\qquad |\Lambda_t-1|\le8\varepsilon_2.
\]
For $t<T$, we have $\Sigma_t^{-128}\le1/8$, so
$(1-2\Sigma_t^{-128})^{-2}\le1+8\Sigma_t^{-128}$.
Thus \eqref{eq:general-active-second-moment} implies
\[
 \mathbb E({\mathsf M}_{t+1}^2\mid\mathcal F_t)
 \le(1+32\Sigma_t^{-4\delta_3}){\mathsf M}_t^2+8{\mathsf M}_t/\pi_{\min}.
\]
We divide by $\prod_{j<\min\{t,T\}}(1+32\Sigma_j^{-4\delta_3})$,
which is at most $e^{32\varepsilon_1}$ by
\eqref{eq:general-error-sum}, and use ${\mathsf M}_t={\mathsf M}_T$ for $t\ge T$.
Since $32\varepsilon_1=32c_{3}(\log D)^{-15\delta_3/8}\le1$ and
$\mathbb E({\mathsf M}_t)=n$, it follows that
\[
 \frac{\operatorname{Var}({\mathsf M}_{n_0})}{n^2}
 \le64\varepsilon_1+\frac{32n_0}{n\pi_{\min}}
 \le \Gamma_{\rm pop}.
\]
Integrability follows from the bounded horizon and positive lower
bounds for $\pi_t,\Lambda_t$.
The assumptions $(\log D)^\epsilon\ge c_{3}^{64}$ and $c_{3}\ge\delta_0^{-1}$,
with $\log\log D\le\sqrt{\log D}$, give
$2^{17}c_{3}^5(\log D)^{c_{3}+2+3\epsilon}D^{-\delta_0}\le1$.
Together with $\varepsilon\ge c_{3}^{-1}(\log D)^{-\epsilon}$ and
$15\delta_3/8=60\epsilon$, this implies
$\Gamma_{\rm pop}\le\varepsilon^3/2^{16}$.
Stopping $({\mathsf M}_t-n)^2$ at its first crossing of
$(\varepsilon n/64)^2$ gives
\[
 \mathbb P\left(\max_{t\le n_0}|{\mathsf M}_t-n|>
                  \varepsilon n/64\right)
 \le 2^{12}\Gamma_{\rm pop}/\varepsilon^2\le\varepsilon/16.
\]
Since $(\log D)^\epsilon\ge c_{3}^{64}$ and
$\varepsilon\ge c_{3}^{-1}(\log D)^{-\epsilon}$, we also have
$\varepsilon_2=c_{3}(\log D)^{\delta_3/8-64}\le\varepsilon/256$.
On the complement of this event,
$|U_t/n-1|\le2\varepsilon/64+8\varepsilon_2\le\varepsilon/16$,
as required.
\end{proof}

\begin{lemma}\label{lem:general-accepted-count-tracking}
The accumulated accepted set satisfies
\[
 \mathbb P\left(\left|\frac{|J^T|}{n\rho(\vb d)}-1\right|
                       >\varepsilon/2\right)\le\varepsilon/8.
\]
\end{lemma}
\begin{proof}
Let $m_t=\mathbb E(|I^t|\mid\mathcal F_t)$, with $m_t=0$ after $T$.
We apply \eqref{eq:general-uniform-acceptance} and
Lemma~\ref{lem:conditional-lll}~(2) at both steps to each eligible atom. The acceptance errors are small relative to $p\sigma$, since
$\Sigma^{-128}/(p\sigma)\le2\Sigma^{-127}$,
$\Sigma^{-128}/q\le\Sigma^{-127}$ and
$\Sigma^{-127}\le\eta\Sigma^{-\delta_3}$. Thus
\begin{equation}\label{eq:general-accepted-conditional-mean}
 \bigl|m_t-z_t|X^t\setminus \mathsf{W}^t|\bigr|
 \le c_{3}\eta\Sigma_t^{-\delta_3}
             z_t|X^t\setminus \mathsf{W}^t|.
\end{equation}
For $t<T$, its relative error satisfies
\[
 c_{3}\eta\Sigma_t^{-\delta_3}
 \le c_{3}e^{\delta_3/2}(\log D)^{-5\delta_3/8}\le\varepsilon/16,
\]
using $\eta=(\log D)^{-4\epsilon}$, $\varepsilon\ge c_{3}^{-1}(\log D)^{-\epsilon}$
and $(\log D)^\epsilon\ge c_{3}^{64}$.
For distinct eligible $u,v$, the inclusion $I^t\subseteq Y^t$
and Lemma~\ref{lem:conditional-lll}~(2) for two radius-$16r$ balls give
$\mathbb P(u,v\in I^t\mid\mathcal F_t)\le6z_t^2$.
Adding the diagonal yields
\[
 \operatorname{Var}(|I^t|\mid\mathcal F_t)
 \le2z_t|X^t\setminus \mathsf{W}^t|+6z_t^2|X^t\setminus \mathsf{W}^t|^2.
\]
Let $\tau$ be the first failure of the bound in
Lemma~\ref{lem:general-active-population-tracking},
or $n_0+1$ if none occurs. The martingale
\[
 {\mathsf N}_s=\sum_{t<s}\vb1_{\{t<\tau\}}(|I^t|-m_t)
\]
has conditional variance per increment at most
$4n\Delta \iota_t+24n^2(\Delta \iota_t)^2$.
Using $\Delta \iota_t\le\eta/D$, $\iota_T\le2\rho$ and $\rho\ge D^{-1}$ gives
\[
 \frac{\mathbb E({\mathsf N}_{n_0}^2)}{(n\rho)^2}
 \le\frac{8}{n\rho}+\frac{48\eta}{D\rho}
 \le16cD^{-\delta_0}+64\eta\le\varepsilon^3/256.
\]
Here $\rho\ge D^{-1}$ follows from $f_0(D^{-1})\le r-1<\log D$,
and \eqref{eq:general-population-room} bounds $n/D$.
The last inequality follows from $c_{3}\ge\max\{c,\delta_0^{-1}\}$,
$(\log D)^\epsilon\ge c_{3}^{64}$ and $\varepsilon\ge c_{3}^{-1}(\log D)^{-\epsilon}$.
Thus $\mathbb P(|{\mathsf N}_{n_0}|>\varepsilon n\rho/4)
\le16\mathbb E({\mathsf N}_{n_0}^2)/(\varepsilon n\rho)^2
\le\varepsilon/16$.
Outside this event and the population failure, also of probability
at most $\varepsilon/16$, we use
\eqref{eq:general-accepted-conditional-mean},
Lemmas~\ref{lem:general-active-population-tracking}, \ref{lem:general-stopping-cost},
and $\eta=(\log D)^{-4\epsilon}\le\varepsilon/16$ to obtain
\[
 (1-\varepsilon/16)^3n\rho
 \le\sum_{t<T}m_t\le(1+\varepsilon/16)^3n\rho.
\]
Their deviation from $n\rho$ is at most $\varepsilon n\rho/4$.
Adding the martingale fluctuation proves the assertion.
\end{proof}

\subsection{\texorpdfstring{Bounds for $|I|$}{\detokenize{Bounds for |I|}}}
\label{subsec:quantitative-assembly}

\begin{lemma}\label{lem:bookkeeping}
After any $T$ rounds, $J^T$ is independent in the original system
and has a maximal independent extension satisfying
\begin{equation}\label{eq:charge}
 |J^T|\le|I|\le\sum_{t<T}|Y^t|+|X^T|.
\end{equation}
\end{lemma}
\begin{proof}
Inductively, each current block is the remainder $S\setminus J^t$
of an original block, and every original remainder contained in
$X^t$ has size at least two and contains a current block.
Independence of $I^t$ prevents an empty next remainder; a singleton
would be deleted through $\partial I^t$. Removing repeated or
redundant blocks preserves both assertions, proving independence
of $J^t$ at every step.
Every unused atom outside $X^T\cup\bigcup_{t<T}Y^t$ was removed
through some $\partial I^t$ and completes an
original block with atoms of $J^T$. Greedily extend $J^T$ inside
$X^T\cup\bigcup_{t<T}(Y^t\setminus I^t)$.
The resulting set is maximal, and disjointness of the samples
proves \eqref{eq:charge}.
\end{proof}

\begin{proof}[Proof of Theorem~\textup{\ref{prop:quantitative-main}}]
We take $C=e^{c_3^{64/\epsilon}}$ and use the conditioning events
in Definitions~\ref{def:good-focus} and~\ref{def:good-construction}.
For a run of the process, let
$\mathsf Z=\sum_{t<T}|Y^t|+|X^T|$.
By \eqref{eq:general-total-cost} and
$8c_3\eta=8c_3(\log D)^{-4\epsilon}\le\varepsilon/4$,
we have $\mathbb E(\mathsf Z)\le(1+\varepsilon/4)n\rho$.
Since $0<\varepsilon\le1$, it follows that
\[
 \mathbb P(\mathsf Z\le(1+\varepsilon)n\rho)
 \ge1-\frac{\mathbb E(\mathsf Z)}{(1+\varepsilon)n\rho}
 \ge\frac{3\varepsilon}{4(1+\varepsilon)}\ge\varepsilon/4.
\]
Lemma~\ref{lem:general-accepted-count-tracking} gives
$\mathbb P(|J^T|<(1-\varepsilon)n\rho)\le\varepsilon/8$.
Subtracting this probability, we obtain
\[
 \mathbb P\bigl(|J^T|\ge(1-\varepsilon)n\rho,
                \ \mathsf Z\le(1+\varepsilon)n\rho\bigr)
 \ge\varepsilon/8.
\]
Thus a run meets both conditions in the selection rule of
Subsection~\ref{subsec:process-focus-parameters} with probability
at least $\varepsilon/8$.

Each run terminates almost surely: $T\le n_0$ by
Subsection~\ref{sec:general-completion}, and every conditioning event
has positive probability. Repeating independent runs therefore
terminates almost surely, with at most $8/\varepsilon$ expected runs.
For the selected run, Lemma~\ref{lem:bookkeeping} gives a maximal
independent extension $I$ of $J^T$ satisfying
\[
 (1-\varepsilon)n\rho\le|J^T|\le|I|\le\mathsf Z
 \le(1+\varepsilon)n\rho.
\]
This completes the proof of Theorem~\ref{prop:quantitative-main}.
\end{proof}

\begin{proof}[Proof of the explicit bound on $K$]
In the deduction of Theorem~\ref{thm:main}, we use $c=1$ and
\[
 D_0=\exp\!\left(c_3^{64/\epsilon}\right),\qquad
 K=\frac{(r-1)^2D_0(\log D_0)^\epsilon}{\log2}.
\]
It therefore suffices to bound $c_3$. We put $R=4r/\delta_0\ge16$
and write $\exp_j$ for the $j$-fold iterated exponential.

We first estimate $b_2$ and $b_3$, which determine $b_5$.
Since $r\le R/8$, $a=2^{10}r^3$ and
$\delta_2^{-1}=256rR\le32R^2$, the definition of $b_2$ gives
\[
 \log b_2
 =r\log\left(\frac{4r(2^{11}r^3+1)}{\delta_2}\right)
 \le r\log(80R^6)\le R^2.
\]
Using $b_2+1\le2b_2$, $2^{10}r^4\le R^4/4$ and
$\log(2^{17}r)\le R^2$ in the definition of $b_3$, we obtain
\[
 \log b_3
 =2^{10}r^4\log\bigl(2^{16}r(b_2+1)\bigr)\le R^6.
\]
Now $\log b_1=4r^3\log2\le R^3$ and
$\delta_3^{-1}=2^{14}rR\le2^{11}R^2$, so
\[
 \log b_5=\frac{16}{\delta_3}(\log b_1+\log b_3)
 \le2^{16}R^8\le R^{12}.
\]

The other constant needed for $c_2$ is $b_6$.
Here $a+1\le3R^3$, $1+\log b_1\le2R^3$,
$\delta_0^{-1}\le R/8$ and $2^{16r^2}\le e^{R^2/4}$ give
\[
 \log b_6=\frac{2^{16r^2}(a+1)(1+\log b_1)}{\delta_0}
 \le R^7e^{R^2/4}\le e^{R^2}.
\]
Since $R^{12}\le e^{R^2}$ and
$\log(2^{16}r/\delta_3)\le R^2$, we can now substitute the bounds
for $b_5,b_6$ directly into $c_2$:
\[
 \log c_2=\log\left(\frac{2^{16}r}{\delta_3}\right)
       +\max\{\log b_5,\log b_6\}
 \le R^2+e^{R^2}\le e^{2R^2}.
\]
Using $512/\delta_3\le2^{20}R^2$ in $c_3=c_2^{512/\delta_3}$
then yields
\[
 \log c_3\le2^{20}R^2e^{2R^2}\le e^{3R^2}.
\]
Finally, $\epsilon^{-1}=32/\delta_3\le2^{16}R^2$ gives
\[
 \log\log D_0=\frac{64}{\epsilon}\log c_3
 \le2^{22}R^2e^{3R^2}\le e^{4R^2}.
\]

To obtain the bound on $K$, we use its formula at the start of
the proof. Since $\epsilon\le1$ and
$\log D_0=c_3^{64/\epsilon}\ge R^2$ (as $c_2\ge R$),
\[
 \log K\le2\log D_0
 \le2\exp_2(4R^2)<\exp_3(R^{16}).
\]
Thus $K\le\exp_4(R^{16})=\exp_4((4r/\delta_0)^{16})$.
This completes the proof of the stated bound on $K$.
\end{proof}

For a $d$-regular $k$-uniform system, we set $d_k=d$ and $d_i=0$
for $i\ne k$. Then
\[
D=d^{1/(k-1)},\qquad \rho=\left(\frac{\log d}{(k-1)d}\right)^{1/(k-1)}.
\]
The two uniform sparsity bounds become
$D^{k-|A|-\delta}$ and $D^{k-1-\delta}$, respectively.
Theorem~\ref{thm:main} applies for every $d>1$.
We absorb $(k-1)^\epsilon$ from
$(\log D)^{-\epsilon}=(k-1)^\epsilon(\log d)^{-\epsilon}$
into its coefficient, and write $r=k-1$.
This proves Theorem~\ref{main}, with
$\epsilon=\min\{\delta,1/2\}/(2^{21}(r+1)^2)$.

\section{Set systems with low maximum degree}\label{sec:maxdegree}

For a nonregular $r$-uniform set system of maximum degree $\triangle$,
we use the uniform sparseness bounds of Section~1 with $\triangle$
in place of $d$. Atom sets are denoted by $X(\mathcal S)$, and
$d_r(A)$ and $d_r(x,y)$ retain their respective codegree and
common-link meanings.

\begin{proposition}\label{prop:regular-extension}
Let $r\ge2$, $0<\delta\le1$ and $\triangle\ge2$, and let
$\mathcal S$ be an $n$-atom $r$-uniform $\delta$-sparse set system
with maximum degree at most $\triangle$. If $r=2$, suppose also that
$\delta\le1-(\log2)/(\log\triangle)$.
For every positive integer $N$, there is a $\triangle$-regular
$\delta$-sparse system $\mathcal S^+$ consisting of atom-disjoint
copies of $\mathcal S$ together with further blocks, such that
\[
 N\le |X(\mathcal S^+)|<N+8rn(\triangle+r)^2.
\]
For fixed $r$, its incidence list can be constructed deterministically
in time polynomial in $n$, $\triangle$ and $N$.
\end{proposition}

\begin{proof}
Let $p$ be the least prime at least $\max\{r,\triangle\}$.
Bertrand's postulate gives $p<2\max\{r,\triangle\}$, and trial
division finds $p$ in time polynomial in $r+\triangle$.
For each integer $0\le i\le\triangle$, let $A_i=\{0,\ldots,i-1\}$
in $\mathbb F_p$; in particular, $A_0=\varnothing$.
We first construct $i$-regular auxiliary systems on a common atom set.

If $r\ge3$, the atom set is $\{0,\ldots,r-1\}\times\mathbb F_p$,
and the blocks are, for $x\in\mathbb F_p$ and $a\in A_i$,
\[
 \{(j,x+ja):0\le j<r\}.
\]
Every atom lies in exactly one block for each $a\in A_i$.
Two atoms in different parts determine $x$ and $a$ uniquely,
since their part indices differ by a nonzero element of $\mathbb F_p$.
Thus distinct blocks meet in at most one atom. In particular, every
proper codegree is at most one and every common $r$-link count is zero.
The common number of atoms in this case is $M_0=rp$.

If $r=2$, the two parts of the auxiliary graph are disjoint copies
of $\mathbb F_p^2$, denoted by subscripts $0$ and $1$. Its edges are
\[
 \bigl\{\{(x,y)_0,(x+a,y+ax)_1\}:
                 x,y\in\mathbb F_p,\ a\in A_i\bigr\}.
\]
For each $a$, the map $(x,y)\mapsto(x+a,y+ax)$ is a bijection,
and these matchings are edge-disjoint, so the graph is $i$-regular.
Two distinct left vertices $(x,y)_0,(u,v)_0$ with a common
neighbor satisfy, for the corresponding $a,b\in A_i$,
\[
 x+a=u+b,\qquad y+ax=v+bu.
\]
If $x=u$, these equations force $y=v$; otherwise they determine
$a$ uniquely. Thus any two distinct left vertices have at most one
common neighbor. The graph has no $4$-cycle, so the same bound
holds for two right vertices; vertices in opposite parts have no
common neighbor. Here $M_0=2p^2$.

Let $M=M_0\lceil N/(nM_0)\rceil$. For each $i$, let
$\mathcal H_i$ be the disjoint union of $M/M_0$ copies of the
corresponding auxiliary system, indexed on $[M]$.
Then $nM\ge N$ and $nM<N+nM_0<N+8rn(\triangle+r)^2$.
Let $\mathcal T_1,\ldots,\mathcal T_M$ be atom-disjoint copies
of $\mathcal S$, and let $v_j$ be the copy of $v$ in $\mathcal T_j$.
For every $v\in X(\mathcal S)$, let $i(v)=\triangle-d_r(v)$,
and place $\mathcal H_{i(v)}$ on $\{v_1,\ldots,v_M\}$.
The resulting system is
\[
 \mathcal S^+=\bigcup_{j=1}^M\mathcal T_j
       \;\cup\!\bigcup_{v\in X(\mathcal S)}\mathcal H_{i(v)}(v).
\]
Every atom has degree $\triangle$.

For $2\le|A|<r$, all blocks containing $A$ lie in a single copy
$\mathcal T_j$ or a single auxiliary system, or there are no such
blocks. Since $\delta\le1$, the preceding construction gives
\[
 d_r(A)\le\max\{1,\triangle^{(r-|A|-\delta)/(r-1)}\}
          =\triangle^{(r-|A|-\delta)/(r-1)}.
\]
For $r\ge3$, a common link has at least two atoms. Such a link
cannot belong both to a block in a copy $\mathcal T_j$ and to
a block in an auxiliary system: the former contains at most one
copy of each original atom. Thus common links come entirely from
one copy or one auxiliary system, and
\[
 d_r(x,y)\le\triangle^{1-\delta/(r-1)}.
\]
For $r=2$, two atoms in one copy have only the common neighbors
already in that copy, and two copies of the same original atom
have at most one common neighbor in the auxiliary graph. For
$x=v_j$ and $y=w_k$ with $v\ne w$ and $j\ne k$, the only possible
common neighbors are $w_j$ and $v_k$. Hence
\[
 d_2(x,y)\le\max\{\triangle^{1-\delta},2\}
             =\triangle^{1-\delta},
\]
where the last equality uses $2\le\triangle^{1-\delta}$.
This proves sparseness. All displayed blocks and copies can be
enumerated by arithmetic modulo $p$; their number is
$nM\triangle/r$. This also proves the stated construction bound.
\end{proof}

\begin{proposition}
Suppose there exists $n_0=n_0(r,\delta)$ such that, for $n\ge n_0$,
every $d$-regular $\delta$-sparse $r$-uniform $n$-atom system
contains an independent set of size at least $f(d,r)n$.
Let $\mathcal S$ be a $\delta$-sparse $r$-uniform $m$-atom system
with maximum degree at most $\triangle\ge2$, where $0<\delta\le1$.
If $r=2$, suppose also that
$\delta\le1-(\log2)/(\log\triangle)$.
Then $\mathcal S$ contains an independent set of size at least
$f(\triangle,r)m$.
\end{proposition}

\begin{proof}
We apply Proposition~\ref{prop:regular-extension} with $N=n_0$.
An independent set $I$ in $\mathcal S^+$ of size at least
$f(\triangle,r)|X(\mathcal S^+)|$ has an intersection of size at
least $f(\triangle,r)m$ with one of the copies of $\mathcal S$,
by averaging. This intersection is independent in that copy.
\end{proof}

For the application to Theorem~\ref{thm:maxdegree}, let
$\delta'=\min\{\delta,1/2\}$ and suppose $\triangle\ge4$.
The original system is $\delta'$-sparse, and
$2\le\triangle^{1-\delta'}$, so the graph exception is excluded.
Applying Theorem~\ref{main} to the regular extension, with uniformity
$r+1$, gives an independent set of size at least
\[
 (1-O((\log\triangle)^{-\epsilon}))n
       \left(\frac{\log\triangle}{r\triangle}\right)^{1/r}.
\]
The constants depend only on $r$ and $\delta$; increasing the error
constant covers the bounded values $\triangle<4$ as well.
This proves Theorem~\ref{thm:maxdegree}.

For its algorithmic assertion we take $N=n$ in
Proposition~\ref{prop:regular-extension} and apply
Corollary~\ref{cor:deterministic-main} to $\mathcal S^+$.
The bound $\triangle\le\binom{n-1}{r}$ and the explicit size bound
on $\mathcal S^+$ show that its incidence list and the resulting
computation have size and running time polynomial in $n$, for fixed
$r$ and $\delta$. A scan of the copies finds the largest intersection
with the constructed independent set. For bounded $\triangle$ a
greedy independent set suffices for the stated error bound.
To the best of our knowledge, the resulting bound improves those
in~\cite{Dhawan,DMV,NV,KMV,RS,VWilson}.

\Needspace{8\baselineskip}
\section{Zero sets of polynomials}\label{sec:polynomial-zero-sets}

In this section, we prove Theorem~\ref{thm:poly-maximal-sets}. The natural set system to consider for a general polynomial $f \in \mathbb Z[x_1,x_2,\dots,x_r]$ is described as follows. Let $X_q = \{a \in \mathbb F_q : f(a,a,\dots,a) \neq 0\}$. Let $\mathcal{S} = \mathcal{S}_f$ be the set system  with atom set $X_q$ whose blocks are the inclusion-minimal subsets $S \subseteq X_q$ for which $f(x_1,x_2,\dots,x_r) = 0$ for some $(x_1,x_2,\dots,x_r) \in S^r$.
Evidently, a maximal independent set in $\mathcal{S}$ is a maximal $f$-free subset of $\mathbb F_q$. Furthermore, if $|X_q| \geq r$, and the affine transformations preserving $f$ up to a nonzero scalar multiple act transitively on $X_q$, then it turns out $\mathcal{S}$ is a $(|X_q|,\vb{d})$-antichain. 
For Theorem \ref{thm:poly-maximal-sets}, we have the additional conditions that $f$ is homogeneous, symmetric, and absolutely irreducible. If $f(1,\ldots,1)\ne0$ in $\mathbb F_q$, then $X_q = \mathbb F_q^{\times}$; otherwise $X_q=\varnothing$ and $\gamma_q(f)=0$. In the nonempty case we will show that $\mathcal{S}$ (or a suitable modification of $\mathcal{S}$) is $\delta$-sparse, so that Theorem \ref{thm:main} can be applied. 

\subsection{The lower bound}

Assume that $f(1,\ldots,1) \ne 0$ in $\mathbb F_q$, and let $d = \deg(f)$ and
$\gamma_2 = (d(r - 1) + 2)^{-1/(r - 1)}$.
We prove the following explicit forms of the lower bound:
\begin{align}
 \gamma_q(f)& \ge 
 \left(\frac{(r - 1)!q}{d}\right)^{\frac{1}{r - 1}}
 -\frac{r}{2} - \frac{(r - 2)!}{d},
 \label{eq:poly-lower-leading}\\
 \gamma_q(f)& \ge 
 \left(\frac{q}{d(r - 1) + 2}\right)^{\frac{1}{r - 1}} = \gamma_2q^{\frac{1}{r - 1}}.
 \label{eq:poly-lower-allq}
\end{align}

\begin{proof}
Let $I$ be any maximal $f$-free set and let $k = |I|$.
Homogeneity gives
\[
 f(z,\ldots,z) = z^{\deg(f)} f(1,\ldots,1).
\]
Thus $0$ belongs to no $f$-free set, while every nonzero
singleton is $f$-free. In particular, $0 \notin I$ and $k \ge 1$.

For each $z \in \mathbb F_q\setminus(I\cup\{0\})$,
maximality supplies a zero of $f$ in $(I\cup\{z\})^r$.
This zero uses $z$, but not in all $r$ positions.
By symmetry, for some $1 \le j \le r - 1$ and a multiset
$\{y_1,\ldots,y_j\}$ of elements of $I$, it has the form
\[
 f(\underbrace{z,\ldots,z}_{r - j},y_1,\ldots,y_j) = 0.
\]
For each such multiset, the polynomial
\[
 T\longmapsto
 f(\underbrace{T,\ldots,T}_{r - j},y_1,\ldots,y_j)
\]
has degree at most $\deg(f)$ and is nonzero: its value at $T = y_1$
is an evaluation of $f$ on $I^r$, hence is nonzero.
It has at most $\deg(f)$ roots. There are
$\binom{k + j - 1}{j}$ multisets of size $j$, so
\begin{equation}\label{eq:poly-lower-count}
 q - 1 - k
 \le \deg(f)\sum_{j = 1}^{r - 1}\binom{k + j - 1}{j} = \deg(f)\left(\binom{k + r - 1}{r - 1} - 1\right).
\end{equation}
The equality follows by summing
$\binom{k+j-1}{j}=\binom{k+j}{j}-\binom{k+j-1}{j-1}$.

Since $\deg(f) \ge 1$, \eqref{eq:poly-lower-count} implies
\[
 q \le \deg(f)\binom{k + r - 1}{r - 1} + k
 \le \frac{\deg(f)}{(r - 1)!}(k + \frac{r}{2})^{r - 1} + k,
\]
since $\prod_{j=1}^{r-1}(k+j)\le(k+r/2)^{r-1}$.
Keeping the first two terms of
$(k+r/2+(r-2)!/\deg(f))^{r-1}$ gives
\[
 \frac{\deg(f)}{(r - 1)!}
 \left(k + \frac{r}{2} + \frac{(r - 2)!}{\deg(f)}\right)^{r - 1}
 \ge \frac{\deg(f)}{(r - 1)!}(k + \frac{r}{2})^{r - 1}
 +(k + \frac{r}{2})^{r - 2}.
\]
Here $(k + r/2)^{r - 2} \ge k$, since $r \ge 3$ and $k \ge 1$.
Taking $(r - 1)$st roots proves
\eqref{eq:poly-lower-leading} for every maximal $I$.

\emph{Proof of \textup{\eqref{eq:poly-lower-allq}}.}
Also $\binom{k + j - 1}{j} \le k^j \le k^{r - 1}$:
every multiset can be represented by an ordered tuple.
Thus \eqref{eq:poly-lower-count} yields
\[
 q \le 1 + k + \deg(f)(r - 1)k^{r - 1}
 \le [\deg(f)(r - 1) + 2]k^{r - 1},
\]
which proves \eqref{eq:poly-lower-allq}.
\end{proof}

\subsection{\texorpdfstring{Set systems for $f$-free sets}{\detokenize{Set systems for f-free sets}}}

For an eligible polynomial $f \in \mathbb Z[x_1,x_2,\dots,x_r]$, if $f(1,1,\dots,1) = 0$ then
the empty set is the unique maximal $f$-free subset of $\mathbb F_q$. In the rest of the proof, we put $d = \deg(f)$ and assume
\begin{equation}\label{eq:polyq-diagonal}
 f(1,\ldots,1) \ne 0.
\end{equation}
We define a set system whose maximal independent sets correspond to
maximal $f$-free subsets of $\mathbb F_q^{\times}$.
This requires the definition of the following multiplicative subgroup of $\mathbb F_q^{\times}$:
\begin{equation*}
 H(f) = \{\zeta \in \mathbb F_q^\times\colon 
 f(\zeta x_1,x_2,\ldots,x_r) = f(x_1,x_2,\ldots,x_r)\}.
 \end{equation*}
Equality in the definition of $H = H(f)$ means equality
of polynomials in $\mathbb F_q[x_1,x_2,\dots,x_r]$. For $f,g \in \mathbb F_q[x_1,x_2,\dots,x_r]$, we write
$\gcd(f,g) = 1$ to mean that $f$ and $g$ have no common nonconstant
factor over $\mathbb F_q$.

\Needspace{8\baselineskip}
\nopagebreak[4]
\begin{lemma}\label{lem:polyq-group}
Let $f \in \mathbb F_q[x_1,\ldots,x_r]$ be eligible, and let $H = H(f)$.
\begin{enumerate}[leftmargin=*,itemsep=.35em]
\item $H$ is a subgroup of $\mathbb F_q^\times$ of order at most $d$.
\item For $z,\iota \in \mathbb F_q^\times$,
\[
 \gcd(f(z,x_2,\ldots,x_r),f(\iota,x_2,\ldots,x_r)) \ne 1
 \quad\Longleftrightarrow\quad z H = \iota H.
\]
\item For $z_1,\ldots,z_r \in H$,
\begin{equation}\label{eq:polyq-invariance}
 f(z_1x_1,\ldots,z_rx_r) = f(x_1,\ldots,x_r).
\end{equation}
\end{enumerate}
\end{lemma}

\begin{proof}[Proof of \textup{(1)}]
We write
\[
 f = \sum_{j = 0}^{d}x_1^jF_j(x_2,\ldots,x_r),
\]
where each nonzero $F_j$ is homogeneous of degree $d - j$.
The polynomial $f$ depends on every variable by symmetry.
The definition of $H$ gives
\[
 H = \{\zeta \in \mathbb F_q^\times\colon 
 \zeta^j = 1\text{ for every }j \ge 1\text{ with }F_j \ne 0\}.
\]
This is an intersection of multiplicative subgroups.
At least one positive $j \le d$ occurs, so $|H| \le j \le d$.

\emph{Proof of \textup{(2)}.}
We have $F_0 \ne 0$: otherwise $x_1$ divides $f$, and absolute
irreducibility would make $f$ a scalar multiple of $x_1$,
contrary to symmetry. Thus $f_z = f(z,x_2,\ldots,x_r)$ has degree $d$.
Homogenizing a factorization of $f_z$ would factor
$f(z x_1,x_2,\ldots,x_r)$, which is absolutely irreducible
because $z \ne 0$. Hence $f_z$ is absolutely irreducible.
The polynomials $f_z$ and $f_\iota$ therefore share a nonconstant
factor exactly when they are scalar multiples. Their homogeneous
parts of degree $d$ are both $F_0$, so that scalar is $1$.
Since the $F_j$ have different homogeneous degrees,
\[
 f_z = f_\iota
 \quad\Longleftrightarrow\quad
 z^j = \iota^j\text{ for every }j \ge 1\text{ with }F_j \ne 0
 \quad\Longleftrightarrow\quad z H = \iota H.
\]

\emph{Proof of \textup{(3)}.}
Symmetry transfers the defining invariance of $H$ from the first
coordinate to every coordinate. Applying these invariances
successively proves \eqref{eq:polyq-invariance}.
\end{proof}

We put $h = |H|$. Equation~\eqref{eq:polyq-invariance} implies that
$HI = \{\zeta x\colon \zeta \in H,\ x \in I\}$ is $f$-free whenever $I$ is
$f$-free: a zero in $(HI)^r$ could be rescaled coordinate by
coordinate to a zero in $I^r$. Consequently every maximal
$f$-free set is a union of complete $H$-cosets.

Let $X = \mathbb F_q^\times$, and let $\mathcal S$ consist of the
inclusion-minimal subsets $B \subseteq X$ for which $f$ vanishes
somewhere on $B^r$. By \eqref{eq:polyq-diagonal}, each block has
size between $2$ and $r$, and minimality makes $\mathcal S$ an
antichain. Homogeneity shows that multiplication by any element
of $X$ preserves $\mathcal S$. This action is transitive on $X$,
so $\mathcal S$ is a $(q - 1,\vb{d})$-antichain.

For a collection $B \subseteq X/H$ of cosets, we write
\[
 \widehat B = \bigcup_{C \in B}C \subseteq X.
\]
The quotient system $\mathcal S/H$ has atom set $X/H$; its blocks
are the inclusion-minimal nonempty collections $B \subseteq X/H$
such that $f$ vanishes somewhere on $(\widehat B)^r$.
We write $\widehat d_i(B)$ and $\widehat d_i(x,y)$ for the
codegrees and common-link sizes in $\mathcal S/H$, defined
as in Section~1. Its degree vector will be denoted by
$\widehat{\vb{d}} = (\widehat d_2,\ldots,\widehat d_r)$.

\Needspace{8\baselineskip}
\nopagebreak[4]
\begin{lemma}\label{lem:polyq-quotient}
The system $\mathcal S/H$ is a $\left((q - 1)/h,\widehat{\vb{d}}\right)$-antichain,
and the following statements hold:
\begin{enumerate}[leftmargin=*,itemsep=.35em]
\item For $2 \leq i \leq r$,
\[
 \widehat d_i = \frac{d_i}{h^{i - 1}}.
\]
\item For $2 \leq i \leq r$, let $B \subseteq X/H$ with $1 \leq |B| < i$, and let $B_0 \subseteq X$ contain one representative from each coset in $B$. Then
\begin{equation}\label{eq:polyq-quotient-codegrees}
 \widehat d_i(B) = \frac{d_i(B_0)}{h^{i - |B|}}.
\end{equation}
\item For $2 \leq i \leq r$ and distinct cosets $z H,\iota H \in X/H$,
\begin{equation}\label{eq:polyq-quotient-common}
 \widehat d_i(z H,\iota H) = \frac{d_i(z,\iota)}{h^{i - 1}}.
\end{equation}
\item If $I$ is a maximal independent set in $\mathcal S/H$, then $\widehat I$ is a maximal $f$-free set in $\mathbb F_q$.
\end{enumerate}
\end{lemma}

\begin{proof}
We first describe the blocks in the two systems.
A single coset $z H$ is $f$-free, since, for all
$\zeta_1,\ldots,\zeta_r \in H$,
\[
 f(z\zeta_1,\ldots,z\zeta_r) = z^df(1,\ldots,1) \ne 0.
\]
A minimal forbidden collection has at most $r$ cosets,
since a zero of $f$ uses at most $r$ of them. Scalar
multiplication acts transitively on $X/H$ and preserves
each rank, proving regularity.

A block of $\mathcal S$ cannot contain two atoms from one coset.
In a zero of $f$ using the elements of such a block, rescale all coordinates
from that coset to one representative. Equation~\eqref{eq:polyq-invariance}
preserves the zero and gives a forbidden proper subset,
contrary to minimality.

Every choice of one representative from each coset of a
block of $\mathcal S/H$ is a block of $\mathcal S$:
rescale the coordinates of a zero to those representatives.
A forbidden proper subset would project to a forbidden
proper collection of cosets, proving minimality.
Conversely, the cosets meeting a block of $\mathcal S$ form a quotient
block: a smaller forbidden collection would give a forbidden
proper subset of the original block.

\emph{Proof of \textup{(1)}.}
We fix an atom $z \in X$. Every quotient block of size $i$
through $z H$ gives exactly $h^{i - 1}$ blocks through $z$,
one for each choice of representatives in the other cosets.
Thus $d_i = h^{i - 1}\widehat d_i$.

\emph{Proof of \textup{(2)}.}
With the representatives $B_0$ fixed, a quotient block
through $B$ gives exactly $h^{i - |B|}$ blocks containing $B_0$.
These account for every original block containing $B_0$,
so \eqref{eq:polyq-quotient-codegrees} follows.

\emph{Proof of \textup{(3)}.}
A common link at the distinct quotient atoms $z H,\iota H$
avoids both cosets. Its $h^{i - 1}$ choices of representatives
are exactly the common links at $z,\iota$ that use these cosets.
Conversely every
common link at $z,\iota$ projects to one of these quotient
common links. This proves \eqref{eq:polyq-quotient-common}.

\emph{Proof of \textup{(4)}.}
Let $I$ be a maximal independent set in $\mathcal S/H$.
The union $\widehat I$ is $f$-free by definition.
For $z \in X\setminus\widehat I$, maximality supplies a zero
using cosets from $I\cup\{z H\}$. Rescale every coordinate
from $z H$ to the single element $z$, leaving the other
coordinates unchanged. This gives a zero in
$(\widehat I\cup\{z\})^r$. Also $0$ cannot be added, since
$f(0,\ldots,0) = 0$. Therefore $\widehat I$ is maximal
$f$-free in $\mathbb F_q$.
\end{proof}

The remaining estimates verify the sparsity bounds for the quotient.
For common links, the original atoms must lie in distinct $H$-cosets.

\Needspace{8\baselineskip}
Let $d\ge1$ and $s\ge1$ be integers. Every nonzero polynomial
$P$ of degree at most $d$ in $s$ variables satisfies
\begin{equation}\label{eq:polyq-zero-count}
 |\{y \in \mathbb F_q^s\colon P(y) = 0\}| \le d q^{s - 1}.
\end{equation}
To prove this, we induct on $s$. We write
$P = \sum_{j = 0}^{\mathfrak{d}}P_jy_s^j$ with $P_{\mathfrak{d}} \ne 0$.
Outside the zero set of $P_{\mathfrak{d}}$, each fibre has at most $\mathfrak{d}$
zeros. By induction, at most $(d - \mathfrak{d})q^{s - 2}$ fibres are
exceptional; allowing all $q$ points in each gives
$\mathfrak{d}q^{s - 1} + (d - \mathfrak{d})q^{s - 1} = dq^{s - 1}$.
For $s=1$, a nonzero polynomial has at most $d$ roots.

If $s \ge 2$ and $P,Q$ are nonzero polynomials of degree
at most $d$ which are coprime over the algebraic closure, then
\begin{equation}\label{eq:polyq-common-zero}
 |\{y \in \mathbb F_q^s\colon P(y) = Q(y) = 0\}|
 \le d^2q^{s - 2}.
\end{equation}
Indeed, the common zero set has dimension at most $s - 2$ and
sum of component degrees at most $d^2$. Theorem~1.1 and Corollary~2.2
of Lachaud and Rolland~\cite{PolyLachaudRolland} therefore give \eqref{eq:polyq-common-zero}.
If either polynomial is a nonzero constant, the common
zero set is empty.

For integers $s\ge2$ and $d\ge1$,
The Lang--Weil estimate~\cite{LangWeil1954}, in the form stated by Poonen~\cite[Theorem~7.7.1]{PolyPoonen}, gives a constant
$\gamma_3(s,d)\ge1$ such that, for every absolutely irreducible
polynomial $P$ of degree at most $d$ in $s$ variables,
\begin{equation}\label{eq:polyq-lang-weil}
 \left|\,|\{y \in \mathbb F_q^s\colon P(y) = 0\}| - q^{s - 1}\,\right|
 \le \gamma_3(s,d)q^{s - \frac{3}{2}}.
\end{equation}
\Needspace{8\baselineskip}
We put
\[
 \gamma_4(r,d) = \gamma_3(r - 1,d) + r^r d^2 + (r - 1)d.
\]

\nopagebreak[4]
\begin{lemma}\label{lem:polyq-degree-estimates}
The following hold:
\begin{enumerate}[leftmargin=*,itemsep=.35em,ref=\arabic*]
\item For $B\subseteq X$ and $1 \leq |B| < i \leq r$,
\begin{equation}\label{eq:polyq-codegrees}
 d_i(B) \leq d i^r q^{i - |B| - 1}.
\end{equation}
\item\label{it:polyq-top-degree}
$\left|d_r-\frac{q^{r-2}}{(r-1)!}\right|
\le\frac{\gamma_4(r,d)}{(r-1)!}q^{r-\frac52}$.
\item For $z,\iota\in X$ with $z H \ne \iota H$,
\begin{equation}\label{eq:polyq-top-common}
 d_r(z,\iota) \leq d^2q^{r - 3}.
\end{equation}
\item For distinct $z,\iota \in X$ and $2 \leq i < r$,
\begin{equation}\label{eq:polyq-lower-common}
 d_i(z,\iota) \leq d r^r q^{i - 2}.
\end{equation}
\end{enumerate}
\end{lemma}

\begin{proof}[Proof of \textup{(1)}]
There is nothing to prove
if $B$ is not $f$-free: no larger minimal forbidden set
contains it. Otherwise, with $1 \le s=|B| < i$, we write
$B = \{\iota_1,\ldots,\iota_s\}$.
Every $i$-element block containing $B$ is the set of coordinate
values of a zero of $f$. Name the other elements
$y_{s + 1},\ldots,y_i$, and for $1\le j\le s$ we put $y_j = \iota_j$.
Thus, for some surjection
$\sigma\colon \{1,\ldots,r\} \to \{1,\ldots,i\}$,
\[
 f(y_{\sigma(1)},\ldots,y_{\sigma(r)}) = 0.
\]
After substitution of the $\iota_j$, this is a nonzero polynomial
in $i - s$ variables of degree at most $d$.
It is nonzero because setting all the remaining variables
equal to $\iota_1$ gives an evaluation on $B^r$, where $f$
does not vanish. There are at most $i^r$ choices of $\sigma$.
Applying \eqref{eq:polyq-zero-count} to each choice, and
allowing ordered assignments to overcount the blocks,
proves \eqref{eq:polyq-codegrees}.

\emph{Proof of \textup{(2)}.}
We fix $z \in X$. Lemma~\ref{lem:polyq-group} and
\eqref{eq:polyq-lang-weil} give
\[
 \left|\,|\{y \in \mathbb F_q^{r - 1}\colon f_z(y) = 0\}| - q^{r - 2}\,\right|
 \le \gamma_3(r - 1,d)q^{r - \frac{5}{2}}.
\]
We next remove solutions whose $r$ coordinate values do not
form a block of size $r$. For every nonsurjective map
$\sigma\colon \{1,\ldots,r\} \to \{1,\ldots,r\}$, we put
\[
 g_\sigma(x_1,\ldots,x_r) = f(x_{\sigma(1)},\ldots,x_{\sigma(r)}).
\]
This polynomial is homogeneous of degree $d$ and nonzero,
since $g_\sigma(1,\ldots,1) = f(1,\ldots,1) \ne 0$.
It omits a variable, whereas $f$ depends on every variable,
so it is not a scalar multiple of $f$.

Specialization at $x_1 = z \ne 0$ is injective on homogeneous
degree-$d$ polynomials: the coefficient of a monomial
$y^\alpha$ is multiplied by $z^{d - |\alpha|}$.
Hence $g_{\sigma,z} = g_\sigma(z,x_2,\ldots,x_r)$ is nonzero.
If it shared a factor with the irreducible polynomial
$f_z$, then $f_z\mid g_{\sigma,z}$; since
$\deg(g_{\sigma,z})\le d=\deg(f_z)$, they would be scalar multiples. Injectivity would
then make $g_\sigma$ a scalar multiple of $f$, a contradiction.
By \eqref{eq:polyq-common-zero}, each $\sigma$ accounts for
at most $d^2q^{r - 3}$ solutions with first coordinate $z$.

We count solutions having a zero coordinate separately. For each
$j \in \{2,\ldots,r\}$, the specialization
\[
 f(z,x_2,\ldots,x_{j - 1},0,x_{j + 1},\ldots,x_r)
\]
is a nonzero polynomial of degree at most $d$ in $r - 2$ variables.
Indeed, the coefficient multipliers $z^{d-|\alpha|}$ are nonzero.
If this specialization were identically zero, then $f(x_1,\ldots,x_{j - 1},0,x_{j + 1},\ldots,x_r)$
would be the zero polynomial, so $x_j$ would divide $f$. Absolute
irreducibility would then force $f$ to be a scalar multiple of $x_j$,
contrary to symmetry. Hence
\eqref{eq:polyq-zero-count} shows that the solutions with at least
one zero coordinate contribute at most
$(r - 1)dq^{r - 3}$.

Every remaining excluded solution satisfies $g_\sigma=0$ for some
nonsurjective $\sigma$: if coordinates repeat, we identify their indices;
otherwise we choose a zero of $f$ using a proper subset of the coordinate
values. Conversely, for distinct nonzero coordinates, such a zero
shows that their set is not a minimal forbidden set. There are at most
$r^r$ maps, so these solutions contribute
at most $r^r d^2q^{r - 3}$.
By symmetry, each remaining block through $z$ gives exactly
$(r - 1)!$ solutions with first coordinate $z$. Consequently
\[
 \left|(r - 1)!d_r - q^{r - 2}\right|
 \le \gamma_3(r - 1,d)q^{r - \frac{5}{2}}
 +\bigl(r^r d^2 + (r - 1)d\bigr)q^{r - 3}
 \le \gamma_4(r,d)q^{r - \frac{5}{2}}.
\]
Dividing by $(r-1)!$ proves part~\textup{(\ref{it:polyq-top-degree})}.

\emph{Proof of \textup{(3)}.}
If $z H \ne \iota H$, Lemma~\ref{lem:polyq-group} says that
$f_z$ and $f_\iota$ are coprime. Every set in $N_r(z,\iota)$,
with its atoms in any order, gives a common zero of these polynomials, so
\eqref{eq:polyq-common-zero} proves
\eqref{eq:polyq-top-common}.

\emph{Proof of \textup{(4)}.}
We have $d_i(z,\iota) \le d_i$, and
\eqref{eq:polyq-codegrees} with $|B|=1$ gives
$d_i \le d i^r q^{i - 2} \le d r^r q^{i - 2}$.
\end{proof}

\subsection{Proof of the upper bound}

\begin{proof}
Assume \eqref{eq:polyq-diagonal}; otherwise homogeneity makes every
constant tuple a zero, so only the empty set is $f$-free.
We fix, before choosing the field-size threshold,
\[
 \delta_0 = \frac{1}{2(r - 2)},\qquad \gamma_5 = 2(r - 1)!\max\{dr^r,d^2\},
\]
and put
\[
 \gamma_6 = \frac{1}{2^{22}r^2(r - 2)}.
\]
Let $\gamma_7 = C(r,\delta_0,\gamma_5)$ be the threshold in
Theorem~\ref{prop:quantitative-main} with its parameter
$c = \gamma_5$. We will prove that some
maximal $f$-free set $I$ satisfies
\begin{equation}\label{eq:polyq-upper}
 |I| \le h\bigl(1 + 3(\log q)^{-\gamma_6}\bigr)
 \bigl((r - 2)(r - 2)!\,q\log q\bigr)^{\frac{1}{r - 1}}.
\end{equation}
For $q \ge 3$, we have $h \le d$ and
$1 + 3(\log q)^{-\gamma_6} \le 4$, so this gives
\[
 |I| \le 
 4d\bigl((r - 2)(r - 2)!\bigr)^{\frac{1}{r - 1}}
 (q\log q)^{\frac{1}{r - 1}}.
\]

For $2\le i<r$, Lemma~\ref{lem:polyq-quotient}~(1),
\eqref{eq:polyq-codegrees} and $1\le h\le d$ give
\[
 \widehat d_i \le \frac{dr^r}{h^{i - 1}}q^{i - 2}
 \le dr^r q^{i - 2}.
\]
For the top rank, Lemma~\ref{lem:polyq-degree-estimates}\,\textup{(\ref{it:polyq-top-degree})} gives
\[
 \left|\frac{(r - 1)!h^{r - 1}\widehat d_r}{q^{r - 2}} - 1\right|
 \le \gamma_4(r,d)q^{-\frac{1}{2}}.
\]
In particular, for $2 \le i < r$,
$\widehat d_i = o_{r,d}(\widehat d_r)$ and
$\widehat d_i^{1/(i - 1)} = o_{r,d}(\widehat d_r^{1/(r - 1)})$.
We choose $\gamma_1 = \gamma_1(r,d) \ge 4$ so that, uniformly over all
such polynomials satisfying \eqref{eq:polyq-diagonal}
and all $q \ge \gamma_1$,
the following hold, with $\kappa = (\log q)^{-\gamma_6}$:
\begin{equation}\label{eq:polyq-threshold-degrees}
 \left(1 - \frac{\kappa}{4}\right)\frac{q^{r - 2}}{(r - 1)!}
 \le d_r \le q^{r - 2},\qquad \log\widehat d_r \ge \frac12\log q.
\end{equation}
For $2 \le i < r$, we also require
\begin{equation*}
 \widehat d_i < \widehat d_r,\qquad \widehat d_i^{\frac{1}{i - 1}} < \widehat d_r^{\frac{1}{r - 1}},\qquad \widehat d_r^{\frac{1}{r - 1}} \ge \gamma_7.
\end{equation*}
This choice is uniform by
\eqref{eq:polyq-codegrees} and Lemma~\ref{lem:polyq-degree-estimates}\,\textup{(\ref{it:polyq-top-degree})},
Lemma~\ref{lem:polyq-quotient}~(1), $1\le h\le d$ and
$q^{-1/2}=o((\log q)^{-\gamma_6})$.

For the quotient degree vector,
\[
 D(\widehat{\vb{d}}) = \max_{2 \le i \le r}\widehat d_i^{\frac{1}{i - 1}} = \widehat d_r^{\frac{1}{r - 1}} = \frac{d_r^{\frac{1}{r - 1}}}h.
\]
Thus $D(\widehat{\vb{d}}) \ge \gamma_7$,
$\log D(\widehat{\vb{d}}) = (\log\widehat d_r)/(r - 1)$, and
$\max_i\widehat d_i = \widehat d_r$.
Since $0 < \kappa \le 1$, \eqref{eq:polyq-threshold-degrees} gives
$d_r \ge q^{r - 2}/(2(r - 1)!)$.
Consequently, for every real $z$ with $0 < z \le r - 1$,
\begin{equation}\label{eq:polyq-power-comparison}
 \gamma_5d_r^{\frac{z}{r - 1}}
 \ge \max\{dr^r,d^2\}\,
 q^{\frac{r - 2}{r - 1}z}.
\end{equation}
The extra constant factor is
$[2(r - 1)!]^{1 - z/(r - 1)} \ge 1$.

For $2 \le s < i \le r$,
\[
 \frac{r - 2}{r - 1}(i - s - \delta_0) - (i - s - 1) = 1 - \frac{i - s + \frac{1}{2}}{r - 1} \ge \frac{1}{2(r - 1)}.
\]
Let $B \subseteq X/H$ have size $s$, and let $B_0$ be a set
of one representative from each coset in $B$.
Equations~\eqref{eq:polyq-codegrees} and
\eqref{eq:polyq-power-comparison} give
$d_i(B_0) \le \gamma_5d_r^{(i - s - \delta_0)/(r - 1)}$.
Dividing by $h^{i - s}$ in \eqref{eq:polyq-quotient-codegrees}
gives every required quotient codegree:
\[
 \widehat d_i(B)
 \le \gamma_5h^{-\delta_0}D(\widehat{\vb{d}})^{i - s - \delta_0}
 \le \gamma_5D(\widehat{\vb{d}})^{i - s - \delta_0}.
\]

For $2 \le i < r$,
\[
 \frac{r - 2}{r - 1}(i - 1 - \delta_0) - (i - 2) = \frac{r - i - \frac{1}{2}}{r - 1} \ge \frac{1}{2(r - 1)},
\]
whereas for the top rank,
\[
 \frac{r - 2}{r - 1}(r - 1 - \delta_0) - (r - 3) = 1 - \frac{1}{2(r - 1)} > 0.
\]
For $2\le i\le r$ and $z,\iota\in X$ with
$z H \ne \iota H$, equations
\eqref{eq:polyq-top-common}--\eqref{eq:polyq-lower-common} and
\eqref{eq:polyq-power-comparison} give
$d_i(z,\iota) \le \gamma_5d_r^{(i - 1 - \delta_0)/(r - 1)}$.
Thus \eqref{eq:polyq-quotient-common} gives
\[
 \widehat d_i(z H,\iota H)
 \le \gamma_5h^{-\delta_0}D(\widehat{\vb{d}})^{i - 1 - \delta_0}
 \le \gamma_5D(\widehat{\vb{d}})^{i - 1 - \delta_0}.
\]
Since $\delta_0/(2^{21}r^2) = \gamma_6$,
Theorem~\ref{prop:quantitative-main} gives
a maximal independent set $I_0 \subseteq X/H$ with
\begin{align*}
 |I_0| &\le 
 \bigl(1 + (\log\widehat d_r)^{-\gamma_6}\bigr)
 \frac{q - 1}{h}\,\rho(\widehat{\vb{d}}),\\
 \sum_{i = 2}^r\widehat d_i\rho(\widehat{\vb{d}})^{i - 1} &= \log D(\widehat{\vb{d}}).
\end{align*}
Keeping the rank-$r$ term yields
\[
 \rho(\widehat{\vb{d}}) \le 
 \left(\frac{\log\widehat d_r}{(r - 1)\widehat d_r}\right)^{\frac{1}{r - 1}}.
\]
We put $I = \widehat{I_0}$. By Lemma~\ref{lem:polyq-quotient}~(4),
$I$ is maximal $f$-free and $|I| = h|I_0|$. Hence
\begin{equation}\label{eq:polyq-before-final}
 |I| \le 
 \bigl(1 + (\log\widehat d_r)^{-\gamma_6}\bigr)(q - 1)
 \left(\frac{\log\widehat d_r}{(r - 1)\widehat d_r}\right)^{\frac{1}{r - 1}}.
\end{equation}

By \eqref{eq:polyq-threshold-degrees},
$(\log\widehat d_r)^{-\gamma_6} \le 2\kappa$ and
$\log\widehat d_r \le \log d_r \le (r - 2)\log q$.
Also $\widehat d_r = d_r/h^{r - 1}$. Therefore
\[
 (q - 1)
 \left(\frac{\log\widehat d_r}{(r - 1)\widehat d_r}\right)^{\frac{1}{r - 1}}
 \le h\left(1 - \frac{\kappa}{4}\right)^{-\frac{1}{r - 1}}
 \bigl((r - 2)(r - 2)!\,q\log q\bigr)^{\frac{1}{r - 1}}.
\]
For $0 \le \kappa \le 1$,
\[
 \left(1 - \frac{\kappa}{4}\right)^{-\frac{1}{r - 1}}
 \le \left(1 - \frac{\kappa}{4}\right)^{-1} \le 1 + \frac{\kappa}{3},\qquad (1 + 2\kappa)(1 + \frac{\kappa}{3}) \le 1 + 3\kappa.
\]
Using $(1+2\kappa)(1+\kappa/3)\le1+3\kappa$ in
\eqref{eq:polyq-before-final} proves \eqref{eq:polyq-upper}.
\end{proof}

\section{Linear equations and Pythagorean triples}\label{sec:Pythag}

In this section, we give almost sharp upper bounds on the size of smallest maximal Pythagorean triple free subsets of $\mathbb F_q$ and 
smallest maximal subsets of finite abelian groups with no non-trivial solutions to linear equations $a_1 + a_2 + \dots + a_k = b_1 + b_2 + \dots + b_{\ell}$.

\subsection{A useful lemma}

The following lemma will be useful throughout this section; it is a specific setting of the parameters in Theorem \ref{thm:main}.
For a set system $\mathcal H$, let $\gamma(\mathcal H)$ denote the minimum size of a maximal independent set in $\mathcal H$.

\begin{lemma}\label{lem:criterion}
Fix $r\ge3$ and $C\ge2$.
Let $\HH$ be a regular antichain on $X$, with block
sizes between $2$ and $r$. Suppose that, for a parameter $M\ge2$,
$M/C\le|X|\le CM$, that for $2\le i<r$,
\[
 C^{-1}M^{r-2}\le d_r\le CM^{r-2},\qquad
 d_i\le CM^{i-2},
\]
that for $2\le |A|<i\le r$,
\[
 d_i(A)\le CM^{i-|A|-1},
\]
and that for distinct $x,y\in X$,
\[
 |N_r(x)\cap N_r(y)|\le CM^{r-3}.
\]
Then, as $M\to\infty$ with $r$ and $C$ fixed,
\[
 \gamma(\HH)\le
 \left(C^{\frac{r}{r-1}}\left(\frac{r-2}{r-1}\right)^{\frac{1}{r-1}}+o(1)\right)
 (M\log M)^{\frac{1}{r-1}}.
\]
\end{lemma}

\begin{proof}
One has $C^{-1/(r-1)}M^{(r-2)/(r-1)}\le D\le CM^{(r-2)/(r-1)}$.
For $i<r$, use $|N_i(x)\cap N_i(y)|\le d_i$.
Before the subtraction of the $\delta$ term, the differences between
the required and the available exponents of $M$ in the mixed
sparseness conditions are, respectively, for $2\le |A|<i\le r$
in the first inequality and $2\le i<r$ in the second,
\[
 \frac{r-1-(i-|A|)}{r-1}\ge\frac{1}{r-1},\qquad
 \frac{r-i}{r-1}\ge\frac{1}{r-1}.
\]
For the intersection of links of rank $r$, the difference is $1$.
We take $\delta=1/(2(r-2))$. The $\delta$ term lowers the required
exponent of $M$ by $1/(2(r-1))$. If $M\ge C^{4(r-1)}$, the
remaining difference contributes a factor at least
$M^{1/(2(r-1))}\ge C^2$. This absorbs both the factor $C$ in
each hypothesis and the factor at most $C$ from the lower bound
on the required power of $D$. Thus $\HH$ is $\delta$-sparse.
In this range, $1\le\log D\le\log M$.
Put $z=(\log D/d_r)^{1/(r-1)}$.
The defining equation $\sum_{i=2}^r d_i\rho^{i-1}=\log D$ gives
$\rho\le z$. The explicit constants in the deduction of
Theorem~\ref{thm:main} permit
$\log\log\log\log K=(8r(r-2))^{16}$ and
$\epsilon=1/(2^{22}r^2(r-2))$. The bounds on $D$ give
$\log D/\log M\to(r-2)/(r-1)$. That theorem therefore gives
\[
 \gamma(\HH)\le(1+K(\log D)^{-\epsilon})|X|z
 =(1+o(1))|X|\left(\frac{(r-2)\log M}{(r-1)d_r}\right)^{\frac{1}{r-1}}.
\]
The claimed bound follows from $|X|\le CM$ and
$d_r\ge C^{-1}M^{r-2}$. The error tends to zero uniformly over
systems satisfying the hypotheses with the fixed values of $r$ and $C$.
\end{proof}

\subsection{Sum-free sets in abelian groups}

\begin{definition}\label{def:kl}
Let $G$ be an abelian group and $k,l\ge1$. A subset $A\subseteq G$
is \emph{$(k,l)$-free} if, for $a_i,b_j\in A$, every equality
\begin{equation}\label{eq:kl}
 a_1+\cdots+a_k=b_1+\cdots+b_l
\end{equation}
has equal multisets on the two sides. Repetitions are allowed.
For finite $G$, denote by $\gamma_{k,l}(G)$ the least cardinality
of a maximal $(k,l)$-free subset.
\end{definition}

For $k\ne l$, this is the usual $(k,l)$-sum-free condition.
For $k=l=h$, it is the generalised Sidon set condition $B_h$. The minimum among maximal
sets is distinct from the maximum cardinality problem. For the
latter, Bajnok and Matzke give an exact formula for cyclic
groups~\cite{BajnokMatzke2018}. Concerning small maximal sets, the
construction of Ruzsa gives a maximal Sidon subset of $[1,N]$ of size
at most $2^{16}(N\log N)^{1/3}$ for $N\ge2$~\cite{Ruzsa}.
Indeed, for $N\ge2^{60}$, his proof permits a prime
$8(N\log N)^{1/3}<p\le18(N\log N)^{1/3}$, with
$q=p^2+p+1$ and $M=\lfloor N/q\rfloor\ge8$.
The two probabilities estimated there are at least $1/(8M)$,
giving a maximal set of size at most $p+1+2N/q\le20(N\log N)^{1/3}$.
For $2\le N<2^{60}$, the bound $2\sqrt N$ for any Sidon set
implies the stated estimate.
The construction of Redman, Rose, and Walker gives at most
$2^{16}(m2^m)^{1/3}$ elements in $(\mathbb Z/2\mathbb Z)^m$ for
$m\ge1$, with the distinct-summand convention~\cite{RedmanRoseWalker}.
For $m\ge32$, their proof permits $T=12$ and the least even integer
$\ell>(2/3)\log_2(12(\log2)m2^m)$, giving size at most
$2^{\ell/2}+2^{m-\ell}\le8(m2^m)^{1/3}$.
For $m<32$, the elementary bound $1+\sqrt{2^{m+1}}$ suffices.
Small maximal sum-free sets in groups are studied by
Giudici and Hart~\cite{GiudiciHart}.

\begin{lemma}\label{lem:lower}
Let $r=k+l\ge3$. If $G$ is a nontrivial finite abelian group of
order $n$ and $\gcd(n,r!)=1$, then
\[
 \gamma_{k,l}(G)\ge\left(\frac{n}{2^r}\right)^{\frac{1}{r-1}}.
\]
\end{lemma}

\begin{proof}
Let $A$ be maximal, note $A\neq \emptyset$. For $x\notin A$, adjoining $x$ produces a nontrivial equality
\eqref{eq:kl}. Its net coefficient $c$ of $x$ is nonzero; otherwise
replacing every occurrence of $x$ by one fixed element of $A$ would
give a nontrivial equality in $A$. Since $0<|c|\le r$,
multiplication by $c$ is an automorphism of $G$. Thus $x$ is
determined by the nonempty set of positions occupied by $x$ and
the entries of $A$ in all remaining positions. The number of
these choices is at most
$\sum_{j=1}^r\binom rj|A|^{r-j}\le(2^r-1)|A|^{r-1}$.
Consequently $n\le|A|+(2^r-1)|A|^{r-1}\le2^r|A|^{r-1}$.
\end{proof}

{\bf Higher Sidon sets.} We first consider the balanced case $a_1 + a_2 + \dots + a_h = b_1 + b_2 + \dots + b_h$. When $h = 2$, 
sets with no non-trivial solutions to the equation are called {\em Sidon sets}. 

\begin{proposition}\label{prop:bh}
Fix $h\ge2$. As $n\to\infty$ through orders of finite abelian
groups $G$ satisfying $\gcd(n,(2h)!)=1$,
\[
 \left(\frac{n}{2^{2h}}\right)^{\frac{1}{2h-1}}\le\gamma_{h,h}(G)
 \le\left(\left(\frac{(2h-2)h!(h-1)!}{2h-1}\right)^{\frac{1}{2h-1}}+o(1)\right)
 (n\log n)^{\frac{1}{2h-1}}.
\]
The error tends to zero uniformly over these groups for fixed $h$.
In particular, these estimates hold for the additive group of
$\F_q$ whenever $\operatorname{char}\F_q>2h$.
\end{proposition}

\begin{proof}
The lower bound is Lemma~\ref{lem:lower}. Put $r=2h$, and let
$\HH$ consist of the inclusion-minimal supports of nontrivial
equalities
\begin{equation}\label{eq:bh}
 a_1+\cdots+a_h=b_1+\cdots+b_h.
\end{equation}
Its independent sets are precisely the $B_h$ subsets of $G$.
Cancelling common terms shows that every block of size $i$ is
witnessed by a relation on distinct elements $z_1,\ldots,z_i$,
with $0<|c_j|\le h$, $\sum_jc_j=0$, and
$\sum_{c_j>0}c_j\le h$:
\begin{equation}\label{eq:patterns}
 \sum_{j=1}^i c_jz_j=0.
\end{equation}
Here the $c_j$ are integers. Conversely, such a relation can be
padded on both sides to give \eqref{eq:bh}. Minimality excludes
zero coefficients. For each $i\le r$ there are at most
$(2h)^i\le r^r$ coefficient vectors, and there are no
blocks of size one or two.

Translations act transitively on $G$ and preserve $\HH$; hence
$\HH$ is regular. Every nonzero coefficient in
\eqref{eq:patterns} acts invertibly on $G$. Fixing $a<i$ distinct
values in such a relation leaves at most $n^{i-a-1}$ choices.
There are at most $r^r$ assignments of the fixed values to
positions. It follows, with $C=(2r)^{2r}$, that for $2\le i\le r$,
and for $2\le |S|<i$ in the second inequality,
\[
 d_i\le Cn^{i-2},\qquad
 d_i(S)\le Cn^{i-|S|-1}.
\]

There are exactly $n^{r-1}$ ordered solutions of \eqref{eq:bh}.
A coincidence between two coordinates, or a relation
\eqref{eq:patterns} on a proper subset of the coordinates, adds
an equation whose coefficient vector is not rationally
proportional to the vector of \eqref{eq:bh}. The two vectors
therefore have a nonzero $2\times2$ minor of absolute value at
most $2h$. This minor acts invertibly on $G$. Fixing the other
coordinates gives at most $n^{r-2}$ solutions. There are at most
$\binom r2+(r+1)^r\le2(r+1)^r$ such equations, since a proper
relation has its coefficients in $\{-h,\ldots,h\}^r$.
After their exclusion the number of ordered solutions with
distinct coordinates and minimal support is between
$n^{r-1}-2(r+1)^rn^{r-2}$ and $n^{r-1}$.
Each block of size $r$ has exactly $2(h!)^2$ ordered
representations. Indeed, if two partitions of its elements into
equal-sum $h$-sets are neither identical nor interchanged, subtracting
their relations and dividing by $2$ gives a nontrivial relation
\eqref{eq:patterns} on a proper subset, contrary to minimality.
The two orders of the sides and the $h!$ orders within each side
are therefore all the representations. Since by regularity
$d_rn$ is $r$ times the number of blocks of size $r$, we obtain
\[
 \frac{n^{r-2}-2(r+1)^rn^{r-3}}{h!(h-1)!}
 \le d_r\le\frac{n^{r-2}}{h!(h-1)!}.
\]
For $n\ge4(r+1)^r$, these bounds imply
$C^{-1}n^{r-2}\le d_r\le Cn^{r-2}$.

Fix distinct $x,y\in G$. A common link of rank $r$ gives, after
ordering its elements, two equations
\[
 u\cdot z+c x=0,\qquad u'\cdot z+c'y=0,
\]
whose coefficients lie in $\{1,-1\}$ and whose full coefficient
sums are zero. If $u'=\lambda u$ over $\mathbb Q$, then
$\lambda=\pm1$ and $c'=\lambda c$; consistency implies $x=y$.
Otherwise a nonzero $2\times2$ minor is $\pm2$ and acts invertibly
on $G$. The two equations then have $n^{r-3}$ solutions.
We fix an ordering of $G$ and list each common link in that order.
There are at most $2^r$ choices of the full coefficient vector
for each equation, giving
$|N_r(x)\cap N_r(y)|\le2^{2r}n^{r-3}\le Cn^{r-3}$.
The proof of Lemma~\ref{lem:criterion}, with $M=n$,
now verifies the hypotheses of Theorem~\ref{thm:main} once
$n\ge\max\{4(r+1)^r,C^{4(r-1)}\}$.
The degree estimates give
$d_r/n^{r-2}\to1/(h!(h-1)!)$ and
$\log D/\log n\to(r-2)/(r-1)$.
Since $\rho\le(\log D/d_r)^{1/(r-1)}$, Theorem~\ref{thm:main}
gives
\[
 \gamma_{h,h}(G)\le
 \left(\left(\frac{(r-2)h!(h-1)!}{r-1}\right)^{\frac{1}{r-1}}+o(1)\right)
 (n\log n)^{\frac{1}{r-1}},
\]
which is the asserted bound because $r=2h$.
\end{proof}

\Needspace{9\baselineskip}

{\bf Even sums and one summand.} Next we consider the equation $a_1 + a_2 + \dots + a_{2h} = b_1$. 

\begin{proposition}\label{prop:even}
Fix $h\ge1$. As $q\to\infty$ through prime powers with
$\operatorname{char}\F_q>2h+1$,
\[
 \left(\frac{q}{2^{2h+1}}\right)^{\frac{1}{2h}}\le\gamma_{2h,1}(\F_q)
 \le\left(\bigl((2h-1)(2h-1)!\bigr)^{\frac{1}{2h}}+o(1)\right)
 (q\log q)^{\frac{1}{2h}}.
\]
Moreover, there is a symmetric maximal $(2h,1)$-free set satisfying
the same upper bound. The error tends to zero uniformly over the
indicated fields for fixed $h$.
\end{proposition}

\begin{proof}
The lower bound follows from Lemma~\ref{lem:lower}. Put $r=2h+1$
and $Y=\F_q^*/\{1,-1\}$. Write $[a]$ for the class of $a$.
Let $\HH$ be the antichain of minimal supports in $Y$ of signed
relations with $\epsilon_j\in\{1,-1\}$,
\begin{equation}\label{eq:signed}
 \epsilon_1z_1+\cdots+\epsilon_rz_r=0,
\end{equation}
where classes may be repeated. The definition is independent of
the representatives. A singleton cannot be a block: its net
coefficient would be a nonzero odd integer (as $r=2h+1$ is odd) of
absolute value at most $r$, hence nonzero in $\F_q$. In a minimal support each net
coefficient is nonzero. Indeed, a class with coefficient zero can
be removed, and the resulting shorter relation can be padded by
opposite pairs from another class to restore length $r$.

Multiplication by $\F_q^*$ acts transitively on $Y$ and preserves
$\HH$. For each $i\le r$, there are at most $(2r)^r$
coefficient vectors, with all nonzero coefficients of absolute
value at most $r$, and at most $r^r$ assignments of fixed classes
to positions. Fixing one representative of each class, the
one-equation count gives, with $C=(2r)^{2r}$, for $2\le i\le r$,
and for $2\le |S|<i$ in the second inequality,
\[
 d_i\le Cq^{i-2},\qquad
 d_i(S)\le Cq^{i-|S|-1}.
\]
Relations on $r$ distinct classes have coefficients $\pm1$.
Coincidences up to sign and relations on proper subsupports each
impose a second independent equation over $\F_q$. Indeed, the
second row has a zero coordinate where the first has coefficient
$\pm1$, and has a nonzero coefficient elsewhere. Solutions with some
$z_j=0$ are excluded in the same way, the second row being the unit
vector $e_j$. Starting with the $q^{r-1}$ solutions of
$z_1+\cdots+z_r=0$, there are at most
$r^2+(2r+1)^r\le2(2r+1)^r$ equations to exclude.
The number of remaining ordered solutions is therefore between
$q^{r-1}-2(2r+1)^rq^{r-2}$ and $q^{r-1}$.
Every block of size $r$ has an all-positive representation
after changing signs of its representatives, and has exactly
$2r!$ ordered all-positive representations. To see this, fix
representatives of its classes. If two sign vectors giving zero
are neither equal nor opposite, their half-sum and half-difference
give relations on complementary nonempty proper subsets.
Since $r$ is odd, one of these subsets has odd size. Padding
that relation by opposite pairs gives \eqref{eq:signed} on a
proper subset, contrary to minimality. Thus there are exactly
two choices of signs and $r!$ orders of the classes.
By regularity, $d_r|Y|$ is $r$ times the number of these blocks,
so
\[
 \frac{q^{r-1}-2(2r+1)^rq^{r-2}}{(r-1)!(q-1)}
 \le d_r\le\frac{q^{r-1}}{(r-1)!(q-1)}.
\]
For $q\ge4(2r+1)^r$, these bounds imply
$C^{-1}q^{r-2}\le d_r\le Cq^{r-2}$.
For two top-rank common-link equations,
proportional coefficient rows have proportionality $\pm1$;
consistency forces the omitted representatives to agree up to
sign. This is impossible at distinct vertices of $Y$. The
remaining pairs of equations have at most $q^{r-3}$ solutions.
After fixing representatives and an ordering of the classes,
there are at most $2^{2r}$ pairs of full coefficient vectors,
so $|N_r(x)\cap N_r(y)|\le2^{2r}q^{r-3}\le Cq^{r-3}$.
The proof of Lemma~\ref{lem:criterion}, with $M=q$,
verifies the hypotheses of Theorem~\ref{thm:main} once
$q\ge\max\{4(2r+1)^r,C^{4(r-1)}\}$.
The degree estimates give $d_r/q^{r-2}\to1/(r-1)!$ and
$\log D/\log q\to(r-2)/(r-1)$.
Since $|Y|=(q-1)/2$ and
$\rho\le(\log D/d_r)^{1/(r-1)}$, Theorem~\ref{thm:main}
gives a maximal independent set $I\subseteq Y$ satisfying
\[
 |I|\le\left(\frac{1}{2}\bigl((r-2)(r-2)!\bigr)^{\frac{1}{r-1}}+o(1)\right)
 (q\log q)^{\frac{1}{r-1}}.
\]

Let $A=\{a\in\F_q^*:[a]\in I\}$. Since $A=-A$, the absence of
\eqref{eq:signed} is equivalent to $(2h,1)$-freeness. It remains
to prove maximality with respect to adjoining a single element.
For $x\in\F_q^*\setminus A$, maximality of $I$ gives a signed
relation in $A\cup\{x,-x\}$. Its net coefficient of $x$ cannot
vanish, since cancellation and padding by opposite pairs from
$A$ would then give a relation in $A$. Cancel opposite occurrences
of $x$ and replace them by opposite pairs from $A$. Reverse all
signs if necessary. The resulting relation has $t$ copies of $x$
with positive sign, where $1\le t\le r-1$; the case $t=r$ is
excluded by the characteristic assumption. Signs on the other
terms may be absorbed into their values in $A$. Moving one such
term to the right yields an equality of $2h$ summands and one
summand in $A\cup\{x\}$. Finally, $0$ cannot be adjoined, since
$2h\cdot0=0$. Thus $A$ is maximal and $|A|=2|I|$.
\end{proof}

\subsection{Pythagorean triples}

For any prime power $q$, a subset $A\subseteq\F_q^*$ is called
\emph{Pythagorean-free} if it contains no three pairwise distinct
elements $x,y,z$ satisfying $x^2+y^2=z^2$. Denote by
$\gamma_{\rm P}(q)$ the least cardinality of a maximal such subset of $\mathbb F_q^{\times}$.

{\bf Fields $\mathbb F_q$ of odd characteristic with $q \neq 3^{2m}$.} The case of finite fields of odd characteristic 
can be treated fairly uniformly via Theorem \ref{thm:main}, except when $q = 3^{2m}$, in which case an 
explicit maximal Pythagorean triple free set of size $2(\sqrt{q} - 1)$ will be constructed. 
The relation to the freeness convention allowing repeated entries in Theorem~\ref{thm:Pythag} is proved separately below.

\begin{theorem}\label{thm:pyth}
For every odd prime power $q$,
\begin{equation}\label{eq:pyth-lower}
 \gamma_{\rm P}(q)\ge
 \begin{cases}
 \displaystyle\Big\lceil\frac{1+\sqrt{8q-7}}{4}\Big\rceil,
     &q\equiv3\pmod4,\\[6pt]
 \displaystyle\Big\lceil\frac{1+\sqrt{3q-2}}{3}\Big\rceil,
     &q\equiv1\pmod4.
 \end{cases}
\end{equation}
Let $p=\operatorname{char}\F_q$, and let $\chi$ be the quadratic
character of $\F_q$. As $q\to\infty$ through odd prime powers,
\begin{equation}\label{eq:pyth-upper}
 \gamma_{\rm P}(q)\le\bigl(c(q)+o(1)\bigr)\sqrt{q\log q},
\end{equation}
where 
\begingroup
\renewcommand{\arraystretch}{1.3}
\begin{center}
\begin{tabular}{|l|l|}
\hline
 \text{\rm Condition}& $c(q)$ \\ \hline 
 $q\equiv3\pmod4,\ \chi(2)=-1$ & $\frac{2}{\sqrt{3}}$ \\
 \hline
 $q\equiv3\pmod4,\ \chi(2)=1$ & $\sqrt3$ \\
 \hline
  $q\equiv1\pmod4,\ p\ne3,\ (\chi(2)=-1\text{ or }p=5)$ & $2$ \\
  \hline
 $q\equiv1\pmod4,\ p\notin\{3,5\},\ \chi(2)=1$ & $3$ \\
 \hline
 $q=3^{2m}$ & $0$. \\
 \hline
\end{tabular}
\end{center}
\endgroup
The rows are exhaustive for odd $q$.
The error tends to zero uniformly over the indicated fields.
In particular, for all $q$, 
$$\gamma_{\rm P}(q)\le(3+o(1))\sqrt{q\log q}.$$
\end{theorem}

\begin{proof}
Let $A$ be maximal and put $s=|A|$. Every $t\notin A$ is a
completion of two distinct elements $a,b\in A$ under one of
\[
 t^2=a^2+b^2,\qquad t^2=a^2-b^2,\qquad t^2=b^2-a^2.
\]
An unordered pair has at most six completions. If $-1$ is not a
square, at most one of the last two right-hand sides is a nonzero
square, so there are at most four completions. Consequently
\[
 q-1\le
 \begin{cases}
 s+4\binom{s}{2}=2s^2-s,&q\equiv3\pmod4,\\
 s+6\binom{s}{2}=3s^2-2s,&q\equiv1\pmod4.
 \end{cases}
\]
Solving these quadratic inequalities proves \eqref{eq:pyth-lower}.

Put $Q=(\F_q^*)^2$, $\Gamma=Q\cap\{1,-1\}$, $g=|\Gamma|$,
and $Y=Q/\Gamma$. Thus $|Y|=(q-1)/(2g)$. Denote the class of
$a\in Q$ by $[a]$. Let $\mathcal S_3$ consist of the triples of
distinct classes admitting representatives $a,b,c\in Q$ with
$a+b=c$. Define
\[
 \mathcal S_2=
 \begin{cases}
 \{\{[a],[2a]\}:a\in Q,\ [a]\ne[2a]\},&2\in Q,\\
 \varnothing,&2\notin Q,
 \end{cases}
\]
and let
\[
 \HH=\mathcal S_2\ \cup\
 \{B\in\mathcal S_3: \text{no member of }\mathcal S_2
                         \text{ is contained in }B\}.
\]

The group $Q$ acts on $Y=Q/\Gamma$ by multiplication, transitively, and
preserves $\mathcal S_2$ and $\mathcal S_3$, hence $\mathcal H$. Each $d_i(\{y\})$
is therefore independent of $y$, so $\mathcal H$ is regular and below
$d_i$ denotes this common value. Moreover, $d_2\le2$.

Put \(\varepsilon=1\) if \(-1\) is a square in \(\F_q\), and \(\varepsilon=-1\) otherwise. We claim that
\begin{equation}\label{eq:square-count}
 |\{a\in Q:1-a\in Q\}|
 =|\{a\in Q:1+a\in Q\}|
 =\frac{q-4-\varepsilon}{4}.
\end{equation}
For the first count, consider \(x^2+y^2=1\). A standard parametrisation from the point \((1,0)\) is
\[
 x=\frac{1-t^2}{1+t^2},\qquad
 y=\frac{2t}{1+t^2}.
\]
If \(-1\) is not a square, this gives \(q+1\) points; if it is a square, the two values with \(1+t^2=0\) are excluded, giving \(q-1\) points. Removing the four points with \(xy=0\), and dividing by four to account for the choices of signs of the square roots, gives \((q-4-\varepsilon)/4\).

For the second count, the equation \(y^2-x^2=1\) is equivalent to
\[
 (y-x)(y+x)=1.
\]
It therefore has \(q-1\) solutions. Removing the two solutions with \(x=0\), and also the two solutions with \(y=0\) if \(-1\) is a square, leaves \(q-3-2\mathbf 1_{\{\varepsilon=1\}}=q-4-\varepsilon\) solutions with \(xy\ne0\). Again, dividing by four gives the claimed count.

Rescaling a fixed square to $1$, we compute its degree in the system of
triples $\{a,b,a+b\}$ of distinct elements of $Q$. A three-element
subset of $Q$ has at most one element which is the sum of the other
two, since $a+b=c$ together with $a+c=b$ or $b+c=a$ would give
$2a=0$ or $2b=0$. Put $N=(q-4-\chi(-1))/4$. By \eqref{eq:square-count},
the number of triples in which $1$ is the sum lies between
$(N-1)/2$ and $N/2$, and the number in which it is a summand
lies between $N-1$ and $N$. The only possible coincidences are
$a=1-a$ in the first count and $b=1$ in the second.
The degree is therefore between $3N/2-3/2$ and $3N/2$.

If $g=2$, each triple of distinct classes in $\mathcal S_3$
lifts to exactly three triples on $Q$ containing a fixed
representative of one of its classes. Triples involving both
$a$ and $-a$ contribute at most three incidences at a fixed vertex
$v$: the possible triples are
$\{v,-v,2v\}$, $\{v,-v,-2v\}$, and $\{v,v/2,-v/2\}$.
The degree in $\mathcal S_3$ is therefore between $q/8-17/8$
and $q/8$ in this case, and between $3q/8-21/8$ and $3q/8$
when $g=1$.
Deleting triples containing a pair $\{[a],[2a]\}$ decreases the
common degree by at most nine. Indeed, for a pair not containing
the given class, the third class is represented by $3a$, $a$, or
$-a$. Each nonzero multiplier determines at most one class $[a]$.
The at most two pairs containing the given class each lie in at
most three triples, since a pair represented by $a,b$ has
completions represented by $a+b$, $a-b$, or $b-a$. Thus
\begin{equation}\label{eq:degrees}
 \alpha q-16\le d_3\le\alpha q,\qquad
 \alpha=\begin{cases}\frac{3}{8},&g=1,\\\frac{1}{8},&g=2.\end{cases}
\end{equation}
These estimates hold for every odd $q$.

Every pair of vertices has at most three completions to a triple.
For a common link at distinct vertices $[v],[w]$, the block
conditions give two linear equations in representatives $a,b$
of the common vertices, with coefficients $\pm1$. Dependent
rows can be consistent only if $w=\pm v$. If $g=1$, then
$w=-v$ is impossible since $v,w\in Q$; if $g=2$, either sign
would imply $[v]=[w]$. Independent rows have determinant $\pm2$,
which is nonzero. Fix representatives of all classes and an
ordering of $Y$. After normalising the coefficient of $v$ or $w$
to $1$, each row has four possible pairs of coefficients of $a,b$.
The sixteen pairs of rows therefore give
$|N_3([v])\cap N_3([w])|\le16$.

Suppose first that $q\ge2^{20}$. By \eqref{eq:degrees},
$q/16\le d_3\le3q/8$, so $D=\sqrt{d_3}$ and
$(1/4)\log q\le\log D\le(1/2)\log q$.
The inequalities $3\le D^{1/2}$ and $16\le D^{3/2}$,
together with $d_2\le2$, verify all the mixed sparseness
conditions with $\delta=1/2$. The equation
$d_2\rho+d_3\rho^2=(1/2)\log d_3$ gives
\[
 \rho\le\sqrt{\frac{\log q}{2(\alpha q-16)}}.
\]
Let $\log\log\log\log K=24^{16}$ and
$\epsilon=1/(9\cdot2^{22})$. Theorem~\ref{thm:main}, with
these explicit constants, supplies a maximal independent set
$I\subseteq Y$ of size at most
$$\bigl(1+2K(\log q)^{-\epsilon}\bigr) \cdot
\frac{\sqrt{q\log q}}{2g\sqrt{2\alpha - 32/q}}.$$
Here $4^\epsilon\le2$ and
$\bigl(1-16/(\alpha q)\bigr)^{-1/2}\le1+128/q
\le1+(1/2)(\log q)^{-\epsilon}$.
Since $(\log q)^{-\epsilon}\le1$ and $K\ge1$, we obtain
\begin{equation}\label{eq:I-size}
 |I|\le\frac{1+4K(\log q)^{-\epsilon}}{2g\sqrt{2\alpha}}
                 \sqrt{q\log q}
 =\begin{cases}
 \displaystyle{\frac{1+4K(\log q)^{-\epsilon}}{\sqrt3}\sqrt{q\log q}},&g=1,\\
 \displaystyle{\frac{1+4K(\log q)^{-\epsilon}}{2}\sqrt{q\log q}},&g=2.
 \end{cases}
\end{equation}
For $3\le q<2^{20}$, the same bound holds for any maximal
independent set in $Y$: its size is at most $q$, whereas
$4K(\log q)^{-\epsilon}\ge K\ge2^{20}$ and
$\sqrt{q\log q}/(2g\sqrt{2\alpha})\ge1$.

Consider $\pi:\F_q^*\to Y$, given by $\pi(x)=[x^2]$.
Its fibres have size $2g$. First suppose that $(p,g)\ne(3,2)$,
and put $A_0=\pi^{-1}(I)$. A Pythagorean triple on three fibres
would give a member of $\mathcal S_3$ contained in $I$, hence
a block of $\HH$ contained in $I$. A triple on exactly two
fibres would give a pair in $\mathcal S_2$ contained in $I$.
A triple on one fibre is impossible if $g=1$; if $g=2$, it
would imply $2a=-a$ for a nonzero $a$, and hence $p=3$.
Thus $A_0$ is Pythagorean-free.

If $\mathcal S_2=\varnothing$, maximality of $I$ implies that
every unselected class is completed to a triple by two classes
of $I$. Every element of the unselected fibre is then blocked
by two elements of $A_0$ (for $g=2$, multiply the relation by
$-1\in Q$ if necessary, so that it involves the given element).
Thus $A_0$ is maximal and has size
$2g|I|$.

Suppose $\mathcal S_2\ne\varnothing$. Multiplication by $2$
induces a permutation $\tau([a])=[2a]$ of $Y$, without fixed
points. By maximality of $I$, every class outside
$I\cup\tau I\cup\tau^{-1}I$ is blocked by a triple with two
classes of $I$. Every element over $\tau I$ is also blocked:
if its square is $2a$ or $-2a$, two opposite roots of $a$ or
$-a$ already lie in $A_0$. These roots $y,-y$ and the element $x$
are pairwise distinct, since $p$ is odd and $x^2=2y^2\ne y^2$.
Choose a maximal Pythagorean-free set $A\supseteq A_0$.
It follows that
\[
 A\setminus A_0\subseteq\pi^{-1}(\tau^{-1}I\setminus\tau I).
\]
For a square $b$ whose class lies in $\tau^{-1}I$, $A$ can
contain at most one of its two roots: otherwise these roots,
together with an element of $A_0$ whose square is $2b$, form a
Pythagorean triple. Each such fibre therefore contributes at
most $g$ new elements. Consequently
\begin{equation*}
 |A|\le2g|I|+g|\tau^{-1}I\setminus\tau I|\le3g|I|.
\end{equation*}
If $p=5$, then $g=2$ and $\tau^2=1$ whenever
$\mathcal S_2\ne\varnothing$, since $\tau^2([a])=[4a]=[-a]=[a]$.
In this case no new elements
are needed and $|A|=2g|I|$.

It remains to consider $p=3$ and $g=2$, equivalently $q=3^{2m}$.
For each $[a]\in I$, choose one square root of $a$ and one
square root of $-a$, and let $A_0$ consist of these choices.
Here $\mathcal S_2=\varnothing$, since $[2a]=[-a]=[a]$.
A triple on three fibres is excluded by independence of $I$.
Two chosen elements in one fibre have opposite squares; a
Pythagorean equation involving them forces the third square
to be $0$, $2a$, or $-2a$, hence zero or in the same class.
There are only two chosen elements in that fibre, so $A_0$
is Pythagorean-free.
Every element of an unselected fibre is blocked by a triple
with two selected fibres (every representative in $Q$ of a class of
$I$ is the square of a chosen element, and the relation may be
multiplied by $-1\in Q$). Within a selected fibre, write the
two chosen roots as $e,f$, where $e^2=a$ and $f^2=-a$.
The missing elements $-e,-f$ are blocked by
\[
 e^2+(-e)^2=f^2,\qquad f^2+(-f)^2=e^2.
\]
Thus $A_0$ is already maximal and has size $2|I|$.
Combining these lifting bounds with \eqref{eq:I-size} gives
\eqref{eq:pyth-upper} with the stated constants when $q\ne3^{2m}$.
For $q=3^{2m}$, Theorem~\ref{thm:pyth-char-three} gives
$\gamma_{\rm P}(q)\le2\sqrt q-2$, proving the row $q=3^{2m}$.
In all rows the $o(1)$ in \eqref{eq:pyth-upper} can be bounded
above by $12K(\log q)^{-\epsilon}+2/\sqrt{\log q}$,
which tends to zero uniformly.
\end{proof}

\begin{proof}[Proof of Theorem~\textup{\ref{thm:Pythag}} in odd characteristic]
Here freeness excludes all zeros of $f(x,y,z)=x^2+y^2-z^2$,
including tuples with repeated entries. We use the notation
$Q,\Gamma,g,Y,\pi,I$ from the proof of Theorem~\ref{thm:pyth}.
When $(p,g)\ne(3,2)$, the full inverse image $A_0=\pi^{-1}(I)$
is free under this convention. Indeed, a zero using three classes
gives a member of $\mathcal S_3$ contained in $I$; a zero using
exactly two classes gives a member of $\mathcal S_2$ contained
in $I$; and a zero using one class is impossible unless $p=3$ and
$g=2$. These statements do not require the three entries to be distinct.

For $x\notin A_0\cup\{0\}$, maximality of $I$ supplies a block
in $I\cup\{\pi(x)\}$. A block of size three gives a zero using
$x$ and two elements of $A_0$, since $A_0$ contains every root
over each selected class. If the block has size two, there is
$y\in A_0$ with either $x^2=2y^2$ or $2x^2=y^2$.
Thus either $(y,y,x)$ or $(x,x,y)$ is a forbidden tuple.
Also $0$ cannot be adjoined, since $f(0,0,0)=0$.
Consequently $A_0$ is already maximal for this convention and,
by \eqref{eq:I-size},
\[
 \gamma_q(f)\le |A_0|=2g|I|
       \le(3+o(1))\sqrt{q\log q}.
\]
For $q=3^{2m}$, parts \textup{(2)} and \textup{(3)} of
Theorem~\ref{thm:pyth-char-three} directly give the required stronger
freeness and maximality, with size $2\sqrt q-2$.

Finally, let $A$ be any maximal $f$-free set and let $s=|A|$.
An element outside $A\cup\{0\}$ either completes two distinct
elements of $A$, or satisfies $x^2=2a^2$ or $2x^2=a^2$ for some
$a\in A$. The first alternative gives at most six elements per
unordered pair, and the other two give at most four per element.
Hence
\[
 q-1\le s+6\binom{s}{2}+4s=3s^2+2s,
 \qquad
 s\ge\frac{\sqrt{3q-2}-1}{3}.
\]
This proves the stated asymptotic lower bound in odd characteristic.
\end{proof}

\begin{remark}\label{rem:three}
The family $\mathcal S_2$ does not detect obstructions within a
single fibre. This matters when $q=3^{2m}$: although
$\mathcal S_2=\varnothing$, each fibre contains Pythagorean
triples. Choosing roots with opposite squares resolves this
obstruction and gives $\gamma_{\rm P}(q)\le(1+o(1))\sqrt{q\log q}$
by the preceding argument. The following construction improves
this to $2\sqrt q-2$, giving the zero in the row $q=3^{2m}$ of
Theorem~\ref{thm:pyth}.
\end{remark}

{\bf Square fields of characteristic three.} Surprisingly, the case of fields of order $q = 3^{2m}$ is exceptional. In this case, there is a simple explicit construction of a maximal Pythagorean free set of size $2(\sqrt{q} - 1)$, as follows: 

\begin{theorem}\label{thm:pyth-char-three}
Let $t\ge1$, let $s=3^t$, and let $q=s^2$.
Write $K=\F_s\subseteq F=\F_q$, and let
$\theta\in F^\times$ be any nonsquare. Define
\begin{equation}\label{eq:A}
A=\{u^3+\theta u:u\in F^\times,\ u^2\in K\}.
\end{equation}
Then:
\begin{enumerate}
\item[\textup{(1)}] $ |A|=2s-2=2\sqrt q-2$.
\item[\textup{(2)}] $x^2+y^2\ne z^2$ for all $x,y,z\in A$,
even when the entries are allowed to coincide.
\item[\textup{(3)}] For every $v\in F\setminus A$, the set
$A\cup\{v\}$ contains three pairwise distinct elements
$x,y,z$ satisfying $x^2+y^2=z^2$.
\end{enumerate}
In particular, if $q = 3^{2m}$, then 
$$\gamma_{\rm P}(q)\le 2(\sqrt{q} - 1).$$
\end{theorem}

\begin{proof}
Every element $a\in K^\times$ is a square in $F$, since
\[
a^{\frac{s^2-1}{2}}=(a^{s-1})^{\frac{s+1}{2}}=1.
\]
Thus $\theta\notin K$, and $\{1,\theta\}$ is a basis of $F$
over $K$. Let
\[
\tau=\theta+\theta^s,\qquad \nu=\theta^{s+1}.
\]
Both $\tau$ and $\nu$ belong to $K$, and
\begin{equation}\label{eq:minpoly}
\theta^2-\tau\theta+\nu=0.
\end{equation}
Moreover,
\[
\nu^{\frac{s-1}{2}}=\theta^{\frac{s^2-1}{2}}=-1,
\]
so $\nu$ is a nonsquare in $K$.

For $a\in K$, define
\[
L(a)=a^3-\nu a,\qquad M(a)=\tau a-a^2.
\]
Using \eqref{eq:minpoly} and $2=-1$, we obtain
\begin{equation}\label{eq:P}
P(a):=a(a+\theta)^2=L(a)+\theta M(a).
\end{equation}
The map $L:K\to K$ is additive. Also,
\[
L(a)=0\ \Longrightarrow\ a(a^2-\nu)=0\
\Longrightarrow\ a=0,
\]
because $\nu$ is a nonsquare in $K$. Hence $L$ is bijective.
Comparison of the coefficients of $1$ in \eqref{eq:P} shows that
$P$ is injective as well.

\medskip
\noindent\textit{Proof of \textup{(1)}.}
For $u\in F$ with $u^2\in K$,
\begin{equation}\label{eq:cubic-square}
(u^3+\theta u)^2=u^2(u^2+\theta)^2=P(u^2).
\end{equation}
Every $a\in K^\times$ has two square roots $u,-u$ in $F$.
Their images $u^3+\theta u$ and $-u^3-\theta u$ are nonzero:
indeed, $u^3+\theta u=u(a+\theta)\ne0$.
By \eqref{eq:cubic-square}, these are exactly the two square roots
of $P(a)$. Therefore the set in \eqref{eq:A} satisfies
\begin{equation}\label{eq:roots}
A=\{x\in F^\times:x^2\in P(K^\times)\}.
\end{equation}
Since $P$ is injective, we have $A=-A$ and
\[
|A|=2|P(K^\times)|=2(s-1).
\]

\medskip
\noindent\textit{Proof of \textup{(2)}.}
Suppose $x,y,z\in A$ satisfy $x^2+y^2=z^2$.
By \eqref{eq:roots}, there are $a,b,c\in K^\times$ such that
\[
x^2=P(a),\qquad y^2=P(b),\qquad z^2=P(c).
\]
Thus $P(a)+P(b)=P(c)$. Comparing the coefficients of $1$
in \eqref{eq:P} gives $L(a)+L(b)=L(c)$.
Since $L$ is additive and injective, $c=a+b$.
Comparison of the coefficients of $\theta$ now gives
\[
0=M(a)+M(b)-M(a+b)=2ab,
\]
contrary to $a,b\ne0$.

\noindent\textit{Proof of \textup{(3)}.}
First let $v\in F^\times\setminus A$, and write
$v^2=a+b\theta$ with $a,b\in K$.
We have $a\ne0$. Otherwise $b\ne0$ and $v^2=b\theta$
would be a nonsquare in $F$, since $b\in K^\times$ is a square
in $F$ whereas $\theta$ is a nonsquare.

Let $\delta=L^{-1}(a)$, so $\delta\ne0$.
There are unique $\lambda,\mu\in K$ satisfying
\[
\lambda-\mu=\delta,\qquad
\lambda+\mu=\tau-\frac{b}{\delta},
\]
because $2\ne0$. In particular, $\lambda\ne\mu$.
By \eqref{eq:P},
\begin{equation}\label{eq:difference}
\begin{aligned}
P(\lambda)-P(\mu)
&=L(\lambda-\mu)
  +\theta(\lambda-\mu)(\tau-\lambda-\mu)\\
&=a+b\theta=v^2.
\end{aligned}
\end{equation}
We must have $\mu\ne0$. If $\mu=0$, then
$\lambda\ne0$, and \eqref{eq:difference} together with
\eqref{eq:roots} would imply $v\in A$.

If $\lambda\ne0$, let $x,y\in A$ satisfy
$x^2=P(\mu)$ and $y^2=P(\lambda)$.
Equation \eqref{eq:difference} gives
\[
x^2+v^2=y^2.
\]
Injectivity of $P$ implies $x^2\ne y^2$, hence $x\ne y$.
Since $v\notin A$, the elements $x,v,y$ are pairwise distinct.

If $\lambda=0$, let $x\in A$ satisfy $x^2=P(\mu)$.
As $A=-A$, both $x$ and $-x$ belong to $A$, and
\[
x^2+(-x)^2=2P(\mu)=-P(\mu)=v^2.
\]
Here $x\ne-x$, since $x\ne0$, and $v\notin A$.
Again the three entries are pairwise distinct.

Finally, if $v=0$, any $x\in A$ gives the three pairwise
distinct elements $0,x,-x$, with
\[
0^2+x^2=(-x)^2.
\]
This proves part~\textup{(3)} and completes the proof.
\end{proof}

{\bf Characteristic two.} The case of characteristic two reduces Pythagorean triple free sets to sum-free sets.

\begin{proposition}\label{prop:char-two}
Let $q=2^m$. Then
\[
 \gamma_{\rm P}(q)=\gamma_{2,1}(\F_q),\qquad
 \gamma_{\rm P}(q)\ge
 \left\lceil\frac{\sqrt{8q-7}-1}{2}\right\rceil.
\]
Moreover, as $m\to\infty$,
\begin{equation}\label{eq:binary}
 \gamma_{\rm P}(2^m)\le
 \begin{cases}
 (2+o(1))\sqrt{2^m},&m\text{ even},\\
 \Bigl(\displaystyle{\frac{3}{\sqrt{2}}+o(1)}\Bigr)\sqrt{2^m},&m\text{ odd}.
 \end{cases}
\end{equation}
In particular, $\gamma_{\rm P}(2^m) \geq 2^{(m-1)/2}$ for $m\ge1$.
\end{proposition}

\begin{proof}
In characteristic two, the Pythagorean equation is equivalent to
$x+y=z$. On nonzero elements such an equality has three distinct
entries, and a sum-free set cannot contain $0$. Thus the two
notions of freeness and maximality coincide. If $A$ is maximal
sum-free of size $s$, every nonzero element outside $A$ is a sum
of two elements of $A$. Hence
\[
 q-1\le s+\binom{s}{2},
\]
giving the lower bound.

Identify $(\F_q,+)$ with $\F_2^m$. Maximal sum-free sets are
precisely complete caps in $\operatorname{PG}(m-1,2)$.
By~\cite{IW}, for $n\ge3$ there are complete caps in
$\operatorname{PG}(n,2)$ of size
$3(2^{n/2}-1)$ for even $n$ and $2^{(n+3)/2}-3$ for odd $n$.
Taking $n=m-1$ gives the bounds $2\sqrt{2^m}-3$ for even
$m\ge4$ and $(3/\sqrt2)\sqrt{2^m}-3$ for odd $m\ge5$,
proving \eqref{eq:binary}. For $m\le3$, the lower
bound is attained by a coset of a subgroup of index two. Its size
$2^{m-1}$ is at most $(3/\sqrt2)\sqrt{2^m}$ for these three values
of $m$. The displayed lower bound and these constructions prove
Theorem~\ref{thm:Pythag} in characteristic two.
\end{proof}

\section{Geometric applications}
\label{sec:pg-circuits}

Throughout this section, $r \ge k \ge 2$ are fixed integers, $q$ is a
prime power, and asymptotic statements refer to $q \to \infty$.
The ambient projective dimension is $r$; the largest block size is
$k + 1$. We put $K_{1} = 2^{k - 1}k!$ and
\[
 X = \PG(r,q),\qquad n = |X| = \frac{q^{r + 1} - 1}{q - 1}.
\]
Here and below, $X$ denotes the point set of the projective space.
Representatives of its points are nonzero vectors of
$\mathbb F_q^{r + 1}$, considered up to nonzero scalar multiplication.
A set of points is \emph{projectively independent} if its vector
representatives are linearly independent.

\subsection{\texorpdfstring{The set system on $\PG(r,q)$}{\detokenize{The set system on PG(r,q)}}}

\begin{definition}
A \emph{projective circuit} is a projectively dependent set of
points all of whose proper subsets are projectively independent.
Let $\mathcal S$ be the set system on $X$ whose blocks are the
projective circuits of sizes $3,\ldots,k + 1$.
\end{definition}

\emph{Remark.} An $i$-point circuit spans an $(i - 2)$-dimensional projective
subspace. The family $\mathcal S$ is an antichain, since a circuit
cannot contain a dependent proper subset. Projective transformations
preserve $\mathcal S$ and act transitively on $X$. Thus $\mathcal S$
is exactly regular in every rank. Lemma~\ref{lem:pg-circuit-counts}
computes the degrees and intersection counts needed to apply
Theorem~\ref{thm:main}. We write
\[
 \vb{d} = (d_2,d_3,\ldots,d_{k + 1}),\qquad d_2 = 0,\qquad D = D(\vb{d}) = \max_{3 \le i \le k + 1}d_i^{\frac{1}{i - 1}}.
\]
The notation $d_i(A)$ and $d_i(x,y)$ has the meaning given in
Section~1. In particular, $d_i(\{x,y\})$ is a pair codegree,
whereas $d_i(x,y)$ is a common-link count.

\subsection{\texorpdfstring{Formulas for $d_i$, $d_i(A)$ and $d_i(x,y)$}{\detokenize{Formulas for d_i, d_i(A) and d_i(x,y)}}}

Background on finite projective geometry is given by Hirschfeld and Thas~\cite{HT}.
Empty products are taken to be $1$.

\begin{lemma}\label{lem:pg-circuit-counts}
For $3 \le i \le k + 1$, the following statements hold.
\begin{enumerate}[leftmargin=*,itemsep=.35em]
\item Every point has rank-$i$ degree
\begin{equation}\label{eq:pg-circuit-degree}
 d_i = \frac{1}{(i - 1)!}
 \prod_{j = 1}^{i - 2}(q^{r + 1} - q^j).
\end{equation}
\item Let $A \subseteq X$ and $2 \leq z = |A| < i$. If $A$ is projectively dependent, then $d_i(A) = 0$. If $A$ is projectively independent, then
\begin{equation}\label{eq:pg-circuit-codegree}
 d_i(A) = \frac{(q - 1)^{z - 1}}{(i - z)!}
 \prod_{j = z}^{i - 2}(q^{r + 1} - q^j).
\end{equation}
\item For distinct $x,y \in X$,
\begin{equation}\label{eq:pg-circuit-common}
 d_i(x,y) = d_i\,
 \frac{(q - 1)^{i - 2} - 1}{n - 1}.
\end{equation}
\end{enumerate}
\end{lemma}

\begin{proof}[Proof of \textup{(1)}]
We fix $x \in X$. We choose an ordered list of $i - 2$ further points
which, together with $x$, is projectively independent. When the
points already chosen have vector span of dimension $j$, the
number of possible next points is
\[
 n - \frac{q^j - 1}{q - 1} = \frac{q^{r + 1} - q^j}{q - 1}.
\]
The resulting $(i - 1)$-set is a projective basis of its span.
A point completes this basis to a circuit precisely when all
its coordinates in this basis are nonzero. There are
$(q - 1)^{i - 2}$ such projective points: there are
$(q - 1)^{i - 1}$ nonzero-coordinate vectors, and each projective
point has $q - 1$ vector representatives.

Each circuit containing $x$ is counted $(i - 1)!$ times, since
any of its other $i - 1$ points can be chosen last and the remaining
$i - 2$ points can be ordered arbitrarily. Consequently
\[
 d_i = \frac{(q - 1)^{i - 2}}{(i - 1)!}
 \prod_{j = 1}^{i - 2}
 \frac{q^{r + 1} - q^j}{q - 1},
\]
which is \eqref{eq:pg-circuit-degree}.

\emph{Proof of \textup{(2)}.}
A circuit has no dependent proper subset, so the assertion is
immediate when $A$ is dependent. Otherwise, we choose $i - z - 1$
ordered points extending $A$ to an independent $(i - 1)$-set,
and then choose the last point with all coordinates nonzero.
Counting these choices gives
\[
 d_i(A) = \frac{(q - 1)^{i - 2}}{(i - z)!}
 \prod_{j = z}^{i - 2}
 \frac{q^{r + 1} - q^j}{q - 1}.
\]
The divisor is $(i - z)!$: there are $i - z$ choices for the last
point and $(i - z - 1)!$ orders for the preceding added points.
Cancelling the powers of $q - 1$ proves
\eqref{eq:pg-circuit-codegree}.

\emph{Proof of \textup{(3)}.}
We fix $x \in X$. Every $B \in N_i(x)$ is an independent $(i - 1)$-set.
Exactly $(q - 1)^{i - 2}$ points complete $B$ to a circuit, one of
which is $x$. Counting the pairs consisting of $B \in N_i(x)$
and another such point gives
\[
 \sum_{y \in X\setminus\{x\}}d_i(x,y) = d_i\bigl((q - 1)^{i - 2} - 1\bigr).
\]
The projective transformations fixing $x$ act transitively on
$X\setminus\{x\}$, so every summand on the left is equal.
Division by $n - 1$ proves \eqref{eq:pg-circuit-common}.
\end{proof}

\Needspace{18\baselineskip}
\subsection{\texorpdfstring{The size $|I|$ of a $k$-complete arc}{\detokenize{The size |I| of a k-complete arc}}}

The counts in Lemma~\ref{lem:pg-circuit-counts} give the following
verification of the hypotheses of Theorem~\ref{thm:main}.

\begin{corollary}\label{cor:pg-circuit-sparsity}
\begin{enumerate}[leftmargin=*,itemsep=.35em]
\item For every fixed
$0 < \delta < (r - k + 1)/((r + 1)(k - 1))$, the set system $\mathcal S$
is a $\delta$-sparse $(n,\vb{d})$-antichain for all sufficiently
large $q$.
\item Every maximal independent set in $\mathcal S$ is a
$k$-complete arc in $\PG(r,q)$.
\end{enumerate}
\end{corollary}

\begin{proof}[Proof of \textup{(1)}]
Since $q \ge 2$ and $r \ge k$,
every factor $1 - q^{j - r - 1}$ with $1 \le j \le k - 1$ is at least
$1/2$. Equation~\eqref{eq:pg-circuit-degree} gives
\[
 d_{k + 1} \ge K_{1}^{-1}q^{(r + 1)(k - 1)},\qquad D \ge K_{1}^{-\frac{1}{k}}q^{\frac{(r + 1)(k - 1)}{k}}.
\]
For $3 \le i \le k + 1$, \eqref{eq:pg-circuit-degree} gives
\begin{equation}\label{eq:pg-circuit-degree-asymptotic}
 d_i = (1 + o(1))\frac{q^{(r + 1)(i - 2)}}{(i - 1)!}.
\end{equation}
For projectively independent $A \subseteq X$ with
$2 \le z = |A| < i$, \eqref{eq:pg-circuit-codegree} gives
\begin{equation*}
 d_i(A) = (1 + o(1))
 \frac{q^{(r + 1)(i - 2) - r(z - 1)}}{(i - z)!}.
\end{equation*}
For distinct $x,y \in X$, \eqref{eq:pg-circuit-common} gives
\begin{equation*}
 d_i(x,y) = (1 + o(1))
 \frac{q^{(r + 2)(i - 2) - r}}{(i - 1)!}.
\end{equation*}
Uniformity over the point sets follows from
\eqref{eq:pg-circuit-degree}--\eqref{eq:pg-circuit-common}.
The exponents $((r + 1)(i - 2))/(i - 1)$ increase strictly with $i$.
Consequently, for all sufficiently large $q$,
\begin{equation}\label{eq:pg-circuit-scale}
 D = d_{k + 1}^{\frac{1}{k}} = (1 + o(1))(k!)^{-\frac{1}{k}}q^{\frac{(r + 1)(k - 1)}{k}}.
\end{equation}

For $3 \le i \le k + 1$ and $2 \le z < i$, the difference between
the powers of $q$ in $D^{i - z}$ and $d_i(A)$ is
\begin{align*}
 \frac{(r + 1)(k - 1)}k(i - z)
 -\bigl((r + 1)(i - 2) - r(z - 1)\bigr)
 & = r + 2 - i - \frac{r - k + 1}{k}(i - z)\\
 & \ge r - k + 1 - \frac{r - k + 1}{k}(k - 1)\\
 & = \frac{r - k + 1}{k}.
\end{align*}
Thus, for every projectively independent $A$ of size $z$ in
this range, \eqref{eq:pg-circuit-codegree} and $0 < i - z - \delta \le k$ give
\[
 \frac{d_i(A)}{D^{i - z - \delta}}
 \le K_{1}
 q^{-\frac{r - k + 1}{k} + \frac{\delta(r + 1)(k - 1)}{k}} = o(1).
\]
The exponent is negative by the assumed upper bound on $\delta$.
For dependent $A$ the codegree is zero. Hence, for
$3 \le i \le k + 1$ and $2 \le |A| < i$, all sufficiently large $q$
satisfy
\begin{equation*}
 d_i(A) \le D^{i - |A| - \delta}.
\end{equation*}

For $3 \le i \le k + 1$, the corresponding difference for common
links is
\begin{align*}
 \frac{(r + 1)(k - 1)}k(i - 1)
 -\bigl((r + 2)(i - 2) - r\bigr)
 & = 2(r + 2 - i) - \frac{r - k + 1}{k}(i - 1)\\
 & \ge 2(r - k + 1) - (r - k + 1)\\
 & = r - k + 1.
\end{align*}
For distinct $x,y \in X$, it follows that
\[
 \frac{d_i(x,y)}{D^{i - 1 - \delta}}
 \le K_{1}
 q^{-(r - k + 1) + \frac{\delta(r + 1)(k - 1)}{k}} = o(1).
\]
Here \eqref{eq:pg-circuit-common}, $n - 1 \ge q^r$ and
$d_i \le q^{(r + 1)(i - 2)}/((i - 1)!)$ bound the numerator by
$q^{(r + 2)(i - 2) - r}$, while $0 < i - 1 - \delta \le k$ bounds
the constant arising from the denominator by $K_{1}$.
Consequently, for $3 \le i \le k + 1$ and distinct $x,y \in X$,
all sufficiently large $q$ satisfy
\begin{equation*}
 d_i(x,y) \le D^{i - 1 - \delta}.
\end{equation*}
Rank two has no blocks. Regularity follows from
\eqref{eq:pg-circuit-degree}, and circuits form an antichain.

\emph{Proof of \textup{(2)}.}
Independence in $\mathcal S$ means that every subset of at most
$k + 1$ points is projectively independent: any dependent set
contains a circuit. A maximal independent set $I$ has at least
$k + 1$ points. Otherwise $I$ is projectively independent, and
its span is a proper subspace of $\PG(r,q)$, so a point
outside this span could be added.

Since $|I| \ge k + 1$, the independence condition is equivalent
to having no $k + 1$ points in a $(k - 1)$-dimensional projective
subspace. Indeed, any dependent subset of size at most $k + 1$
can be extended within $I$ to a dependent $(k + 1)$-subset.
Finally, if $x \notin I$, maximality supplies a circuit in
$I\cup\{x\}$. Extending this circuit to a $(k + 1)$-subset of
$I\cup\{x\}$ shows that $x$ cannot be added while preserving
the geometric condition. This proves the required maximality.
\end{proof}

\begin{proof}[Proof of Theorem~\textup{\ref{thm:k-complete}}]
We fix $\delta = (r - k + 1)/(2(r + 1)(k - 1))$.
Corollary~\ref{cor:pg-circuit-sparsity}~(1) lets us apply
Theorem~\ref{thm:main}, with upper rank
$k + 1$ and $n$ atoms. It produces a maximal independent set
$I$ in $\mathcal S$ satisfying
\begin{equation}\label{eq:pg-circuit-main-application}
 |I| = (1 + o_{r,k}(1))n\rho(\vb{d}),\qquad \sum_{i = 3}^{k + 1}d_i\rho(\vb{d})^{i - 1} = \log D.
\end{equation}
We compute this quantity explicitly.

We put $\kappa = \left((\log D)/d_{k + 1}\right)^{1/k}$. Equations
\eqref{eq:pg-circuit-degree-asymptotic} and
\eqref{eq:pg-circuit-scale} give
\[
 \kappa \sim \left(
 \frac{(r + 1)(k - 1)(k - 1)!\log q}
 {q^{(r + 1)(k - 1)}}
 \right)^{\frac{1}{k}}.
\]
We put $K_{2} = (r + 1)K_{1}$.
Equation~\eqref{eq:pg-circuit-degree} gives
$D \le q^{((r + 1)(k - 1))/k}$ and
\[
 \kappa \le K_{2}^{\frac{1}{k}}q^{-\frac{(r + 1)(k - 1)}{k}}(\log q)^{\frac{1}{k}}.
\]
For each $3 \le i \le k$, we consequently have
\[
 d_i \kappa^{i - 1} \le K_{2} q^{(r + 1)(i - 2) - \frac{(r + 1)(k - 1)(i - 1)}{k}} (\log q)^{\frac{i - 1}{k}} = K_{2} q^{\frac{(r + 1)(i - k - 1)}{k}}(\log q)^{\frac{i - 1}{k}} = o(1).
\]
When $k = 2$, the sum over these ranks is empty and equals zero.
Since the sum in \eqref{eq:pg-circuit-main-application} has nonnegative
terms, $\rho(\vb{d}) \le \kappa$. Therefore
\[
 0 \le 1 - \left(\frac{\rho(\vb{d})}{\kappa}\right)^k = \frac{\displaystyle\sum_{i = 3}^{k} d_i\rho(\vb{d})^{i - 1}}{\log D} \le \frac{\displaystyle\sum_{i = 3}^{k}d_i \kappa^{i - 1}}{\log D} = o(1).
\]
It follows that $\rho(\vb{d}) \sim \kappa$. Substituting this into
\eqref{eq:pg-circuit-main-application}, and using $n \sim q^r$,
gives
\begin{align*}
 |I|
 & = (1 + o_{r,k}(1))n\rho(\vb{d})\\
 & = (1 + o_{r,k}(1))q^r
 \left(
 \frac{(r + 1)(k - 1)(k - 1)!\log q}
 {q^{(r + 1)(k - 1)}}
 \right)^{\frac{1}{k}}\\
 & = (1 + o_{r,k}(1))
 \left((r + 1)(k - 1)(k - 1)!\,
 q^{rk - (r + 1)(k - 1)}\log q\right)^{\frac{1}{k}}\\
 & = (1 + o_{r,k}(1))
 \left((r + 1)(k - 1)(k - 1)!\,
 q^{r - k + 1}\log q\right)^{\frac{1}{k}},
\end{align*}
proving \eqref{eq:pg-circuit-size}. By
Corollary~\ref{cor:pg-circuit-sparsity}~(2), this $I$ is a
$k$-complete arc.
\end{proof}

\subsection{Collinear subsets and partial ovoids}

We verify the parameters for Theorems~\ref{thm:plane-collinear}
and~\ref{thm:classical-geometries}.

\begin{proof}[Proof of Theorem~\textup{\ref{thm:plane-collinear}}]
Let $\mathcal S$ consist of the collinear $(k+1)$-subsets of
$\PG(2,q)$. Each point lies on $q+1$ lines, so the number of atoms
and the common degree are
\[
 n=q^2+q+1,\qquad d=(q+1)\binom qk.
\]
For a set $A$ with $2\le|A|\le k$, its codegree is zero unless
$A$ is collinear, and in all cases
\[
 d_{k+1}(A)\le\binom{q+1-|A|}{k+1-|A|}.
\]
For two distinct points $x,y$, a common link consists of $k$
points on the line $xy$, excluding $x$ and $y$. Consequently
\[
 |N_{k+1}(x)\cap N_{k+1}(y)|=\binom{q-1}{k}.
\]
Since $d\sim q^{k+1}/k!$, any fixed $0<\delta<1/(k+1)$
satisfies both inequalities in Definition~\ref{def:sparse} for
all sufficiently large $q$. Indeed, writing $a=|A|$, the difference
between the exponents of $q$ on their right and left sides is
$(k+1-a-(k+1)\delta)/k>0$ for the codegree bound, and
$1-(k+1)\delta/k>0$ for the common-link bound.
Theorem~\ref{main} therefore gives a maximal independent set with
\[
 |I|=(1+o_k(1))n\left(\frac{\log d}{kd}\right)^{1/k}
   =(1+o_k(1))
       \bigl((k+1)(k-1)!q^{k-1}\log q\bigr)^{1/k},
\]
which is \eqref{eq:plane-collinear-size}. Independence and maximality
are precisely the asserted geometric properties.
\end{proof}

\begin{proof}[Proof of Theorem~\textup{\ref{thm:classical-geometries}}]
We first consider the five polar-space rows. The point and residue
counts used below are the standard counts for classical polar spaces;
see~\cite{HT,Thas,deBKMS}. Write $Q$ for the size of the coordinate
field, $m$ for the rank, and $e$ for the parameter in the following table:
\[
\begin{array}{c|ccc}
\mathcal Q&Q&m&e\\\hline
Q(2n,q)&q&n&1\\
Q^+(2n+1,q)&q&n+1&0\\
Q^-(2n+1,q)&q&n&2\\
H(2n,q^2)&q^2&n&3/2\\
H(2n-1,q^2)&q^2&n&1/2
\end{array}
\]
Let $v_0=0$ and, for $j\ge1$, let
\[
 v_j=\frac{(Q^j-1)(Q^{j+e-1}+1)}{Q-1}.
\]
The collinearity graph has $v_m$ vertices and degree $d=Qv_{m-1}$.
For distinct vertices $x,y$, its common-neighbor counts are
\[
 |N(x)\cap N(y)|=
 \begin{cases}
 Q-1+Q^2v_{m-2},&x\sim y,\\
 v_{m-1},&x\not\sim y.
 \end{cases}
\]
To see the first formula, the line $xy$ contributes $Q-1$ points
other than its endpoints, and each point of its perpendicular
quotient contributes $Q^2$ further points. For the second formula,
the span of noncollinear $x,y$ is a nondegenerate hyperbolic plane;
its perpendicular space is a polar space of rank $m-1$ with the
same parameter $e$. The degree formula follows similarly by taking
the perpendicular quotient at a point, with $Q$ lifts per quotient
point. Thus these are counts of common neighbors, as required by
Definition~\ref{def:sparse}, and not counts of lines.

For the three quadratic rows, respectively, let
$a=2n-2,2n-1,2n-1$. Then $d\sim q^a$, $v_m/d\sim q$, and,
for all sufficiently large $q$,
\[
 \max_{x\ne y}|N(x)\cap N(y)|\le8q^{a-1}
       \le d^{1-1/(2a)}.
\]
These give exactly the three values of $\delta$ in the table.
For the two Hermitian rows let $a=4n-3,4n-5$, respectively.
Here $d\sim q^a$, $v_m/d\sim q^2$, and
\[
 \max_{x\ne y}|N(x)\cap N(y)|\le8q^{\max\{a-2,2\}}
       \le d^{1-1/(4n-3)}
\]
for all sufficiently large $q$. In the second row with $n=2$,
the last exponent of $q$ is $3(1-1/5)>2$; in all other Hermitian
cases it exceeds $a-2$ and $2$ as well. The graph case of
Theorem~\ref{main}, namely $(1+o(1))v_m(\log d)/d$, now gives
the five stated leading constants.

For a generalized quadrangle $GQ(q,t)$, the corresponding counts are
\[
 v=(q+1)(qt+1),\qquad d=q(t+1),\qquad
 |N(x)\cap N(y)|=
 \begin{cases}q-1,&x\sim y,\\t+1,&x\not\sim y.\end{cases}
\]
The first common-neighbor count is given by the line through $x,y$;
the second follows from the unique point collinear with $y$ on each
of the $t+1$ lines through $x$.
Higman's inequalities give $\sqrt q\le t\le q^2$ for thick
generalized quadrangles~\cite{Higman1,PT}. Hence, for large $q$,
\[
 \max\{q-1,t+1\}\le d^{3/4},\qquad
 \frac vd\sim q,\qquad \log d\le(3+o(1))\log q.
\]
Indeed, $(q-1)^4\le q^3(t+1)^3$ follows from $t\ge\sqrt q$,
and $(t+1)^4\le q^3(t+1)^3$ follows from $t+1\le q^3$.
These inequalities verify $\delta=1/4$, and Theorem~\ref{main}
gives $\gamma(GQ(q,t))\le(3+o(1))q\log q$, as required.
\end{proof}

\section{Construction of random Steiner systems}\label{sec:Steiner}
 
We construct Steiner triple systems in which every maximal independent
set has size at least $(1-o(1))\sqrt{n\log n}$. Before doing so, we give a construction showing that the common-link condition in the definition of
$\delta$-sparseness is necessary for Theorem~\ref{main}.

\subsection{Blowups}\label{sec:blowups}

We show that without condition (2) of $\delta$-sparseness (the common link condition), condition (1) can hold and yet the conclusion of the main theorem does not hold. 

\begin{proposition}\label{prop:blowup}
Fix an integer $r\ge2$ and a real number
\[
 0<\delta\le \frac{1}{2r-1}.
\]
There is an infinite family of  $(r+1)$-uniform set systems
$\mathcal S_q$ with $n_q$ atoms and degree $d_q$ such that the following hold, with degrees and links taken in
$\mathcal S_q$:

\begin{enumerate}
\item for every $A\subseteq X(\mathcal S_q)$ with
      $2\le |A|\le r$,
      \[
       d_{r+1}(A)
       \le d_q^{(r+1-|A|-\delta)/r};
      \]
\item there are distinct atoms $u,v$ for which
      \[
       d_{r+1}(u,v)=d_q,
      \]
      so the common-link condition
      $d_{r+1}(u,v)\le d_q^{1-\delta/r}$ fails;
\item every maximal independent set in $\mathcal S_q$ has at least $rq$
      atoms, whereas
      \[
       n_q\left(\frac{\log d_q}{r d_q}\right)^{1/r}
       =O_r\bigl(q^{1/r}(\log q)^{1/r}\bigr).
      \]
\end{enumerate}
\end{proposition}

\medskip

\begin{proof}
We give the construction as a blow-up of a fixed Steiner
$(r+1)$-tuple system. Let $\mathcal T$ be a Steiner $(r+1)$-tuple system on a
fixed point set $P$ of size $m$. Thus every pair of points of $P$ lies
in exactly one block of $\mathcal T$. Every point of $\mathcal T$ has the same degree
\[
 d_0=\frac{m-1}{r}.
\]
The trivial choice $P=[r+1]$ and $\mathcal T=\{P\}$ is already sufficient.

Let $q\ge2$ be an integer, and write $\Zq=\mathbb Z/q\mathbb Z$. Replace
each point $p\in P$ by the set
\[
 \{p\}\times\Zq\times\Zq.
\]
Hence
\[
 X(\mathcal S_q)=P\times\Zq\times\Zq,
 \qquad n_q=mq^2.
\]
For every block $B\in \mathcal T$, every choice
$(x_p)_{p\in B}\in\Zq^B$ satisfying
\[
 \sum_{p\in B}x_p=0,
\]
and every choice $(a_p)_{p\in B}\in\Zq^B$, put
\[
 \{(p,x_p,a_p):p\in B\}
\]
in $\mathcal S_q$. This defines a simple $(r+1)$-uniform set system.

\medskip
\noindent\textbf{Degree.}
Fix $(p,x,a)\in X(\mathcal S_q)$. There are $d_0$ choices for the block
$B\in \mathcal T$ containing $p$. Once $B$ is fixed, the remaining $r$
first coordinates have $q^{r-1}$ choices satisfying their one linear
sum condition, and the remaining $r$ second coordinates are arbitrary.
Therefore
\[
 d_q=d_0q^{r-1}q^r=d_0q^{2r-1}.
\]
In particular, $\mathcal S_q$ is regular.

\medskip
\noindent\textbf{Ordinary codegrees.}
Let $A\subseteq X(\mathcal S_q)$ have size $j$, where $2\le j\le r$. If two
atoms of $A$ project to the same point of $P$, then no block contains
$A$. Otherwise the projections of the atoms of $A$ are distinct.
Since $\mathcal T$ is a Steiner system, at most one block of $\mathcal T$ contains
all these projected points. For that block, the remaining
$r+1-j$ first coordinates have
$q^{r-j}$ choices satisfying the sum condition, while the remaining
$r+1-j$ second coordinates are arbitrary. Hence
\[
 d_{r+1}(A)\le q^{r-j}q^{r+1-j}
 =q^{2r-2j+1}.
\]
Now
\[
 \frac{(2r-1)(r+1-j-\delta)}{r}-(2r-2j+1)
 =\frac{j-1-(2r-1)\delta}{r}\ge0,
\]
because $j\ge2$ and $\delta\le1/(2r-1)$. Since $d_0\ge1$, it follows
that
\[
 d_{r+1}(A)
 \le q^{2r-2j+1}
 \le (d_0q^{2r-1})^{(r+1-j-\delta)/r}
 =d_q^{(r+1-j-\delta)/r}.
\]
Thus the ordinary codegree condition holds.

\medskip
\noindent\textbf{Failure of the common-link condition.}
For $p\in P$ and $x\in\Zq$, define the fibre
\[
 F_{p,x}=\{(p,x,a):a\in\Zq\}.
\]
Choose distinct atoms
\[
 u=(p,x,a),\qquad v=(p,x,a')
\]
in the same fibre. The block condition does not involve the last
coordinate. Therefore replacing $u$ by $v$ in any block containing $u$
produces a block with exactly the same remaining $r$ atoms. Hence
\[
 N_{r+1}(u)=N_{r+1}(v)
 \qquad\text{and consequently}\qquad
 d_{r+1}(u,v)=d_q.
\]
Since $d_q>1$ and $\delta>0$,
\[
 d_q>d_q^{1-\delta/r},
\]
so the common-link condition fails as strongly as possible.

\medskip
\noindent\textbf{Maximal independent sets.}
Let $I$ be a maximal independent set in $\mathcal S_q$. We first show that
$I$ is a union of fibres. Suppose that $u\in I\cap F_{p,x}$ and that
$v\in F_{p,x}\setminus I$. If $I\cup\{v\}$ contained a block, then
replacing $v$ in that block by $u$ would produce a block contained in
$I$, a contradiction. Thus $I\cup\{v\}$ is independent, contrary to
maximality. Hence
\[
 I\cap F_{p,x}\ne\varnothing
 \quad\Longrightarrow\quad
 F_{p,x}\subseteq I.
\]
Every block of $\mathcal S_q$ meets exactly $r+1$ distinct fibres. If $I$ were
the union of fewer than $r$ fibres, then adjoining any vertex outside
$I$ would produce a set meeting at most $r$ fibres and hence still
containing no block. This would contradict maximality. Therefore $I$
contains at least $r$ fibres, and
\[
 |I|\ge rq.
\]

Finally,
\[
 n_q\left(\frac{\log d_q}{r d_q}\right)^{1/r}
 =mq^2
   \left(
    \frac{\log(d_0q^{2r-1})}{r d_0q^{2r-1}}
   \right)^{1/r} 
 =O_{r,\mathcal T}\bigl(q^{1/r}(\log q)^{1/r}\bigr)
 =o(q).
\]
Since every maximal independent set has size at least $rq$, the ratio
in the proposition tends to infinity.
\end{proof}

\medskip

When $r=1$, take the one-block Steiner $2$-tuple system and the same
construction. The resulting graph has atom set
$[2]\times\Zq\times\Zq$ and blocks
\[
 \{(1,x,a),(2,-x,b)\},
 \qquad x,a,b\in\Zq.
\]
It is the disjoint union of $q$ copies of $K_{q,q}$. Thus it is
$q$-regular, the ordinary codegree condition is vacuous, and two
atoms on the same side of one component have $q$ common neighbours.
Every maximal independent set has size $q^2$, while the scale in the
graph case is
\[
 2q^2\frac{\log q}{q}=2q\log q=o(q^2).
\]
So the common-neighbour condition is necessary in the graph case as
well.

\bigskip

We remark that we do not have a construction that works for fixed $\delta$ for all $r \geq 2$.

\subsection{The antichain condition}\label{subsec:antichain}

We give a counterexample to the mixed-system main theorem
(Theorem~\ref{thm:main}) when only the antichain hypothesis
is removed. The examples have upper rank three, are regular in
each block size, and satisfy both sparsity inequalities with the
fixed parameter $\delta=1/4$. Nevertheless, every maximal independent
set is larger than the size predicted by the theorem by a factor
that tends to infinity.

For a set system $\cS$ on $X$, let
\[
 N_i(A)=\{S\setminus A:A\subseteq S\in\cS,\ |S|=i\},
 \qquad d_i(A)=|N_i(A)|.
\]
For distinct atoms $u,v$, we write
\[
 d_i(u,v)=|N_i(\{u\})\cap N_i(\{v\})|.
\]
Thus $d_3(\{u,v\})$ counts triples containing both atoms, whereas
$d_3(u,v)$ counts pairs that complete each of the two atoms to a
triple. For upper rank three, the two sparsity conditions require that,
for all distinct $u,v$,
\begin{equation}\label{eq:sparsity}
 d_3(\{u,v\})\le D^{1-\delta},\qquad
 d_2(u,v)\le D^{1-\delta},\qquad
 d_3(u,v)\le D^{2-\delta},
\end{equation}
where $D=\max\{d_2,\sqrt{d_3}\}$ for a regular system with
degree vector $\vb{d}=(d_2,d_3)$. We define $\rho=\rho(\vb{d})>0$ by
\begin{equation}\label{eq:rho}
 d_2\rho+d_3\rho^2=\log D.
\end{equation}
All logarithms are natural.

\begin{proposition}\label{prop:counterexample}
Let $q=2^h\ge64$ and $s=q/2-1$. There is a set system $\cS_q$
with distinct blocks of sizes two and three such that
\begin{enumerate}
\item $|X|=n=q^2(2q+1)$, every atom belongs to exactly $d_2=q$
two-element blocks and $d_3=q^2s$ three-element blocks, and
$D=q\sqrt{s}$;
\item the inequalities \eqref{eq:sparsity} hold with $\delta=1/4$;
\item $\cS_q$ is not an antichain, and every maximal independent
set $I$ satisfies
\[
 |I|\ge\frac{q^2}{2},\qquad
 \frac{|I|}{n\rho}\ge
 \frac{q\sqrt{s}}{2(2q+1)\sqrt{\log(q\sqrt{s})}}
 \ge\frac{\sqrt q}{6\sqrt{6\log q}}.
\]
\end{enumerate}
In particular, the conclusion of Theorem~\ref{thm:main} is false without the
antichain condition, even with $r=3$ and $\delta=1/4$ fixed.
\end{proposition}

\medskip
\noindent\emph{The construction.}
Let $\mathcal P=\F_q^2$, and let $\cL$ consist of its nonvertical
affine lines:
\[
 \ell_{a,b}=\{(x,ax+b):x\in\F_q\},\qquad
 \cL_a=\{\ell_{a,b}:b\in\F_q\},\qquad
 \cL=\bigcup_{a\in\F_q}\cL_a.
\]
There are $q^2$ such lines. Each contains $q$ points, and each
point lies on exactly one line of each slope, hence on $q$ lines
of $\cL$. Two distinct lines of different slopes meet in exactly
one point. Two distinct points lie on a unique line of $\cL$
if their first coordinates differ, and on no line of $\cL$
otherwise. These facts follow directly from the equations
$y=ax+b$.

Since $\F_q$ has characteristic two, the pairs $\{a,a+1\}$ partition
the slopes into $q/2$ pairs. Let $A\subseteq\F_q$ contain one
representative from each such pair.
Let $[q]=\{1,\ldots,q\}$ be partitioned into two sets $A_0,A_1$
of size $q/2$. For $i,j\in[q]$, the notation $i\sim j$ means
that $i\ne j$ and $i,j$ belong to the same $A_e$. Thus each
$i$ has exactly $s$ indices $j$ with $i\sim j$.

The atoms are
\[
 X=\cL\ \sqcup\ \{L(P,i):P\in\mathcal P,\ i\in[q]\}
          \ \sqcup\ \{R(P,i):P\in\mathcal P,\ i\in[q]\}.
\]
Here the lines are themselves atoms, and $L(P,i),R(P,i)$ are
distinct formal symbols. In particular, the points $P$ are indices
for these symbols, not additional atoms.

The two-element blocks are
\[
 \cS_2=\bigl\{\{\ell,\ell'\}:a\in A,\ \ell\in\cL_a,\ \ell'\in\cL_{a+1}\bigr\}
 \cup\bigl\{\{L(P,i),R(P,j)\}:P\in\mathcal P,\ i,j\in[q]\bigr\}.
\]
For $\ell\in\cL$ and $i\in[q]$, let
\[
 L_i(\ell)=\{L(P,i):P\in\ell\},\qquad
 R_i(\ell)=\{R(P,i):P\in\ell\}.
\]
The three-element blocks are
\begin{align*}
 \cS_3^{(1)}&=\bigl\{\{\ell,L(P,i),R(P,j)\}:
                 \ell\in\cL,\ P\in\ell,\ i\sim j\bigr\},\\
 \cS_3^{(2)}&=\bigcup_{\substack{\ell\in\cL\\i\in[q]}}
       \left(\binom{L_i(\ell)}3\cup\binom{R_i(\ell)}3\right),
 \\
 \cS_3&=\cS_3^{(1)}\cup\cS_3^{(2)}.
\end{align*}
Finally, let $\cS_q=\cS_2\cup\cS_3$. The two types of triples
are disjoint, since the first contains a line-atom and the second
does not. A triple in $\cS_3^{(2)}$ determines its label, its
choice of $L$ or $R$, and the unique line through its three
distinct point indices. Thus no triple is counted twice.

\begin{proof}[Proof of Proposition~\textup{\ref{prop:counterexample}}]
\emph{Proof of \textup{(1)}.}
The number of atoms is $q^2+2q^3=q^2(2q+1)$.
The graph whose edges are $\cS_2$ is the disjoint union of
$q/2$ copies of $K_{q,q}$ on paired parallel classes, and $q^2$
copies of $K_{q,q}$ on the two sets
$\{L(P,i):i\in[q]\}$ and $\{R(P,j):j\in[q]\}$.
Consequently every atom has degree $d_2=q$.

For a line-atom $\ell$, a triple containing it is determined by
a point $P\in\ell$, an index $i\in[q]$, and one of the $s$
indices $j\sim i$. Its degree in $\cS_3$ is therefore
\[
 d_3(\{\ell\})=q\cdot q\cdot s=q^2s.
\]
For an atom $L(P,i)$, there are $qs$ triples in $\cS_3^{(1)}$:
one may use any of the $q$ nonvertical lines through $P$ and any
$j\sim i$. There are $q\binom{q-1}{2}$ triples in
$\cS_3^{(2)}$ containing this atom. Since $(q-2)/2=s$, its total
degree is
\[
 qs+q\binom{q-1}{2}=qs+q(q-1)s=q^2s.
\]
The same count applies to $R(P,i)$. Hence $d_3=q^2s$ and,
as $s\ge1$, $D=\max\{q,q\sqrt{s}\}=q\sqrt{s}$.

\medskip
\noindent\emph{Proof of \textup{(2)}: pairs contained in triples.}
The only nonzero possibilities for $d_3(\{u,v\})$ are as follows;
the corresponding statements with $L$ and $R$ interchanged also hold:
\begin{enumerate}
\item If $u=\ell$ and $v=L(P,i)$ with $P\in\ell$, the other atom
is $R(P,j)$ for one of the $s$ indices $j\sim i$.
\item If $u=L(P,i)$ and $v=R(P,j)$ with $i\sim j$, the other
atom is one of the $q$ lines through $P$.
\item If $u=L(P,i)$ and $v=L(Q,i)$, where $P\ne Q$ lie on a
line of $\cL$, the other atom is $L(U,i)$ for one of the $q-2$
points $U$ on that line other than $P,Q$.
\end{enumerate}
Every other pair is contained in no triple. Thus, for distinct $u,v$,
\begin{equation}\label{eq:pair-degree}
 d_3(\{u,v\})\le q.
\end{equation}

\medskip
\noindent\emph{Proof of \textup{(2)}: intersections of links.}
For distinct $u,v$, the disjoint union of copies of $K_{q,q}$ described above gives
\begin{equation}
 d_2(u,v)\le q.
\end{equation}
We next enumerate $N_3(\{u\})\cap N_3(\{v\})$.
\begin{enumerate}
\item For distinct line-atoms $\ell,\ell'$, this intersection is
empty if the lines are parallel. Otherwise they meet at a point
$P$, and their common link consists precisely of
\[
 \{\{L(P,i),R(P,j)\}:i\sim j\},
\]
which has $qs$ members.
\item For $u=L(P,i)$ and $v=L(P,k)$ with $i\ne k$, the common
link is empty when $i,k$ belong to different halves. If they
belong to the same half, its members, for $P\in\ell$ and
$j\in A_e\setminus\{i,k\}$, are
\[
 \{\ell,R(P,j)\}.
\]
where $i,k\in A_e$. There are $q(q/2-2)=q(s-1)$ such pairs.
\item For $u=L(P,i)$ and $v=L(Q,i)$ with $P\ne Q$, the common
link is empty if $P,Q$ lie on no nonvertical line. Otherwise,
if $\ell$ is their common nonvertical line, its members are
\[
 \{L(U,i),L(V,i)\},\qquad
 \{U,V\}\in\binom{\ell\setminus\{P,Q\}}2.
\]
There are $\binom{q-2}{2}$ such pairs.
\end{enumerate}
The second and third cases also apply with $R$ in place of $L$.
All other intersections are empty. Indeed, a pair in the link
of a line-atom has one $L$-atom and one $R$-atom; a pair in the
link of an $L$-atom has either one line-atom and one $R$-atom,
or two $L$-atoms with the same label. The analogous description
holds for an $R$-atom. For two $L$-atoms with different point
indices and different labels, neither kind of pair can coincide.
It follows that
\begin{equation}\label{eq:three-common}
 d_3(u,v)\le
 \max\left\{qs,q(s-1),\binom{q-2}{2}\right\}=qs.
\end{equation}

Since $q\ge64$ and $s=q/2-1\ge q/4$, we have
$s^3\ge q^3/64\ge q^2$. Therefore
\[
 q\le q^{3/4}s^{3/8}=D^{3/4}.
\]
Also $s\le q$ gives
\[
 \frac{qs}{D^{7/4}}=\frac{s^{1/8}}{q^{3/4}}
       \le q^{-5/8}\le1.
\]
Together with \eqref{eq:pair-degree}--\eqref{eq:three-common},
these are exactly the inequalities \eqref{eq:sparsity} for
$\delta=1/4$.

\medskip
\noindent\emph{Proof of \textup{(3)}.}
For every triple in $\cS_3^{(1)}$, we have
\[
 \{L(P,i),R(P,j)\}\subsetneq\{\ell,L(P,i),R(P,j)\},
 \qquad \{L(P,i),R(P,j)\}\in\cS_2.
\]
Thus $\cS_q$ is not an antichain. Moreover, a set avoids all
blocks of $\cS_q$ if and only if it avoids all blocks of
$\cS_2\cup\cS_3^{(2)}$: the triples in $\cS_3^{(1)}$ add no
restriction to independence. The maximal independent sets of
these two systems are consequently identical.

No triple of $\cS_3^{(2)}$ contains a line-atom. Hence the
restriction of a maximal independent set $I$ to each copy of
$K_{q,q}$ on $\cL_a\cup\cL_{a+1}$ must itself be maximal
independent in that graph. An independent set in $K_{q,q}$ lies
entirely in one of its two parts, and it is maximal precisely
when it is one whole part. Therefore, for every $a\in A$,
\[
 |I\cap(\cL_a\cup\cL_{a+1})|=q,
 \qquad |I\cap\cL|=\frac{q^2}{2},
\]
which proves $|I|\ge q^2/2$.

Equation~\eqref{eq:rho} becomes
\[
 q\rho+q^2s\rho^2=\log(q\sqrt{s}).
\]
Both terms on the left are nonnegative, so
\[
 \rho\le\frac{\sqrt{\log(q\sqrt{s})}}{q\sqrt{s}},\qquad
 n\rho\le\frac{q(2q+1)}{\sqrt{s}}\sqrt{\log(q\sqrt{s})}.
\]
It follows that
\[
 \frac{|I|}{n\rho}\ge
 \frac{q\sqrt{s}}{2(2q+1)\sqrt{\log(q\sqrt{s})}}
 \ge\frac{\sqrt q}{6\sqrt{6\log q}},
\]
where the last inequality uses $s\ge q/4$, $2q+1\le3q$, and
$\log(q\sqrt{s})\le\frac32\log q$.
This lower bound tends to infinity as $q\to\infty$ through powers
of two. For any fixed $K$ and $\epsilon>0$, all sufficiently
large members of the construction therefore have no maximal
independent set satisfying
\[
 |I|\le\bigl(1+K(\log D)^{-\epsilon}\bigr)n\rho.
\]
Thus even the upper bound in the asserted conclusion fails.
\end{proof}

\medskip
Deleting the triples in $\cS_3^{(1)}$ leaves all independent sets
unchanged, but changes the three-element degrees: a line-atom then
has degree zero, while an auxiliary atom has degree
$q\binom{q-1}{2}$. Thus passing to the inclusion-minimal blocks
does not preserve the degree vector used to define $D$ and $\rho$.
The construction makes this obstruction explicit while satisfying
all the remaining hypotheses of the main theorem.

\subsection{Random systems}

To make a Steiner triple system in which the smallest maximal independent sets have size asymptotic to $\sqrt{n\log n}$, we consider the construction in~\cite{EustisV}, where we fill
the lines of a projective plane with independently relabelled Steiner triple systems.  The following more precise statement is useful.

\begin{theorem}
\label{thm:sts-plane-lower-bound}
Let $q\to\infty$ through the sequence $q=2^{2a+1}$, where $a$ is a
positive integer, and put $n=q^2+q+1$.  On every line of
$\PG(2,q)$, place a uniformly random relabelling of an arbitrary
fixed Steiner triple system on $q+1$ points, independently on different
lines.  Their union $\mathcal S$ is a Steiner triple system on $n$ points,
and, with probability tending to one,
\[
 \gamma(\mathcal S)>
 \left\lfloor\sqrt{n(\log n-3\log\log n)}\right\rfloor.
\]
In particular, for infinitely many $n$ there is a Steiner triple system
on $n$ points in which every maximal independent set has size at least
$(1-o(1))\sqrt{n\log n}$.
\end{theorem}

We first describe the construction and prove the two probabilistic
estimates needed in its analysis.  In particular, no independence
assumption will be made about the events that individual points remain
uncovered.

\subsection{The projective plane and its local triple systems}

Take the points of the plane to be the one-dimensional subspaces of
$\mathbb F_q^3$, and its lines to be the two-dimensional subspaces.
There are
\[
 \frac{q^3-1}{q-1}=q^2+q+1=n
\]
points, and each line contains $(q^2-1)/(q-1)=q+1$ points.
Each point lies on $q+1$ lines: the lines through a fixed point
correspond to the one-dimensional subspaces of the two-dimensional
quotient by that point.  Counting incidences consequently shows that
there are also $n$ lines.  Two distinct points span a unique line, and
two distinct lines meet in exactly one point.  Write $X$ for the point
set and $\mathcal L$ for the line set.

For completeness, a Steiner triple system on $q+1$ points can be
constructed explicitly.  Since $q=2^{2a+1}$, we may write $q+1=3u$,
where $u$ is odd. On $\mathbb Z_u\times\mathbb Z_3$, we take, for each
$x\in\mathbb Z_u$, the triple
\[
 \{(x,0),(x,1),(x,2)\},
\]
and, for distinct $x,y\in\mathbb Z_u$ and each $i\in\mathbb Z_3$,
the triple
\[
 \left\{(x,i),(y,i),
 \left(\frac{x+y}{2},i+1\right)\right\}
 ,
\]
where $\{x,y\}$ is unordered, division by $2$ is in $\mathbb Z_u$,
and the second coordinate is read modulo $3$.
A pair in one level determines its unique displayed triple.
A pair in different levels can be written uniquely as
$(x,i),(z,i+1)$.  If $x=z$, its unique triple is vertical;
otherwise its unique triple is obtained by taking $y=2z-x$.
This proves that every pair lies in exactly one triple.

For each $\ell\in\mathcal L$, fix any Steiner triple system
$\mathcal T_\ell$ on $q+1$ abstract points.  Choose a uniformly random
bijection from its points to $\ell$, independently for different
$\ell$, and transfer its triples to $\ell$ using this bijection.
Let $\mathcal S$ be the union of all the transferred triples.
Every pair of points of $X$ lies on a unique projective line and in a
unique triple of the system placed on that line.  Thus $\mathcal S$
is a Steiner triple system.  Although two lines share a point, their
random bijections, and hence their local random triple systems, are
independent.

\subsection{An elementary avoidance inequality}

We give a short proof of the probabilistic inequality used to control
large intersections with a line.

\begin{lemma}
\label{lem:sts-plane-avoidance}
Let $R$ be a random subset of a finite set, obtained by selecting its
elements independently.  Let $E_1,\ldots,E_m$ be nonempty fixed subsets,
and put $A_i=\{E_i\subseteq R\}$.  Define
\[
 \mu=\sum_{i=1}^m\mathbb P(A_i),
 \qquad
 \Delta=\sum_{\substack{i<j\\E_i\cap E_j\ne\varnothing}}
              \mathbb P(A_i\cap A_j).
\]
Then
\[
 \mathbb P\left(\bigcap_{i=1}^m\overline{A_i}\right)
 \le \exp(-\mu+\Delta).
\]
\end{lemma}

\begin{proof}
By continuity, it suffices to consider selection probabilities strictly
between zero and one.
An event is increasing if changing a coordinate from unselected to
selected cannot make it fail.  We first verify that two increasing
events $A,B$ in independent Bernoulli coordinates satisfy
\[
 \mathbb P(A\cap B)\ge\mathbb P(A)\mathbb P(B).
\]
Induct on the number of coordinates.  Condition on the last coordinate,
which equals one with probability $p$, and write $a_0,a_1$ and $b_0,b_1$
for the corresponding conditional probabilities of $A$ and $B$.
The induction hypothesis applies within each conditional space.
Moreover, $a_1\ge a_0$ and $b_1\ge b_0$.  Subtracting the product of
the two unconditional probabilities from the resulting lower bound
on $\mathbb P(A\cap B)$ leaves
\[
 p(1-p)(a_1-a_0)(b_1-b_0)\ge0.
\]
This proves the assertion.  Taking the complement of one event shows
that an increasing event and a decreasing event satisfy the reverse
inequality.

For each $i$, let $F_i$ be the event that no earlier $A_j$ whose
coordinate set is disjoint from $E_i$ occurs, and let $G_i$ be the
event that no earlier $A_j$ whose coordinate set meets $E_i$ occurs.
Thus $F_i\cap G_i=\bigcap_{j<i}\overline{A_j}$.
If this event has probability zero for some $i$, the desired conclusion
is immediate.  Otherwise put
\[
 \theta_i=\mathbb P(A_i\mid F_i\cap G_i).
\]
Since $F_i$ depends only on coordinates outside $E_i$, it is independent
of $A_i$.  Also, $F_i$ is decreasing and $A_i\cap A_j$ is increasing,
so the inequality just proved gives
\[
 \mathbb P(A_i\cap A_j\mid F_i)\le\mathbb P(A_i\cap A_j).
\]
Consequently, the union bound gives
\begin{align*}
 \theta_i
 &\ge \mathbb P(A_i\cap G_i\mid F_i)\\
 &\ge \mathbb P(A_i\mid F_i)
       -\sum_{\substack{j<i\\E_j\cap E_i\ne\varnothing}}
          \mathbb P(A_i\cap A_j\mid F_i)\\
 &\ge \mathbb P(A_i)
       -\sum_{\substack{j<i\\E_j\cap E_i\ne\varnothing}}
          \mathbb P(A_i\cap A_j).
\end{align*}
The first inequality follows by dividing
$\mathbb P(A_i\cap G_i\mid F_i)$ by $\mathbb P(G_i\mid F_i)\le1$.
Multiplying the successive conditional avoidance probabilities now gives
\[
 \mathbb P\left(\bigcap_{i=1}^m\overline{A_i}\right)
 =\prod_{i=1}^m(1-\theta_i)
 \le \exp\left(-\sum_{i=1}^m\theta_i\right)
 \le \exp(-\mu+\Delta).
\]
\end{proof}

\begin{lemma}
\label{lem:sts-plane-local}
Put $h=\lceil q^{2/5}\rceil$.
There is an absolute constant $c>0$ such that, for all sufficiently
large $q$, a uniformly random subset of cardinality at least $h$ in
any Steiner triple system on $q+1$ points is independent with
probability at most $\exp(-c q^{1/5})$.
\end{lemma}

\begin{proof}
First select each point independently with probability
\[
 p=\frac{h}{2(q+1)},
\]
and let $R$ be the selected set.  Apply
Lemma~\ref{lem:sts-plane-avoidance} to the events that the individual
triples are contained in $R$.
A Steiner triple system on $q+1$ points has $q(q+1)/6$ triples,
and each point lies in $q/2$ triples.  Therefore
\[
 \mu=\frac{q(q+1)}6p^3
     =\frac{qh^3}{48(q+1)^2}
     =\Theta(q^{1/5}).
\]
Two distinct triples meet in at most one point, so the number of
unordered intersecting pairs of triples is
$(q+1)\binom{q/2}{2}$.  Their unions have five points.  Thus
\[
 \Delta=(q+1)\binom{q/2}{2}p^5,
 \qquad
 \frac{\Delta}{\mu}
 =\frac{3(q-2)h^2}{16(q+1)^2}
 =O(q^{-1/5}).
\]
It follows that
\[
 \mathbb P(R\text{ is independent})\le\exp(-c_0q^{1/5})
\]
for an absolute constant $c_0>0$ and all sufficiently large $q$.

Let $a_j$ be the probability that a uniform $j$-subset is independent.
The sequence $a_j$ is nonincreasing: deleting a uniformly random point
from a uniform $(j+1)$-subset gives a uniform $j$-subset, and preserves
independence.  Since $\mathbb E(|R|)=h/2$, Markov's inequality gives
$\mathbb P(|R|\le h)\ge1/2$.  Conditional on its size, $R$ is a uniform
subset.  Hence
\[
 \mathbb P(R\text{ is independent})
 =\sum_j\mathbb P(|R|=j)a_j
 \ge \frac12a_h.
\]
We conclude that $a_h\le2\exp(-c_0q^{1/5})$.
Decreasing the constant in the exponent absorbs the factor $2$.
Monotonicity gives the same bound for every cardinality at least $h$.
All estimates are uniform over the choice of the Steiner triple system.
\end{proof}

\subsection{A covering inequality}

The next lemma handles the dependence between the events that different
points are covered.  The sets from which the random subsets are chosen
are allowed to overlap arbitrarily.

\begin{lemma}
\label{lem:sts-plane-cover}
Let $B_1,\ldots,B_m$ be independent random sets, where $B_j$ is a
uniform $b_j$-subset of a fixed finite set $A_j$.
For a fixed target set $W$ and each $v\in W$, let
\[
 p_v=\mathbb P\left(v\notin\bigcup_{j=1}^mB_j\right)
 .
\]
Then
\[
 \mathbb P\left(W\subseteq\bigcup_{j=1}^mB_j\right)
 \le\prod_{v\in W}(1-p_v)
 \le\exp\left(-\sum_{v\in W}p_v\right).
\]
\end{lemma}

\begin{proof}
We first describe a coupling for one uniform $b$-subset $B$ of an
$m$-element set $A$, where $v\in A$ and $b<m$.
If $v\notin B$, leave $B$ unchanged.  If $v\in B$, remove $v$ and
replace it by a uniformly chosen point of $A\setminus B$.
Call the resulting set $B'$.  For any fixed $b$-subset
$T\subseteq A\setminus\{v\}$,
\[
 \mathbb P(B'=T)
 =\frac1{\binom mb}
  +\frac{b}{(m-b)\binom mb}
 =\frac1{\binom{m-1}b}.
\]
The first term comes from $B=T$.  The second comes from the $b$
possible sets $(T\setminus\{w\})\cup\{v\}$, with $w\in T$, followed
by choosing $w$ as the replacement.  Thus $B'$ has exactly the law
of $B$ conditioned to avoid $v$.  Moreover, no point other than $v$
that belonged to $B$ is removed.

Now fix $v\in W$.  If $p_v=1$, complete coverage is impossible.
If $p_v=0$, the factor $1-p_v$ is one and $v$ can be omitted from
the target set.  Otherwise apply the preceding replacement independently
to all $B_j$ whose ground sets contain $v$.  The result has the exact
joint law of the $B_j$ conditioned on $v$ being missed by all of them:
this conditioning factors over the independent random sets.
No other target point loses coverage in the coupling.

Let $H$ be the event that every point of $W\setminus\{v\}$ is covered.
The coupling proves that
\[
 \mathbb P(H\mid v\text{ is missed})\ge\mathbb P(H).
\]
Therefore
\[
 \mathbb P(W\text{ is covered})
 =\mathbb P(H)-p_v\mathbb P(H\mid v\text{ is missed}) 
 \le(1-p_v)\mathbb P(H).
\]
Induction on $|W|$ proves the product bound.  The exponential bound
follows from $1-x\le e^{-x}$.
\end{proof}

\subsection{Bounding the probability for a fixed candidate set}

Fix a deterministic set $S\subseteq X$, with
\[
 s=|S|\le\sqrt{n\log n}.
\]
For each line $\ell$, put $t_\ell=|S\cap\ell|$, and call $\ell$
heavy if $t_\ell\ge h=\lceil q^{2/5}\rceil$.
The heavy lines depend only on $S$ and the fixed projective plane.

\paragraph{Many heavy lines.}
Suppose that there are at least $q^{9/10}$ heavy lines.
For $S$ to be independent in $\mathcal S$, its intersection with
every such line must be independent in that line's random local
system.  Under a uniform relabelling, a fixed $t_\ell$-subset of the
line corresponds to a uniform $t_\ell$-subset of the fixed local
Steiner triple system.  Lemma~\ref{lem:sts-plane-local} therefore
applies.  Since different lines are filled independently,
\begin{equation}
 \label{eq:sts-plane-heavy}
 \mathbb P(S\text{ is maximal independent})
 \le\mathbb P(S\text{ is independent})
 \le\exp(-c q^{11/10}).
\end{equation}

\paragraph{Few heavy lines.}
Suppose instead that there are fewer than $q^{9/10}$ heavy lines.
Remove $S$ and every point on a heavy line, and write
\[
 W=X\setminus\left(S\cup\bigcup_{\ell\text{ heavy}}\ell\right),
 \qquad N=|W|.
\]
Because a line has $q+1$ points,
\begin{equation}
 \label{eq:sts-plane-size-W}
 N\ge n-s-(q+1)q^{9/10}
   =n\bigl(1-O(q^{-1/10})\bigr).
\end{equation}
The estimate is uniform over all the candidate sets under consideration.
Every line through a point of $W$ is light.

For a light line $\ell$, let $C_\ell$ be the set of points of
$\ell\setminus S$ that complete some pair of $S\cap\ell$ to a triple
of the local system.  Put
\[
 b_\ell=\binom{t_\ell}{2},
 \qquad m_\ell=q+1-t_\ell.
\]
Each pair has exactly one completion, so $|C_\ell|\le b_\ell$.
Some different pairs may have the same completion, and some may have
their completion in $S$; neither possibility causes a problem.
Since $\ell$ is light,
\begin{equation}
 \label{eq:sts-plane-small-fraction}
 \frac{b_\ell}{m_\ell}=O(q^{-1/5}),
\end{equation}
uniformly in $\ell$ and $S$.  In particular, $b_\ell<m_\ell$ for
all sufficiently large $q$.

The law of $C_\ell$ is invariant under every permutation of
$\ell\setminus S$: such a permutation, extended to fix $S\cap\ell$,
preserves the law of the uniformly relabelled local system.
Consequently, conditional on $|C_\ell|=j$, the set $C_\ell$ is a
uniform $j$-subset of $\ell\setminus S$.
Enlarge it to a set $B_\ell$ of size exactly $b_\ell$ by choosing
$b_\ell-j$ further points uniformly from its complement.
Then $B_\ell$ is a uniform $b_\ell$-subset of $\ell\setminus S$.
Indeed, for any fixed $b_\ell$-subset $T$ of that set,
\[
 \mathbb P(B_\ell=T\mid |C_\ell|=j)
 =\frac{\binom{b_\ell}{j}}
 {\binom{m_\ell}{j}\binom{m_\ell-j}{b_\ell-j}}
 =\frac1{\binom{m_\ell}{b_\ell}}.
\]
Use independent auxiliary choices on different lines.  The sets
$B_\ell$ are then independent, and $C_\ell\subseteq B_\ell$ for
every light line.  No conditioning on the independence of $S$ has
been used here.

If $S$ is maximal independent, every point $v\notin S$ completes
some pair of $S$ to a triple.  For $v\in W$, that triple lies on a
light line through $v$.  Thus maximality implies
\[
 W\subseteq\bigcup_{\ell\text{ light}}C_\ell
 \subseteq\bigcup_{\ell\text{ light}}B_\ell.
\]
For $v\in W$, let $p_v$ be the probability that $v$ is missed by
all these enlarged sets.  Their independence gives
\[
 p_v=\prod_{\ell\ni v}
       \left(1-\frac{b_\ell}{m_\ell}\right),
\]
where all the lines in the product are light.
By Lemma~\ref{lem:sts-plane-cover},
\begin{equation}
 \label{eq:sts-plane-cover-probability}
 \mathbb P(S\text{ is maximal independent})
 \le\exp(-\mu_S),
 \qquad \mu_S=\sum_{v\in W}p_v.
\end{equation}

It remains to obtain a lower bound on $\mu_S$ with the correct leading
constant.  By~\eqref{eq:sts-plane-small-fraction} and
$-\log(1-z)\le z/(1-z)$ for $0\le z<1$,
\begin{align*}
 \sum_{v\in W}-\log p_v
 &=\sum_{\ell\text{ light}}|W\cap\ell|
       \left[-\log\left(1-\frac{b_\ell}{m_\ell}\right)\right]\\
 &\le\bigl(1+O(q^{-1/5})\bigr)
       \sum_{\ell\text{ light}}|W\cap\ell|
             \frac{b_\ell}{m_\ell}\\
 &\le\bigl(1+O(q^{-1/5})\bigr)
       \sum_{\ell\text{ light}}b_\ell\\
 &\le\bigl(1+O(q^{-1/5})\bigr)\binom{s}{2}.
\end{align*}
The third line uses $|W\cap\ell|\le m_\ell$.
The last line uses the exact identity
\begin{equation*}
 \sum_{\ell\in\mathcal L}\binom{|S\cap\ell|}{2}
 =\binom{s}{2},
\end{equation*}
which holds because every pair of $S$ lies on a unique projective line.

The arithmetic--geometric mean inequality now gives
\begin{equation}
 \label{eq:sts-plane-mu-raw}
 \mu_S\ge
 N\exp\left(-\frac{(1+O(q^{-1/5}))\binom{s}{2}}{N}\right).
\end{equation}
To make the uniformity of the error explicit, use
\eqref{eq:sts-plane-size-W} and $s^2/n\le\log n$.
Replacing $N$ by $n$ in the denominator of the exponent changes the
resulting upper bound on that exponent by at most
$O(q^{-1/10}\log n)=o(1)$.  The factor $1+O(q^{-1/5})$ has an even
smaller effect.  Also $\log(N/n)=O(q^{-1/10})$.
Consequently there is an absolute constant $C>0$ such that
\eqref{eq:sts-plane-mu-raw} implies, uniformly over all the candidate
sets in the few-heavy-lines case,
\begin{equation}
 \label{eq:sts-plane-mu-uniform}
 \mu_S\ge
 n\exp\left(-\frac{s^2}{2n}-Cq^{-1/10}\log q\right)
 =(1-o(1))n\exp\left(-\frac{s^2}{2n}\right).
\end{equation}
In particular, the $o(1)$ is independent of the choice of $S$.

\subsection{The final union bound}

\begin{proof}[Proof of Theorem~\textup{\ref{thm:sts-plane-lower-bound}}]
Put
\[
 s_*=
 \left\lfloor\sqrt{n(\log n-3\log\log n)}\right\rfloor,
 \qquad A_n=\sqrt n(\log n)^{3/2}.
\]
These quantities are well defined for all sufficiently large $n$.
If $s\le s_*$, then
\[
 n\exp\left(-\frac{s^2}{2n}\right)
 \ge n\exp\left(-\frac{\log n-3\log\log n}{2}\right)
 =A_n.
\]
Thus~\eqref{eq:sts-plane-cover-probability} and
\eqref{eq:sts-plane-mu-uniform} show that every fixed candidate of
size at most $s_*$ with fewer than $q^{9/10}$ heavy lines is maximal
independent with probability at most
\[
 \exp\bigl(-(1-o(1))A_n\bigr),
\]
with an error uniform over all such candidates.  Candidates with at
least $q^{9/10}$ heavy lines satisfy~\eqref{eq:sts-plane-heavy}.

The total number $M_n$ of candidates of size at most $s_*$ satisfies
\[
 M_n=\sum_{j=0}^{s_*}\binom nj
 \le(s_*+1)\left(\frac{en}{s_*}\right)^{s_*},
\]
since $s_*<n/2$ for large $n$.  Moreover,
\[
 s_*=(1+o(1))\sqrt{n\log n},
 \qquad
 \log\left(\frac{en}{s_*}\right)
 =\left(\frac12+o(1)\right)\log n.
\]
It follows that
\[
 \log M_n\le\left(\frac12+o(1)\right)A_n.
\]
Finally, $n=q^2+q+1$ implies
$A_n=O(q(\log q)^{3/2})=o(q^{11/10})$.
Let $\mathcal A$ be the event that some maximal independent set has
size at most $s_*$. The union bound gives
\[
 \mathbb P(\mathcal A)\le M_n\left\{e^{-cq^{11/10}}
             +e^{-(1-o(1))A_n}\right\}\longrightarrow0.
\]
This proves the asserted probability statement and hence the existence
of the required Steiner triple systems.
\end{proof}

\section{Random regular set systems}\label{sec:rrh}

We use $\gamma(\mathcal S)$ for the minimum size of a maximal
independent set, as in Section~\ref{sec:Pythag}.
The following result shows that the leading constant in the upper bound
for $\gamma(\mathcal S)$ in Theorem \ref{main} cannot be decreased, even after conditioning on
$\delta$-sparseness.  All logarithms in this section are natural.

\begin{theorem}
\label{thm:random-regular-delta-sparse}
Fix an integer $r\ge1$, real numbers $0<\delta\le1$ and
$0<\alpha<1$, and let $d$ be sufficiently large in terms of
$r,\delta,\alpha$.  Let $n\to\infty$ through integers satisfying
$(r+1)\mid nd$, put $X=[n]$, and let $\mathcal S$ be chosen uniformly from the simple
$d$-regular $(r+1)$-uniform antichains on $X$ which are
$\delta$-sparse.  Then asymptotically almost surely
\[
 \gamma(\mathcal S)
 >
 (1-\alpha)n
 \left(\frac{\log d}{rd}\right)^{1/r}.
\]
Equivalently,
\[
 \lim_{d\to\infty}\lim_{n\to\infty}
 \mathbb P\left(
 \gamma(\mathcal S)
 \le
 (1-\alpha)n
 \left(\frac{\log d}{rd}\right)^{1/r}
 \right)
 =0,
\]
where the inner probability is taken in the conditioned uniform model
described above.
\end{theorem}

\subsection{The configuration model}

Put
\[
 k=r+1,\qquad M=nd,\qquad Q=\frac{M}{k}.
\]
For every $v\in X$, let $P_v$ be a set of $d$ distinguishable
points, and put
\[
 P=\bigcup_{v\in X}P_v.
\]
Thus $|P|=M$.  Also take $Q$ disjoint right buckets
$C_1,\ldots,C_Q$, each containing $k$ distinguishable positions.
Choose a uniformly random perfect matching between the $M$ points of
$P$ and the $M$ right-bucket positions.  Each right bucket gives one
$(r+1)$-multiblock, whose atoms are the left buckets containing its
matched points.

Let $\mathcal D_\delta$ be the event that the resulting set system is
simple and $\delta$-sparse. Thus, on $\mathcal D_\delta$, for $2\le|A|\le r$,
\begin{equation}
 d_{r+1}(A)
 \le
 d^{(r+1-|A|-\delta)/r}
,
 \label{eq:random-regular-proper-codegree}
\end{equation}
and, for distinct $x,y$,
\begin{equation}
 d_{r+1}(x,y)
 \le
 d^{1-\delta/r}
.
 \label{eq:random-regular-common-link}
\end{equation}
Here degrees and links are taken in $\mathcal S$, and the first
family of inequalities is empty when $r=1$.

We first prove that $\mathcal D_\delta$ has probability bounded away
from zero for every fixed $r,d,\delta$ in the stated range.

\begin{lemma}
\label{lem:short-incidence-cycles}
For fixed $r$ and $d$, let $\mathsf C_\ell$ be the number of cycles of length
$2\ell$ in the random bipartite incidence multigraph of the
configuration, with two parallel incidence edges regarded as a cycle
of length two.  For every fixed $L$, the random variables
\[
 \mathsf C_1,\ldots,\mathsf C_L
\]
converge jointly to independent Poisson random variables with means
\begin{equation}
 \lambda_\ell
 =
 \frac{(r(d-1))^\ell}{2\ell},
 \qquad
 1\le\ell\le L.
 \label{eq:short-cycle-means}
\end{equation}
Consequently,
\begin{equation}
 \mathbb P(\mathsf C_1=\mathsf C_2=\mathsf C_3=\mathsf C_4=0)
 \longrightarrow
 \exp\left(
 -\sum_{\ell=1}^4\frac{(r(d-1))^\ell}{2\ell}
 \right)
 >0.
 \label{eq:positive-short-cycle-avoidance}
\end{equation}
\end{lemma}

\begin{proof}
For a fixed $\ell$, choose the $\ell$ distinct left vertices and the
$\ell$ distinct right vertices of an oriented cycle, choose the two
incident points at every chosen left vertex and the two positions at
every chosen right vertex, and then divide by the $2\ell$ choices of
a starting point and orientation.  This gives
\begin{equation}
 \mathbb E(\mathsf C_\ell)
 =
 \frac{
 (n)_\ell(Q)_\ell
 \bigl(d(d-1)k(k-1)\bigr)^\ell
 }{
 2\ell\,(M)_{2\ell}
 }.
 \label{eq:expected-short-cycles}
\end{equation}
Since $Q=M/k$, $M=nd$, and $k-1=r$, the right side of
\eqref{eq:expected-short-cycles} tends to the value in
\eqref{eq:short-cycle-means}.

More generally, fix nonnegative integers $m_1,\ldots,m_L$ and expand
the joint factorial moment
\[
 \mathbb E\left(\prod_{\ell=1}^L(\mathsf C_\ell)_{m_\ell}\right).
\]
The contribution from collections of pairwise vertex-disjoint cycles
is asymptotic to
\[
 \prod_{\ell=1}^L\lambda_\ell^{m_\ell}.
\]
Every other intersection pattern contains a connected component with
at least two cycles.  If such a pattern uses $v$ vertices of the
incidence graph and $e$ prescribed matching edges, then $e\ge v+1$.
There are $O(n^v)$ embeddings of the pattern, while the probability
of its prescribed matching edges is $O(M^{-e})=O(n^{-e})$.  Hence
every overlapping pattern contributes $O(n^{-1})$.  There are only
finitely many patterns for fixed $m_1,\ldots,m_L$, so
\[
 \mathbb E\left(\prod_{\ell=1}^L(\mathsf C_\ell)_{m_\ell}\right)
 \longrightarrow
 \prod_{\ell=1}^L\lambda_\ell^{m_\ell}.
\]
The method of moments gives the asserted joint independent Poisson
limit, and \eqref{eq:positive-short-cycle-avoidance} follows.
\end{proof}

\begin{lemma}
\label{lem:short-cycles-imply-sparse}
Let $\mathcal G$ be the event
\[
 \mathsf C_1=\mathsf C_2=\mathsf C_3=\mathsf C_4=0.
\]
Then, for $0<\delta\le1$,
\[
 \mathcal G\subseteq\mathcal D_\delta.
\]
In particular,
\begin{equation}
 \liminf_{n\to\infty}\mathbb P(\mathcal D_\delta)
 \ge
 \exp\left(
 -\sum_{\ell=1}^4\frac{(r(d-1))^\ell}{2\ell}
 \right)
 >0.
 \label{eq:sparse-conditioning-positive}
\end{equation}
\end{lemma}

\begin{proof}
The absence of a cycle of length two says that no right bucket is
matched twice to the same left bucket.  Thus every block has
$r+1$ distinct atoms.  The absence of a cycle of length four says
that two distinct blocks have at most one common atom.  Hence the
resulting set system is simple and linear.

Suppose first that $r\ge2$. For $2\le|A|\le r$, linearity gives
\[
 d_{r+1}(A)\le1
.
\]
Since $0<\delta\le1$,
\[
 r+1-|A|-\delta\ge0,
\]
and therefore \eqref{eq:random-regular-proper-codegree} holds.  If
$x\ne y$ and
\[
 B\in N_{r+1}(x)\cap N_{r+1}(y),
\]
then the two distinct blocks $B\cup\{x\}$ and $B\cup\{y\}$ share the
$r$-set $B$.  This contradicts linearity, since $r\ge2$.  Thus
\[
 d_{r+1}(x,y)=0,
\]
and \eqref{eq:random-regular-common-link} holds.

Now suppose that $r=1$, so that the resulting object is a graph.
The proper-codegree condition is empty.  If two distinct vertices
$x,y$ had two distinct common neighbors, the graph would contain a
$4$-cycle, which corresponds to a cycle of length eight in its
incidence graph.  Hence
\[
 d_{r+1}(x,y)\le1.
\]
Since $\delta\le1$, we have $1\le d^{1-\delta}$, and the common-link
condition follows.  Equation~\eqref{eq:sparse-conditioning-positive}
now follows from Lemma~\ref{lem:short-incidence-cycles}.
\end{proof}

\subsection{A necessary covering condition for maximality}

Fix a set $S\subseteq X$ of size $m$, and put
\[
 x=\frac{m}{n}.
\]
Let
\[
 A=\bigcup_{v\in S}P_v,
 \qquad
 B=P\setminus A,
 \qquad
 a=|A|=dm,
 \qquad
 b=|B|=d(n-m).
\]
A point $z\in B$ is called \emph{$S$-good} if the right bucket
matched to $z$ has its other $r$ positions matched to points of $A$.
Let
\[
 B_S=\{z\in B:z\text{ is $S$-good}\},
 \qquad
 L=L_S=|B_S|.
\]
We say that $S$ \emph{covers its complement} if, for every $v\in X\setminus S$,
\begin{equation}
 B_S\cap P_v\ne\varnothing
.
 \label{eq:covering-condition}
\end{equation}

If the resulting set system is simple and $S$ is a maximal
independent set, then $S$ covers its complement.  Indeed, for every
$v\notin S$, maximality gives a block contained in $S\cup\{v\}$.
Since $S$ is independent, that block contains $v$ and exactly $r$
atoms of $S$, so one point of $P_v$ is $S$-good.

Conditional on $L=\ell$, the set $B_S$ is a uniformly random
$\ell$-subset of $B$.  This follows from the invariance of the
configuration model under arbitrary permutations of the points of
$B$.

\begin{lemma}
\label{lem:conditional-covering-probability}
For $1\le\ell<b$,
\begin{equation}
 \mathbb P\bigl(S\text{ covers its complement}\mid L=\ell\bigr)
 \le
 3\sqrt b\,
 \exp\left(
 -(n-m)\left(1-\frac{\ell}{b}\right)^d
 \right).
 \label{eq:conditional-covering-bound}
\end{equation}
\end{lemma}

\begin{proof}
Put $p=\ell/b$, and let $\widetilde B\subseteq B$ be formed by
including every point independently with probability $p$.  Conditional
on $|\widetilde B|=\ell$, the set $\widetilde B$ is a uniformly random
$\ell$-subset of $B$, and therefore has the same distribution as
$B_S$ conditional on $L=\ell$.

For every $v\notin S$,
\[
 \mathbb P(\widetilde B\cap P_v=\varnothing)=(1-p)^d.
\]
The buckets $P_v$, $v\notin S$, are disjoint, so these events are
independent.  Consequently,
\begin{equation}
 \mathbb P\bigl(
 \widetilde B\cap P_v\ne\varnothing
 \text{ for every }v\notin S
 \bigr)
 =
 \bigl(1-(1-p)^d\bigr)^{n-m}
 \le
 \exp\bigl(-(n-m)(1-p)^d\bigr).
 \label{eq:independent-covering-bound}
\end{equation}
Stirling's inequalities give, at the mean of a binomial random
variable,
\begin{equation}
 \mathbb P(\operatorname{Bin}(b,p)=\ell)
 =
 \binom b\ell p^\ell(1-p)^{b-\ell}
 \ge
 \frac{1}{3\sqrt b}.
 \label{eq:binomial-mass-lower}
\end{equation}
Dividing \eqref{eq:independent-covering-bound} by
\eqref{eq:binomial-mass-lower} proves
\eqref{eq:conditional-covering-bound}.
\end{proof}

\subsection{An upper tail for the number of good points}

For an integer $u\ge0$, write
\[
 (z)_u=z(z-1)\cdots(z-u+1).
\]

\begin{lemma}
\label{lem:good-point-factorial-moment}
For every integer $1\le t\le Q$,
\begin{equation}
 \mathbb E\left(\binom Lt\right)
 =
 \binom Qt k^t
 \frac{(b)_t(a)_{rt}}{(M)_{kt}}.
 \label{eq:good-point-factorial-moment}
\end{equation}
For $t>Q$, $\mathbb E(\binom Lt)=0$, since at most one point is good in each right bucket.
Moreover, if $1\le kt\le M/2$, then
\begin{equation}
 \mathbb E\left(\binom Lt\right)
 \le
 \frac{(bx^r)^t}{t!}
 \exp\left(\frac{2k^2t^2}{M}\right).
 \label{eq:good-point-factorial-upper}
\end{equation}
\end{lemma}

\begin{proof}
Choose the $t$ right buckets containing the selected good points and,
in each chosen bucket, choose the unique position matched to a point
of $B$.  The selected positions must receive $t$ distinct points of
$B$, while the other $rt$ positions must receive distinct points of
$A$.  This gives \eqref{eq:good-point-factorial-moment}.

For $0\le u\le M/2$,
\begin{equation*}
 (M)_u
 =
 M^u\prod_{j=0}^{u-1}\left(1-\frac{j}{M}\right)
 \ge
 M^u\exp\left(-\frac{2u^2}{M}\right).
\end{equation*}
Also,
\[
 \binom Qt k^t
 \le
 \frac{M^t}{t!},
 \qquad
 (b)_t\le b^t,
 \qquad
 (a)_{rt}\le a^{rt}.
\]
Substitution into \eqref{eq:good-point-factorial-moment}, using
$a/M=x$ and $k=r+1$, gives
\eqref{eq:good-point-factorial-upper}.
\end{proof}

Put
\begin{equation}
 x_0
 =
 (1-\alpha)
 \left(\frac{\log d}{rd}\right)^{1/r}.
 \label{eq:random-regular-x0}
\end{equation}
Choose a constant $\eta>0$ so small that
\begin{equation}
 \beta=(1+4\eta)(1-\alpha)^r<1.
 \label{eq:random-regular-beta}
\end{equation}
For a fixed $S$ with $|S|=m\le x_0n$, define
\begin{equation*}
 \ell_0
 =
 \left\lceil(1+\eta)bx_0^r\right\rceil.
\end{equation*}

\begin{lemma}
\label{lem:good-point-upper-tail}
For all sufficiently large $d$, uniformly over all
$S\subseteq X$ with $|S|\le x_0n$,
\begin{equation}
 \mathbb P(L\ge\ell_0)
 \le
 \exp(-c_1n\log d),
 \label{eq:good-point-tail}
\end{equation}
where $c_1=c_1(r,\alpha,\eta)>0$.
\end{lemma}

\begin{proof}
For integers $1\le t\le\ell_0/2$,
\begin{equation*}
 \binom{\ell_0}{t}
 \ge
 \frac{\ell_0^t}{t!}
 \exp\left(-\frac{2t^2}{\ell_0}\right).
\end{equation*}
Since $x\le x_0$ and
$\ell_0\ge(1+\eta)bx_0^r$, Lemma~\ref{lem:good-point-factorial-moment}
and Markov's inequality give
\begin{equation}
 \mathbb P(L\ge\ell_0)
 \le
 \frac{\mathbb E\left(\binom Lt\right)}{\binom{\ell_0}{t}} 
 \le
 \exp\left(
 -t\log(1+\eta)
 +\frac{2k^2t^2}{M}
 +\frac{2t^2}{\ell_0}
 \right).
 \label{eq:good-point-tail-intermediate}
\end{equation}
Choose
\[
 t=\lfloor c\eta\ell_0\rfloor,
\]
where $c=c(r)>0$ is sufficiently small.  Since
\[
 \frac{\ell_0}{M}
 =
 O_{r,\alpha,\eta}\left(\frac{\log d}{d}\right),
\]
we may first choose $c$ sufficiently small and then choose $d$
sufficiently large so that the exponent in
\eqref{eq:good-point-tail-intermediate} is at most
\[
 -c_2\eta^2\ell_0
\]
for some $c_2=c_2(r)>0$.

For sufficiently large $d$, we have $x_0\le1/2$.  Since
$m\le x_0n$,
\[
 \ell_0\ge(1+\eta)d(n-m)x_0^r
 \ge\frac{(1+\eta)(1-\alpha)^r}{2r}\,n\log d.
\]
This proves \eqref{eq:good-point-tail}.
\end{proof}

\subsection{The probability that a fixed set covers}

\begin{lemma}
\label{lem:fixed-set-covering}
For all sufficiently large $d$, uniformly over all
$S\subseteq X$ with $|S|\le x_0n$,
\begin{equation}
 \mathbb P(S\text{ covers its complement})
 \le
 \exp(-c_1n\log d)
 +
 3\sqrt{nd}\,
 \exp\left(-\frac{n}{2}d^{-\beta/r}\right),
 \label{eq:fixed-set-covering-bound}
\end{equation}
where $\beta<1$ is defined in \eqref{eq:random-regular-beta}.
\end{lemma}

\begin{proof}
If $L<\ell_0$, then, for all sufficiently large $n$,
\[
 \frac{L}{b}
 \le
 \frac{\ell_0}{b}
 \le
 (1+2\eta)x_0^r.
\]
The inequality valid for $0\le z\le1/2$,
\[
 \log(1-z)\ge-z-2z^2,
\]
therefore gives
\begin{equation*}
 \left(1-\frac{L}{b}\right)^d
 \ge
 \exp\left(
 -d(1+2\eta)x_0^r
 -2d(1+2\eta)^2x_0^{2r}
 \right) 
 \ge
 d^{-\beta/r},
\end{equation*}
for all sufficiently large $d$, by
\eqref{eq:random-regular-x0} and
\eqref{eq:random-regular-beta}.  Also $n-m\ge n/2$ and $b\le nd$.
Lemma~\ref{lem:conditional-covering-probability} now gives
\[
 \mathbb P(S\text{ covers its complement}\mid L<\ell_0)
 \le
 3\sqrt{nd}\,
 \exp\left(-\frac{n}{2}d^{-\beta/r}\right).
\]
Combining this with Lemma~\ref{lem:good-point-upper-tail} proves
\eqref{eq:fixed-set-covering-bound}.
\end{proof}

\subsection{The union bound}

Let $\mathsf Z$ be the number of sets $S\subseteq X$ such that
\[
 |S|\le x_0n
\]
and $S$ covers its complement.  Let
\[
 \mathsf H(x)=-x\log x-(1-x)\log(1-x)
\]
denote the binary entropy function.  For $x_0\le1/2$,
\begin{equation}
 \sum_{m\le x_0n}\binom nm
 \le
 (n+1)\exp(n\mathsf H(x_0)).
 \label{eq:small-set-count}
\end{equation}
Moreover,
\begin{equation*}
 \mathsf H(x_0)
 \le
 x_0\log\left(\frac{e}{x_0}\right)
 =
 O_{r,\alpha}\left(
 d^{-1/r}(\log d)^{1+1/r}
 \right).
\end{equation*}
Since $\beta<1$,
\begin{equation*}
 d^{-1/r}(\log d)^{1+1/r}
 =
 o(d^{-\beta/r}).
\end{equation*}
Choose $d$ sufficiently large that
\begin{equation}
 \mathsf H(x_0)
 \le
 \min\left\{
 \frac{c_1\log d}{2},
 \frac{d^{-\beta/r}}{4}
 \right\}.
 \label{eq:entropy-dominated}
\end{equation}
Equations
\eqref{eq:fixed-set-covering-bound},
\eqref{eq:small-set-count}, and
\eqref{eq:entropy-dominated} give
\begin{equation*}
 \mathbb E(\mathsf Z)
 \le
 (n+1)e^{n\mathsf H(x_0)}
 \left[
 e^{-c_1n\log d}
 +
 3\sqrt{nd}\,
 e^{-\frac{n}{2}d^{-\beta/r}}
 \right] 
 =o(1).
\end{equation*}
Therefore
\begin{equation}
 \mathbb P(\mathsf Z>0)=o(1).
 \label{eq:no-small-covering-set}
\end{equation}

\subsection{Conditioning on simplicity and \texorpdfstring{$\delta$}{delta}-sparseness}

\begin{proof}[Proof of Theorem~\textup{\ref{thm:random-regular-delta-sparse}}]
Let $\mathcal E$ be the event that the resulting simple set system
contains a maximal independent set of size at most $x_0n$.  On
$\mathcal D_\delta$, every such set covers its complement, by
\eqref{eq:covering-condition}.  Hence
\[
 \mathcal E\cap\mathcal D_\delta
 \subseteq
 \{\mathsf Z>0\}.
\]
By \eqref{eq:no-small-covering-set} and
\eqref{eq:sparse-conditioning-positive},
\begin{equation}
 \mathbb P(\mathcal E\mid\mathcal D_\delta)
 =\frac{\mathbb P(\mathcal E\cap\mathcal D_\delta)}{\mathbb P(\mathcal D_\delta)}
 \le\frac{\mathbb P(\mathsf Z>0)}{\mathbb P(\mathcal D_\delta)}=o(1).
 \label{eq:conditioned-small-maximal-probability}
\end{equation}

It remains only to identify the conditioned distribution.  Every
simple labeled $d$-regular $(r+1)$-uniform set system with $Q$ blocks
is represented by exactly
\[
 (d!)^n Q!(k!)^Q
\]
perfect matchings in the model above: order the $d$ incidences at
each left vertex, assign the blocks to the labeled right buckets,
and order the $k$ positions in every right bucket.  This number is
independent of the set system.  Consequently, conditional on
$\mathcal D_\delta$, the configuration model is uniform over the
simple $d$-regular $\delta$-sparse $(r+1)$-uniform set systems on
$[n]$.  Equation~\eqref{eq:conditioned-small-maximal-probability}
therefore proves the theorem.
\end{proof}

\section{Deterministic algorithm for maximal independent sets}\label{sec:alg}

In this section we give a deterministic algorithm for the conclusion of
Theorem~\ref{thm:main}. We retain the residual set systems and the
antichain update from its proof, but change the sampling probabilities
and determine the choices by the deterministic local-lemma algorithm of
Harris~\cite{HarrisDetLLL}, based on the resampling algorithm of
Moser and Tardos~\cite{MoserTardos}.

The process used for Theorem~\ref{thm:main} has the fixed parameter
$\eta=(\log D)^{-\delta_3/8}$ and, by
Subsection~\ref{sec:general-completion}, terminates after
$T\le c_3(\log D)^{1+\delta_3/8}$ rounds.
For the algorithm in this section, we instead define
\begin{equation}\label{eq:alg-sampling}
 \overline{\eta}_t=\Sigma_t^{-\delta_3/4},\qquad
 q_t=\Sigma_t^{-1+\delta_2},\qquad
 p_t=\overline{\eta}_t\Sigma_t^{-\delta_2}.
\end{equation}
For each $x\in X^t\setminus\mathsf W^t$, let $\xi_x$ be independent
Bernoulli$(q_t)$ variables, and let
$Z^t=\{x\in X^t\setminus\mathsf W^t:\xi_x=1\}$.
For a fixed $Z^t$, let $(\chi_x:x\in Z^t)$ be independent
Bernoulli$(p_t)$ variables, and let $Y^t=\{x\in Z^t:\chi_x=1\}$.
The algorithm
deterministically assigns values to these variables while preserving
the inequalities proved using their product distributions.
Thus $q_tp_t=\overline{\eta}_t/\Sigma_t$, and the one-round estimates
are applied with $\overline{\eta}_t$ in place of $\eta$.
We prove
\[
 \overline{T}\le\frac{2^{32r^2}}{\delta_3}D^{\delta_3/4}
\]
for the number of rounds in Lemma~\ref{lem:alg-tracking}.
The bound is a power of $D$, in place of a power of $\log D$.
At time zero the ratio of the two sampling parameters is
$\overline{\eta}_0/\eta
=\Sigma_0^{-\delta_3/4}(\log D)^{\delta_3/8}$.
These smaller probabilities permit polynomial approximations
whose degrees are independent of $n$ and $D$.

We also replace the minimum survival probability in
\eqref{eq:degree-survival-probabilities} by an explicitly calculated
number. Let $\sigma_t,z_t$ be given by
\eqref{eq:general-acceptance-reference}, with $\eta=\overline{\eta}_t$,
and define
\begin{align}
 e_t&=4b_3\bigl(\overline{\eta}_t\Sigma_t^{-\delta_3}
                               +\Sigma_t^{-4\delta_3}\bigr),
       \label{eq:alg-survival-error}\\
 \overline{\omega}_t&=\exp\left(-\sum_{i=2}^r d_{i,t}z_t^{i-1}-2e_t\right).
       \label{eq:alg-survival-reference}
\end{align}

The main additional difficulty is to compute conditional expectations
for acceptance and survival: an indicator such as $J_x$ or $\psi_x$
may depend on every block containing the atom and hence on a number
of variables that grows with the input. We bound each required product
of these indicators above and below by explicit polynomials whose
expected difference is at most $\Sigma_t^{-A}$, for a fixed exponent
$A$. Their conditional expectations can then be evaluated by the
monomial formula in Lemma~\ref{lem:alg-local-polynomial-bounds}.
The weighted selection result in Proposition~\ref{prop:alg-weighted-selection}
enforces the degree and intersection inequalities while controlling
the two quantities $|J^{\overline{T}}|$ and
$\sum_{t<\overline{T}}|Y^t|+|X^{\overline{T}}|$ that bound the final independent set.

We use the incidence-list representation of a set system, and write
\begin{equation*}
 N=N(\mathcal S)=n+\sum_{S\in\mathcal S}|S|,
 \qquad n=|X(\mathcal S)|.
\end{equation*}
All blocks have size at most $r$, and repeated blocks are discarded.
Consequently,
\begin{equation}\label{eq:alg-input-n}
 n\le N\le n+\sum_{i=2}^ri\binom ni\le r^2n^r.
\end{equation}
The constants $\delta_1,\delta_2,\delta_3$, the functions
$D(\vb d),\rho(\vb d)$, and $\varepsilon(\triangle)$ retain their
definitions from Sections~\ref{sec:process}--\ref{sec:completing-process}.

\begin{theorem}\label{thm:deterministic-algorithm}
Let $r\ge2$, $c\ge1$ and $0<\delta_0\le1/2$.  There are constants
$\overline{c}_4=\overline{c}_4(r,\delta_0,c)$ and $\overline{C}=\overline{C}(r,\delta_0,c)$ with the
following property. If $\mathcal S$ satisfies the hypotheses of
Theorem~\ref{prop:quantitative-main} and $D\ge \overline{C}$, then a
deterministic algorithm, using at most $n^{\overline{c}_4}$ bit operations, produces
a maximal independent set $I$ such that
\[
 (1-\varepsilon(\triangle))n\rho(\vb d)
 \le |I|\le(1+\varepsilon(\triangle))n\rho(\vb d).
\]
Here $\varepsilon(\triangle)=(\log\triangle)^{-\delta_3/32}$.
\end{theorem}

\begin{corollary}\label{cor:deterministic-main}
For fixed $\delta>0$ and $r\ge2$, the maximal independent set in
Theorem~\ref{thm:main} can be constructed deterministically in time
$n^{\overline{c}_5(\delta,r)}$. Here one may take
$\overline{c}_5(\delta,r)=\overline{c}_4(r,\min\{\delta,1/2\},1)+2$,
where $\overline{c}_4$ is the exponent in Theorem~\ref{thm:deterministic-algorithm}.
\end{corollary}

\begin{proof}
Let $\delta_0=\min\{\delta,1/2\}$, $c=1$, and
$\overline{D}_0=\overline{C}(r,\delta_0,1)$. Theorem~\ref{thm:deterministic-algorithm}
applies when $D\ge \overline{D}_0$. For $1<D<\overline{D}_0$, we construct a maximal
independent set by considering the atoms in a fixed order and adding
an atom precisely when its addition preserves independence. The bounded-degree calculation in the deduction
of Theorem~\ref{thm:main} from Theorem~\ref{prop:quantitative-main}
applies with the common choice
\[
 \overline{K}(\delta,r)=\frac{(r-1)^2\overline{D}_0(\log \overline{D}_0)^{\delta_3/32}}{\log2}.
\]
Thus taking the existential constant $K$ to be $\overline{K}$ gives the bound
in Theorem~\ref{thm:main} in both ranges. Testing whether a new
atom can be added requires a scan of the incidence list, so the
greedy algorithm uses at most $nN$ such elementary incidence
operations. Equation~\eqref{eq:alg-input-n} and the exponent
specified in the corollary bound both algorithms by
$n^{\overline{c}_5(\delta,r)}$.
\end{proof}

All constants in this section are independent of $n$ and $D$.
We place a bar over a constant or process parameter when its choice
differs from the corresponding choice in
Sections~\ref{sec:process}--\ref{sec:completing-process}.
In particular, $\overline{T}_{\mathrm{init}}$ and $\overline{T}$ are the
initial-phase endpoint and stopping time, and $\overline{C}$ is the
threshold for $D$. The constant $c$ multiplying the right sides of
\eqref{eq:initialP} and \eqref{eq:initialCommon} is unchanged and
remains unbarred. The constants $\delta_i,a,b_1,b_2,b_3,c_1,c_3$
that we use unchanged retain their original notation.
The auxiliary constants $\overline{a}_1,\ldots,\overline{a}_6$ and
$\overline{s}$ are new constants for the polynomial and computation
bounds, with no counterparts in the earlier process.
The new constants $\overline{c}_4(r,\delta_0,c)$ and
$\overline{c}_5(\delta,r)$ are the running-time exponents in the theorem
and corollary, respectively. We state the arguments of auxiliary
constants when they depend on additional integers.
Throughout this section $\varepsilon$ abbreviates
$(\log\triangle)^{-\delta_3/32}$.

\subsection{Moments and polynomial bounds}

We first bound even moments and approximate products of acceptance
and survival indicators by polynomials of bounded degree.
The explicit polynomial formulas permit the calculation of conditional
expectations after any subset of the independent variables has been fixed.

\begin{proposition}\label{prop:alg-fixed-moments}
Let $m\ge1$ be an integer.
\begin{enumerate}
\item Let $F$ be a sum of indicators with read bound $k$ in
Definition~\ref{def:coordinate-supports}, and let $\mu=\mathbb E(F)$.
Then
\begin{equation}\label{eq:alg-read-moment}
 \mathbb E\bigl(|F-\mu|^{2m}\bigr)
 \le 2m!(4k\mu)^m+2(2m)!(4k)^{2m}.
\end{equation}
\item Let $F$ be a polynomial of degree at most $s\ge1$ with
nonnegative coefficients, as in Definition~\ref{def:polynomial}.
Let $\mu=\mathbb E(F)$, $E_* =\max_{1\le j\le s}E_j$, and
\[
 \overline{a}_1(s,m)=\bigl[8m(2ms+1)^{s+1}\bigr]^{4m}.
\]
Then
\begin{equation}\label{eq:alg-polynomial-moment}
 \mathbb E\bigl(|F-\mu|^{2m}\bigr)
 \le \overline{a}_1(s,m)\bigl((\mu E_*)^m+E_*^{2m}\bigr).
\end{equation}
\end{enumerate}
\end{proposition}

\begin{proof}[Proof of \textup{(1)}]
Proposition~\ref{prop:readk} gives, for $x>0$,
\[
 \mathbb P(|F-\mu|\ge x)
 \le2\exp\left(-\frac{x^2}{2k(\mu+x)}\right)
 \le2e^{-x^2/(4k\mu)}+2e^{-x/(4k)}.
\]
For $\mu>0$, integration gives
\begin{align*}
 \mathbb E\bigl(|F-\mu|^{2m}\bigr)
 &=\int_0^\infty 2mx^{2m-1}\mathbb P(|F-\mu|>x)\,dx\\
 &\le4m\int_0^\infty x^{2m-1}
                  \bigl(e^{-x^2/(4k\mu)}+e^{-x/(4k)}\bigr)\,dx\\
 &=2m!(4k\mu)^m+2(2m)!(4k)^{2m}.
\end{align*}
If $\mu=0$, then $F=0$ almost surely.

\emph{Proof of \textup{(2)}.}
We write $U_A=\prod_{x\in A}U_x$ and
$F=z_{\varnothing}+\sum_{1\le|A|\le s}z_AU_A$.
The centered moment has the exact expansion
\begin{equation}\label{eq:alg-centered-expansion}
 \mathbb E\bigl((F-\mu)^{2m}\bigr)
 =\sum_{A_1,\ldots,A_{2m}}
 \left(\prod_{j=1}^{2m}z_{A_j}\right)
 \mathbb E\left(\prod_{j=1}^{2m}
                      (U_{A_j}-\mathbb E(U_{A_j}))\right),
\end{equation}
where $1\le|A_j|\le s$. For each tuple, let the graph on $[2m]$
have edge set $\{\{j,k\}:j\ne k,\ A_j\cap A_k\ne\varnothing\}$. 
If this graph has an isolated vertex, independence and
$\mathbb E(U_{A_j}-\mathbb E(U_{A_j}))=0$ make its summand zero.
Otherwise the graph has at most $m$ connected components.

For a component with index set $C$ and $|C|=q$, the product is
\[
 \prod_{j\in C}(U_{A_j}-\mathbb E(U_{A_j}))
 =\sum_{D\subseteq C}(-1)^{q-|D|}
     \left(\prod_{j\in D}U_{A_j}\right)
     \left(\prod_{j\in C\setminus D}\mathbb E(U_{A_j})\right).
\]
For each $D\subseteq C$, independence of the Bernoulli variables gives
\[
 \mathbb E\left(\prod_{j\in D}U_{A_j}\right)
       \prod_{j\in C\setminus D}\mathbb E(U_{A_j})
 \le \mathbb E\left(\prod_{j\in C}U_{A_j}\right).
\]
Indeed, every probability on the right occurs at least once on the
left, and additional occurrences multiply the left side by factors
in $[0,1]$.
By connectedness, the indices of each component admit an order
$j_1,\ldots,j_q$ such that
$A_{j_\ell}\cap\bigcup_{k<\ell}A_{j_k}\ne\varnothing$ for $2\le\ell\le q$.
We sum over the sets in this order. The sum for the first set is at most $\mu$.
If $R$ is the preceding union, then $|R|\le2ms$, and the sum for
the next set is bounded by
\[
 \sum_{\substack{A:A\cap R\ne\varnothing\\1\le|A|\le s}}
 z_A\prod_{x\in A\setminus R}\mathbb P(U_x=1)
 \le\sum_{\substack{T\subseteq R\\1\le|T|\le s}}
                         \mathbb E(\partial_TF)
 \le(2ms+1)^sE_*.
\]
There are at most $(2m)^{2m}$ partitions of $[2m]$ and at most
$(2m)!$ orders. Thus \eqref{eq:alg-centered-expansion} is at most
\[
 2^{2m}(2m)^{2m}(2m)!(2ms+1)^{2ms}
                 \sum_{j=1}^m\mu^jE_*^{2m-j}.
\]
If $\mu\ge E_*$, each term of the last sum is at most
$(\mu E_*)^m$; if $\mu<E_*$, each is at most $E_*^{2m}$.
The specified value of $\overline{a}_1(s,m)$ dominates the displayed
coefficient multiplied by $m$.
\end{proof}

\begin{lemma}\label{lem:alg-union-moment}
Let $s,L\ge1$ be integers, $c_0\ge1$, and $h\ge1$.
Let $\mathcal H$ be a family of distinct nonempty sets of size at
most $s$ such that, for $0\le j\le i\le s$ and every $j$-set $A$,
\begin{equation}\label{eq:alg-extension-count}
 |\{B\in\mathcal H:|B|=i,\ A\subseteq B\}|\le c_0h^{i-j}.
\end{equation}
Let $(U_x)$ be independent Bernoulli variables with
$\mathbb P(U_x=1)=p=\eta/h$, where $0<\eta\le1/2$, and put
$K=\sum_{B\in\mathcal H}\prod_{x\in B}U_x$. Let
\[
 \overline{a}_2(s,c_0,L)=c_0^L(s+1)^{2L}(sL)^{sL},\qquad
 \nu_s(L)=\min\left\{v:\sum_{i=1}^s\binom vi\ge L\right\}.
\]
Then
\begin{equation}\label{eq:alg-factorial-moment}
 \mathbb E\left(\binom KL\right)
       \le \overline{a}_2(s,c_0,L)\eta^{\nu_s(L)}.
\end{equation}
\end{lemma}

\begin{proof}
The definition of $K$ gives
\begin{equation}\label{eq:alg-factorial-expansion}
 \mathbb E\left(\binom KL\right)
 =\frac1{L!}\sum_{\substack{B_1,\ldots,B_L\in\mathcal H\\
                             B_j\ne B_k\ (j\ne k)}}
                              p^{|B_1\cup\cdots\cup B_L|}.
\end{equation}
For a tuple in this sum put $R_0=\varnothing$,
$R_\ell=B_1\cup\cdots\cup B_\ell$, $i_\ell=|B_\ell|$, and
$j_\ell=|B_\ell\cap R_{\ell-1}|$.
We first consider tuples for which the integers
$(i_\ell,j_\ell)_{\ell=1}^L$ have specified values.
There are at most $\binom{sL}{j_\ell}$ choices for
$B_\ell\cap R_{\ell-1}$. After this intersection is fixed,
\eqref{eq:alg-extension-count} gives at most
$c_0h^{i_\ell-j_\ell}$ choices for $B_\ell$.
Since $|R_L|=\sum_\ell(i_\ell-j_\ell)$, the contribution for
these integers in \eqref{eq:alg-factorial-expansion} is at most
\[
 \frac{c_0^L}{L!}\prod_{\ell=1}^L\binom{sL}{j_\ell}
                       (hp)^{|R_L|}
 \le c_0^L(sL)^{sL}\eta^{|R_L|}.
\]
The $L$ distinct sets $B_\ell$ are nonempty subsets of $R_L$ of
size at most $s$, so $|R_L|\ge\nu_s(L)$.
There are at most $(s+1)^{2L}$ choices of the integers
$(i_\ell,j_\ell)_{\ell=1}^L$. Summing the last bound proves
\eqref{eq:alg-factorial-moment}.
\end{proof}

We use the Bonferroni inequalities in the following pointwise
form; see Feller~\cite[Chapter IV, Section 5(c)]{Feller}.
For indicators $B_1,\ldots,B_M$, let $K_B=\sum_{a=1}^MB_a$ and
\begin{equation}\label{eq:alg-bonferroni}
 B^{[j]}=\sum_{\substack{A\subseteq[M]\\|A|\le j}}
                   (-1)^{|A|}\prod_{a\in A}B_a
         =\sum_{\ell=0}^j(-1)^\ell\binom{K_B}{\ell}.
\end{equation}
For every integer $L\ge0$,
\begin{equation}\label{eq:alg-bonferroni-gap}
 B^{[2L+1]}\le\prod_{a=1}^M(1-B_a)\le B^{[2L]},\qquad
 B^{[2L]}-B^{[2L+1]}=\binom{K_B}{2L+1}.
\end{equation}
The equality follows by subtracting the sums in
\eqref{eq:alg-bonferroni}: their only different term is the term
with $\ell=2L+1$. The inequalities also follow directly: if
$K_B=0$, all three quantities equal one; if $K_B\ge1$, then
\[
 \sum_{\ell=0}^j(-1)^\ell\binom{K_B}{\ell}
       =(-1)^j\binom{K_B-1}{j},
 \qquad \prod_{a=1}^M(1-B_a)=0.
\]
Here a binomial coefficient is zero when its lower index exceeds
its nonnegative upper index.

\begin{lemma}\label{lem:alg-local-polynomial-bounds}
Let $A>0$ and let $m\ge1$ be an integer. Suppose $Z$ satisfies
\eqref{eq:general-constant-focused-completions} and
\eqref{eq:general-focused-degrees}, and
$E_i,\widehat E_i\le\Sigma^{i-1}$. Let
\[
 \overline{\eta}=\Sigma^{-\delta_3/4},\qquad h=\Sigma^{\delta_2},\qquad
 p=\overline{\eta}/h.
\]
Let $(\chi_x:x\in Z)$ be independent Bernoulli$(p)$ variables,
and set $\chi_x=0$ for $x\notin Z$, as in
Definition~\ref{def:degree-original-variables}.
Independently of these variables, let
$(\mathsf w_x:x\in X\setminus\mathsf W)$ be independent Bernoulli
variables with specified probabilities; set $\mathsf w_x=1$ on
$\mathsf W$. Thus $\chi_x=1$ means $x\in Y$, whereas
$\psi_x(1-\mathsf w_x)=1$ means $x\in X'\setminus\mathsf W'$.
Let $\mathsf{F}$ be a product of at most $m$ indicators of the following types:
\begin{enumerate}[label=\textup{(\roman*)}]
\item $J_x$, $\psi_x$, or $\psi_x(1-\mathsf w_x)$;
\item $\mathsf Z_B^{\deg}$ or $\widehat{\mathsf Z}_B^{\deg}$
from Definition~\ref{def:degree-sums};
\item $\mathsf Z_B^{\rm cov}$ from
Definition~\ref{def:coverage-random-variables}.
\end{enumerate}
One may also replace specified variables $\chi_x$ by zero in
these formulas, as in Definitions~\ref{def:degree-modified-variables}
and \ref{def:coverage-random-variables}.
There are constants $\overline{s}(A,m,r,\delta_0)$ and
$\overline{a}_3(A,m,r,\delta_0)$, specified in
\eqref{eq:alg-explicit-polynomial-constants}, such that, if
$\Sigma\ge \overline{s}(A,m,r,\delta_0)$, the following statements hold.
\begin{enumerate}[label=\textup{(\arabic*)},ref=\arabic*]
\item There are multilinear polynomials $\mathsf{F}^-,\mathsf{F}^+$ satisfying
\begin{equation}\label{eq:alg-sandwich}
 \mathsf{F}^-\le \mathsf{F}\le \mathsf{F}^+\quad\text{on }\{0,1\}^{X\times\{1,2\}},
 \qquad \mathbb E(\mathsf{F}^+-\mathsf{F}^-)\le\Sigma^{-A}.
\end{equation}
\item The degrees of $\mathsf{F}^-,\mathsf{F}^+$ are at most $\overline{a}_3$,
and they involve only variables occurring in the factors of $\mathsf{F}$.
\item Their descriptions and conditional expectations can be computed
using at most $N^{\overline{a}_3}$ arithmetic operations on the
supplied Bernoulli probabilities.
\end{enumerate}
\end{lemma}

\begin{proof}[Proof of \textup{(1)}]
We give both the formulas and the truncation orders. For $V\subseteq X$, we define
\begin{equation}\label{eq:alg-zero-substitution}
 \chi_x^V=\begin{cases}0,&x\in V,\\\chi_x,&x\notin V,\end{cases}
 \qquad
 J_x^V=\chi_x^V\prod_{x\in S\in\mathcal S}
                 \left(1-\prod_{y\in S\setminus\{x\}}\chi_y^V\right),
 \qquad J_B^V=\prod_{x\in B}J_x^V.
\end{equation}
These are functions of the original independent variables.
They do not require conditioning their joint law. For example,
$J_x^{\{v\}}=J_x^v$ in \eqref{eq:degree-modified-J}, and
\eqref{eq:coverage-variable-formula} becomes
\begin{equation}\label{eq:alg-coverage-explicit}
 \mathsf Z_B^{\rm cov}
       =1-\prod_{x\in C_B\cap Z}(1-J_x^{F_0\cup B}).
\end{equation}
In particular, the phrase ``setting variables to zero'' means the
literal substitution in \eqref{eq:alg-zero-substitution}.

For comparison, in the experiment without substitutions the
indicator of $Q\subseteq I$ and $R\subseteq X'$ is exactly
\begin{equation}\label{eq:alg-joint-explicit}
 \prod_{x\in Q}J_x\prod_{v\in R}\psi_v
 =\left(\prod_{x\in Q}J_x\right)
  \left(\prod_{v\in R}(1-\chi_v)\right)
  \prod_{v\in R}\prod_{\substack{v\in S\in\mathcal S\\
                                  S\setminus\{v\}\subseteq Z}}
                       (1-J_{S\setminus\{v\}}).
\end{equation}
This is an equality of functions for every assignment of the
variables, by \eqref{eq:degree-original-variables}. For instance,
if $Q\cap R\ne\varnothing$, the right side contains
$\chi_x(1-\chi_x)=0$; in that case the event on the left is empty.
Formula \eqref{eq:alg-joint-explicit}, rather than any independence
assumption about $J_x$ or $\psi_v$, is used below.

Let $s=r-1$ and
\[
 t=\left\lceil m2^{r+4}r^2(rb_2+1)\right\rceil,
 \qquad u=t2^t,\qquad c_0=2^{t+2}tr(rb_2+3).
\]
We express each factor of $\mathsf{F}$ using
\eqref{eq:degree-original-variables}, \eqref{eq:degree-modified-R},
\eqref{eq:degree-restricted-psi}, \eqref{eq:degree-restricted-R}, and
\eqref{eq:alg-coverage-explicit}.
If $P_1,\ldots,P_k$ are the products for which a factor $1-P_a$
occurs in this expression, we use the identity
\[
 \prod_{a=1}^k(1-P_a)
       =\sum_{D\subseteq[k]}(-1)^{|D|}\prod_{a\in D}P_a.
\]
Thus $\mathsf{F}=\sum_{\ell=1}^{M_0}e_\ell T_\ell$, where
$M_0\le2^t$, $e_\ell\in\{-1,1\}$, and each $T_\ell$ has the form
\begin{equation}\label{eq:alg-basic-term}
 T=\left(\prod_{a=1}^dJ_{Q_a}^{V_a}\right)
    \left(\prod_{x\in R}(1-\chi_x)\right)
    \left(\prod_{x\in W_0}(1-\mathsf w_x)\right)
    \prod_{(a,B)\in\mathcal A}(1-J_B^{V_a}).
\end{equation}
Here $d,|R|,|W_0|,|\bigcup_aQ_a|\le t$;
$\mathcal A$ is a finite indexed collection, and each $B$ is
nonempty with $|B|\le s$. Factors with $B=\varnothing$ equal
$1-J_\varnothing^{V_a}=0$ and are discarded together with their
whole term. Repeated identical factors are retained just once.
There are at most $t$ families of sets $B$: subfamilies of links
$N_i(v,Z)$, the singleton families in
\eqref{eq:alg-coverage-explicit}, and the completion families
$\mathcal C_Z(F)$ with $|F|\ge2$. The last have at most $rb_2+1$
members by \eqref{eq:general-constant-focused-completions}.
The bound on $t$ follows by counting the at most $2^r$ subsets
$F\subseteq B$ in \eqref{eq:degree-hat-indicator}, at most $r$
atoms per completion, and at most $m$ original factors.

We first approximate the last product in
\eqref{eq:alg-basic-term}. Let $Q=\bigcup_aQ_a$ and
$K=\sum_{(a,B)\in\mathcal A}J_B^{V_a}$. If $Q\nsubseteq Z$, some $Q_a\cap V_a\ne\varnothing$, or
$Q\cap R\ne\varnothing$, then $T$ is identically zero and is
omitted. Otherwise $\mathbb P(\chi_x=1\text{ for all }x\in Q)>0$.
Let
\[
 \mathcal D_Q=\{B\setminus Q:(a,B)\in\mathcal A,\ B\nsubseteq Q\},
 \qquad K_0=\sum_{F\in\mathcal D_Q}\prod_{x\in F}\chi_x,
\]
where $\mathcal D_Q$ contains each distinct set just once.
Conditional on $\chi_x=1$ for all $x\in Q$, the variables in this
sum are products of independent Bernoulli$(p)$ variables outside $Q$. There are at most $u$ indexed sets $B\subseteq Q$ and at
most $u$ indexed sets giving any one nonempty difference
$B\setminus Q$. Therefore
\begin{equation}\label{eq:alg-multiplicity-bound}
 K\le u+uK_0,\qquad
 \binom Kj\le[u(j+1)]^j
 \sum_{\ell=\max\{0,\lceil j/u\rceil-1\}}^j\binom{K_0}{\ell}.
\end{equation}
For the second inequality, a set of $j$ indexed members with
$J_B^{V_a}=1$ determines $\ell$ distinct members of $\mathcal D_Q$
whose products equal one. There are $\binom{K_0}{\ell}$ possible
collections of these members. For each collection, the number of
available indexed members is at most $u(\ell+1)\le u(j+1)$, so the
number of selections is at most $[u(j+1)]^j$.
The inequality $j\le u(\ell+1)$ gives
$\ell\ge\lceil j/u\rceil-1$, proving the stated summation range.

The family of these distinct nonempty differences satisfies
\eqref{eq:alg-extension-count} with the displayed $c_0$.
For each specified set $B\cap Q$, of which there are at most $2^t$,
we count the possible original blocks. For a link
$B=S\setminus\{v\}$ with $|S|=i$, if only $v$ has been specified,
\eqref{eq:general-focused-degrees} gives at most $3h^{i-1}$ choices
of $S$. If a set of at least two but fewer than $i$ atoms of $S$
has been specified, \eqref{eq:general-constant-focused-completions}
gives the required bound on its completions. If all atoms have
been specified, there is at most one choice of $S$. For the singleton families,
$|C_B\cap Z|\le3rh$, since
$C_B\subseteq\bigcup_{x\in B}N_2(x,X\setminus\mathsf W)$.
These bounds, summed over the at most $t$ families, give $c_0$.

We define the integers
\begin{equation}\label{eq:alg-truncation-orders}
 v=\left\lceil\frac{4(A+4)}{\delta_3}\right\rceil,
 \qquad j=2u\bigl((v+1)^r+2\bigr)+1,
 \qquad j_1=2(v+1)^r+1,
 \qquad q_1=t+rj.
\end{equation}
Both $j$ and $j_1$ are odd. For $k\in\{j-1,j\}$, let
\[
 T^{[k]}=
 \left(\prod_{a=1}^dJ_{Q_a}^{V_a}\right)
 \left(\prod_{x\in R}(1-\chi_x)\right)
 \left(\prod_{x\in W_0}(1-\mathsf w_x)\right)
 \sum_{\substack{\mathcal A'\subseteq\mathcal A\\|\mathcal A'|\le k}}
       (-1)^{|\mathcal A'|}\prod_{(a,B)\in\mathcal A'}J_B^{V_a}.
\]
Here subsets of $\mathcal A$ refer to its indexed members.
Equation~\eqref{eq:alg-bonferroni-gap} gives
\[
 T^{[j]}\le T\le T^{[j-1]},\qquad
 0\le T^{[j-1]}-T^{[j]}
       \le\left(\prod_{x\in Q}\chi_x\right)\binom Kj.
\]
Equations \eqref{eq:alg-multiplicity-bound} and
\eqref{eq:alg-factorial-moment} show that its expectation is at most
\begin{equation}\label{eq:alg-outer-error}
 E_0\overline{\eta}^{v+1},\qquad
 E_0=[u(j+1)]^j
       \sum_{\ell=\lceil j/u\rceil-1}^j \overline{a}_2(s,c_0,\ell).
\end{equation}
Here $\nu_s(\ell)\ge v+1$ for every index in this sum, because
$\ell>(v+1)^r\ge\sum_{i=1}^s\binom vi$.

Every term of these two sums contains a product
$\prod_aJ_{Q'_a}^{V_a}$ with $Q'=\bigcup_aQ'_a$ and $|Q'|\le q_1$.
If $Q'\nsubseteq Z$ or some $Q'_a\cap V_a\ne\varnothing$,
this product is identically zero by \eqref{eq:alg-zero-substitution}. Otherwise let
\[
 \mathcal D=\{S\setminus Q':S\in\mathcal S,\ S\subseteq Z,\
       S\cap Q'_a\ne\varnothing,\ S\cap V_a=\varnothing
                        \text{ for at least one }a\}.
\]
With $\chi_F=\prod_{x\in F}\chi_x$, the exact identity is
\begin{equation}\label{eq:alg-acceptance-expansion}
 \prod_aJ_{Q'_a}^{V_a}
         =\chi_{Q'}\prod_{F\in\mathcal D}(1-\chi_F).
\end{equation}
If $\varnothing\in\mathcal D$, there are $a$ and $S\in\mathcal S$
such that $S\subseteq Q'$, $S\cap Q'_a\ne\varnothing$, and
$S\cap V_a=\varnothing$. Conditional on $\chi_{Q'}=1$, every atom
of $S$ is selected, so $J_x^{V_a}=0$ for every $x\in S\cap Q'_a$. Algebraically,
the factor $1-\chi_\varnothing=1-1$ in
\eqref{eq:alg-acceptance-expansion} is zero. The product in
\eqref{eq:alg-acceptance-expansion} is then identically zero,
and its contribution is omitted.

If $\varnothing\notin\mathcal D$, the same block count as above
gives \eqref{eq:alg-extension-count} for $\mathcal D$, with
$c_*=2^{q_1+2}q_1r(rb_2+3)$ in place of $c_0$.
For $k\in\{j_1-1,j_1\}$, let
\[
 A^{[k]}=\chi_{Q'}
       \sum_{\substack{\mathcal D'\subseteq\mathcal D\\|\mathcal D'|\le k}}
              (-1)^{|\mathcal D'|}\prod_{F\in\mathcal D'}\chi_F.
\]
Equations~\eqref{eq:alg-bonferroni-gap} and
\eqref{eq:alg-acceptance-expansion} imply
\[
 A^{[j_1]}\le\prod_aJ_{Q'_a}^{V_a}\le A^{[j_1-1]}.
\]
Conditional on $\chi_{Q'}=1$, the variables outside $Q'$ are
independent Bernoulli$(p)$ variables. Lemma~\ref{lem:alg-union-moment}
therefore gives
\[
 \mathbb E\bigl(A^{[j_1-1]}-A^{[j_1]}\mid\chi_{Q'}=1\bigr)
       \le\overline{a}_2(s,c_*,j_1)\overline{\eta}^{v+1}.
\]
In a term $\lambda A$ with $A=\prod_aJ_{Q'_a}^{V_a}$, the upper bound
is $\lambda A^{[j_1-1]}$ if $\lambda\ge0$ and $\lambda A^{[j_1]}$ if $\lambda<0$.
The lower bound uses the other sum in each case. This gives upper
and lower bounds for every term of $T^{[j-1]}$ and $T^{[j]}$.

We next bound the sum of the expectations of the factors that
multiply these errors. Under the original independent distribution
of the variables $\chi_x$, let $Y=\{x:\chi_x=1\}$ and define
\[
 K_\chi=\sum_{(a,B)\in\mathcal A}\prod_{x\in B}\chi_x
            =|\{(a,B)\in\mathcal A:B\subseteq Y\}|.
\]
Thus $K_\chi$ counts the indexed pairs $(a,B)$ for which every atom
of $B$ is sampled. Different indices with the same set $B$ are counted
separately. It is a polynomial of degree at most $s$ with nonnegative
coefficients, and $K\le K_\chi$ because $J_B^{V_a}\le\prod_{x\in B}\chi_x$.
The extension counts give every derivative expectation, including
its mean, at most $rc_0$. Proposition~\ref{prop:positive-moments}
therefore gives
\begin{equation}\label{eq:alg-inner-weight-sum}
 \sum_{\ell=0}^j\mathbb E\left(\binom{K_\chi}{\ell}\right)
 \le(j+1)\mathbb E\bigl((1+K_\chi)^j\bigr)
 \le(j+1)\bigl[2rc_0(sj+1)^s\bigr]^j.
\end{equation}
The initial factor $\chi_Q$ and the factors $1-\chi_x$ can only
decrease this bound. Independent factors $1-\mathsf w_x$ are
multiplied exactly, without truncation. Combining
\eqref{eq:alg-outer-error} and \eqref{eq:alg-inner-weight-sum},
gives upper and lower polynomials $T_\ell^+,T_\ell^-$ for each
initial term. If $\mathsf{F}=\sum_\ell e_\ell T_\ell$, where
$e_\ell\in\{-1,1\}$, define
\[
 \mathsf{F}^+=\sum_{e_\ell=1}T_\ell^+-\sum_{e_\ell=-1}T_\ell^-,\qquad
 \mathsf{F}^-=\sum_{e_\ell=1}T_\ell^--\sum_{e_\ell=-1}T_\ell^+.
\]
Then $\mathsf{F}^-\le \mathsf{F}\le \mathsf{F}^+$, and summing over the at most $2^t$
terms bounds $\mathbb E(\mathsf{F}^+-\mathsf{F}^-)$ by
$B_0\overline{\eta}^{v+1}$, where
\begin{equation}\label{eq:alg-explicit-polynomial-constants}
 \begin{aligned}
 B_0&=2^{t+2}\left(E_0+(j+1)\bigl[2rc_0(sj+1)^s\bigr]^j
                       \overline{a}_2(s,c_*,j_1)\right),\\
 \overline{s}(A,m,r,\delta_0)&=\max\{2,B_0^{1/4}\},\\
 \overline{a}_3(A,m,r,\delta_0)&=2^5r(t+j+j_1+1)^2.
 \end{aligned}
\end{equation}
Since $\delta_3(v+1)/4>A+4$ and $\Sigma\ge B_0^{1/4}$,
$B_0\overline{\eta}^{v+1}\le\Sigma^{-A}$.

\medskip
\noindent\emph{Proof of \textup{(2)}.}
The constructed polynomials have degree at most
$q_1+sj_1+2t\le\overline{a}_3$. All monomials are products of
variables that occur in the factors of $\mathsf{F}$, so no other variable is introduced.

\medskip
\noindent\emph{Proof of \textup{(3)}.}
The coefficients of $\mathsf{F}^-$ and $\mathsf{F}^+$ are computed
from the displayed formulas for $T^{[k]}$ and $A^{[k]}$. Each sum
ranges over subsets of at most $j$ or $j_1$ members of a family
containing at most $rN^r$ sets. Distributing these
finite sums, then replacing every positive power of a Boolean
variable by its first power, takes at most $N^{\overline{a}_3}$
arithmetic operations. For conditional expectations, we write
$X_{(x,1)}=\chi_x$ and $X_{(x,2)}=\mathsf w_x$. If $F$ indexes a monomial, let
$\mathcal A_0=\{X_y=\alpha_y\text{ for every }y\in A_0\}$ be a
specified event with positive probability. Independence gives
\begin{equation}\label{eq:alg-monomial-oracle}
 \mathbb E\left(\prod_{y\in F}X_y\,\middle|\,\mathcal A_0\right)
 =\prod_{y\in F\cap A_0}\alpha_y
       \prod_{y\in F\setminus A_0}\mathbb P(X_y=1).
\end{equation}
For a polynomial $\sum_F a_F\prod_{y\in F}X_y$, its conditional
expectation is therefore
\[
 \sum_F a_F\prod_{y\in F\cap A_0}\alpha_y
             \prod_{y\in F\setminus A_0}\mathbb P(X_y=1).
\]
Only variables appearing in the original factors occur in this expression.
\end{proof}

\subsection{Deterministic selection with local conditions}

We state the deterministic selection result that simultaneously
excludes specified events and bounds a weighted sum of nonnegative
functions. We then construct the polynomial upper bounds for event
indicators needed to apply this result to degree deviations.

For an event $B$ determined by Bernoulli variables indexed by a
finite set $V$, we use a polynomial $S_B$ such that, for every
assignment $\alpha\in\{0,1\}^V$,
\[
 S_B(\alpha)\ge\vb1_B(\alpha).
\]
Thus $S_B$ is nonnegative on all assignments and at least one on $B$.
Its coefficients need not be nonnegative.

The following is the weighted form of Harris's deterministic
local-lemma algorithm~\cite[Theorems 3.4 and 3.6(1)]{HarrisDetLLL}.
It uses a deterministic resampling table for the
Moser--Tardos algorithm~\cite{MoserTardos}.
Related derandomization results appear in~\cite{HarrisPolynomial,HarrisJuntas}.

\begin{proposition}\label{prop:alg-weighted-selection}
Let $r\ge2$ and $d\ge1$ be integers, let $a=2^{10}r^3$ and
$b=2^9r^3$, and suppose $N\ge\max\{4,r\}$ and $4\le\Sigma\le N^d$.
Let $(X_x:x\in V)$ be independent Bernoulli variables, and let
$\mathcal B$ be a family of events determined by these variables.
Suppose the following conditions hold.
\begin{enumerate}[label=\textup{(\arabic*)},ref=\arabic*]
\item $|V|\le N^d$ and $|\mathcal B|\le N^d$.
\item For every $B\in\mathcal B$, a set $T_B\subseteq V$ determines $B$, and
\[
 |\{B'\in\mathcal B\setminus\{B\}:T_{B'}\cap T_B\ne\varnothing\}|
       \le\Sigma^b.
\]
\item For each $B\in\mathcal B$, a supplied polynomial $S_B$ depends
only on $(X_x:x\in T_B)$ and satisfies, for every
$\alpha\in\{0,1\}^V$,
\[
 S_B(\alpha)\ge\vb1_B(\alpha),\qquad
 \mathbb E(S_B)\le\Sigma^{-a}.
\]
\item For $1\le j\le M\le N^d$, a supplied function
$\mathsf{F}_j:\{0,1\}^V\to\mathbb R_{\ge0}$ depends only on
$(X_x:x\in T_j)$ for a specified set $T_j\subseteq V$, and
\[
 |\{B\in\mathcal B:T_B\cap T_j\ne\varnothing\}|\le\Sigma^{b+1}.
\]
Each supplied weight $w_j$ is a nonnegative rational number.
\item The Bernoulli probabilities are rational, and their descriptions,
the weights, and the descriptions of the events and functions each
have at most $N^d$ bits. Each event $B$ can be tested in at most
$N^d$ bit operations. After any subset of the variables has been
fixed to values of positive joint probability, each conditional
expectation of $S_B$ or $\mathsf{F}_j$ can be computed in at most
$N^d$ bit operations.
\end{enumerate}
Let
\[
 \overline{a}_4(d,r)=2^{20}(d+1)^2(a+b+2^9)^2.
\]
A deterministic algorithm using at most $N^{\overline{a}_4(d,r)}$
bit operations finds an assignment $\alpha\in\{0,1\}^V$ outside
every event in $\mathcal B$ and satisfying
\begin{equation}\label{eq:alg-weighted-selection}
 \sum_{j=1}^M w_j\mathsf{F}_j(\alpha)
 \le(1+\Sigma^{-2^7})\sum_{j=1}^M w_j\mathbb E(\mathsf{F}_j).
\end{equation}
\end{proposition}

\begin{proof}
We apply the cited deterministic local-lemma theorem with
probability bound $p=\Sigma^{-a}$, dependency bound
$d_0=\Sigma^b+1$, and exponent $1/2$ in its hypothesis
$epd_0^{1+1/2}<1$. Here this hypothesis follows from
\[
 e\Sigma^{-a}(\Sigma^b+1)^{3/2}
 \le e2^{3/2}\Sigma^{-a+3b/2}<1.
\]
For each function $\mathsf{F}_j$, the proof of the cited theorem bounds the
sum of the weights of the graphs associated with $\mathsf{F}_j$ by
\[
 \exp(e\Sigma^{-a}\Sigma^{b+1})
 \le1+2e\Sigma^{-a+b+1}\le1+\Sigma^{-2^8}.
\]
The first bound is the estimate in the proof of
\cite[Theorem 3.6(1)]{HarrisDetLLL}.
With accuracy $\Sigma^{-2^8}$ in
\cite[Theorem 3.4]{HarrisDetLLL}, the two multiplicative bounds satisfy
$(1+\Sigma^{-2^8})(1+\Sigma^{-2^8})\le1+\Sigma^{-2^7}$. Since
$\Sigma^{2^8}\le N^{2^8d}$, its accuracy requirement is bounded
by a fixed power of $N$.

We justify its application to the stated functions.
Let $(X_{x,k}:k\ge1)$ be independent copies of each input variable
$X_x$, with all these copies mutually independent. When a variable
is resampled, the algorithm replaces it by the next unused copy in
its row. This array is the resampling table in the cited theorem.
Each directed acyclic graph $G$ in that proof has a final vertex
$0$ associated with $\mathsf{F}_j$; every other vertex $v$ is labelled by a
bad event $B_v$. Vertices whose functions involve a common variable
are ordered by the directed edges. For a variable $x$ occurring at
$v$, its table index is
\[
 k_G(x,v)=1+|\{u:u\text{ precedes }v\text{ and its function involves }x\}|.
\]
Let $X_{G,v}$ be the vector of entries $X_{x,k_G(x,v)}$ used at $v$.
The associated nonnegative random variable is
\[
 \mathsf V_G=\mathsf{F}_j(X_{G,0})\prod_{v\in V(G)\setminus\{0\}}
                              S_{B_v}(X_{G,v}),
 \qquad
 \mathbb E(\mathsf V_G)=\mathbb E(\mathsf{F}_j)
       \prod_{v\in V(G)\setminus\{0\}}\mathbb E(S_{B_v}).
\]
The vectors $X_{G,v}$ use disjoint table entries, so their
independence proves the expectation identity. If the resampling tests
recorded by $G$ occur, each corresponding factor $S_{B_v}$ is at
least one. Consequently, on the event that all tests represented by $G$
occur, $\mathsf{F}_j(X_{G,0})\le \mathsf V_G$. Summing these bounds and fixing the table entries
by conditional expectations is exactly the proof of the cited
theorem: at a vertex labelled by $B_v$, the factor $S_{B_v}$
is nonnegative everywhere and is at least one when $B_v$ occurs.
Its expectation is bounded by $\Sigma^{-a}$, as required in the
graph-weight estimate. This proves \eqref{eq:alg-weighted-selection} for
nonnegative functions as well as indicators.

We record a numerical exponent for this use of the algorithm.
In the notation of the cited proof, we take the exponent parameter
to be $1/4$ and set $p_* =\Sigma^{-a}$ and $q_* =p_*^{3/4}$. The symmetric local-lemma estimates give
\[
 e q_*(\Sigma^b+1)<1,\qquad
 e q_*\Sigma^{b+1}<\log2,\qquad
 W_{1/4}\le8N^d.
\]
Here $W_{1/4}$ is the sum of the graph weights defined before
\cite[Lemma 3.1]{HarrisDetLLL}, including graphs whose final vertex is associated with one of
$\mathsf{F}_1,\ldots,\mathsf{F}_M$. Let $\theta=N^{-2^{11}d}$ be the threshold used to
define the collections $\mathcal F_\theta$ and
$\mathcal F^{\rm low}_\theta$ in that proof.
The proof of \cite[Proposition 3.3]{HarrisDetLLL} gives
\[
 w(\mathcal F^{\rm low}_\theta)
 \le\theta^{1/4}W_{1/4}\le\tfrac12\Sigma^{-2^8},
 \qquad
 |\mathcal F_\theta|
 \le W_{1/4}(\theta^{-2}+\theta^{-1})\le N^{2^{13}d}.
\]
Since $q_*<e^{-1}$, each graph in this collection has at most
$1+2\log(1/\theta)\le N^{2^4d}$ vertices. The procedure in
\cite[Appendix B, Algorithm 2]{HarrisDetLLL} constructs this
collection in at most $N^{2^4d}$ iterations. Each iteration tests at most
$N^{2^{14}d}$ pairs. A graph is represented by its labels and
adjacency matrix; occurrences of the same label are canonically
ordered by their directed edges. Membership tests and comparison
with the stored graphs therefore use at most $N^{2^{14}d}$ bit
operations per pair. Thus enumeration uses at most $N^{2^{16}d}$
bit operations. The resampling table has at most $N^{2^5d}$ entries.
For each entry and each of its two possible values, we compute the
conditional expectations of the at most $N^{2^{13}d}$ products $\mathsf V_G$.
Each product requires at most $N^{2^4d}$ of the conditional
expectations supplied in the hypotheses. Even adding rational terms
without cancelling common factors gives at most $N^{2^{14}d}$ bits
per intermediate number. Adding or multiplying two integers with at most
$N^{2^{14}d}$ bits requires at most $2N^{2^{14}d}$ or $4N^{2^{15}d}$ bit
operations, respectively. These bounds and the preceding counts give at
most $N^{2^{16}d}$ bit operations for the table calculation. Testing and performing the
subsequent resamplings obeys the same bound. Both bounds, their
sum, and the input calculations are less than $N^{\overline{a}_4(d,r)}$
for the value specified in the statement.

\end{proof}

\begin{lemma}\label{lem:alg-moment-estimator}
Let $d,m\ge1$ be integers and $A_1>0$. Suppose
$F=\sum_{j=1}^{M}\mathsf{F}_j$, where $M\le\Sigma^d$ and the $\mathsf{F}_j$
are indicators from Lemma~\ref{lem:alg-local-polynomial-bounds}.
Assume $\Sigma^{-d}\le\Delta\le\Sigma^d$ and
\begin{equation}\label{eq:alg-normalized-moment}
 \frac{\mathbb E\bigl(|F-\mathbb E(F)|^{2m}\bigr)}{\Delta^{2m}}
 \le\Sigma^{-A_1-2}.
\end{equation}
Let
\[
 b_m=4m(d+1),\qquad A_2=A_1+b_m+2md+8,
 \qquad A_3=A_2+d+4.
\]
If $\Sigma\ge\max\{2^{8m+8},\overline{s}(A_3,2m,r,\delta_0)\}$,
the following statements hold.
\begin{enumerate}[label=\textup{(\arabic*)},ref=\arabic*]
\item There is a computable rational number $\widetilde\mu$ with
\[
 |\widetilde\mu-\mathbb E(F)|\le\Sigma^{-A_2}\le\Delta/8.
\]
\item There is a computable polynomial $S^+$ satisfying
\[
 S^+\ge\vb1_{\{|F-\widetilde\mu|>3\Delta/4\}},\qquad
 \mathbb E(S^+)\le\Sigma^{-A_1}.
\]
\item If $|F-\widetilde\mu|\le3\Delta/4$, then
$|F-\mathbb E(F)|\le\Delta$.
\end{enumerate}
\end{lemma}

\begin{proof}[Proof of \textup{(1)}]
For each $\mathsf{F}_j$, let $\mathsf{F}_j^-,\mathsf{F}_j^+$ be the polynomials in
\eqref{eq:alg-sandwich} with exponent $A_3$. Then
\[
 \sum_{j=1}^M\mathbb E(\mathsf{F}_j^-)
 \le\mathbb E(F)\le\sum_{j=1}^M\mathbb E(\mathsf{F}_j^+),\qquad
 \sum_{j=1}^M\mathbb E(\mathsf{F}_j^+-\mathsf{F}_j^-)
       \le\Sigma^{d-A_3}=\Sigma^{-A_2-4}.
\]
Let $\widetilde\mu$ be a rational approximation to the midpoint of
this interval, with error at most $\Sigma^{-A_2}/2$. It satisfies
\[
 |\widetilde\mu-\mathbb E(F)|
       \le\tfrac12\Sigma^{-A_2-4}+\tfrac12\Sigma^{-A_2}
       \le\Sigma^{-A_2}\le\Delta/8.
\]
\medskip
\noindent\emph{Proof of \textup{(2)}.}
We define
\begin{equation*}
 S=\left(\frac{F-\widetilde\mu}{3\Delta/4}\right)^{2m}
 =\left(\frac4{3\Delta}\right)^{2m}
  \sum_{k=0}^{2m}\binom{2m}k(-\widetilde\mu)^{2m-k}
        \sum_{j_1,\ldots,j_k=1}^M \mathsf{F}_{j_1}\cdots \mathsf{F}_{j_k}.
\end{equation*}
The inner sum is one when $k=0$. We index its terms by
$\ell=(k,j_1,\ldots,j_k)$ and write
\[
 B_\ell=\prod_{u=1}^k\mathsf F_{j_u},\qquad
 c_\ell=\left(\frac4{3\Delta}\right)^{2m}
                \binom{2m}k(-\widetilde\mu)^{2m-k}.
\]
Thus $S=\sum_\ell c_\ell B_\ell$, each $B_\ell$ is a product of at
most $2m$ of the original indicators, and
\[
 \sum_\ell|c_\ell|
 \le\left(\frac{4(M+|\widetilde\mu|)}{3\Delta}\right)^{2m}
 \le(4\Sigma^{2d})^{2m}\le\Sigma^{b_m}.
\]
The elementary inequality
$|x+y|^{2m}\le2^{2m-1}(|x|^{2m}+|y|^{2m})$ and
\eqref{eq:alg-normalized-moment} give
\[
 \mathbb E(S)
 \le(4/3)^{2m}2^{2m-1}
       \bigl(\Sigma^{-A_1-2}
                   +\Sigma^{-2m(A_2-d)}\bigr)
 \le\tfrac12\Sigma^{-A_1}.
\]
For each $B_\ell$, let $B_\ell^-,B_\ell^+$ be the bounds in
Lemma~\ref{lem:alg-local-polynomial-bounds} with exponent $A_3$,
and let
\begin{equation}\label{eq:alg-signed-upper}
 S^+=\sum_{c_\ell\ge0}c_\ell B_\ell^+
                    +\sum_{c_\ell<0}c_\ell B_\ell^-.
\end{equation}
Then $S^+\ge S\ge0$ pointwise, and $S^+\ge1$ whenever
$|F-\widetilde\mu|>3\Delta/4$. Moreover,
\[
 \mathbb E(S^+)\le\mathbb E(S)+\Sigma^{b_m-A_3}
                         \le\Sigma^{-A_1}.
\]
Formula~\eqref{eq:alg-monomial-oracle} computes its conditional
expectations.

\medskip
\noindent\emph{Proof of \textup{(3)}.}
By part~\textup{(1)},
$|F-\mathbb E(F)|\le3\Delta/4+\Delta/8<\Delta$ on the
complement of the bad event.
\end{proof}

\subsection{The modified process and its local conditions}

We now define the varying sampling probabilities and the deterministic
reference parameters for the algorithm. The lemmas in this subsection
show that the modified focus and construction conditions preserve the
required degree and intersection bounds through the stopping time.

We retain the sets $X^t,Z^t,Y^t,I^t,J^t,\mathsf W^t$ and the
update of Subsection~\ref{subsec:actual-process}, which retains the
distinct inclusion-minimal residual blocks.
The sampling probabilities are those in \eqref{eq:alg-sampling}.
Let $a=2^{10}r^3$, $b=2^9r^3$, and
\begin{equation}\label{eq:alg-orders}
 m_0=\left\lceil\frac{4a+4b+128}{\delta_3}\right\rceil,
 \qquad
 A_0=2^7(m_0+1)(r+1)[16r(r-1)+1]+2^4a+2^{10}.
\end{equation}

The value $\overline{\eta}_t$ in \eqref{eq:alg-sampling} is not the fixed
$\eta=(\log D)^{-\delta_3/8}$ in Subsection~\ref{subsec:actual-process}.
For $D>\exp(4^{8/\delta_3})$ and $t<\overline{T}$, we have
$\Sigma_t>\sqrt{\log D}$ and hence
\begin{equation*}
 \Sigma_t^{-\delta_3/2}\le\overline{\eta}_t
 <(\log D)^{-\delta_3/8}<\frac14
.
\end{equation*}
The one-round results cited below assume a numerical value
$0<\eta\le1/4$; in their formulas we substitute $\eta=\overline{\eta}_t$.
We do not change the fixed value used in the proof of
Theorem~\ref{thm:main}.

With $\overline{\omega}_t$ defined in
\eqref{eq:alg-survival-reference}, the reference parameters satisfy
\begin{align}
 d_{i,t+1}&=\overline{\omega}_t^{i-1}
       \sum_{j=i}^r\binom{j-1}{i-1}z_t^{j-i}d_{j,t},
       \label{eq:alg-degree-reference}\\
 \pi_{t+1}&=\overline{\omega}_t\pi_t,\qquad
 \iota_{t+1}=\iota_t+\pi_tz_t,\qquad
 \pi_0=1,\quad\iota_0=0.
       \notag
\end{align}
The two deterministic stopping indices are
\begin{equation}\label{eq:alg-stopping}
 \overline{T}_{\mathrm{init}}=\min\{t:f_0(\iota_t)\ge2^{4r}\},\qquad
 \overline{T}=\min\{t:\Sigma_t\le\sqrt{\log D}\},
\end{equation}
where $f_0,g_i,g$ are the functions preceding
Lemma~\ref{lem:general-generating}. The initial phase here ends at
$\overline{T}_{\mathrm{init}}$, in place of $T_{\mathrm{init}}$ in
\eqref{eq:general-initial-period-length}.

Let $b_1,b_2,b_3,c_1,c_3$ be the constants prescribed in
Lemmas~\ref{lem:good-focus-existence},
\ref{lem:general-error-closure} and
\ref{lem:general-stopping-cost}. We define
\begin{equation}\label{eq:alg-large-constant}
 \overline{c}_2=\exp\left(
 \bigl[2^{16}r(b_1+b_2+b_3+c_1+c+\delta_3^{-1}+1)\bigr]^{2^{12}r^4}
                         \right).
\end{equation}
Let
\[
 \overline{c}_1=\exp(\overline{c}_2^{\,2}).
\]
These constants replace $c_1$ and $c_2$, respectively, in the
degree-error estimates. Both depend only on $r,\delta_0,c$.

For each intersection family, we define its parameters by
\begin{equation}\label{eq:alg-intersection-types}
 (\mathcal B,m,\zeta)=
 \begin{cases}
 (N_i(A),i-|A|,\delta_5(i)),&2\le|A|<i\le r,\\
 (N_i(x,y),i-1,\delta_4(i)),&x\ne y, 2\le i\le r.
 \end{cases}
\end{equation}
For $0\le t<\overline{T}_{\mathrm{init}}$, the initial bounds and recurrence are
\begin{equation}\label{eq:alg-initial-envelope}
 |\mathcal B|\le \overline{k}_tD^{m-\zeta-\delta_0/4},\qquad
 \overline{k}_0=c,\qquad \overline{k}_{t+1}=(1+\overline{\eta}_t)\overline{k}_t.
\end{equation}
This replaces \eqref{eq:general-initial-envelope} only in the
algorithmic process. For $t\ge \overline{T}_{\mathrm{init}}$ we use
\eqref{eq:general-codegree-barriers}--\eqref{eq:general-common-barriers}.

For each family $\mathcal B$ in \eqref{eq:alg-intersection-types},
let $\mathcal Y(\mathcal B)$ be the following indexed family of
polynomials from \eqref{eq:general-birth-witness-polynomials}:
\[
 \mathcal Y(N_i(A))=\bigl(\mathsf Y_{i,A;k}\bigr)_{i<k\le r},
 \qquad
 \mathcal Y(N_i(x,y))=
 \bigl(\mathsf Y_{i,x,y;k,l,h}\bigr)_{
   \substack{i\le k,l\le r,\ (k,l)\ne(i,i)\\
             0\le h\le\min\{k-i,l-i\}}}.
\]
An empty indexed family has sum zero. For
$\mathsf Y\in\mathcal Y(\mathcal B)$, let $s,m,\iota$ have the
values in \eqref{eq:birth-polynomial-parameters}, and let
$\zeta=\delta_5(i)$ when $\mathcal B=N_i(A)$ and
$\zeta=\delta_4(i)$ when $\mathcal B=N_i(x,y)$. We define
\begin{equation*}
 h_t=\begin{cases}D,&t<\overline{T}_{\mathrm{init}},\\\Sigma_t,&t\ge \overline{T}_{\mathrm{init}},\end{cases}
 \qquad
 M_t(\mathsf Y)=
 \begin{cases}
 \overline{k}_tD^{m-\zeta-\delta_0/4-\delta_1},&t<\overline{T}_{\mathrm{init}},\\
 \iota,&t\ge \overline{T}_{\mathrm{init}},
 \end{cases}
\end{equation*}
Here $\iota=\Sigma_t^{i-|A|-\delta_5(k)}$ for
$\mathsf Y_{i,A;k}$ and
$\iota=\Sigma_t^{i-1-\delta_4(i+1)}$ for
$\mathsf Y_{i,x,y;k,l,h}$, as in
\eqref{eq:birth-polynomial-parameters}. These values are distinct
from the time parameter $\iota_t$.

In a single round below, we write
$\overline{\eta}=\overline{\eta}_t$ and
$\overline{\omega}=\overline{\omega}_t$ when the round index is omitted.
The state variables and the quantities defined from them by the same
formulas as in the earlier process retain their original notation.

\begin{definition}\label{def:alg-good-round}
A set $Z\subseteq X^t\setminus\mathsf W^t$ is a \emph{good focus}
in round $t$ of the algorithm if it satisfies
\eqref{gf:small-counts}--\eqref{gf:coverage} of
Definition~\ref{def:good-focus}, with $\overline{T}_{\mathrm{init}}$ in place of
$T_{\mathrm{init}}$, and, for $|A|\le s$, the derivative bounds
\begin{equation}\label{eq:alg-focus-derivatives}
 (\partial_A\mathsf Y)(Z)
 \le \overline{c}_2 M_t(\mathsf Y)h_t^{(s-|A|)\delta_2}
.
\end{equation}
These bounds replace \eqref{gf:new-counts}. Every subsequent
reference to a good focus in this section means a set satisfying
these modified conditions.

For this fixed good focus, let $s_x^+$ be the upper estimate for $s_x$
computed by Lemma~\ref{lem:alg-local-polynomial-bounds}; it satisfies
\begin{equation}\label{eq:alg-upper-survival}
 s_x\le s_x^+\le s_x+\Sigma_t^{-A_0}.
\end{equation}
For $x\in X^t\setminus\mathsf W^t$, let $\mathsf w_x$ be independent
Bernoulli variables, independent also of all $\chi_y$, with
\begin{equation}\label{eq:alg-retention}
 \mathbb P(\mathsf w_x=0)=\overline{\omega}_t/s_x^+.
\end{equation}
The construction is required to satisfy the following inequalities;
a local failure is the event that one of these inequalities fails.
\begin{enumerate}
\item For $F\in\{F_i(v),\widehat F_i(v)\}$ under the law
$\mathbb P_{Z,v}^0$ of \eqref{eq:degree-conditioned-law},
\begin{equation}\label{eq:alg-degree-test}
 |F-\widetilde\mu_F|\le3\Delta_i/4,
 \qquad \Delta_i=\Sigma_t^{i-1-2\delta_3}.
\end{equation}
Here Lemma~\ref{lem:alg-moment-estimator} computes a rational
number $\widetilde\mu_F$ satisfying
\[
 \bigl|\widetilde\mu_F-\mathbb E_{Z,v}^0(F)\bigr|
           \le\min\{\Delta_i/8,\Sigma_t^{-A_0}\}.
\]
Consequently \eqref{eq:alg-degree-test} implies
$|F-\mathbb E_{Z,v}^0(F)|\le\Delta_i$.
\item If $t\ge \overline{T}_{\mathrm{init}}$ and the family in
\eqref{eq:alg-intersection-types} satisfies
$d_{2,t}\ge\Sigma_t/2$ and
$|\mathcal B|\ge\frac12\Sigma_t^{m-\zeta}$, as in
\eqref{eq:coverage-family-activation}, then, with $q_0=1-\overline{\omega}_t$,
\begin{equation}\label{eq:alg-coverage-test}
 \sum_{B\in\mathcal B}\mathsf Z_B^{\rm cov}
                    \ge(m-\delta_1/8)q_0|\mathcal B|.
\end{equation}
\item For each family $\mathcal B$ in
\eqref{eq:alg-intersection-types}, we require
\begin{equation}\label{eq:alg-birth-test}
 \sum_{\mathsf Y\in\mathcal Y(\mathcal B)}\mathsf Y(\chi)
 \le\begin{cases}
 \overline{\eta}_t \overline{k}_tD^{m-\zeta-\delta_0/4},&t<\overline{T}_{\mathrm{init}},\\
 (\delta_1/64)\overline{\eta}_t\Sigma_t^{m-\zeta},&t\ge \overline{T}_{\mathrm{init}}.
 \end{cases}
\end{equation}
Each polynomial is evaluated at the sampling variables
$\chi_x=\vb1_{\{x\in Y^t\}}$. By
Lemma~\ref{lem:general-births}, if the specified atoms $A$, or
$x,y$, survive, then the corresponding next family
$\mathcal B'=N_i'(A)$, or $N_i'(x,y)$, satisfies
\[
 |\mathcal B'\setminus\mathcal B|
       \le\sum_{\mathsf Y\in\mathcal Y(\mathcal B)}\mathsf Y(\chi).
\]
\end{enumerate}
\end{definition}

\begin{lemma}\label{lem:alg-local-estimators}
Suppose the current degree errors have the bounds
\begin{equation}\label{eq:alg-degree-invariants}
 E_i\le \overline{c}_1^{2-i}\Sigma^{i-1-\delta_3},\qquad
 \widehat E_i\le(2\overline{c}_1)^{r-i+1}\Sigma^{i-1-\delta_3}.
\end{equation}
Suppose the appropriate
intersection bounds in \eqref{eq:alg-initial-envelope} or
\eqref{eq:general-codegree-barriers}--\eqref{eq:general-common-barriers}
hold. During the initial phase also suppose
\[
 e^{-2^{16r^2}}D\le\Sigma\le2^{4r^2}D,
 \qquad \overline{k}_t\le c\exp(2^{16r^2}).
\]
Then each local failure in Definition~\ref{def:alg-good-round}
has a polynomial $S_B\ge\vb1_B$ whose degree is bounded
independently of $n$. After grouping the events as below, these
polynomials satisfy
\begin{equation}\label{eq:alg-local-allowances}
 \mathbb E(S_B)\le\Sigma^{-a},\qquad
 d_{\rm dep}\le\Sigma^b.
\end{equation}
The sets of variables that determine them satisfy
Lemma~\ref{lem:conditioning-locality} with
\eqref{eq:indicator-locality-parameters}. Their conditional
expectations are computable using at most $N^{\overline{a}_5}$ arithmetic operations, where $\overline{a}_5=2^4(r+a+b+m_0+A_0+\overline{a}_3(A_0,2m_0,r,\delta_0)+1)^2$.
\end{lemma}

\begin{proof}[Proof for the focus tests]
We give the moment estimates for all the types of focus tests.
In each application of Proposition~\ref{prop:alg-fixed-moments}, the
power is $2m_0$ from \eqref{eq:alg-orders}. We impose the condition
\begin{equation}\label{eq:alg-moment-constant-threshold}
 \Sigma^{m_0\delta_3-2a-32}\ge
 \bigl(\overline{a}_1(r,m_0)+2(2m_0)!4^{2m_0}\bigr)
       \bigl[4\overline{c}_2^5(2rm_0+1)^r\bigr]^{2m_0}.
\end{equation}
The exponent on the left is positive by \eqref{eq:alg-orders}.
Since $\Sigma>\sqrt{\log D}$ before stopping, this condition follows from
\[
 (\log\overline{C})^{(m_0\delta_3-2a-32)/2}\ge
 \bigl(\overline{a}_1(r,m_0)+2(2m_0)!4^{2m_0}\bigr)
       \bigl[4\overline{c}_2^5(2rm_0+1)^r\bigr]^{2m_0}.
\]
We impose this explicit inequality on $\overline{C}$.

For $d_i(v,Z)$, the proof of Lemma~\ref{lem:general-focus-estimates}
gives $\mu\le2\Sigma^{(i-1)\delta_2}$ and $E_*\le1$. With
$\Delta=\Sigma^{(i-1)\delta_2-4\delta_3}$, it follows that
\begin{equation}\label{eq:alg-focus-degree-ratio}
 \frac{\mu E_*}{\Delta^2}+\frac{E_*^2}{\Delta^2}
                \le 4\Sigma^{-\delta_2+8\delta_3}.
\end{equation}
For $M_{j,k}(v),\widehat M_{j,k}(v)$, write $s=j-k$ and
$B=\Sigma^{k-1+s\delta_2}$. The same proof gives
$\mu\le2^{r+1}B$ and $E_*\le2^rB\Sigma^{-\delta_5(j)-\delta_2}$.
For $\Delta=B\Sigma^{-4\delta_3}$ the ratio in
\eqref{eq:alg-focus-degree-ratio} is
at most $2^{2r+2}\Sigma^{-\delta_5(j)-\delta_2+8\delta_3}$.
Both $\delta_2-8\delta_3$ and
$\delta_5(j)+\delta_2-8\delta_3$ are at least $\delta_3$.

For the polynomials in \eqref{eq:degree-occurrence-polynomial}
and \eqref{eq:degree-direct-occurrence-polynomial}, the bounds
\eqref{eq:degree-occurrence-derivatives} and the derivative calculation
in Lemma~\ref{lem:degree-indicator-read}(2) apply. If $B$ is the
corresponding upper bound for the mean, then $\mu\le B$ and $E_*\le\overline{c}_2 B\Sigma^{-\delta_2}$.
Proposition~\ref{prop:alg-fixed-moments}(2), with deviation $B$, replaces the applications of
Corollary~\ref{cor:polynomial-envelope}(2) in that proof. Equations
\eqref{eq:degree-J-variable-count} and
\eqref{eq:degree-read-union-tail} then imply
\eqref{eq:general-degree-read}.

For $|C_x\cap Z|$, the expectation and deviation are those in
\eqref{eq:coverage-candidate-bernstein}, and
$q|C_x|\ge\Sigma^{\delta_2}/4$. Its read bound is one, so
\eqref{eq:alg-read-moment} gives a ratio
at most $\overline{c}_2\Sigma^{-\delta_2}$.
We next treat the polynomials used to bound the number of
coverage indicators depending on one variable. Let
$U=X\setminus\mathsf W$, and let $(\xi_x:x\in U)$ be independent
Bernoulli$(q)$ variables, where $q=\Sigma^{-1+\delta_2}$ and
$Z=\{x\in U:\xi_x=1\}$. With $\mathcal B,F_0,\ell,C_B$ as in
\eqref{eq:coverage-candidate-definitions}, for
$u\in X\setminus(\mathsf W\cup F_0)$ and $w\in U$ let
\[
 m_u(w)=|\{B\in\mathcal B:u\notin B,\ w\in C_B\}|.
\]
For $2\le s\le r$, the required polynomial is
\[
 \mathsf Q_{u,s}(\xi)
 =\sum_{f\in N_s(u)\cap\binom{U}{s-1}}
       \left(\sum_{w\in f}m_u(w)\right)\prod_{x\in f}\xi_x.
\]
The terms with zero coefficient may be omitted, so this is the same
polynomial as in the proof of Lemma~\ref{lem:focus-coverage}.
For $s=2$, we consider the normalized polynomial
\[
 \frac{\mathsf Q_{u,2}}{\ell+r|\mathcal B|\Sigma^{2\delta_2-\delta_1}}.
\]
For $3\le s\le r$, we consider
\[
 \frac{\mathsf Q_{u,s}}{(s-1)\ell}.
\]
Their coefficients are at most one. The first has degree one and
mean at most $\Sigma^{-\delta_2}$. The second has degree $s-1$;
every proper derivative, including the polynomial itself, has
expectation at most $\Sigma^{-(\delta_1-r\delta_2)}$, by the
calculation preceding \eqref{eq:coverage-weighted-polynomial-probability}.

Let $\Phi(\xi)=\sum_{F\in\binom{U}{k}}z_F\prod_{x\in F}\xi_x$
be one of these normalized polynomials, where $k=1$ or $k=s-1$,
respectively. Let
\[
 \ell_0=\left\lceil\frac{2a+32}{\delta_2}\right\rceil,
 \qquad \lambda_k=\sum_{j=0}^{k-1}\binom{k(\ell_0-1)}j.
\]
Let $\mathcal A_{\ell_0}$ consist of the ordered tuples
$\mathbf F=(F_1,\ldots,F_{\ell_0})$ of $k$-subsets of $U$ for which
$F_j\nsubseteq\bigcup_{u<j}F_u$ for every $1\le j\le\ell_0$.
We define
\[
 \mathsf V_{\ell_0}(\xi)
 =\sum_{\mathbf F\in\mathcal A_{\ell_0}}
       \prod_{j=1}^{\ell_0}\left(z_{F_j}\prod_{x\in F_j}\xi_x\right),
 \qquad S_B=\frac{\mathsf V_{\ell_0}}{\lambda_k^{\ell_0}}.
\]
Here $B=\{\Phi>b_2\}$. The polynomial $\mathsf V_{\ell_0}$ is the
polynomial denoted $W_m$ in the proof of
Corollary~\ref{cor:polynomial-envelope}(3), with $m=\ell_0$ and degree $k$.
It is unrelated to the set $\mathsf W^t$ of atoms excluded from sampling.
The value of $b_2$ in Definition~\ref{def:good-focus} satisfies
$b_2\ge(r\ell_0)^r\ge(k\ell_0)^k$.
The proof of that corollary therefore gives $S_B\ge\vb1_B$ and
$\mathbb E(S_B)\le\Sigma^{-\delta_2\ell_0}\le\Sigma^{-2a-32}$.
For each completion-count polynomial in \eqref{gf:small-counts},
we define $\mathsf V_{\ell_0}$ and $S_B$ by the same formulas,
with $k$ equal to its degree. The derivative bounds
\eqref{eq:focus-completion-derivative-mean} give
$S_B\ge\vb1_B$ and $\mathbb E(S_B)\le\Sigma^{-2a-32}$.
Thus each of these focus conditions has a polynomial upper bound
for its failure indicator, of degree independent of $n$.

For a derivative of $\mathsf Y$ of remaining degree $j>0$,
\eqref{eq:birth-derivative-mean} gives the upper bound
$B=(\overline{c}_2/4)M_t(\mathsf Y)h_t^{j\delta_2}$ for its mean, and every further derivative
has expectation at most $\overline{c}_2 Bh_t^{-\delta_2}$.
For $t<\overline{T}_{\mathrm{init}}$, the same assertion follows by using
\eqref{eq:alg-initial-envelope} in the enumeration in
Lemma~\ref{lem:focus-birth-derivatives}: one prescribed proper
codegree or common link has exponent at least $\zeta+\delta_1$
and the factor $D^{-\delta_0/4}$ in
\eqref{eq:alg-initial-envelope} multiplies its upper bound.
All remaining counts use $d_i(v,X)\le2^{4r^3+1}D^{i-1}$ or the trivial
completion bound. This is exactly the division into cases in the
proofs of Lemmas~\ref{lem:focus-birth-derivatives} and
\ref{lem:initial-focus-estimates}. Their coefficients are bounded by
\[
 \bigl[2^{16}r(b_1+b_2+b_3+1)\exp(2^{18r^2})\bigr]^{2^{10}r^2}
                 \le\overline{c}_2/2,
\]
by \eqref{eq:alg-large-constant}. Equation
\eqref{eq:alg-polynomial-moment} proves
\eqref{eq:alg-focus-derivatives}; the highest derivatives are
deterministic.

\emph{Proof for the construction tests.}
For $F_i(v),\widehat F_i(v)$, Lemma~\ref{lem:degree-indicator-identities}(1)
and \eqref{eq:general-degree-read} give the read bound
$k\le b_3\Sigma^{i-1-\delta_1+3r\delta_2}$.
Equation~\eqref{eq:general-oriented-mean-upper} gives
$\mu\le2^{r+2}\Sigma^{i-1}$. Hence, with $\Delta_i$ in
\eqref{eq:alg-degree-test},
\begin{equation}\label{eq:alg-degree-moment-ratio}
 \frac{k\mu}{\Delta_i^2}+\frac{k^2}{\Delta_i^2}
 \le(2^{r+3}b_3+2b_3^2)\Sigma^{-\kappa_1},\qquad
 \kappa_1=\delta_1-3r\delta_2-4\delta_3>\delta_3.
\end{equation}
For $F=\sum_{B\in\mathcal B}\mathsf Z_B^{\rm cov}$, the sets of variables
in Definition~\ref{def:coverage-random-variables} and
\eqref{eq:coverage-focus-read} give
\[
 k\le rb_1b_3^2|\mathcal B|\Sigma^{2\delta_2-\delta_1},
 \qquad \mu\le|\mathcal B|.
\]
Lemma~\ref{lem:general-coverage}, specifically
\eqref{eq:general-coverage-mean}, places the mean above the
threshold in \eqref{eq:alg-coverage-test} by at least
$(\delta_1/32)q_0|\mathcal B|$ when
$\Sigma^{\delta_3}\ge128\overline{c}_2/\delta_1$.
We ensure this by imposing
$\log\overline{C}\ge(128\overline{c}_2/\delta_1)^{2/\delta_3}$.
Here $\overline{\eta}/8\le q_0\le2\overline{\eta}$, by the argument for
\eqref{eq:general-graph-hazard} and \eqref{eq:alg-survival-reference}.
For $\Delta=(\delta_1/64)q_0|\mathcal B|$,
\begin{equation}\label{eq:alg-coverage-moment-ratio}
 \frac{k\mu}{\Delta^2}+\frac{k^2}{\Delta^2}
 \le\overline{c}_2\Sigma^{-\kappa_2},\qquad
 \kappa_2=\delta_1-2\delta_2-\delta_3/2>\delta_3.
\end{equation}
The change in the survival reference contributes
at most $\overline{c}_2\Sigma^{-\delta_3}q_0|\mathcal B|$ to the mean comparison,
which is smaller than the displayed gap.

Proposition~\ref{prop:alg-fixed-moments}(1),
\eqref{eq:alg-orders}, \eqref{eq:alg-degree-moment-ratio}, and
\eqref{eq:alg-coverage-moment-ratio} imply
\[
 \mathbb E\bigl(|F-\mu|^{2m_0}\bigr)/\Delta^{2m_0}
                  \le\Sigma^{-2a-32}.
\]
We apply Lemma~\ref{lem:alg-moment-estimator} to these sums.
By Lemma~\ref{lem:conditioning-locality}(1), their number of
summands is at most $\Sigma^d$, where $d=r[16r(r-1)+1]$.
In the proof of that lemma, the sum of the absolute coefficients is
at most $\Sigma^{4m_0(d+1)}$. If each indicator product is replaced
by a polynomial with expected error at most $\Sigma^{-A_0}$,
the total error in \eqref{eq:alg-signed-upper} is therefore at most
\[
 \Sigma^{4m_0(d+1)-A_0}\le\Sigma^{-2a-32},
\]
by the definition of $A_0$ in \eqref{eq:alg-orders}.

For $\mathsf Y\in\mathcal Y(\mathcal B)$ of degree $s$, fixing a focus
satisfying \eqref{eq:alg-focus-derivatives} gives, for $0\le j\le s$,
\begin{equation*}
 E_j\le\overline{c}_2^2M_t(\mathsf Y)\overline{\eta}^{s-j}.
\end{equation*}
Indeed, independence gives
\[
 \mathbb E(\partial_A\mathsf Y(\chi))
   =p^{s-j}(\partial_A\mathsf Y)(Z)
   \le\overline{c}_2 M_t(\mathsf Y)\overline{\eta}^{s-j}
                          (h_t/\Sigma)^{(s-j)\delta_2},\qquad |A|=j.
\]
In the initial phase, $D/\Sigma\le\exp(2^{16r^2})$, and
$\exp(r\delta_2 2^{16r^2})\le\overline{c}_2$.
The ratios of these bounds to each assigned fraction of the
threshold in \eqref{eq:alg-birth-test} are
at most $\overline{c}_2^4h_t^{-\delta_1+\delta_3/4}$. Proposition~\ref{prop:positive-moments},
with power $2m_0$, therefore gives the polynomial
$S_B=(\mathsf Y/u)^{2m_0}\ge\vb1_{\{\mathsf Y>u\}}$
with expectation at most $\Sigma^{-2a-16}$,
where $u$ is the threshold in \eqref{eq:alg-birth-test} divided
by $2r^3$. The indices $k,l,h$ in
\eqref{eq:general-birth-witness-polynomials} give at most $2r^3$
polynomials for each family $\mathcal B$. In particular,
\[
 \mathbb E\left((\mathsf Y/u)^{2m_0}\right)
 \le\left(\sum_{j=0}^s\binom{2sm_0}j\frac{E_j}{u}\right)^{2m_0}
 \le\left[\overline{c}_2^4(2rm_0+1)^r
                     h_t^{-\delta_1+\delta_3/4}\right]^{2m_0}
 \le\Sigma^{-2a-16},
\]
where the final inequality follows from
\eqref{eq:alg-moment-constant-threshold}. Here
$h_t\ge2^{-4r^2}\Sigma$, so replacing $h_t$ by $\Sigma$ in the
preceding upper bound multiplies its bracket by at most $\overline{c}_2$.

Finally, the grouping in the proofs of
Lemmas~\ref{lem:good-focus-existence} and
\ref{lem:general-degree-round} associates each collection of
elementary failures with its specified atom $v$, set $A$, or pair
$x,y$. For a collection $\mathcal A$, we define its bad event and
its polynomial by
\[
 B_{\mathcal A}=\bigcup_{B\in\mathcal A}B,\qquad
 S_{B_{\mathcal A}}=\sum_{B\in\mathcal A}S_B.
\]
There are at most $\Sigma^b$ terms in this sum. It is nonnegative
and at least one on $B_{\mathcal A}$, and
\[
 \mathbb E(S_{B_{\mathcal A}})
   \le\Sigma^b\Sigma^{-2a-16}\le\Sigma^{-a}.
\] Lemma~\ref{lem:conditioning-locality}(2),(3), with
\eqref{eq:indicator-locality-parameters}, proves
\eqref{eq:alg-local-allowances}. Formula
\eqref{eq:alg-monomial-oracle} supplies the required conditional
expectations.
\end{proof}

\begin{lemma}\label{lem:alg-tracking}
With $\overline{c}_2$ defined in \eqref{eq:alg-large-constant},
the threshold $\overline{C}$ can be chosen so that every sequence of rounds satisfying
Definition~\ref{def:alg-good-round} has the following properties:
\begin{enumerate}[label=\textup{(\arabic*)},ref=\arabic*]
\item The retention probabilities in \eqref{eq:alg-retention} lie
in $[0,1]$, and \eqref{eq:alg-degree-invariants} holds.
\item Bounds \eqref{eq:alg-initial-envelope} hold through $\overline{T}_{\mathrm{init}}$.
For $t\ge \overline{T}_{\mathrm{init}}$, the bounds
\eqref{eq:general-codegree-barriers}--\eqref{eq:general-common-barriers}
hold, and $d_{2,t}\ge\Sigma_t/2$.
\item The reference trajectory satisfies the following bounds.
\begin{enumerate}[label=\textup{(3.\arabic*)},ref=3.\arabic*,leftmargin=3.3em]
\item\label{it:alg-initial-length}
$\overline{T}_{\mathrm{init}}\le2^{20r^2}D^{\delta_3/4}$.
\item\label{it:alg-initial-sampling}
$\sum_{t<\overline{T}_{\mathrm{init}}}\overline{\eta}_t\le2^{16r^2}$.
\item\label{it:alg-initial-scale} For $0\le t\le\overline{T}_{\mathrm{init}}$,
$e^{-2^{16r^2}}D\le\Sigma_t\le2^{4r^2}D$.
\item\label{it:alg-scale-ratios} For $\overline{T}_{\mathrm{init}}\le t<\overline{T}$,
$e^{-2\overline{\eta}_t}\le\Sigma_{t+1}/\Sigma_t\le e^{-\overline{\eta}_t/8}$.
\item\label{it:alg-horizon}
$\overline{T}\le\frac{2^{32r^2}}{\delta_3}D^{\delta_3/4}$.
\end{enumerate}
\end{enumerate}
\end{lemma}

\begin{proof}
We use induction on $t$, assuming that all asserted bounds hold
through the beginning of round $t$. The proofs of the three parts
are interdependent, so we indicate the assertion at each stage.

\medskip
\noindent\emph{Proof of \textup{(1)}.}
Under
\eqref{eq:alg-degree-invariants},
$\mathsf E_{\rm acc}\le2\overline{\eta}\Sigma^{-\delta_3}$.
For $x\in X$ and $y\in X\setminus\mathsf W$,
Lemma~\ref{lem:general-survival-estimates}(3) and
\eqref{eq:general-survival-log-reference} give
\begin{equation}\label{eq:alg-equalization-bounds}
 \begin{aligned}
 \overline{\omega}&\le s_x\le\overline{\omega}e^{3e_t},\\
 \overline{\omega}(1-8\Sigma^{-A_0})
 &\le\mathbb E\bigl(\psi_y(1-\mathsf w_y)\bigr)\le\overline{\omega}.
 \end{aligned}
\end{equation}
Here $s_x\ge1/8$ by \eqref{eq:general-survival-comparison}.
This proves the first assertion about probabilities.

When we apply \eqref{eq:general-marked-mixed} with survival reference
$\overline{\omega}$, each factor $\overline{\omega}/s_x$ is replaced by
$\overline{\omega}/s_x^+$. A product of at most $r$ such factors changes
by a relative amount at most $8r\Sigma^{-A_0}$. Repeating the subtraction in the
proof of Lemma~\ref{lem:general-degree-round}, using
Lemma~\ref{lem:general-cleanup} and
\eqref{eq:general-active-conditional-mean}--\eqref{eq:general-full-conditional-mean}, gives
\begin{align}
 E_i'&\le\overline{\omega}^{i-1}\sum_{j=i}^r\binom{j-1}{i-1}z^{j-i}E_j
 +\overline{c}_2\bigl(\mathsf E_{\rm acc}\mathcal B_i+\Sigma^{i-1-2\delta_3}\bigr),
 \label{eq:alg-active-error}\\
 \widehat E_i'&\le\overline{\omega}^{i-1}(1+\overline{c}_2\overline{\eta}\Sigma^{-\delta_3})
             \sum_{j=i}^r\binom{j-1}{i-1}z^{j-i}\widehat E_j
 +\overline{c}_2\Sigma^{i-1-\delta_3}\bigl(\overline{\eta}+\Sigma^{-\delta_3}\bigr).
 \label{eq:alg-full-error}
\end{align}
The additional full-degree term follows from
$s_x/\overline{\omega}\le e^{3e_t}$. The deviation is bounded by
\eqref{eq:alg-degree-test}; the conversion from the unconditional
mean under $\mathbb P_{Z,v}^0$ to the mean conditional on $\psi_v=1$
is the calculation in
\eqref{eq:general-degree-unconditioned-tail}--\eqref{eq:general-degree-conditioned-tail}.
Its tail bound is supplied here by
\eqref{eq:alg-degree-moment-ratio} and Markov's inequality. The
$\Sigma^{-A_0}$ terms are smaller than
$\Sigma^{i-1-2\delta_3}$.

\medskip
\noindent\emph{Proof of \textup{(3.1)--(3.3)}.}
The algebra of Lemma~\ref{lem:general-generating} applies exactly
to \eqref{eq:alg-degree-reference}. In particular,
\begin{equation}\label{eq:alg-generating}
 d_{i,t}=\pi_t^{i-1}g_i(\iota_t),\quad
 \Sigma_t=\pi_tg(\iota_t),\quad
 \Delta\iota_t=\frac{\overline{\eta}_t\sigma_t}{g(\iota_t)},\quad
 \log\pi_t+f_0(\iota_t)=-2\sum_{s<t}e_s.
\end{equation}
For $0\le t\le\overline{T}_{\mathrm{init}}$, the coefficient estimates in
the proof of Lemma~\ref{lem:general-scale} give
\[
 D\le g(\iota_t)\le2^{4r^2}D,\qquad
 \iota_{\overline{T}_{\mathrm{init}}}\le(2^{4r}+1)/D.
\]
Since $\sigma_t\ge e^{-1/2}$, the third identity in
\eqref{eq:alg-generating} yields
\[
 \sum_{t<\overline{T}_{\mathrm{init}}}\overline{\eta}_t
 =\sum_{t<\overline{T}_{\mathrm{init}}}\frac{g(\iota_t)}{\sigma_t}\Delta\iota_t
 \le e^{1/2}2^{4r^2}(2^{4r}+1)\le2^{16r^2-1}.
\]
Therefore $\overline{\omega}_t\ge e^{-2\overline{\eta}_t}$ implies
$\pi_t\ge\exp(-2^{16r^2})$. Together with $\pi_t\le1$ and
$\Sigma_t=\pi_tg(\iota_t)$, this proves the scale bounds in
Lemma~\ref{lem:alg-tracking}\,\textup{(\ref{it:alg-initial-scale})}. Also
\[
 \overline{T}_{\mathrm{init}}=\sum_{t<\overline{T}_{\mathrm{init}}}\overline{\eta}_t\Sigma_t^{\delta_3/4}
 \le2^{16r^2}(2^{4r^2}D)^{\delta_3/4}
 \le2^{20r^2}D^{\delta_3/4},
\]
since $\delta_3/4<1$.

\medskip
\noindent\emph{Proof of \textup{(2)} for $t\le \overline{T}_{\mathrm{init}}$.}
Initially \eqref{eq:alg-initial-envelope} follows from
\eqref{eq:initialP}--\eqref{eq:initialCommon}, because
$\zeta<\delta_0/8$. Each update adds at most the first line of
\eqref{eq:alg-birth-test}. Consequently,
$|\mathcal B'|\le(1+\overline{\eta}_t)\overline{k}_tD^{m-\zeta-\delta_0/4}
=\overline{k}_{t+1}D^{m-\zeta-\delta_0/4}$, and
$\overline{k}_t\le c\exp(\sum_{s<\overline{T}_{\mathrm{init}}}\overline{\eta}_s)\le c\exp(2^{16r^2})$.
At $\overline{T}_{\mathrm{init}}$, the additional factor $D^{-\delta_0/4}$ places every
count below $\frac12\Sigma_{\overline{T}_{\mathrm{init}}}^{m-\zeta}$.

\medskip
\noindent\emph{Proof of \textup{(1)}: the degree bounds.}
For the initial normalized degree error
$e_t^\circ=\max_i(E_{i,t}+\widehat E_{i,t})/\Sigma_t^{i-1}$,
\eqref{eq:alg-active-error}--\eqref{eq:alg-full-error} give
\[
 e_{t+1}^\circ\le(1+\overline{c}_2^{3/2}\overline{\eta}_t)e_t^\circ
       +\overline{c}_2^{3/2}\overline{\eta}_tD^{-\delta_3}+\overline{c}_2^{3/2}D^{-2\delta_3}.
\]
Here the factors from comparing $D$ and $\Sigma_t$ are bounded
by $\exp(2r\,2^{16r^2})\le\overline{c}_2^{1/2}$, by
\eqref{eq:alg-large-constant}. Also
$\overline{c}_2^{3/2}2^{16r^2}\le\overline{c}_2^2/4$.
For $\max_i E_{i,t}/\Sigma_t^{i-1}$, the same recurrence holds
without the term $\overline{c}_2^{3/2}\overline{\eta}_tD^{-\delta_3}$. Hence 
$$\frac{E_{i,\overline{T}_{\mathrm{init}}}}{\Sigma_{\overline{T}_{\mathrm{init}}}^{i-1}}
\le \overline{c}_1^{1/2}D^{-7\delta_3/4} \quad \quad \mbox{ and } \quad \quad \frac{\widehat E_{i,\overline{T}_{\mathrm{init}}}}{\Sigma_{\overline{T}_{\mathrm{init}}}^{i-1}}\le \overline{c}_1^{1/2}D^{-\delta_3}.$$
Both this coefficient and $c_1$ are at most $\overline{c}_1$.
We also impose $\overline{C}^{3\delta_3/4}\ge \overline{c}_1^{2r}$. The displayed estimates then establish
\eqref{eq:alg-degree-invariants} at $\overline{T}_{\mathrm{init}}$.
The same bounds, with their partial sums, hold before $\overline{T}_{\mathrm{init}}$.

For $t\ge \overline{T}_{\mathrm{init}}$, the proof of Lemma~\ref{lem:general-scale} uses
$e_t/\overline{\eta}_t\le8b_3\Sigma_t^{-\delta_3}$, and proves Lemma~\ref{lem:alg-tracking}\,\textup{(\ref{it:alg-scale-ratios})},
$d_{2,t}\ge\Sigma_t/2$, and
$\Sigma_{t+1}/(\overline{\omega}_t\Sigma_t)\le e^{\overline{\eta}_t/2}$.
The definition \eqref{eq:alg-large-constant} gives
\begin{align*}
 \overline{c}_2^2\sum_{j>i}\overline{c}_1^{2-j}
    &\le\frac{\delta_3}{128}\overline{c}_1^{2-i},\\
 \overline{c}_2^2\left(1+\sum_{j>i}(2\overline{c}_1)^{r-j+1}\right)
    &\le\frac{\delta_3}{128}(2\overline{c}_1)^{r-i+1}.
\end{align*}
Here $\overline{c}_2^2$ bounds the coefficients in
\eqref{eq:alg-active-error}--\eqref{eq:alg-full-error} after using
\eqref{eq:general-born-mass}. The term containing $E_{i,t}$ in
\eqref{eq:alg-active-error}, after division by
$\Sigma_{t+1}^{i-1-\delta_3}$, has multiplier at most
\[
 \overline{\omega}_t^{i-1}
       \left(\frac{\Sigma_t}{\Sigma_{t+1}}\right)^{i-1-\delta_3}
       \le\overline{\omega}_t^{\delta_3}\le1-\delta_3\overline{\eta}_t/16.
\]
After the same division for the remaining terms in
\eqref{eq:alg-active-error} and \eqref{eq:alg-full-error}, the
additional error is bounded by
\[
 \overline{c}_2^2(\overline{\eta}_t^2+\Sigma_t^{-\delta_3})
       \le(\delta_3/128)\overline{\eta}_t,
\]
provided $\Sigma_t^{\delta_3/4}\ge256\overline{c}_2^2/\delta_3$, by
\eqref{eq:alg-sampling}. We ensure this by imposing
$\log\overline{C}\ge(256\overline{c}_2^2/\delta_3)^{8/\delta_3}$. This proves
\eqref{eq:alg-degree-invariants} by the same inequalities as
\eqref{eq:general-active-error-update}--\eqref{eq:general-full-error-update}.

\medskip
\noindent\emph{Proof of \textup{(2)} for $t\ge \overline{T}_{\mathrm{init}}$.}
For intersections, we apply the estimates in the proof of
Lemma~\ref{lem:general-overlap-round}. If
$|\mathcal B|\ge\frac12\Sigma^{m-\zeta}$, the coverage bound
\eqref{eq:alg-coverage-test} is stronger than the bound
$[1-(m-\delta_1/4)q_0]|\mathcal B|$ for the surviving old members
used there. The second line of \eqref{eq:alg-birth-test} is at most
half the upper bound $\delta_1\overline{\eta}\Sigma^{m-\zeta}/32$ in that
proof. For smaller families its final two paragraphs apply
unchanged. This proves the intersection induction.

\medskip
\noindent\emph{Proof of \textup{(3.5)}.}
For $t\ge \overline{T}_{\mathrm{init}}$, Lemma~\ref{lem:alg-tracking}\,\textup{(\ref{it:alg-scale-ratios})} gives
\[
 \Sigma_t^{\delta_3/4}-\Sigma_{t+1}^{\delta_3/4}
 \ge\Sigma_t^{\delta_3/4}
                (1-e^{-\delta_3\overline{\eta}_t/32})\ge\delta_3/64.
\]
Thus
\[
 \overline{T}-\overline{T}_{\mathrm{init}}\le\frac{64}{\delta_3}\Sigma_{\overline{T}_{\mathrm{init}}}^{\delta_3/4}
       \le\frac{2^{4r^2+6}}{\delta_3}D^{\delta_3/4}.
\]
Adding the bound on $\overline{T}_{\mathrm{init}}$ in
Lemma~\ref{lem:alg-tracking}\,\textup{(\ref{it:alg-initial-length})} and using
$2^{20r^2}+2^{4r^2+6}/\delta_3\le2^{32r^2}/\delta_3$ proves
Lemma~\ref{lem:alg-tracking}\,\textup{(\ref{it:alg-horizon})} and completes the simultaneous induction.
\end{proof}

\begin{lemma}\label{lem:alg-reference-sums}
For $\alpha>\delta_3/4$, let
\[
 \overline{c}_3(\alpha)=\exp\bigl((\alpha+1)2^{16r^2}\bigr)
                         +\frac{32}{1-2^{\delta_3/4-\alpha}}.
\]
For $D\ge \overline{C}$, the following statements hold.
\begin{enumerate}[label=\textup{(\arabic*)},ref=\arabic*]
\item For every $\alpha>\delta_3/4$,
\begin{equation}\label{eq:alg-scale-sum}
 \sum_{t<\overline{T}}\Sigma_t^{-\alpha}
 \le \overline{c}_3(\alpha)
       \left(D^{\delta_3/4-\alpha}+(\log D)^{(\delta_3/4-\alpha)/2}\right).
\end{equation}
\item The sampling parameters satisfy the following bounds.
\begin{enumerate}[label=\textup{(2.\arabic*)},ref=2.\arabic*,leftmargin=3.3em]
\item\label{it:alg-max-eta}
$\max_{t<\overline{T}}\overline{\eta}_t
       \le(\log D)^{-\delta_3/8}$.
\item\label{it:alg-eta-sum}
$\sum_{t<\overline{T}}\overline{\eta}_t\Sigma_t^{-\delta_3}
       \le2\overline{c}_3(5\delta_3/4)(\log D)^{-\delta_3/2}$.
\item\label{it:alg-fourth-power-sum}
$\sum_{t<\overline{T}}\Sigma_t^{-4\delta_3}
       \le2\overline{c}_3(4\delta_3)(\log D)^{-15\delta_3/8}$.
\end{enumerate}
\item At time $\overline{T}$,
\begin{equation}\label{eq:alg-terminal-reference}
 \left|\frac{\iota_{\overline{T}}}{\rho}-1\right|
       \le c_3\frac{\log\log D}{\log D},\qquad
 \frac{\pi_{\overline{T}}}{\rho}\le4(\log D)^{-1/2}.
\end{equation}
\item The reference populations satisfy the following bounds.
\begin{enumerate}[label=\textup{(4.\arabic*)},ref=4.\arabic*,leftmargin=3.3em]
\item\label{it:alg-population-lower} For $0\le t\le\overline{T}$,
$\pi_t\ge D^{-1}(\log D)^{-c_3-1}$.
\item\label{it:alg-population-room}
$\sum_{t<\overline{T}}\frac1{n\pi_t}
       \le D^{\delta_3/4-\delta_0}(\log D)^{c_3+2}$.
\end{enumerate}
\end{enumerate}
\end{lemma}

\begin{proof}[Proof of \textup{(1)}]
Lemma~\ref{lem:alg-tracking}\,\textup{(\ref{it:alg-initial-length})} and \textup{(\ref{it:alg-initial-scale})} give
\[
 \sum_{t<\overline{T}_{\mathrm{init}}}\Sigma_t^{-\alpha}
 \le2^{20r^2}\exp(\alpha2^{16r^2})D^{\delta_3/4-\alpha}
 \le\exp\bigl((\alpha+1)2^{16r^2}\bigr)D^{\delta_3/4-\alpha}.
\]
For $t\ge \overline{T}_{\mathrm{init}}$, let $M_j=2^j\sqrt{\log D}$, $j\ge0$.
While $M_j\le\Sigma_t<2M_j$, the bounds in
Lemma~\ref{lem:alg-tracking}\,\textup{(\ref{it:alg-scale-ratios})} give
$\log(\Sigma_t/\Sigma_{t+1})\ge(2M_j)^{-\delta_3/4}/8$.
Thus at most $1+8(2M_j)^{\delta_3/4}\log2\le16(2M_j)^{\delta_3/4}$
rounds have their starting scale in this interval. Their total
contribution is at most $16\,2^{\delta_3/4} M_j^{\delta_3/4-\alpha}$.
Since $\delta_3/4<1$, summing over $j$ gives the second term in
\eqref{eq:alg-scale-sum}.

\medskip
\noindent\emph{Proof of \textup{(\ref{it:alg-max-eta})}.}
Since $\Sigma_t>\sqrt{\log D}$ for $t<\overline{T}$,
$\overline{\eta}_t=\Sigma_t^{-\delta_3/4}\le(\log D)^{-\delta_3/8}$,
which proves \textup{(2.1)}.

\medskip
\noindent\emph{Proof of \textup{(2.2)} and \textup{(2.3)}.}
Substitution of $\alpha=5\delta_3/4$ and $\alpha=4\delta_3$,
respectively, in part~\textup{(1)}, together with
$D^{\delta_3/4-\alpha}\le(\log D)^{(\delta_3/4-\alpha)/2}$, gives
the two assertions.

\medskip
\noindent\emph{Proof of \textup{(3)}.}
For $\alpha=5\delta_3/4$, $4\delta_3$ or $2^6$, the displayed
formula for $\overline{c}_3(\alpha)$ gives
\[
 \overline{c}_3(\alpha)
 \le\exp(2^{16r^2+7})+128/\delta_3
 \le\frac{\overline{c}_2^2}{2^{10}(b_3+c+1)}.
\]
The last inequality follows from \eqref{eq:alg-large-constant}.
Equations~\eqref{eq:alg-survival-error} and \eqref{eq:alg-generating},
together with Lemma~\ref{lem:alg-reference-sums}\,\textup{(\ref{it:alg-eta-sum})} and \textup{(\ref{it:alg-fourth-power-sum})} yield
\begin{equation}\label{eq:alg-reference-log-bound}
 |\log\pi_t+f_0(\iota_t)|
       \le \overline{c}_2^2(\log D)^{-\delta_3/2}\le1,
\end{equation}
where the last inequality is ensured by
$\overline{C}\ge\exp(\overline{c}_2^{4/\delta_3})$.
The stopping rule and Lemma~\ref{lem:alg-tracking}\,\textup{(\ref{it:alg-scale-ratios})} give
$e^{-2\overline{\eta}_{\overline{T}-1}}\sqrt{\log D}\le\Sigma_{\overline{T}}\le\sqrt{\log D}$.
The proof of Lemma~\ref{lem:general-stopping-cost}, using
\eqref{eq:alg-reference-log-bound}, now gives the explicit bounds
\[
 |f_0(\iota_{\overline{T}})-\log D|\le c_3\log\log D,
 \qquad
 \left|\frac{\iota_{\overline{T}}}{\rho}-1\right|
       \le c_3\frac{\log\log D}{\log D},
 \qquad \pi_{\overline{T}}/\rho\le4(\log D)^{-1/2}.
\]
These are precisely the coefficient calculation in
\eqref{eq:general-terminal-hazard} and its two following
inequalities; they use no constancy assumption on $\overline{\eta}_t$.

\medskip
\noindent\emph{Proof of \textup{(4.1)}.}
Since $\pi_t\ge\pi_{\overline{T}}$,
\eqref{eq:alg-reference-log-bound} also gives
$\pi_t\ge e^{-1}D^{-1}(\log D)^{-c_3}
\ge D^{-1}(\log D)^{-c_3-1}$ for $\log D\ge e$.

\medskip
\noindent\emph{Proof of \textup{(\ref{it:alg-population-room})}.}
Equation~\eqref{eq:general-population-room} and
Lemma~\ref{lem:alg-tracking}\,\textup{(\ref{it:alg-horizon})} give
\[
 \sum_{t<\overline{T}}\frac1{n\pi_t}
 \le\frac{2^{32r^2+1}c}{\delta_3}
       D^{\delta_3/4-\delta_0}(\log D)^{c_3+1}
 \le D^{\delta_3/4-\delta_0}(\log D)^{c_3+2},
\]
where the last inequality holds once
$\log \overline{C}\ge2^{32r^2+1}c/\delta_3$.
\end{proof}

\subsection{Bounds for the accepted set and the final extension}

We control $|J^{\overline{T}}|$ from below and
$S_{\overline{T}}+|X^{\overline{T}}|$, where $S_{\overline{T}}=\sum_{t<\overline{T}}|Y^t|$, from above by a single
explicit nonnegative function $\Phi_t$. The inequality $\Phi_t<1$
will imply both required bounds at time $\overline{T}$, and we show how to select
each round so that $\Phi_{t+1}\le\Phi_t$.

Let $M=n\rho$, $\tau_t=\iota_t/\rho$, $\tau_* =\tau_{\overline{T}}$, and
$S_t=\sum_{s<t}|Y^s|$. For $t\le\overline{T}$, we define
\begin{equation}\label{eq:alg-normalized-accepted}
 \begin{aligned}
 P_t&=\frac{|X^t\setminus\mathsf W^t|}{n\pi_t},\qquad
 Q_t=\frac{|\mathsf W^t|}{n\pi_t},\qquad
 b_t=\frac{\pi_{\overline{T}}+\sum_{s=t}^{\overline{T}-1}\pi_sq_sp_s}{\rho},\\
 L_t&=\frac{|J^t|}{M}+(\tau_*-\tau_t)P_t,\qquad
 C_t=\frac{S_t}{M}+b_tP_t+\frac{\pi_{\overline{T}}}{\rho}Q_t.
 \end{aligned}
\end{equation}
Let $r_t=\overline{\eta}_t\Sigma_t^{-\delta_3}
+\Sigma_t^{-4\delta_3}+\Sigma_t^{-2^6}$. We also define
\begin{align}
 R_t&=\frac{\overline{c}_2^8}{\varepsilon^2}
  \sum_{s=t}^{\overline{T}-1}
  \left(r_s+\frac1{n\pi_s}+(\Delta\tau_s)^2+
                                      \frac{\Delta\tau_s}{M}\right),
 \label{eq:alg-reserve}\\
 \Phi_t&=\frac{C_t}{1+\varepsilon}
 +\frac{16}{\varepsilon^2}
            \bigl((P_t-1)^2+(L_t-\tau_*)^2\bigr)+R_t.
 \label{eq:alg-global-potential}
\end{align}
All coefficients involving the reference parameters are deterministic
and computable in advance. The function $\Phi_t$ is nonnegative;
its three nonnegative summands give
\[
 \Phi_t<1\quad\Longrightarrow\quad
 C_t<1+\varepsilon,\qquad
 |P_t-1|<\varepsilon/4,\qquad
 |L_t-\tau_*|<\varepsilon/4.
\]
At time $\overline{T}$, the definitions in \eqref{eq:alg-normalized-accepted}
turn these inequalities into bounds on $|J^{\overline{T}}|$ and $S_{\overline{T}}+|X^{\overline{T}}|$.
The sum $R_t$ in \eqref{eq:alg-reserve} is included so that its
decrease bounds the sum of the errors arising in round $t$.

\begin{lemma}\label{lem:alg-one-round-moments}
Suppose the local hypotheses in Lemma~\ref{lem:alg-tracking} hold,
and let $Z$ be a good focus at time $t$ in the sense of
Definition~\ref{def:alg-good-round}. Let $h=|Z|$ and
$a_t=p_t\sigma_t/M$. The expectations below are over the independent
variables $\chi_x,\mathsf w_x$ with $Z$ fixed, before conditioning
on the inequalities in Definition~\ref{def:alg-good-round}. Then the following bounds hold.
\begin{enumerate}[label=\textup{(\arabic*)},ref=\arabic*]
\item The expected normalized population satisfies
\begin{equation}\label{eq:alg-population-first}
 |\mathbb E(P_{t+1}\mid Z)-P_t|
       \le\overline{c}_2\Sigma_t^{-A_0}P_t.
\end{equation}
\item Its variance satisfies
\begin{equation}\label{eq:alg-population-second}
 \operatorname{Var}(P_{t+1}\mid Z)
       \le\overline{c}_2\left(\Sigma_t^{-4\delta_3}P_t^2
                                      +\frac{P_t}{n\pi_t}\right).
\end{equation}
\item The expected number of accepted atoms satisfies
\begin{equation}\label{eq:alg-accepted-first}
 \left|\mathbb E\left(\frac{|I^t|}{M}\,\middle|\,Z\right)-a_th\right|
       \le\overline{c}_2 r_ta_th.
\end{equation}
\item Its normalized variance satisfies
\begin{equation}\label{eq:alg-accepted-second}
 \operatorname{Var}\left(\frac{|I^t|}{M}\,\middle|\,Z\right)
       \le\overline{c}_2\left(\frac{a_th}{M}+(a_th)^2\right).
\end{equation}
\item The expected numbers of surviving and sampled atoms satisfy
\begin{equation}\label{eq:alg-full-first}
 \mathbb E(|X^{t+1}|\mid Z)
       \le\overline{\omega}_t(1+\overline{c}_2 r_t)|X^t|,
 \qquad \mathbb E(|Y^t|\mid Z)=p_th.
\end{equation}
\end{enumerate}
\end{lemma}

\begin{proof}[Proof of \textup{(1)}]
Let $U=X^t\setminus\mathsf W^t$ and
$A_x=\psi_x(1-\mathsf w_x)$ for $x\in U$. All expectations in
this proof are conditional on the fixed focus $Z$.
Since $P_{t+1}=(n\pi_{t+1})^{-1}\sum_{x\in U}A_x$ and
$\pi_{t+1}=\overline{\omega}_t\pi_t$, equation~\eqref{eq:alg-equalization-bounds}
gives
\[
 (1-8\Sigma_t^{-A_0})P_t
       \le\mathbb E(P_{t+1}\mid Z)\le P_t.
\]
This proves \eqref{eq:alg-population-first}.

\medskip
\noindent\emph{Proof of \textup{(2)}.}
For distinct $x,y\in U$, Lemma~\ref{lem:general-survival-estimates}(1),
specifically \eqref{eq:general-survival-joint}, and independence
of $\mathsf w_x,\mathsf w_y$ give
\[
 \mathbb E(A_xA_y\mid Z)
       \le(1+4b_3\Sigma_t^{-4\delta_3})\overline{\omega}_t^2.
\]
Also $\mathbb E(A_x\mid Z)\mathbb E(A_y\mid Z)
\ge\overline{\omega}_t^2(1-16\Sigma_t^{-A_0})$.
Summing the covariances over unordered pairs $\{x,y\}\in\binom U2$, we obtain
\begin{align*}
 \operatorname{Var}(P_{t+1}\mid Z)
 &=\frac1{(n\pi_{t+1})^2}
   \biggl(\sum_{x\in U}\operatorname{Var}(A_x\mid Z)
       +2\sum_{\{x,y\}\in\binom U2}\operatorname{Cov}(A_x,A_y\mid Z)\biggr)\\
 &\le\frac{P_t}{n\pi_t\overline{\omega}_t}
      +(4b_3\Sigma_t^{-4\delta_3}+16\Sigma_t^{-A_0})P_t^2.
\end{align*}
The bounds $\overline{\omega}_t\ge1/2$ and $A_0>4\delta_3$ imply
\eqref{eq:alg-population-second}. This is also the calculation in
\eqref{eq:general-active-second-moment} of
Lemma~\ref{lem:general-active-population-tracking}.

\medskip
\noindent\emph{Proof of \textup{(3)}.}
Equation~\eqref{eq:general-uniform-acceptance} in
Lemma~\ref{lem:general-survival-estimates}(3), together with
$\mathsf E_{\rm acc}\le2\overline{\eta}_t\Sigma_t^{-\delta_3}$, proves
\eqref{eq:alg-accepted-first}.

\medskip
\noindent\emph{Proof of \textup{(4)}.}
Since $I^t\subseteq Y^t$,
\[
 \mathbb E\bigl(|I^t|^2\bigr)\le p_th+p_t^2h^2.
\]
By \eqref{eq:general-sigma-bounds}, $\sigma_t\ge1/2$, so
\[
 \operatorname{Var}(|I^t|/M\mid Z)
 \le\frac{p_th+p_t^2h^2}{M^2}
 \le\frac{2a_th}{M}+4(a_th)^2.
\]
This proves \eqref{eq:alg-accepted-second}.

\Needspace{6\baselineskip}
\medskip
\noindent\emph{Proof of \textup{(5)}.}\nopagebreak[4]
\begin{align*}
 \mathbb E(|X^{t+1}|\mid Z)&=\sum_{x\in X^t}s_x
 \le\overline{\omega}_te^{3e_t}|X^t|
 \le\overline{\omega}_t(1+\overline{c}_2 r_t)|X^t|,\\
 \mathbb E(|Y^t|\mid Z)&=\sum_{x\in Z}p_t=p_th,
\end{align*}
which proves \eqref{eq:alg-full-first}.
\end{proof}

\Needspace{8\baselineskip}
\begin{lemma}\label{lem:alg-focus-objective}
Suppose $\Phi_t<1$. Let $U=|X^t\setminus\mathsf W^t|$ and
$a_t=p_t\sigma_t/M$. For $h\ge0$, let
\begin{align*}
 A_t(h)&=S_t/M+b_{t+1}P_t+(\pi_{\overline{T}}/\rho)Q_t+p_th/M
                  +\overline{c}_2^2r_t(1+C_t),\\
 B_t(h)&=r_t\bigl(1+(a_th)^2\bigr)+(n\pi_t)^{-1}+a_th/M+(a_th)^2.
\end{align*}
Then
\begin{enumerate}[label=\textup{(\arabic*)},ref=\arabic*]
\item The following nonnegative quadratic polynomial in $h$ is an
upper bound for $\mathbb E(\Phi_{t+1}\mid Z)$ whenever $Z$ is a good focus as defined in
Definition~\ref{def:alg-good-round} and $|Z|=h$:
\begin{equation}\label{eq:alg-focus-objective}
 G_t(h)=R_{t+1}+\frac{A_t(h)}{1+\varepsilon}
 +\frac{16}{\varepsilon^2}\left[(P_t-1)^2+
 \bigl(L_t-\tau_*+a_t(h-q_tU)\bigr)^2+\overline{c}_2^2B_t(h)\right].
\end{equation}
\item When $\xi_x=\vb1_{\{x\in Z\}}$ are independent Bernoulli$(q_t)$ variables for $x\in X^t\setminus\mathsf W^t$,
\begin{equation}\label{eq:alg-focus-objective-mean}
 \mathbb E\bigl(G_t(|Z|)\bigr)
 \le\Phi_t-R_t+R_{t+1}
 +\frac{\overline{c}_2^4}{\varepsilon^2}
       \left(r_t+\frac1{n\pi_t}+(\Delta\tau_t)^2
                                      +\frac{\Delta\tau_t}{M}\right).
\end{equation}
\end{enumerate}
\end{lemma}

\begin{proof}[Proof of \textup{(1)}]
The assumption $\Phi_t<1$ implies
\begin{equation*}
 |P_t-1|<\varepsilon/4,\qquad
 |L_t-\tau_*|<\varepsilon/4,\qquad C_t<1+\varepsilon.
\end{equation*}
In particular, $P_t\le2$ and $|L_t|\le\tau_*+\varepsilon/4$. The identities
\begin{align*}
 L_{t+1}&=|J^t|/M+|I^t|/M
                         +(\tau_*-\tau_{t+1})P_{t+1},\\
 C_{t+1}&=S_t/M+|Y^t|/M
             +(b_{t+1}-\pi_{\overline{T}}/\rho)P_{t+1}
             +\frac{\pi_{\overline{T}}|X^{t+1}|}{\rho n\pi_{t+1}}
\end{align*}
reduce the calculation to Lemma~\ref{lem:alg-one-round-moments}.
For $L_{t+1}$, its mean differs from
$L_t+a_t(h-q_tU)$ by at most $4\overline{c}_2(r_ta_th+\Sigma_t^{-A_0})$.
Writing $v_t=\tau_*-\tau_{t+1}$, where $0\le v_t\le2$, we have
\begin{align*}
 \operatorname{Var}(L_{t+1}\mid Z)
 &\le2\operatorname{Var}(|I^t|/M\mid Z)
          +2v_t^2\operatorname{Var}(P_{t+1}\mid Z)\\
 &\le2\operatorname{Var}(|I^t|/M\mid Z)
          +8\operatorname{Var}(P_{t+1}\mid Z).
\end{align*}
We apply \eqref{eq:alg-population-second} and
\eqref{eq:alg-accepted-second} to the right side.
If $\ell_t=L_t-\tau_*+a_t(h-q_tU)$ and
$d_t=\mathbb E(L_{t+1}\mid Z)-L_t-a_t(h-q_tU)$, then
\[
 |(\ell_t+d_t)^2-\ell_t^2|
       \le2|\ell_t||d_t|+d_t^2
       \le\overline{c}_2^2r_t(1+(a_th)^2).
\]
For $C_{t+1}$, equation~\eqref{eq:alg-full-first} and
$(\pi_{\overline{T}}/\rho)Q_t\le C_t$ give
\[
 \mathbb E(C_{t+1}\mid Z)
 \le S_t/M+b_{t+1}P_t+(\pi_{\overline{T}}/\rho)Q_t+p_th/M
                          +\overline{c}_2^2r_t(1+C_t).
\]
Together with
$\mathbb E((L_{t+1}-\tau_*)^2\mid Z)
=\operatorname{Var}(L_{t+1}\mid Z)
 +(\mathbb E(L_{t+1}\mid Z)-\tau_*)^2$, these inequalities prove
\eqref{eq:alg-focus-objective}. Every displayed term in that
polynomial is nonnegative for $h\ge0$.

\medskip
\noindent\emph{Proof of \textup{(2)}.}
For independent Bernoulli$(q_t)$ variables $\xi_x$ and
$h=\sum_{x\in X^t\setminus\mathsf W^t}\xi_x$,
\[
 \mathbb E(h)=q_tU,\quad\operatorname{Var}(h)\le q_tU,\quad
 a_tq_tU=\Delta\tau_tP_t,\quad
 a_t^2q_tU\le\frac{\Delta\tau_tP_t}{M}.
\]
Since $b_t-b_{t+1}=\pi_tq_tp_t/\rho$, the linear expression satisfies
\[
 \mathbb E\bigl(S_t/M+b_{t+1}P_t+(\pi_{\overline{T}}/\rho)Q_t+p_th/M\bigr)=C_t.
\]
The quadratic expression satisfies
\[
 \mathbb E\bigl((L_t-\tau_*+a_t(h-q_tU))^2\bigr)
       =(L_t-\tau_*)^2+a_t^2\operatorname{Var}(h).
\]
Substitution of these identities and $P_t\le2$ into
\eqref{eq:alg-focus-objective} gives
\eqref{eq:alg-focus-objective-mean}.
\end{proof}

\begin{proposition}\label{prop:alg-one-round}
Assume the local hypotheses in Lemma~\ref{lem:alg-tracking} and
$\Phi_t<1$. A deterministic polynomial-time algorithm
chooses a focus and construction satisfying
Definition~\ref{def:alg-good-round} and
\begin{equation}\label{eq:alg-potential-decrease}
 \Phi_{t+1}\le\Phi_t.
\end{equation}
\end{proposition}

\begin{proof}
We apply Proposition~\ref{prop:alg-weighted-selection} first to
the variables $\xi_x=\vb1_{\{x\in Z\}}$, and then, with $Z$
fixed, to the variables $\chi_x$ and $\mathsf w_x$.
Lemma~\ref{lem:alg-local-estimators} supplies the bad events and
the polynomials that bound their indicators.

We write $G_t(h)=\gamma+\gamma_1h+\gamma_2h^2$, with
\[
 \gamma=G_t(0),\quad
 \gamma_2=\tfrac12(G_t(2)-2G_t(1)+G_t(0)),\quad
 \gamma_1=G_t(1)-G_t(0)-\gamma_2.
\]
Let $<$ be a fixed linear order on $U=X^t\setminus\mathsf W^t$.
Since $\xi_x^2=\xi_x$, the exact expansion is
\begin{equation}\label{eq:alg-focus-polynomial-expansion}
 G_t\left(\sum_{x\in U}\xi_x\right)
 =\gamma+(\gamma_1+\gamma_2)\sum_{x\in U}\xi_x
                    +2\gamma_2\sum_{\substack{x,y\in U\\x<y}}
                                                \xi_x\xi_y.
\end{equation}
For each negative coefficient $\alpha_F$, the identity
$\alpha_F\xi_F=|\alpha_F|(1-\xi_F)-|\alpha_F|$ gives a
nonnegative multiple of $1-\xi_F$ and a constant.
Thus \eqref{eq:alg-focus-polynomial-expansion} is a constant plus
a nonnegative weighted sum of indicators, each using at most two
of the variables $\xi_x$. Lemma~\ref{lem:conditioning-locality}(2)
bounds the number of bad-event sets met by each such set of
variables by $\Sigma_t^{b+1}$.

The sum of the absolute values of the nonconstant coefficients is
at most $\overline{c}_2^3\varepsilon^{-2}\Sigma_t^2$. This follows from
\eqref{eq:alg-focus-objective}, $P_t\le2$, $C_t<1+\varepsilon$,
$|L_t-\tau_*|<\varepsilon/4$, and
\begin{equation*}
 \frac{|U|}{M}=\frac{\pi_tP_t}{\rho}\le2\Sigma_t.
\end{equation*}
The last inequality follows from $\Sigma_t=\pi_tg(\iota_t)$ and
$\rho g(\iota_t)\ge\rho D\ge1$, as in
Lemma~\ref{lem:general-accepted-count-tracking}.
For example, $|U|^2a_t^2\le4\Sigma_t^2$, and every quadratic
coefficient in \eqref{eq:alg-focus-objective} is at most
$64\overline{c}_2^2\varepsilon^{-2}a_t^2$ in absolute value.
The linear coefficients satisfy the same total bound with
$|U|/M$ in place of $|U|^2/M^2$.
Proposition~\ref{prop:alg-weighted-selection} therefore gives a
focus satisfying all focus conditions and
\begin{equation}\label{eq:alg-selected-focus}
 G_t(|Z|)\le\mathbb E\bigl(G_t(|Z|)\bigr)
                   +2\overline{c}_2^3\varepsilon^{-2}\Sigma_t^{-2^6}.
\end{equation}

For the selected $Z$, the probabilities are given by
\eqref{eq:alg-retention}. The four quantities in $\Phi_{t+1}$ are
\begin{equation}\label{eq:alg-next-potential-functions}
 \begin{alignedat}{2}
 P_{t+1}&=\frac{\sum_{x\in U}\psi_x(1-\mathsf w_x)}{n\pi_{t+1}},
 &\quad Q_{t+1}&=\frac{\sum_{x\in X^t}\psi_x}{n\pi_{t+1}}-P_{t+1},\\
 L_{t+1}&=\frac{|J^t|+\sum_{x\in Z}J_x}{M}
              +(\tau_*-\tau_{t+1})P_{t+1},
 &\quad C_{t+1}&=\frac{S_t+\sum_{x\in Z}\chi_x}{M}
              +b_{t+1}P_{t+1}+\frac{\pi_{\overline{T}}}{\rho}Q_{t+1}.
 \end{alignedat}
\end{equation}
We substitute \eqref{eq:alg-next-potential-functions} into
\eqref{eq:alg-global-potential} and expand every square by
$(\sum_x A_x)^2=\sum_{x,y}A_xA_y$.
The resulting expression has the form
\[
 \Phi_{t+1}=c_0+\sum_{\ell=1}^{M_1}c_\ell \mathsf{F}_\ell,
\]
where each $\mathsf{F}_\ell$ is either $\chi_x$ or a product of at most two
of $J_x,\psi_x,\psi_x(1-\mathsf w_x)$.
For each product, let $\mathsf{F}_\ell^-,\mathsf{F}_\ell^+$ be the polynomials in
Lemma~\ref{lem:alg-local-polynomial-bounds} with exponent $A_0$.
For $\mathsf{F}_\ell=\chi_x$ we take $\mathsf{F}_\ell^-=\mathsf{F}_\ell^+=\chi_x$.
The polynomial
\[
 c_0-\sum_{c_\ell<0}|c_\ell|
       +\sum_{c_\ell\ge0}c_\ell \mathsf{F}_\ell^+
       +\sum_{c_\ell<0}|c_\ell|(1-\mathsf{F}_\ell^-)
\]
is at least $\Phi_{t+1}$ for every assignment. Its nonconstant
terms are nonnegative, and its expected excess above
$\mathbb E(\Phi_{t+1}\mid Z)$ is at most
$\Sigma_t^{-A_0}\sum_\ell|c_\ell|$.

The total absolute coefficient sum is again at most
$\overline{c}_2^3\varepsilon^{-2}\Sigma_t^2$. In detail, the sums for
squares of the population, accepted count, and their products are
bounded using
\[
 \left(\frac{|U|}{n\pi_{t+1}}\right)^2\le16,
 \qquad \left(\frac{|U|}{M}\right)^2\le4\Sigma_t^2,
 \qquad 2ab\le a^2+b^2.
\]
The coefficient sum from the full population is at most
\[
 \frac{\pi_{\overline{T}}|X^t|}{\rho n\pi_{t+1}}
 \le2\frac{\pi_{\overline{T}}}{\rho}(P_t+Q_t)
 \le8+2C_t.
\]
Here $\overline{\omega}_t\ge1/2$, $\pi_{\overline{T}}/\rho\le2$,
$P_t\le2$, and $(\pi_{\overline{T}}/\rho)Q_t\le C_t$.
Each summand is determined by variables in the union of at most
two of the balls in Lemma~\ref{lem:conditioning-locality}; the
number of bad-event sets it meets is at most $\Sigma_t^{b+1}$.
Applying Proposition~\ref{prop:alg-weighted-selection} gives
\begin{equation}\label{eq:alg-selected-construction}
 \Phi_{t+1}\le\mathbb E(\Phi_{t+1}\mid Z)
       +4\overline{c}_2^3\varepsilon^{-2}\Sigma_t^{-2^6}
 \le G_t(|Z|)+4\overline{c}_2^3\varepsilon^{-2}\Sigma_t^{-2^6}.
\end{equation}
The polynomial approximation error is at most
$\overline{c}_2^3\varepsilon^{-2}\Sigma_t^{2-A_0}$.
Since $A_0>2^7$ and $2^7-2\ge2^6$, this error and the error from
\eqref{eq:alg-weighted-selection} are included in the preceding bound.
Let
\[
 d_t=r_t+\frac1{n\pi_t}+(\Delta\tau_t)^2+\frac{\Delta\tau_t}{M}.
\]
Then $R_t-R_{t+1}=\overline{c}_2^8\varepsilon^{-2}d_t$ and
$d_t\ge\Sigma_t^{-2^6}$.
Equations~\eqref{eq:alg-focus-objective-mean},
\eqref{eq:alg-selected-focus}, and
\eqref{eq:alg-selected-construction} therefore give
\[
 \Phi_{t+1}-\Phi_t
 \le\frac{\overline{c}_2^4-\overline{c}_2^8}{\varepsilon^2}d_t
       +\frac{6\overline{c}_2^3}{\varepsilon^2}\Sigma_t^{-2^6}
 \le-\varepsilon^{-2}d_t\le0,
\]
where the penultimate inequality uses
$\overline{c}_2^8\ge\overline{c}_2^4+6\overline{c}_2^3+1$.
This proves \eqref{eq:alg-potential-decrease}.
\end{proof}

\begin{lemma}\label{lem:alg-potential-endpoints}
The threshold $\overline{C}$ can be chosen so that, for $D\ge \overline{C}$,
\begin{enumerate}[label=\textup{(\arabic*)},ref=\arabic*]
\item $\Phi_0<1-\varepsilon/2$.
\item Every sequence satisfying \eqref{eq:alg-potential-decrease}
through time $\overline{T}$ then satisfies
\begin{equation}\label{eq:alg-final-counts}
 |J^{\overline{T}}|\ge(1-\varepsilon)n\rho,
 \qquad S_{\overline{T}}+|X^{\overline{T}}|\le(1+\varepsilon)n\rho.
\end{equation}
\end{enumerate}
\end{lemma}

\begin{proof}[Proof of \textup{(1)}]
Equations~\eqref{eq:alg-generating} and
\eqref{eq:alg-terminal-reference}, together with $D\rho\ge1$,
give the following bounds:
\begin{align*}
 \sum_{t<\overline{T}}\Delta\tau_t&=\tau_*,&
 |\tau_*-1|&\le c_3\frac{\log\log D}{\log D},\\
 \max_{t<\overline{T}}\Delta\tau_t&\le\max_{t<\overline{T}}\overline{\eta}_t,&
 \sum_{t<\overline{T}}(\Delta\tau_t)^2&\le2\max_{t<\overline{T}}\overline{\eta}_t.
\end{align*}
The last inequality follows from $\tau_*\le2$, which holds once
$c_3\log\log D\le\log D$.
The bounds on $\overline{c}_3$ in the proof of
Lemma~\ref{lem:alg-reference-sums} give
\[
 2^7\overline{c}_2^8\bigl(1+c+\overline{c}_3(5\delta_3/4)
                  +\overline{c}_3(4\delta_3)
                  +\overline{c}_3(2^6)\bigr)\le\overline{c}_2^{10}.
\]
By Lemma~\ref{lem:alg-reference-sums}\,\textup{(1)--(4)}
and $M=n\rho\ge(2c)^{-1}D^{\delta_0}$,
\begin{equation}\label{eq:alg-initial-reserve-bound}
 R_0\le\frac{\overline{c}_2^{10}}{\varepsilon^2}
 \left((\log D)^{-\delta_3/8}
      +D^{\delta_3/4-\delta_0}(\log D)^{c_3+2}
      +D^{-\delta_0}\right).
\end{equation}
At time zero $P_0=1$, $Q_0=S_0=0$, and $L_0=\tau_*$. Also
\[
 b_0=\frac{\pi_{\overline{T}}}{\rho}
                  +\sum_{t<\overline{T}}\frac{\Delta\tau_t}{\sigma_t}.
\]
By \eqref{eq:general-sigma-bounds},
$0\le\sigma_t^{-1}-1\le4\overline{\eta}_t$ when $\overline{\eta}_t\le1/4$.
Consequently
\begin{equation}\label{eq:alg-initial-cost-bound}
 |b_0-1|\le8(\log D)^{-\delta_3/8}
       +c_3\frac{\log\log D}{\log D}+4(\log D)^{-1/2}.
\end{equation}
We require $\overline{C}$ to satisfy the following conditions for every $D\ge \overline{C}$: the right sides of
\eqref{eq:alg-initial-reserve-bound} and
\eqref{eq:alg-initial-cost-bound} are at most $\varepsilon/16$,
$c_3\log\log D/\log D\le\varepsilon/8$, and
$0<\varepsilon\le1/4$.
These are explicit inequalities in $D$ and the specified constants.
Such a finite $\overline{C}$ exists because
$\log D\le\log\triangle\le(r-1)\log D$,
$\delta_3/4<\delta_0$, and the ratio of each required error bound
to $((r-1)\log D)^{-\delta_3/32}$ tends to zero as $D$ tends to infinity.
For example, the logarithmic term in
\eqref{eq:alg-initial-reserve-bound}, after division by
$\varepsilon$, is at most
$\overline{c}_2^{10}(r-1)^{3\delta_3/32}(\log D)^{-\delta_3/32}$,
which decreases to zero.
It follows that
\[
 \Phi_0=\frac{b_0}{1+\varepsilon}+R_0
 \le\frac{1+\varepsilon/16}{1+\varepsilon}
                                      +\varepsilon/16
 <1-\varepsilon/2.
\]
\medskip
\noindent\emph{Proof of \textup{(2)}.}
At time $\overline{T}$, $R_{\overline{T}}=0$ and $\Phi_{\overline{T}}\le\Phi_0<1$.
Equations~\eqref{eq:alg-global-potential} and
\eqref{eq:alg-normalized-accepted} imply
\[
 S_{\overline{T}}+|X^{\overline{T}}|<(1+\varepsilon)M,\qquad
 \left|\frac{|J^{\overline{T}}|}{M}-\tau_*\right|<\varepsilon/4.
\]
Since $|\tau_*-1|\le\varepsilon/8$, these imply
\eqref{eq:alg-final-counts}.
\end{proof}

\subsection{The algorithm and its running time}

We combine the deterministic choices for each round with a final
greedy extension to a maximal independent set. The proof specifies
the numerical precision and bounds the total number of bit operations
by a fixed power of $n$.

\begin{proof}[Proof of Theorem~\textup{\ref{thm:deterministic-algorithm}}]
We first describe the construction for the exact stopping indices
in \eqref{eq:alg-stopping}. At the end of the proof we give a finite
enumeration that includes these indices and avoids deciding equality
with an irrational threshold. For prescribed indices, we compute the
reference parameters in
\eqref{eq:alg-sampling},
\eqref{eq:alg-survival-error}--\eqref{eq:alg-survival-reference},
\eqref{eq:alg-degree-reference}, and \eqref{eq:alg-stopping}, the coefficients
in \eqref{eq:alg-normalized-accepted}, and the sums $R_t$ in
\eqref{eq:alg-reserve}. The initial sets are
$J^0=\mathsf W^0=\varnothing$. For $0\le t<\overline{T}$, the sets at time
$t+1$ are determined as follows.
\begin{enumerate}
\item We construct the polynomials for the focus conditions in
Lemma~\ref{lem:alg-local-estimators}. Applied to these polynomials
and $G_t(|Z|)$ in \eqref{eq:alg-focus-objective},
Proposition~\ref{prop:alg-weighted-selection} gives a focus $Z^t$
satisfying \eqref{eq:alg-selected-focus} and all focus conditions.
\item For this $Z^t$, we compute the numbers $s_x^+$ in
\eqref{eq:alg-upper-survival}. We define the independent variables
$\chi_x,\mathsf w_x$ by
$\mathbb P(\chi_x=1)=p_t$ for $x\in Z^t$, $\chi_x=0$ otherwise,
and $\mathbb P(\mathsf w_x=0)=\overline{\omega}_t/s_x^+$ for
$x\in X^t\setminus\mathsf W^t$.
\item The construction polynomials in
Lemma~\ref{lem:alg-local-estimators} and the upper polynomial for
$\Phi_{t+1}$ in the proof of Proposition~\ref{prop:alg-one-round}
satisfy the hypotheses of Proposition~\ref{prop:alg-weighted-selection}.
That proposition determines values of $\chi_x,\mathsf w_x$ satisfying
\eqref{eq:alg-selected-construction} and all construction conditions.
\item We define $Y^t=\{x\in Z^t:\chi_x=1\}$ and use the formulas
of Subsection~\ref{subsec:actual-process} for
$I^t,J^{t+1},X^{t+1},\mathsf W^{t+1}$ and the residual blocks.
Among those residual blocks, we retain each distinct
inclusion-minimal member once.
\end{enumerate}
Lemmas~\ref{lem:alg-tracking} and
\ref{lem:alg-potential-endpoints}, and
Proposition~\ref{prop:alg-one-round}, show inductively that every
step is available and that \eqref{eq:alg-final-counts} holds.
At time $\overline{T}$, let $v_1,\ldots,v_q$ be the atoms of
\[
 X^{\overline{T}}\cup\bigcup_{t<\overline{T}}(Y^t\setminus I^t)
\]
in a fixed order. Let $I_0=J^{\overline{T}}$ and, for $1\le j\le q$, let
\[
 I_j=\begin{cases}
 I_{j-1}\cup\{v_j\},&\text{if no }S\in\mathcal S
                           \text{ satisfies }S\subseteq I_{j-1}\cup\{v_j\},\\
 I_{j-1},&\text{otherwise}.
 \end{cases}
\]
We take $I=I_q$. Lemma~\ref{lem:bookkeeping}, including its
assertions about the residual antichain and the atoms that remain
available for this extension, applies to these updates. It gives a maximal
independent set $I$ with
$|J^{\overline{T}}|\le|I|\le S_{\overline{T}}+|X^{\overline{T}}|$, proving the required bounds.

We give separate names to the exponents in the computation.
Let
\begin{align*}
 d_*&=\left\lceil2^7\bigl(r+a+b+A_0+\overline{a}_5
       +\overline{a}_3(A_0,2m_0,r,\delta_0)
       +\log_2(\overline{c}_2+1)+1\bigr)^2\right\rceil,\\
 \overline{a}_6&=\overline{a}_4(d_*,r)+d_*+4,\\
 \overline{c}_4&=(r+2)(\overline{a}_6+3)+2.
\end{align*}
This last $\overline{c}_4$ is the exponent in the theorem.
The integers in \eqref{eq:alg-orders} and
\eqref{eq:alg-truncation-orders} are independent of $n$.
The descriptions and arithmetic-operation bounds for their
polynomials are given in
Lemmas~\ref{lem:alg-local-polynomial-bounds} and
\ref{lem:alg-local-estimators}. Lemma~\ref{lem:conditioning-locality}
bounds the number of events and the number of variables occurring
in each event. These counts are at most $N^{d_*}$.
Once the numerical precision below has been specified,
Proposition~\ref{prop:alg-weighted-selection} requires at most
$N^{\overline{a}_4(d_*,r)}$ bit operations per application.
Since $D\le n$, Lemma~\ref{lem:alg-tracking}\,\textup{(\ref{it:alg-horizon})} gives $\overline{T}\le N^2$ once
$N\ge2^{32r^2}/\delta_3$. Thus
all rounds and the final greedy extension require at most
$N^{\overline{a}_6}$ bit operations, apart from the precision
calculations bounded next.

We represent each real reference probability or parameter $u$ by
rational endpoints $u^-,u^+$ satisfying
\[
 u^-\le u\le u^+.
\]
Each elementary evaluation uses an interval of width at most
$N^{-A_1}$, where the integer $A_1$ is specified below. Subsequent
arithmetic uses rational upper and lower bounds. We first bound
the widths of the resulting intervals.
The reference degrees are evaluated from
$d_{i,t}=\pi_t^{i-1}g_i(\iota_t)$ in \eqref{eq:alg-generating}.
The first step starts at the exact value $\iota_0=0$ and gives
$\iota_1=\overline{\eta}_0\sigma_0/g(0)\ge N^{-3}$.
For $t\ge1$, a nonzero polynomial $g_i$ has nonnegative
coefficients and degree at most $r-i$, so
\[
 \frac{\iota_tg_i'(\iota_t)}{g_i(\iota_t)}\le r-i.
\]
For $t\ge1$, let $u=\log\pi_t$ and $v=\log\iota_t$.
The recurrence for these two logarithms is the map
\[
 (u,v)\longmapsto
 \left(u+\log\overline{\omega}_t,\,
 v+\log\left(1+\frac{\overline{\eta}_t\sigma_t}{e^vg(e^v)}\right)\right),
\]
where $\Sigma_t=e^ug(e^v)$, $d_{i,t}=e^{(i-1)u}g_i(e^v)$,
and $\overline{\eta}_t,\sigma_t,\overline{\omega}_t$ are evaluated from these values.
Differentiating this map, using the last bound for each nonzero
$g_i$, gives the following upper bound for the maximum absolute
row sum of its derivative matrix:
\[
 1+\overline{c}_2\left(\overline{\eta}_t+\frac{\Delta\iota_t}{\iota_t}\right).
\]
Polynomials $g_i$ that are identically zero are kept equal to zero;
no fractional power is differentiated at zero.
Moreover, $g$ is nondecreasing, $\Sigma_{t+1}\ge\overline{\omega}_t\Sigma_t$,
and $\sigma_t^{-1}\le e^{2\overline{\eta}_t}$. Hence
\[
 \frac{\Delta\iota_{t+1}}{\Delta\iota_t}
 \le\frac{\overline{\eta}_{t+1}}{\overline{\eta}_t}\frac1{\sigma_t}
 \le e^{(2+\delta_3/2)\overline{\eta}_t}\le3.
\]
Thus $0\le\Delta\iota_t/\iota_t\le3$ for $t\ge1$.
For $0\le x\le3$, $x\le4\log(1+x)$, and therefore
\[
 \sum_{1\le t<\overline{T}}\frac{\Delta\iota_t}{\iota_t}
 \le4\log\frac{\iota_{\overline{T}}}{\iota_1}\le16\log N,
 \qquad
 \sum_{t<\overline{T}}\overline{\eta}_t\le2^{16r^2+4}\log N.
\]
The second inequality follows by summing
Lemma~\ref{lem:alg-tracking}\,\textup{(\ref{it:alg-scale-ratios})} after $\overline{T}_{\mathrm{init}}$ and using
Lemma~\ref{lem:alg-tracking}\,\textup{(\ref{it:alg-initial-sampling})} before $\overline{T}_{\mathrm{init}}$.
The product of the derivative bounds is consequently at most
\[
 N^{\overline{c}_2(2^{16r^2+4}+16)}\le N^{\overline{c}_2^2},
\]
by \eqref{eq:alg-large-constant}.

Let
\[
 A_1=\left\lceil2^4(\overline{a}_6+\overline{c}_2^2+d_*+A_0+1)\right\rceil.
\]
We use intervals of width at most $N^{-A_1}$ for the elementary
evaluations. The denominators have lower bounds that are inverse
powers of $N$, by Lemma~\ref{lem:alg-reference-sums}\,\textup{(\ref{it:alg-population-lower})}, $s_x\ge1/8$,
and $\rho\ge D^{-1}\ge N^{-1}$. The coefficient bounds for the
polynomials are at most $N^{d_*}$ by their explicit expansions.
The preceding sensitivity estimate then bounds the total propagated
error by $N^{-A_1+\overline{a}_6+\overline{c}_2^2+2d_*}$.
In particular it is smaller than $N^{-A_0-2}$ and than
$\varepsilon^{-2}\Sigma_t^{-2^6}$ in every round.
We include $\overline{C}\ge2^{A_1}$ among the defining conditions on $\overline{C}$.
Then the bit lengths and all elementary interval computations are
bounded by $N^{d_*}$, as assumed above. Products of rational
probabilities in \eqref{eq:alg-monomial-oracle} can be evaluated
exactly; no independence approximation is made in these expectations.

For each computed upper estimate $s_x^+$, let $\widehat s_x^+$
be its rational upper endpoint, and let $\widehat p_x$ be a
rational lower approximation to $\overline{\omega}_t/\widehat s_x^+$.
We denote the error bound just obtained by
\[
 \epsilon_{\rm num}=N^{-A_1+\overline{a}_6+\overline{c}_2^2+2d_*}.
\]
These choices satisfy
\[
 s_x^+\le\widehat s_x^+\le s_x^++\epsilon_{\rm num},\qquad
 0\le\frac{\overline{\omega}_t}{\widehat s_x^+}-\widehat p_x\le \epsilon_{\rm num}.
\]
The algorithm uses $\mathbb P(\mathsf w_x=0)=\widehat p_x$.
Since $s_x^+\ge1/8$ and $\overline{\omega}_t\le1$, we have
\[
 0\le\frac{\overline{\omega}_t}{s_x^+}-\widehat p_x
 \le\frac{\overline{\omega}_t(\widehat s_x^+-s_x^+)}{s_x^+\widehat s_x^+}+\epsilon_{\rm num}
 \le2^7\epsilon_{\rm num}.
\]
For a degree test, the rational approximation
$\widetilde\mu_F$ is the one specified after
\eqref{eq:alg-degree-test}. We take its rational deviation
threshold $\widetilde\Delta_i=(3\Delta_i/4)^-$. Its error is
at most $\Delta_i/16$, so
$11\Delta_i/16\le\widetilde\Delta_i\le3\Delta_i/4$. Then
\[
 |F-\widetilde\mu_F|\le\widetilde\Delta_i
 \quad\Longrightarrow\quad
 |F-\mathbb E_{Z,v}^0(F)|\le3\Delta_i/4+\Delta_i/8<\Delta_i.
\]
For an upper test $F\le u$, the rational threshold lies in
$[u-\epsilon_{\rm num},u]$; for a lower test $F\ge l$, it lies in $[l,l+\epsilon_{\rm num}]$.
These choices imply the required exact inequalities. The bound
$2^7\epsilon_{\rm num}$ is less than one eighth of each permitted deviation, so
the moment estimates in \eqref{eq:alg-focus-degree-ratio},
\eqref{eq:alg-degree-moment-ratio}, and
\eqref{eq:alg-coverage-moment-ratio} still give
$\mathbb E(S_B)\le\Sigma_t^{-a}$ after rounding.

The total additional error in the bound for $\Phi_{t+1}$ is at
most $\varepsilon^{-2}\Sigma_t^{-2^6}$. With
\[
 d_t=r_t+\frac1{n\pi_t}+(\Delta\tau_t)^2+\frac{\Delta\tau_t}{M},
\]
the last inequality in the proof of
Proposition~\ref{prop:alg-one-round} becomes
\[
 \Phi_{t+1}-\Phi_t
 \le\frac{\overline{c}_2^4+6\overline{c}_2^3+1-\overline{c}_2^8}{\varepsilon^2}d_t
 \le0,
\]
because $d_t\ge\Sigma_t^{-2^6}$.

We now specify how the stopping indices are implemented. The intervals
for the reference trajectory are computed before any of the actual
sets are selected. Let
\[
 \widehat{\overline{T}}_{\mathrm{init}}
   =\min\{t:(f_0(\iota_t))^->2^{4r}\},\qquad
 \widehat{\overline{T}}
   =\min\{t:\Sigma_t^+<(\sqrt{\log D})^-\}.
\]
These computable indices need not equal the exact indices in
\eqref{eq:alg-stopping}. We use them to give a finite list containing
the exact pair, rather than running an additional round with the
old terminal coefficients.

To bound the list, write
$\beta_i=d_{i,t}^{1/(i-1)}/\Sigma_t$, so that $\sum_i\beta_i=1$.
The generating identity gives, in the neighborhood of the initial
threshold,
\[
 f_0(\iota_{t+1})-f_0(\iota_t)
   =\sum_{i=2}^r(\overline{\eta}_t\sigma_t\beta_i)^{i-1}
   \ge\left(\frac{\overline{\eta}_t\sigma_t}{r-1}\right)^{r-1}
   \ge N^{-2r}.
\]
Here $\Sigma_t\le N^2$, $\overline{\eta}_t\ge N^{-1}$ and
$\sigma_t\ge1/2$; the lower bound on $\overline C$ ensures
$N\ge2r$. Near the terminal threshold,
Lemma~\ref{lem:alg-tracking}\,\textup{(\ref{it:alg-scale-ratios})}
similarly gives
\[
 \Sigma_t-\Sigma_{t+1}
 \ge\Sigma_t(1-e^{-\overline{\eta}_t/8})
 \ge\frac{\Sigma_t\overline{\eta}_t}{16}\ge N^{-2}.
\]
These reference estimates remain valid for the single further
evaluation required to find the displayed indices: the scale is still
at least $e^{-1}\sqrt{\log D}$ and the sampling parameter is at most
$1/4$. The additional evaluation changes the sensitivity product by
at most $N$, and the explicit constant in
\eqref{eq:alg-large-constant} satisfies
\[
 N^{\overline c_2(2^{16r^2+4}+16)+1}\le N^{\overline c_2^2}.
\]
Thus the same propagated error bound $\epsilon_{\rm num}$ applies.
Since $4\epsilon_{\rm num}<N^{-2r}$, two consecutive reference
values cannot both have intervals meeting the same threshold.
It follows, including when equality holds at an exact threshold, that
\[
 \overline{T}_{\mathrm{init}}
 \le\widehat{\overline{T}}_{\mathrm{init}}
 \le\overline{T}_{\mathrm{init}}+1,
 \qquad
 \overline{T}\le\widehat{\overline{T}}\le\overline{T}+1.
\]

We run the construction separately for each pair of candidate indices
$t_{\mathrm{init}}\in\{\widehat{\overline{T}}_{\mathrm{init}}-1,
\widehat{\overline{T}}_{\mathrm{init}}\}$ and
$t_{\mathrm{end}}\in\{\widehat{\overline{T}}-1,
\widehat{\overline{T}}\}$ with $0\le t_{\mathrm{init}}\le t_{\mathrm{end}}$.
In each run the actual sets start again from
$J^0=\mathsf W^0=\varnothing$; the phase change occurs at
$t_{\mathrm{init}}$ and the terminal index is $t_{\mathrm{end}}$.
All terminal quantities are recomputed for this index: in particular,
$\tau_*=\tau_{t_{\mathrm{end}}}$, the definition of $b_t$ uses
$\pi_{t_{\mathrm{end}}}$ and a sum ending at $t_{\mathrm{end}}-1$,
and the sum defining $R_t$ ends at $t_{\mathrm{end}}-1$.
Thus no run uses coefficients for one endpoint while stopping at another.
There are at most four runs, and one uses exactly
$(t_{\mathrm{init}},t_{\mathrm{end}})
=(\overline{T}_{\mathrm{init}},\overline{T})$.

A run is discarded if a required rational operation is undefined,
a proposed probability lies outside $[0,1]$, or it fails to return an output within
$N^{\overline a_6+3}$ bit operations. The preceding bounds guarantee
that the run with the exact pair finishes within this bound and
satisfies Lemmas~\ref{lem:alg-tracking}--\ref{lem:alg-potential-endpoints}
with the numerical errors already estimated above. We do not assume
these conclusions for the other runs. Instead, every returned set is
checked directly for independence and maximality by scanning the
original incidence list and testing the addition of each omitted atom.

We also check its size by rational inequalities. Let
$\ell=(1-\varepsilon)M$ and $u=(1+\varepsilon)M$, and compute rational
endpoints with
\[
 0\le\ell^+-\ell\le\frac{\varepsilon M}{16},\qquad
 0\le u-u^-\le\frac{\varepsilon M}{16}.
\]
The accepted output is the first maximal independent set $I$ satisfying
$\ell^+\le|I|\le u^-$.
Such an output exists. In the run with the exact pair,
$\Phi_{\overline T}\le\Phi_0<1-\varepsilon/2$ and the proof of
Lemma~\ref{lem:alg-potential-endpoints} gives the stronger inequalities
\[
 \left(1-\frac{3\varepsilon}{8}\right)M
 \le |I|
 <(1+\varepsilon)(1-\varepsilon/2)M
 \le\left(1+\frac{\varepsilon}{2}\right)M.
\]
These imply both rational size tests. The precision required for the
tests is covered by $A_1$, since $M\ge n/N$ and
$\varepsilon\ge N^{-1}$ for the stated range of $N$.
Every accepted output therefore satisfies the original theorem,
regardless of which candidate run produced it.

The four bounded runs, their reference computations, and the direct
output checks take at most $5N^{\overline a_6+3}$ bit operations.
The specified lower bounds on $\overline C$ already ensure $n\ge3$.
Finally, $N\le r^2n^r\le n^{r+2}$ by \eqref{eq:alg-input-n}, so
$5N^{\overline a_6+3}\le n^{(r+2)(\overline a_6+3)+2}$.
This is the stated bound $n^{\overline c_4(r,\delta_0,c)}$.
\end{proof}

\section*{Concluding Remarks}

$\bullet$ Theorem~\ref{main} gives, for fixed $r\ge1$ and $\delta>0$,
a maximal independent set of size asymptotic to
$n((\log d)/(rd))^{1/r}$ in any $d$-regular $\delta$-sparse
$(r+1)$-uniform set system with $n$ atoms, as $d\to\infty$. Together with Theorem~\ref{thm:random-regular-delta-sparse}, this determines the asymptotic minimum size of a maximal independent set in almost all such systems: first $n\to\infty$ through admissible orders with $d$ fixed, and then $d\to\infty$. 
A salient open question is to determine for each $r \geq 1$ whether there exists a fixed $c > 1$ there is a polynomial time algorithm for finding an independent set of size $c((\log d)/rd)^{1/r} \cdot n$ even in the setting of linear $(r + 1)$-uniform $d$-regular hypergraphs and random regular hypergraphs. The graph case $r = 1$ is a well-known open question of Karp~\cite{Karp1976}.

$\bullet$ We left the open problem of determining geometric and algebraic instances where the logarithmic term supplied by Theorem \ref{thm:main} is necessary. The following question was posed in~\cite[Problem 2]{SchillewaertV}; we state it as a conjecture:

\begin{conjecture*}\nonumber
Every maximal partial ovoid in the elliptic quadric $Q^-(5,q)$ has size at least $q\omega(q)$ for some function $\omega(q)\to\infty$ as $q\to\infty$.  
\end{conjecture*}

$\bullet$ It is not clear what the most general 
theorem for maximal $f$-free sets over $\mathbb F_q$ is. Theorem \ref{thm:poly-maximal-sets} requires $f$ to be symmetric and homogeneous, which seem to be  rather strong conditions, given that the main theorem works for $f(x,y,z) = x^2 + y^2 - z^2$ and various linear polynomials. It may be checked that Theorem \ref{thm:main} may be applied to many other polynomials (non-homogeneous and asymmetric). A salient open conjecture is as follows:

\begin{conjecture*}
For each $r \geq 3$, there exists a symmetric homogeneous polynomial $f \in \mathbb Z[x_1,x_2,\dots,x_r]$ that is absolutely irreducible over $\mathbb Q$, satisfies $f(1,\ldots,1)\ne0$ over $\mathbb Z$, and has
$\{\gamma_q(f)q^{-1/(r - 1)}\}$ unbounded over primes $q$. 
\end{conjecture*}

$\bullet$ Another very interesting question is to determine, for each fixed $r\ge3$, whether there exists a polynomial $f \in \mathbb Z[x_1,x_2,\dots,x_r]$ with $f(x,x,\ldots,x)\not\equiv0$ such that every $f$-free subset of $\mathbb F_q$ has size at most $q^{1/(r - 1) + o(1)}$ as $q \rightarrow \infty$ through primes for which the diagonal polynomial remains nonzero.

$\bullet$ The results given by Theorem \ref{thm:steiner-maximal} for Steiner systems can be extended to combinatorial designs, by simply applying Theorem \ref{thm:main} with different parameters.

$\bullet$ A very interesting open problem arises for maximum independent sets in set systems which are not necessarily $\delta$-sparse, but where every block contains few subsets that violate the first $\delta$-sparseness condition. Some special cases were raised by Cooper and Mubayi~\cite{CooperMubayi}. Specifically, suppose we have a triple system $H$ of maximum degree $d$ such that in every block of $H$, at most one pair is contained in another block of $H$. Is it true that $H$ contains an independent set of size $n\omega(d)/\sqrt{d}$ where $\omega(d) \rightarrow \infty$? 

$\bullet$ The following was raised in~\cite{EustisV} (see also~\cite{KMV}): 

\begin{conjecture*}
Every $n$-vertex Steiner triple system has an independent set of size $(1 - o(1))\sqrt{3n\log n}$ as $n \rightarrow \infty$.
\end{conjecture*}

\section*{Acknowledgement}

This work was completed during multiple visits of the first author to the University of California, San Diego (UCSD), 2010--2026. We would like to thank UCSD for the conducive work environment. 

\section*{AI Declaration}

The theorems and their proofs were obtained by the authors before the advent of AI, and without the aid of AI\@. We used ChatGPT to check the extensive list of constants and the manuscript for typographical errors.
The authors take full responsibility for the mathematical content of the paper.

\end{document}